\documentclass[12pt]{article}
\usepackage{amsfonts, amsmath, amssymb}
\usepackage{tikz}

\usepackage{amsthm}
\usepackage{hyperref}
\newtheoremstyle{mythmstyle}%
{}{}%
{\itshape}{}%
{\bfseries}{}%
{ }%
{\thmname{#1}\thmnumber{ #2}\thmnote{ (#3)}}
\theoremstyle{mythmstyle}

\newtheorem{lemma}[subsubsection]{Lemma}
\newtheorem{proposition}[subsubsection]{Proposition}
\newtheoremstyle{mydefstyle}%
{}{}%
{}{}%
{\bfseries}{}%
{ }%
{\thmname{#1}\thmnumber{ #2}\thmnote{ (#3)}}
\theoremstyle{mydefstyle}
\newtheorem{definition}[subsubsection]{Definition}
\newtheorem{remark}[subsubsection]{Remark}
\newtheoremstyle{sssectionstyle}%
{}{}%
{}{}%
{\bfseries}{}%
{ }%
{\thmnumber{#2}\quad\thmnote{ #3}}
\theoremstyle{sssectionstyle}
\newtheorem{sssection}[subsubsection]{}
\newcommand{\tocexclude}[1]{
	\addtocontents{toc}{\protect\setcounter{tocdepth}{-1}} #1 \addtocontents{toc}{\protect\setcounter{tocdepth}{3}}
}
\newcommand{\sssectionnum}
{\stepcounter{subsubsection}{\bf \thesubsubsection}}
\newcommand{\apx}[1]{\underset{#1}{\approx} }
\newcommand{\Sim}[1]{\underset{#1}{\sim} }
\newcommand{\mt}[1]{$\left( \mathbb M^t \right)^{#1}$}
\newcommand{\simt}{$\underset{t}{\sim}$ \ }
\newcommand{\simHt}{$\underset{H_t}{\sim}$ \ }
\newcommand{\simBt}{$\underset{B_t}{\sim}$ \ }
\newcommand{\diamondMove}[1]{
	\begin{tikzpicture}[scale=.3]
		\draw (1,0)--(2,1)--(1,2)--(0,1)--(1,0);
		\node at (1,1) {$#1$};
	\end{tikzpicture}
}
\newcommand{\diamondmove}[1]{
	\begin{tikzpicture}[scale=.25]
		\draw (1,0)--(2,1)--(1,2)--(0,1)--(1,0);
		\node at (1,1) {$#1$};
	\end{tikzpicture}
}
\begin{document}
	
	\begin{center}
		
		{\Large On the Isoperimetric functions of a class of Artin Groups} 
		
		Arye Juh\'{a}sz,\\
		Department of Mathematics\\
		The Technion -- Israel Institute of Technology\\
		Haifa 32000, ISRAEL
		
	\end{center}
	
	\vspace{10pt}
	
	{\bf Abstract}
	
	{\it
		\noindent
		We find a polynomial ($n^6$) isoperimetric function for Artin groups, the defining graph of which contains no edges labelled by 3. This in particular shows that even Artin groups have solvable word problem.
		We use small cancellation theory of relative extended presentations.
	}

	
	\section*{Introduction}
	\addtocontents{toc}{\protect\contentsline{section}
		{Introduction}
		{\thepage}{}
	}
	Let $m,n \in {\mathbb N}$ ($\mathbb N$ - the natural numbers) 
	$m \geq 0$, $n \geq 1$ and let $\Gamma$ be a simple graph without loops, 
	with vertex set $V =\{v_1 , \ldots , v_n\}$ and edge set 
	$E = \{e_1, \ldots , e_m\}$. Label the edges by natural numbers via a labelling function $\lambda : E \rightarrow {\mathbb N}\setminus \{ 1 \}$. 
	Denote $\lambda (e) = n_{ij}$, where the edge $e$ connects 
	$v_i$ to $v_j$, $i \neq j$ ($n_{ij}=n_{ji}$), and let $n_{ij}=0$ if $v_i$ and $v_j$ are not connected by an edge.
	For every such labelled graph correspond a group presentation 
	\[\tag{\rm I}
	{\cal P}(\Gamma) = \langle X \ | \  {\cal R} \rangle
	\]
	such that 
	
	$X= \{x_v, \ v \in V \}$, ${\cal R} = \{R_e \ , \ e \in E\}$, denote $x_{v_i}$ by $x_i$,
	
	$R_e=(x_ix_j)^{\frac{1}{2}n_{ij}} \left(x_i^{-1}x_j^{-1}\right)^{\frac{1}{2}n_{ij}}$, if $n_{ij}$ is even and 
	
	$R_e=(x_ix_j)^{\frac{n_{ij}-1}{2}} x_i  
	\left(x_i^{-1}x_j^{-1}\right)^{\frac{n_{ij}-1}{2}}x_j^{-1}$, if $n_{ij}$ is odd.

	The group presented
	by ${\cal P}(\Gamma)$ is denoted  by $A(\Gamma)$ 
	and is called the
	\emph{Artin group defined by $\Gamma$}.
	
	\vspace{10pt}
	
	\noindent {\bf Examples 1}
	\begin{enumerate}
		\item $n_{ij}=0$ for every $i,j$, $1 \leq i < j \leq n$.\\
		Then $A(\Gamma)$ is the free group on $X$.
		\item $n_{ij}=2$ for every $i,j$, $1 \leq i < j \leq n$.\\ 
		Then $A(\Gamma)$ is the free abelian group on $X$.
		\item $n_{ij}\in \{0,2\}$ for every $i,j$, $1 \leq i < j \leq n$.\\ 
		Then $A(\Gamma)$ is called \emph{right angled Artin group}.
		\item $n_{ij} \in 2{\mathbb N}$ for every $i,j$, $1 \leq i < j \leq n$. \\
		Then $A(\Gamma)$ is called an \emph{even Artin group}.
		\item $n_{ij} \neq 2$ for every $i$, $j$, $1 \leq i < j \leq n$.\\
		Then $A(\Gamma)$ is called of \emph{large type}.
	\end{enumerate}

	To each Artin group $A(\Gamma)$ there corresponds a Coxeter group $W_A$, obtained by adding the relators $x_i^2$, $i=1, \ldots , n$. An Artin group is called of 
	\emph{spherical type} if $W_A$ is a finite group. The graphs of the spherical type Artin groups are classically known. An Artin group is said to be of \emph{FC type} if it belongs to the smallest class of Artin groups which are closed under amalgamation along standard parabolic subgroups and contain all the spherical type Artin groups.
	Standard parabolic subgroups are the subgroups which are generated by subsets of $X$.
	
	We recall isoperimetric functions. Let $G$ be a group presented by
	${\cal P} = \langle X \ | \ {\cal R} \rangle$. 
	Let $W \in F(X)$, $F(X)$ the free group freely generated by $X$, $W$ cyclically reduced (i.e. $WW$ contains neither $xx^{-1}$ nor $x^{-1}x$, $x \in X$). Then $W$ represents 1 if and only if
	
	\[\tag{\rm II}
	W = C_1 \cdots C_{k_w} \ , \ k_w \geq 1 \ , \ C_i=f_iR_i^{\varepsilon_i}f_i^{-1}\ , \]
	\[ f_i \in F(X) \ , \ 
	R_i \in {\cal R} \ \mbox{and} \ \varepsilon \in \{1,-1\}
	\]
	
	A function 
	$g:{\mathbb N} \rightarrow {\mathbb N}$ is an \emph{isoperimetric function} for the presentation ${\cal P} = \langle X  \ | \ {\cal R} \rangle$ if for every word $W$ which represents 1 of the group $G$, $k_W$ in ({\rm II}) satisfies $k_W \leq g(|W|)$.
	Here $|W|$ denotes the length of $W$ in $F(X)$.

	The isoperimetric functions are known for spherical type, FC-type, large type, right angled Artin groups and some other classes of Artin groups. See \cite{3},\cite{8} and \cite{9}.
	
	\vspace{10pt}

	\subsection*{Main Theorem}
	\addtocontents{toc}{\protect\contentsline{subsection}
		{Main Theorem}
		{\thepage}{}
	}
	
	\noindent {\it Let $A$ be an Artin group, $A = A(\Gamma)$. If $n_{ij}\neq 3$ 
		for every $n_{ij}$, $1 \leq i < j \leq n$, then $f(n) = n^6$ 
		is an isoperimetric function for ${\cal P}(\Gamma)$.}
	
	We believe that "$6$" can be replaced by "$2$".
	
	\vspace{10pt}
	
	\noindent {\bf Corollary}\\
	{\it Let $A(\Gamma)$ be an Artin group. If $n_{ij} \neq 3$ for every $i,j$, $1 \leq i < j \leq n$, then $A$ has solvable word problem. In particular, even Artin groups have solvable word problem.}
	
	\vspace{10pt}
	
	The  Corollary  has been proved recently also by Ruben Blasco-Garcia, Maria  Cumplido  and Rose Morris-Wright in arXiv:2204.03523, using a completely  different  method.

	Our basic method is small cancellation theory of extended (relative) presentations via van Kampen and Howie diagrams, see 1.2.2 and 1.2.3.
	(Recall that a van Kampen $\cal R$-diagram over $F(X)$ is a connected, simply connected labelled planar 2-complex, labelled by elements of $F(X)$ such that the labels of the boundary of the 2-cells are elements of $\cal R$. We say that the corresponding 2 cell (region) \emph{realises} the element of $\cal R$.)
	
	\vspace{10pt}
	
	This  work consists of the details of my talk   I  gave in GAGTA  2019.

	\subsection*{Outline of the proof of the Main Theorem}
	\addtocontents{toc}{\protect\contentsline{subsection}
		{Outline of the proof of the Main Theorem.}
		{\thepage}{}
	}
	
	In order to find an isoperimetric function for the presentation $\cal P$ we have to find an upper bound on the number $k_w$ in ({\rm II}) in terms of the length $|W|$ of $W$, for every $W$ representing 1 in the group defined by $\cal P$. By the basic theorem of van Kampen diagrams, for every cyclically reduced  $W$ representing 1 there is a van Kampen diagram $M$ with boundary label $W$ such that the number of regions in $M$ is at most $k_W$.
	See \cite[Ch. V]{14}.
	So our problem can be considered as a counting problem in van Kampen diagrams and this is our approach. We use the following simple basic principles, together with the well known result in Proposition 1 below. Thus, let $S$ be a finite set. In order to count the elements of $S$ we subdivide $S$ into subsets $S_1, \ldots S_m$ such that we know $m$ and know the maximal possible number $n$ of elements in $S_i$, $ i = 1, \ldots , m$. Then $S$ has at most $mn$ elements. In our case Propositon 1 provides the first approximation to the numbers $m$ and $n$, in terms of the length of $W$.

	Before we go on we introduce a few notations. For a set $S$ denote by $|S|$ the number of elements  of $S$.
	For a word $W$ in $F(X)$ denote by $|W|$ the number of letters from $X \cup X^{-1}$ (with multiplicity) that occur in $W$. 
	This is the length of $W$.
	For $t \in X$ denote by $|W|_t$ the number of occurences of $t^{\pm 1}$ in $W$. 
	We denote by $||W||$ the syllable length (free product length) of $W$
	and $||W||_t$ denotes the number of disjoint occurences of $t^{\alpha_i}$ in the word $W$, $t \in X$,
	$0 \neq \alpha \in {\mathbb Z}$.
	Thus, $|W| = \sum\limits_{t \in X} |W|_t$ and $||W|| = \sum\limits_{t \in X}||W||_t$.
	We denote by $Supp(W)$ the set of the letters of $X$ contained in $W$ and $W^{-1}$.
	
	\ \\ 
	\noindent We turn now to diagrams. In general we follow \cite[Ch. V]{14} for notations and definitions. \\
	For a region $D$ in an $\cal R$-diagram over $F$ corresponding to the presentation given by ({\rm I}) we denote by $n(D)$ the number $n_{ij}$, where $D$ realises $R_e$, $\lambda(e)=n_{ij}$.
	Let $M$ be a van Kampen $\cal R$-diagram over $F$.
	We denote by $\Phi :M \rightarrow F(X)$ the labelling function of the diagram $M$ over $F(X)$.
	If $D$ is a region in a diagram, we shall write $Supp(D)$ for $Supp(\Phi(\partial D))$,
	where $\partial D$ denotes the boundary of $D$.
	If $Supp(R) = \{a,t\} \subseteq X$, $ R \in {\cal R}$, we shall denote $R$ by $R(a,t)$, when convenient.
	
	Denote by $Reg (M)$ the set of the regions of $M$. Let $\{D_1,\ldots,D_k\}$ be a set of regions of $M$. \emph{The subdiagram generated} by \emph{$\{D_1,\ldots,D_k\}$} is the smallest subdiagram which contains $\{D_1,\ldots,D_k\}$. Notation $\langle D_1,\ldots,D_k \rangle$. Thus $Reg (\langle D_1,\ldots,D_k \rangle)=\{D_1,\ldots,D_k\}$. Concerning Artin groups with 
	$n_{ij} \neq 3$ we subdivide $Reg(M)$ into $Reg(M)=Reg_2(M) \cup Reg_{4^+}(M)$,
	where $Reg_2(M)$ is the set of all the regions which realise a commuting relation and $Reg_{4^+}(M)$ is the set of the regions which realise a relation $R$ with length $\geq 8$. ($n_{ij} \geq 4$). Denote by $T(M)$ the set of all $x \in X$ such that $x$ occurs in a relation realised by a region in $Reg_{4^+}(M)$.
	\par \ \par
	\noindent {\bf Definition 1:}
	\addtocontents{toc}{\protect\contentsline{subsection}
		{Definition 1\ (C(4)\&T(4))}
		{\thepage}{}
	}
	Recall from \cite[p. 240]{14} that a presentation $\langle X|\cal{R}\rangle$ satisfies the condition $C(p), p\geq 1$, if no element of $\cal{R}$ is the product of fewer than p pieces. $R$ satisfies the condition $T(q), q\geq 1$ if the following holds: 
	Let $3\leq h<q$. Suppose $R_1,\ldots,R_h$ are elements of $\cal{R}$ with no successive elements $R_i$ and $R_{i+1}$ which form a reducing pair (i.e. $R_i=R^{-1}_{i+1}$). Then at least one of the products $R_1 R_2\cdots R_{h-1}R_h, R_h R_1$ is reduced without cancellation. It follows that if a van Kampen 
	diagram $M$ satisfies the condition C(4) \& T(4) then every inner vertex has valency at least 4 and every region 
	which has no edges on the boundary $\partial M$ of $M$,
	has at least 4 neighbouring regions.\\ One of the key results on which the proof of the Main Theorem lies is the following Proposition \cite[Ch. V]{14}.
	
	\vspace{20pt}
	
	\noindent {\bf Proposition 1}
	\addtocontents{toc}{\protect\contentsline{subsection}
		{Proposition 1\ (Properties  of C(4)\&T(4) maps needed)}
		{\thepage}{}
	}
	{\it Let $M$ be a reduced van Kampen $\cal R$-diagram over $F(X)$. 
		Let $W$ be a boundary label of $M$, cyclically reduced. 
		Assume $M$ satisfies the condition} C(4) \&  T(4). {\it 
		Suppose that each piece has length 1.
		Then each of the following holds: 
		\begin{itemize}
			\item[(a)] For each region $\Delta$ of $M$ we have $|\partial \Delta | \leq |W|$
			\item[(b)] The number $|M|$ of regions $\Delta$ in $M$ is at most $|W|^2$
	\end{itemize}}
	
	Hence, if our diagram $M$ would be a C(4) \& T(4) diagram, then 
	by part (b) of the Proposition the number of the regions in $M$ would be bounded by $|W|^2$.
	Being this  not the case, we first modify our approach. The core idea is 
	not to try to estimate the number of all the regions via a result like Proposition 1, but to use that result for the estimation of a part of the regions only, namely, 
	$|Reg_{4^+}(M)|$ and then relate the number of the remaining regions to $|Reg_{4^+}(M)|$.

	\vspace{10pt}
	
	In this direction we have
	
	\noindent {\bf Theorem 1:} 
	\addtocontents{toc}{\protect\contentsline{subsection}
		{Theorem 1\ (upper bound on the number of non-commuting relations)}
		{\thepage}{}
	}
	{\it Let $M$ be a van Kampen $\cal R$-diagram, $\cal R$ given by ({\rm I}). Assume 
		\[ \tag{\rm III}
		{\cal R}=\{R_e \ | \ e \in E \} \   , \ \lambda(e) \neq 3 \ \mbox{for every} \ e \in E
		\]
		Then there is a van Kampen diagram $M'$ with the same boundary label as $M$ such that $|Reg_{4^+}(M')| \leq |\partial M|^4$
	}
	
	\vspace{10pt}
	
	\noindent To prove,
	we proceed as follows. Let $W$ be the word mentioned in ({\rm II}).
	
	\noindent First we "cover" $Reg_{4^+}(M)$ with subsets $Reg_{4^+, t}(M)$, $ t \in T(M)$, defined by $Reg_{4^+, t}(M) = \{ D \in Reg_{4^+}(M), \ t \in Supp(D) \}$. Clearly,
	\[\tag{1} Reg_{4^+}(M) = \bigcup\limits_{t \in T(M)} Reg_{4^+, t}(M) \ , \
	|Reg_{4^+}(M)| \leq \sum\limits_{t \in T(M)} |Reg_{4^+,t} (M)|.
	\]
	
	Hence, it is enough to estimate  $Reg_{4^+,t}(M)$. We define in the First Construction below, an equivalence relation  on $Reg_{4^+, t}(M)$ which has the following properties.
	\begin{itemize}
		\item[($i$)] Let $[D_1], \ldots , [D_{r_t}]$ be the equivalence classes of $Reg_{4^+, t}(M)$. Then $\Delta (D) := Int\left( \cup\overline{E} \ | \ E \in [D] \right)$ is connected and simply connected. When $D$ is clear from the context we shall write $\Delta$ for $\Delta(D)$.
		$Int$ stands for Interior and
		$\overline{E}$ denotes the closure of $E$ in ${\mathbb E}^2$.
		\item[($ii$)] $|[D]| \leq |\partial \Delta (D)|^2$
	\end{itemize}
	
	Thus,
	\[\tag{2}
	\mbox{$|Reg_{4^+, t}(M)| \leq r_tp_t$, where $p_t = \max|[D_i]|$}, \ i = 1, \ldots r_t
	\]

	Hence, it is enough to estimate $r_t$ and $p_t$ in terms of $|W|$. 
	
	\vspace{10pt}
	
	\noindent {\bf Proposition 2}
	\addtocontents{toc}{\protect\contentsline{subsection}
		{Proposition 2\ (The map $N_t, t\in X$)}
		{\thepage}{}
	} 
	{\it For every $t \in T(M)$ there is a map $N_t$ constructed out of $M$ with the following properties:}
	{\it
		\begin{itemize}
			\item[(a)] $|\partial N_t| \leq |W|_t$
			\item[(b)] $|Reg(N_t)| = r_t$
			\item[(c)] For every $\Delta (D)$ $N_t$ contains a region $\tilde{\Delta}$ with $|\partial \Delta|_t = |\partial \tilde{\Delta}|$
			\item[(d)] $N_t$ satisfies the small cancellation condition  C(4) $\&$ T(M)(4) 
		\end{itemize}
	}
	Consequently, by Proposition 1 \  $r_t \leq \left(|W|_t\right)^2$ and $|\Delta| \leq \left( |W|_{t_1}+ |W|_{t_2}\right)^2 \leq |W|^2$,
	where $Supp(\Delta) = \{t_1, \ t_2 \} \subseteq T(M)$. Thus,
	
	\[\tag{3}
	r_t \leq (|W|_t)^2 \ \mbox{and} \  p_t \leq |W|^2.
	\]
	
	Hence, by (1), (2) and (3) 
	\[|Reg_{4^+}(M)| \leq \left(\sum\limits_{t \in T} |W|_t^2\right) |W|^2 \leq \left(\sum\limits_{t \in T}| W_t|\right)^2 |W|^2 \leq |W|^2 \cdot |W|^2,\]
	i.e.
	$|Reg_{4^+,}(M)| \leq |W|^2\cdot |W|^2=|W|^4$, proving Theorem 1.
	
	\vspace{10pt}
	
	For the proof of Theorem 1,Propositions 1,2 and the Main Theorem we carry out three constructions of diagrams, and a variant of the first construction.

	\subsection*{The first construction: $\Delta(D)$ and $\mathbb M$}
	\addtocontents{toc}{\protect\contentsline{subsection}
		{The first construction: $\Delta(D)$ and $\mathbb M$. }
		{\thepage}{}
	} 
	{\bf Definition 2:}
	\addtocontents{toc}{\protect\contentsline{subsection}
		{Definition 2 (Equivalence)}
		{\thepage}{}
	} 
	Let $D_1$ and $D_2$ be regions of $M$. Say that $D_1$ and $D_2$ are \emph{friends} if $D_1$ and $D_2$ share a common edge and $Supp(D_1)=Supp(D_2)$.
	(We consider $D_1$ as a friend of itself.)
	The transitive closure of this relation is an equivalence relation "$\sim$" on the set of regions of $M$. For a region $D$ of $M$ denote by $[D]_M$ the equivalence class of $D$ in $M$ and define 
	$\Delta(D) = Int(\bigcup\limits_{E \in [D]_M} \overline{E})$, as before. We shall call $[D]_M$ an "Equivalence class" if $D\in Reg_4^+(M)$.

	We would like to produce from $M$ a new diagram  with the same boundary label as $M$ which is easier to handle. The idea is to use $\Delta(D)$ in place of $D$. To this end we have to show that $\Delta(D)$ is homeomorphic to the open unit disc.
	Then the totality of $\Delta(D)$ for $D\in Reg_{4^+}(M)$, together with $Reg_2(M)$ generates a diagram which we denote by $\mathbb M$. Indeed, $\mathbb M$
	is a diagram, but it is not better to our purpose than $M$.
	However, there is a variant of this construction which does work. (See Section 2.)
	
	\

	\vspace{10pt}
	

	\subsection*{The second construction $M^t$ and $\mathbb{M}^t$, $t \in T$}
	\addtocontents{toc}{\protect\contentsline{subsection}
		{The second construction $M^t$ and $\mathbb{M}^t$}
		{\thepage}{}
	} 
	Recall that a relative presentation ${\mathbb P} = \langle H, Y  \ | \ {\cal R}\rangle$ consists of a group $H$, a set $Y$ with $Y \cap H =\emptyset$ and a set of words $\cal R$ in 
	$H * \langle Y \rangle \setminus H$. The group presented by $\mathbb P$ is isomorphic to
	$H * \langle Y \rangle / N$, where $N$ is the normal closure of $\cal R$ in $H * F(Y)$.
	(See {\cite{7} and \cite{12}.)

		Consider now a van Kampen $\cal R$-diagram $M$ with connected interior $\cal R$, as in {\rm (III)}.
		For $t \in T$ let $X_t = X \setminus \{t\}$,  let $H_t = \langle X_t \rangle$, and let
		${\mathbb P}_t =\langle H_t, t \ | \ {\cal R}_t \rangle$,
		where ${\cal R}_t$ is the set of all the relations in $\cal R$ which contain $t$. We define $M^t$ as follows: 
		\begin{enumerate}
			\item Each boundary edge of each region $D$ in $M$ which is labelled by words from $H_t$ is shrunk to a point so that if $Supp (D) = \{t,a\}$, $a \in X_t$ and
			
			$\Phi (\partial D) =R_e$ then the region obtained is a $\lambda (e)$-gon denoted by $D^t$, with edges labelled by $t^{\pm 1}$ ($t$-edges) and
			\item  the corners between two consequtive $t$-edges are labelled by $a^{\pm 1}$, accordingly. See Figure 1.
		\end{enumerate}

		\vspace{10pt}
		
		\begin{minipage}{2in}
			\begin{center}
				\begin{tikzpicture}[scale=.3]
					\draw (2,0)--(5,0)--(7,3)--(7,6)--(5,9)--(2,9)--(0,6)--(0,3)--(2,0);
					
					\draw (3.7,.3)--(4,0)--(3.7,-.3);
					\node at (4,.7) {$a$};
					
					\draw (3.7,9.3)--(4,9)--(3.7,8.7);
					\node at (4,8.5) {$a$};
					
					\draw (-.3,4.7)--(0,5)--(.3,4.7);
					\node at (.5,5) {$a$};
					
					\draw (6.7,4.7)--(7,5)--(7.3,4.7);
					\node at (6.5,5) {$a$};
					
					\draw (.5,7.2)--(.8,7.2)--(.8,6.9);
					\node at (1.5,7.5) {$t$};
					
					\draw (5.7,7.7)--(5.7,8)--(6,8);
					\node at (5,7.5) {$t$};
					
					\draw (1,1)--(1,1.5)--(1.5,1.5);
					\node at (1.5,2) {$t$};
					
					\draw (5.7,1.5)--(6,1.5)--(6,1.2);
					\node at (5.5,2) {$t$};

				\end{tikzpicture}
			\end{center}
		\end{minipage}$\bf \longrightarrow$ \hspace{10pt}\begin{minipage}{2in}
			\begin{center}
				\begin{tikzpicture}[scale=.3]
					\draw (0,5)--(3,10)--(8,2)--(4,-4)--(0,5);
					
					\draw (1.3,1.3)--(1.6,1.6)--(1.9,1.3);
					\node at (2,2) {$t$};
					\draw (6.2,-.3)--(6.5,-.3)--(6.6,-.5);
					\node at (5.9,0) {$t$};
					\draw (5.4,5.6)--(5.6,6)--(6.1,5.6);
					\node at (5,6) {$t$};
					\draw (1,7)--(1.3,7)--(1.3,6.7);
					\node at (2,7) {$t$};
					
					\node[right] at (0,5) {$a$};
					
					\node[below] at (3,9.5) {$a$};
					
					\node[left] at (8,2) {$a^{-1}$};
					
					\node[above] at (4.2,-3) {$a^{-1}$};

				\end{tikzpicture}
			\end{center}
		\end{minipage}
		
		\begin{center}
			Figure 1
		\end{center}
		
		We denote the diagram obtained, by $M^t$. It is called the 
		\emph{Howie diagram corresponding to ${\mathbb P}_t$}, see \cite{7}, \cite{12} and Remark 1.2.3. 
		We denote by ${\mathbb M}^t$ the diagram with regions $(\Delta(D))^t$, $D \in Reg_{4^+,t}(M)$.
		We denote by $D^t$ the image of $D \in Reg(M)$ in $M^t$ and by $(\Delta (D))^t$ the image of $\Delta (D)$ in ${\mathbb M}^t$.
		
		Notice that $\Delta(D^t) = (\Delta (D))^t$.
		
		\vspace{10pt}
		
		\noindent {\bf Remark} Let $M$ be a van Kampen diagram with connected interior. Notice that while $M$ has connected interior, the interior of ${\mathbb M}^t$ may be not connected. Also, if $M$ contains no vertices with valency 1, it is possible that ${\mathbb M}^t$ contains. However, since we use diagrams in order to count the regions in them, these changes are not relevant to the problem we consider.
		
		\vspace{10pt}
		We need a few more definitions:
		
		\noindent {\bf Definition 3:} (one-layer maps)
		\addtocontents{toc}{\protect\contentsline{subsection}
			{Definition 3 (one-layer maps)}
			{\thepage}{}
		} 
		\begin {center}
		\begin{tikzpicture} [scale = 1.2]
			\draw (-4.5,2.5) node [label={left:$v$},label={right:$D_1$}] (v1) {} .. controls (-3.5,5) and (0.5,5) .. (1.5,2.5) 
			node[pos=.17 ] (u1) {} 
			node[pos=.33 ] (u2) {}
			node[pos=.47 ] (u3) {}
			node[pos=.6 ] (u4) {}
			node[pos=.67,rotate =-20,label={above:$\mu$} ] (arr1) 
			{\tikz{\draw (-3.2,7.65) -- (-3,7.5) -- (-3.2,7.35);}}
			node[pos=.72 ] (u5) {}
			node[pos=.84 ] (u6) {}
			node [inner sep=0mm,label={right:$u$},label={left:$D_m$}] (v2) {};
			\draw (v1.center)  .. controls (-3.5,0) and (0.5,0) .. (v2.center) 
			node[pos=.17 ] (l1) {} 
			node[pos=.33 ] (l2) {}
			node[pos=.47 ] (l3) {}
			node[pos=.6 ] (l4) {}
			node[pos=.67,rotate = 20,label={below:$\nu$} ] (arr2) 
			{\tikz{\draw (-3.2,7.65) -- (-3,7.5) -- (-3.2,7.35);}}
			node[pos=.72 ] (l5) {}
			node[pos=.84 ] (l6) {};
			
			\draw (u1.center) -- (l1.center);
			\draw (u2.center) -- (l2.center);
			\draw (u3.center) -- (l3.center);
			\draw (u4.center) -- (l4.center);
			\draw (u5.center) -- (l5.center);
			\draw (u6.center) -- (l6.center);
			\draw (u4.center) -- (l3.center);

		\end{tikzpicture}\\
		Figure 2
	\end{center}
	\begin{itemize}
		\item [1)]\emph{Simple one-layer map. See Fig. 2}\\ Let $M$ be a connected, simply connected map with connected interior and with simple boundary, $m=|M|\geq 1$. Call $M$ a \emph{a simple one-layer map} if $m=1$, or $m\geq 2$ and the following hold
		\begin{itemize}
			\item [(i)] There are exactly 2 regions which we denote by $D_1$ and $D_m$ respectively, with exactly one neighbour in $M$, which we denote by $D_2$ and $D_{m-1}$ respectively. $Reg(M)\setminus\{D_1,D_m\}$ can be numbered such that $Reg(M)=\{D_1,\ldots,D_m\}$, and each region $D_i,i\in\{2,\ldots,m-1\}$ has exactly two neighbours in $M$, namely, $D_{i-1}$ and $D_{i+1}$.
			\item [(ii)]$\partial D_1\cap \partial M$ and $\partial D_m\cap \partial M$ contain vertices $u(M)$ and $v(M)$ respectively, such that $v(M)\mu u(M)\nu^{-1}$ is a boundary cycle of $M$ such that 
			$\mu=\mu_1\ldots\mu_m,\,\   \nu=\nu_1\ldots\nu_m,\,\ \mu_i=\partial D_i\cap \mu$ and $\nu_i=\partial D_i\cap \nu$. $\mu_i$ and $\nu_i$ are connected. 
			We have $M=\langle D_1,\ldots ,D_M\rangle$. When dealing with one-layer maps we shall mean that the regions $D_1,\ldots,D_m$ occur in this order. See Fig 2. We call $\mu$ and $\nu$ the \emph{sides}  of $M$.
		\end{itemize}
		\begin{center}
			\begin{tikzpicture} [scale = 0.7]
				\draw (0.45,2.5) node (v1) {} .. controls (1.45,5) and (5.45,5) .. (6.45,2.5) 
				node[pos=.15 ] (u1) {} 
				node[pos=.3 ] (u2) {}
				node[pos=.45 ] (u3) {}
				node[pos=.6 ] (u4) {}	
				node[pos=.75 ] (u5) {}
				node[pos=.87 ] (u6) {}
				node(v2) {};
				\draw (v1.center)  .. controls (1.45,0) and (5.45,0) .. (v2.center) 
				node[pos=.15 ] (l1) {} 
				node[pos=.3 ] (l2) {}
				node[pos=.45 ] (l3) {}
				node[pos=.6 ] (l4) {}
				node[pos=.75 ] (l5) {}
				node[pos=.87 ] (l6) {};
				
				\draw (u1.center) -- (l1.center);
				\draw (u2.center) -- (l2.center);
				\draw (u3.center) -- (l3.center);
				\draw (u4.center) -- (l4.center);
				\draw (u5.center) -- (l5.center);
				\draw (u6.center) -- (l6.center);

				\draw (-9.5,2.5) node (v3) {} .. controls (-8.5,5) and (-4.5,5) .. (-3.5,2.5) 
				node[pos=.15 ] (u1) {} 
				node[pos=.3 ] (u2) {}
				node[pos=.45 ] (u3) {}
				node[pos=.6 ] (u4) {}	
				node[pos=.75 ] (u5) {}
				node[pos=.87 ] (u6) {}
				node (v4) {};
				\draw (v3.center)  .. controls (-8.5,0) and (-4.5,0) .. (v4.center) 
				node[pos=.15 ] (l1) {} 
				node[pos=.3 ] (l2) {}
				node[pos=.45 ] (l3) {}
				node[pos=.6 ] (l4) {}
				node[pos=.75 ] (l5) {}
				node[pos=.87 ] (l6) {};
				
				\draw (u1.center) -- (l1.center);
				\draw (u2.center) -- (l2.center);
				\draw (u3.center) -- (l3.center);
				\draw (u4.center) -- (l4.center);
				\draw (u5.center) -- (l5.center);
				\draw (u6.center) -- (l6.center);
				\filldraw (v4) circle [radius=0.15];
				\filldraw (v1) circle [radius=0.15];
				\filldraw (v2) circle [radius=0.15];
				\draw (v4.center) -- (v1.center) ;
				\node (v5) at (10,2.5) {};
				\draw  (v5.center) -- (v2.center);
				\node at (-2.5,3) {$u(M_{i-1})$};
				\node at (0,3) {$v(M_i)$};
				\node at (7,3) {$u(M_i)$};
			\end{tikzpicture}\\
			Figure 3\\
		\end{center}
		\item [2)] 
		\emph{One-layer maps} \\ Let $M$ be a connected simply connected map and let $Int(M)=M_1\dot{\cup}\cdots\dot{\cup}M_r$,\ $r\geq 1$. Call $M$ {\it one-layer} if $M_i$ are simple one-layer maps and $v(M_i)$ is connected by a simple path to $u(M_{i-1})$. See Fig. 3. 
		We denote by $OL(\mu,\nu)$ the collection of all the one-layer maps, with sides $\mu$ and $\nu$.
	\end{itemize}

	
	
	\noindent {\bf Definition 4:} (a-bands)
	\addtocontents{toc}{\protect\contentsline{subsection}
		{Definition 4 (a-bands)}
		{\thepage}{}
	} 
	Let $D_1, D_2 \in Reg_2(M)$ and let $Supp(D_1)=\{a,x\}$, $Supp(D_2)=\{a,y\}$
	($x=y$ is not excluded). Say that $D_1$ and $D_2$ are \emph{a-neighbours} if
	$\theta := \partial D_1 \cap \partial D_2$ contains a subpath $\theta_0$ with label $a^{\pm 1}$. (We consider $D_i$ as a neighbour of itself, $i=1,2$.) 
	The transitive closure of "a-neighbourhood" is  an equivalence relation.
	For $D \in Reg_2(M)$ with $Supp(D)=\{a,x\}$ denote by $\hat{L}_a(D)$ the a-equivalence class of $D$. Let $L_a(D)$ be the subdiagram of $M$ generated by 
	$\hat{L}_a(D)$. Call $L_a(D)$ the \emph{a-band generated by $D$}. Denote by $Str(M)$ the collection of all the bands in $M$. Clearly
	\[
	Reg_2(M)=\bigcup\limits_{L \in Str(M)} Reg(L)
	\]
	
	\vspace{20pt}
	
	\noindent {\bf Definition 5:} (adequate bands)
	\addtocontents{toc}{\protect\contentsline{subsection}
		{Definition 5 (adequate bands)}
		{\thepage}{}
	} 
	\begin{enumerate}
		\item Let $L_a \in Str(M)$. Call $L_a$ \emph{adequate} if
		$L_a$ is a simply connected one-layer diagram with cyclically reduced boundary label $L_a=\langle D_1,\ldots , D_m \rangle$,
		$m \geq 1$, $D_i \in Reg_2(M)$ and 
		$\partial D_i \cap \partial D_{i+1}$ has label $a$ or label $a^{-1}$, $i=1, \ldots , m-1$. In particular $\partial L_a$ is a simple closed curve. 
		In Section 5 we show that if $M$ is Minimal then every band is an adequate band.
		\item Say that \emph{$Str(M)$ is adequate} if
		\begin{itemize}
			\item[(i)] $L$ is adequate for every $L \in Str(M)$ and 
			\item[(ii)] if $L_a$ and $L_b$ are different bands and $L_a \cap L_b$ contains a region $D$ then $\{D\} = L_a \cap L_b$.
		\end{itemize}
	\end{enumerate}
	
	
	\noindent {\bf Definition 6:} ($\mathcal{M}_3(W)$, Minimal diagrams)\\
	\addtocontents{toc}{\protect\contentsline{subsection}
		{Definition 6 ($\mathcal{M}_3(W)$, Minimal diagrams)}
		{\thepage}{}
	} 
	Let $W$ be a cyclically reduced word which represents 1 in $A(\Gamma)$. Assume that $W\neq 1$ in $A(\Gamma_2)$. Here $\Gamma_2$ is the graph 
	obtained from $\gamma$ by removing all edges not labelled by 2.
	Denote by ${\cal M}(W)$ the set of all the van Kampen diagrams with boundary label $W$. By \cite[Ch. V]{14} ${\cal M}(W) \neq \emptyset$.
	Let ${\cal M}_1(W)$ denote the set of all the diagrams in ${\cal M}(W)$ for which $|Reg_{4^+}(M)|$ is minimal possible. Clearly, ${\cal M}_1(W) \neq \emptyset$. Let ${\cal M}_2(W)$ be the set of the diagrams in ${\cal M}_1(W)$ which contain minimal number of Equivalence classes $[D]_M$ of regions in $Reg_{4^+}(M)$. Clearly, ${\cal M}_2(W) \neq \emptyset$. Finally, let ${\cal M}_3(W)$ be the set of the diagrams in ${\cal M}_2(W)$ with adequate bands.
	Say that $M$ is \emph{Minimal} if $M \in {\cal M}_2(W)$ and $|Reg_2(M)|$ is minimal possible.
	It follows that if M is Minimal then $M\in\mathcal{M}_3(W)$. 
	Hence $\mathcal{M}_3(W)\neq \emptyset$. 
	
	\noindent {\bf Remark} 
	
	\noindent We do not require that the diagrams in $\mathcal{M}(W)$ be
	reduced. It follows from the definition of $\mathcal{M}_2(W)$ that
	\begin{equation*}\tag {*}
		\left.
		\begin{minipage}{4in}
			{\it $M\in\mathcal{M}_2(W)$ contains neither reducing pairs $(D_1,D_2)$ with $D_i\in Reg_{4^+}( M)\, i=1,2$, nor contains $\mathbb M$ reducing pairs $(\Delta_1,\Delta_2)$ with $\Delta_i\in Reg_{4^+}(\mathbb M)$}
		\end{minipage}
		\right\}
	\end{equation*}
	We call diagrams satisfying (*) \quad {\it $4^+$-reduced}.\\
	Thus, $M\in \mathcal{M}_3(W)$ may contain only reducing pairs $(D_1,D_2),\quad D_i\in Reg_2(M)$. However, if $M$ is Minimal then $M$ is reduced.
	
	\subsection*{The third construction $(\widetilde M^t)$ and  $((\widetilde{\mathbb M^t}))$ }
	\addtocontents{toc}{\protect\contentsline{subsection}
		{The third construction $(\widetilde M^t)$ and  $((\widetilde{\mathbb M^t}))$}
		{\thepage}{}
	} 
	Given $M^t$ or ${\mathbb M}^t$, we shrink every bigon in them to an edge labelled 
	$t^{\pm 1}$. Denote the obtained diagram by $\widetilde{M^t}$ and
	$\widetilde{{\mathbb M}^t}$, respectively.   This alters the diagram since we ignore all the commuting relators, which contain $t$, but it alters neither the number of regions $\Delta$ in $Reg_{4^+,t}(\mathbb M)$ with $\Phi (\partial \Delta) = R_e$, $\lambda(e) \geq 4$, nor the number of regions inside $\Delta = \Delta(D)$, which we want to count first. 
	See Sect. 5.
	
	\vspace{20pt}
	
	\noindent {\bf Proposition 3} 
	\addtocontents{toc}{\protect\contentsline{subsection}
		{Proposition 3 (Properties of $\widetilde{\mathbb M}^t$)}
		{\thepage}{}
	} 
	{\it Let $M\in \mathcal{M}_3(W)$ be a van Kampen $\cal R$-diagram over $F(X)$, where $\cal R$ is given in ({\rm III}). 
		
		Then each of the following holds: Let $t \in T(M)$.
		\begin{itemize}
			\item[(a)] $\widetilde{({\mathbb M}^t)}$ satisfies the condition C(4) \& T(4)
			\item[(b)] For each region $\widetilde{\Delta}$ of $\widetilde{({\mathbb M}^t)}$ we have $|\partial \widetilde{\Delta} | \leq |W|_t$
			\item[(c)] The number of regions $\widetilde{\Delta}$ in $\widetilde{({\mathbb M^t)}}$ is at most $|W|^2_t$
	\end{itemize}}
	
	\vspace{10pt}
	
	\noindent {\bf Corollary 1:} 
	\[
	\sum\limits_{\widetilde{\Delta} \in Reg\left(\widetilde{{\mathbb M^t}}\right)}|\partial \widetilde{\Delta} | \leq |W|^3
	\]

	\vspace{20pt}
	
	\noindent {\bf Proof of Corollary 1}
	
	By part (b) of the Proposition 
	$|\partial \widetilde{\Delta} | \leq |W|_t$ and by part (c) of the Proposition the number of regions in $\widetilde{({\mathbb M}^t)}$ is bounded from above by $(|W|_t)^2$.
	The result follows.
	
	\vspace{10pt}
	
	To complete the proof of the Main Theorem we need to estimate $|Reg_2(M)|$ in terms of $|W|$.

	\vspace{10pt}
	
	\noindent {\bf Theorem 2}
	\addtocontents{toc}{\protect\contentsline{subsection}
		{Theorem 2 (upper bound on the number of commuting relations)}
		{\thepage}{}
	} 
	{\it Let $M$ be a Minimal van Kampen $\cal R$-diagram with cyclically reduced boundary label $W$. Then $|Reg_2(M)| \leq |W|^6$.}
	
	\vspace{10pt}
	
	\noindent{ Proof of Theorem 2, relying on Corollary 1 of Proposition 3 and on Proposition 4 below.}

	\noindent {\bf Proposition 4:} {\it Let  notation and assumptions be as in Theorem 2. Then 
		$Str(M)$ is adequate.}

	\vspace{20pt}
	
	\noindent {\bf Corollary 2:} {\it Let $l=\left(\sum\limits_{\Delta \in Reg_{4^+}({\mathbb M})}|\partial \Delta|\right) + |\partial M|$. Then
		\begin{itemize}
			\item[(a)] $|Str(M)| \leq l$( hence by Corollary 1
			$|Str(M)| \leq |W|^3+|W|\sim |W|^3$)
			\item[(b)] For every $L \in Str(M)$ $|L| \leq l $ holds.
	\end{itemize}}
	
	\vspace{10pt}
	
	\noindent {\bf Proof of Corollary 2 }
	\begin{itemize}
		\item[(a)] Each band $L_a$ starts and ends either on the boundary of a region in $Reg_{4^+}({\mathbb M})$ or on $\partial M$. Hence $|Str(M)| \leq l$.
		\item[(b)] Each non-boundary region $D$ in $L_a$ belongs to an $L_b$, $b \neq a$ and since $Str(M)$ is
		adequate, $D$ belongs exactly to one such $L_b$. 
		On the other hand, the number of boundary regions in $L_{a}$ is at most $|W|$.
		Hence $|L_a| \leq |Str(M)| \leq l$. 
	\end{itemize}
	
	\vspace{10pt}
	
	The proof of Theorem 2  now is immediate:
	
	\vspace{10pt}
	
	$Reg_2(M)$ is the totality of all the regions in all the bands, hence
	$|Reg_2(M)| \leq |Str(M)|\cdot |L| \leq l^2$, by parts (a) and (b) of Corollary 2. Hence 
	$|Reg_2(M)| \leq |W|^6$, by Corollary 1. \hfill $\Box$

	\subsection*{Proof of the Main Theorem }
	\addtocontents{toc}{\protect\contentsline{subsection}
		{Proof of the Main Theorem }
		{\thepage}{}
	} 
	Clearly $|Reg(M)|=|Reg_{4^+}(M)|+|Reg_2(M)|$. 
	If 
	$Reg_{4^+}(M) \neq \emptyset$ then 
	the result follows from Theorem 1 and Theorem 2.
	If 
	$|Reg_{4^+}(M)| =0$ then $A(\Gamma)$ is a right angled Artin group for which $f(n) = n^2$. \hfill $\Box$

	Thus, the proof of the Main Theorem is reduced to the proofs of  Propositions 1-4, which we carry out by Theorems A, B and C below and some results on the way to their proofs. In the rest of the work we concentrate on their proofs. We assume that $\cal R$ is given by ({\rm III}).
	
	\vspace{10pt}
	
	\noindent {\bf Theorem A} \\
	\addtocontents{toc}{\protect\contentsline{subsection}
		{Theorem A\ ($\widetilde{\mathbb M}^t$ satisfies C(4)\&T(4))}
		{\thepage}{}
	} 
	{\it Let $M$ be a van Kampen $\cal R$-diagram over $F$ with cyclically reduced boundary label W. 
		If $M\in \mathcal{M}_3(W)$ then $\widetilde{\mathbb M}^t$ satisfies the small cancellation condition C(4) \& T(4), for every $t \in T(M)$}.
	
	\vspace{20pt}
	
	\noindent {\bf Theorem B} \\
	\addtocontents{toc}{\protect\contentsline{subsection}
		{Theorem B\ ($Supp(W)\supseteq\cup Supp(D),\ D\in Reg(M)$}
		{\thepage}{}
	}
	{\it Let $M$ be a van Kampen $\cal R$-diagram over $F$ with cyclically reduced boundary label $W$. If $M\in \mathcal{M}_3(W)$ then for every region $D$ in $Reg(M)$
		realising a relation $R(a,t)$, $||W||_t \geq n(a,t)$ and $||W||_a \geq n(a,t)$, $a$ and $t$ in $T$.}
	
	\vspace{10pt}
	
	Theorem B has a few easy-to-prove consequences which we mention below without proof. (We do not use them in this work.)

	\vspace{10pt}

	\noindent {\bf Corollary 3}\\ {\it Let $A(\Gamma)$ be an Artin group such that $n_{ij} \neq 3$, for all $1 \leq i < j \leq n$. Let $W$ be a cyclically reduced word in $F(X)$. Assume that there is an $a \in X$ such that $W=W_1a^{\alpha}W_2a^{\beta}$, reduced as written and such that $ a \not\in (Supp(W_1) \cup Supp(W_2))$. Then the following hold:
		\begin{enumerate}
			\item If $\alpha + \beta \neq 0$ then $W \neq 1$ in $A$.
			\item If $\alpha + \beta =0$ then $W =1$ only if $W_1W_2 =1$ in $A$.
			
	\end{enumerate}}
	
	\vspace{10pt}
	
	\noindent {\bf Corollary 4}\\  {\it
		Let $A(\Gamma)$ and $W$ be as in the previous Corollary. Assume that there is a letter $a \in X$ such that $W=W_1a^{\alpha}W_2a^{\beta}W_3a^{\gamma}$, $ a \not\in (Supp(W_1) \cup Supp(W_2) \cup Supp(W_3))$,
		\begin{enumerate}
			\item If $\alpha + \beta + \gamma \neq 0$ then $W \neq 1$ in $A$.
			\item If $\alpha + \beta + \gamma = 0$ then $W=1$ in $A$ if and only if $W_1W_2W_3 =1$ in $A$.
		\end{enumerate}
	} 
	\vspace{10pt}
	
	Part 2 of Corollary 5 below has the "flavour of Tits' Conjecture".
	
	\vspace{10pt}
	
	\noindent {\bf Corollary 5}\\  {\it
		Let $W$ be a cyclically reduced cyclic word in which $||W||_t \leq 3$ for every $ t \in Supp(W)$. Assume $W=1$ in $A$. 
		\begin{enumerate}
			\item Let $\alpha_i(t)$, $ i \in \{1,2,3\}$ be the exponents of $t$ in $W$, $ t \in Supp(W)$. Then $\sum\limits_{i=1}^3 \alpha_i(t) =0$, for every $ t \in Supp(W)$.
			\item All the relations in the right hand side of ({\rm II}) (for $W$) are commuting relations. In particular, $W=1$ in $A(\Gamma_2)$, 
			where $\Gamma_2$ is obtained from $\Gamma$ by the removal of all edges which are not labelled by 2.
			(It is a right angled Artin group.)
		\end{enumerate}
	}
	
	\vspace{10pt}
	
	\noindent {\bf  Theorem C (Greendlinger's Lemma)}
	\addtocontents{toc}{\protect\contentsline{subsection}
		{Theorem C\ (generalised Greendlinger's Lemma)}
		{\thepage}{}
	}
	{\it
		Let $A(\Gamma)$ be an Artin group given by ({\rm III}). 
		Let $M_0\in\mathcal{M}_3(W)$ be a connected, simply connected  Minimal ${\cal R}$-diagram over $F$ with cyclically reduced boundary label $W, W\neq 1$ in $A(\Gamma_2)$. 
		
		\noindent Let $W=W_1W^{-1}_2$ be reduced as written. Then there is a van Kampen diagram $M,\ M\in \mathcal{M}_3(W)$, with boundary cycle $\omega$, such that $\omega$ decomposes into $\omega = \omega_1 \omega^{-1}_2$,\quad $\omega_i$ labelled with $W_i$, $i=1,2$ and one of the following holds:
		\begin{itemize}
			\item [(a)] There is a boundary region $K \in Reg_{4^+}({\mathbb M})$ with $\partial K \cap  \omega_i$ connected, $||\partial K \cap \omega_i|| \geq n(K)$, for $i=1$ or $i=2$,
			where $Reg_{4^+}({\mathbb M})=\{K \in Reg({\mathbb M}) \ | \ n(K) \geq 4 \}$.
			Moreover, $|\xi|\leq |\partial K\cap\omega_i|$ where $\xi$ is the complement of $\partial K\cap \omega_i$ on $\partial K$. 
			\item [(b)] There are boundary regions $D_1,D_2\in Reg_2(M)$ with $\partial D_i\cap \omega_j$ connected, $|\partial D_i\cap \omega_j|\geq 2$, for $j=1$ or for $j=2$.
			\item [(c)] There is a region $D\in Reg_2(M)$ with $\partial D\cap \omega_j$ connected, $|\partial D\cap \omega_j|\geq 2$ and there is a region $K\in Reg_{4^+}(\mathbb M)$ with $\partial K\cap\omega_j$ connected, $\|\partial D\cap \omega_j\|\geq 2$, for $j=1$ or for $j=2$. 
		\end{itemize}
	}

	\vspace{10pt}
	
	\subsubsection* {We give below a short description of what is done in the work.}
	\addtocontents{toc}{\protect\contentsline{subsection}
		{Description of the main steps in the work}
		{\thepage}{}
	}
	The goals of the work are: 
	\begin{itemize}
		\item[1] to estimate the number of regions in $Reg_{4^+}(M)$
		\item[2] to estimate the number of regions in $Reg_2(M)$.
	\end{itemize}  
	The idea for 1) is  to count the number of regions in $\widetilde{\mathbb{M}}^t, t\in T(M)$ where due to theorem A the condition C(4)\&T(4) is satisfied, hence we may do this, and then to show that the number of regions in $\widetilde{\mathbb{M}}^t$  
	and $Reg_{4^+,t}{\mathbb{M}}$ is the same, via a natural mapping $\psi_t$ of diagrams which sends a region 
	$\Delta$ in ${\mathbb{M}}_t$ to a uniquely defined region $\widetilde{\Delta}$  in $\widetilde{\mathbb{M}^t}$. One of the problems is that $\widetilde{\mathbb{M}}^t$ is obtained from ${\mathbb{M}}_t$ by shrinking edges and identifying vertices hence in principle $\widetilde{\Delta}$ 
	may shrink to an edge or a vertex.
	Thus our task is to show that $\psi_t$  (and $\psi^{-1}_t$) do not cause deformation and collapse of regions. One of the classical methods to avoid deformations in diagrams 
	is by Greendlinger's Lemma. 
	It roughly  says that a diagram with at least two regions has at least two Greendlinger regions. See def. 1.3.1
	In general terms it shows that if there is a kind of deformation in certain subdiagrams then certain events are simultaneously unavoidable and forbidden, which of course is absurd. 
	Hence no deformation occurs. 
	For example consider a reduced diagram in which every region has boundary label $U^n$, $U$ cyclically reduced, not proper power, $n\geq 7$.
	Suppose that the boundary of a region E is not simple closed. (This is the deformation which we would like to avoid) Then there is a loop in the boundary of the region which surrounds a disc Q, the boundary of which is labelled by a subword of the cyclic word $U^n$. 
	Suppose that in Q we have Greendlinger's Lemma. Then Q has a boundary region D with $\partial D\cap\partial Q$ connected, such that $|\partial D\cap\partial Q|\geq 3|U|$. Hence D and E have necessarily a common (unavoidable) piece of length at least $3|U|$. But such a word is forbidden, because due to $U$ being not a proper power, this would mean that D cancels E, violating that the diagram is reduced. (Though we do not use this example we use it's underlying idea. 
	See Proofs of 2.1.7 and 2.1.8). 
	Coming back to the context of our work, 
	notice that by Theorem C every bigonal ${\cal R}$-diagram $M$, i.e. $\partial M=\omega_1\omega_2^{-1}$
	has a generalised Greendlinger region on one of the sides $\omega_i$ of the bigon. 
	Since $\widetilde{\mathbb M^t}$ does satisfy C(4)\& T(4) hence has Greendlinger regions. 
	The natural candidates for Greendlinger's regions in $M$ are the images (by $\Psi_t^{-1}$) of the Greendlinger regions of $\widetilde{\mathbb{M}^t}$. The problem is that while Greendlinger regions are boundary regions, our candidates in $M$ need not be. 
	So we need a technique  to "move" inner regions to the boundary. More precisely, to show that there is a diagram with the same boundary label as M in which  the candidates are boundary regions.
	We do this by observing that since the generators of the second homotopy group of the defining complex are prisms, we can replace one half of the prism by its complement by a rotation of the prism. 
	In M this has the effect of moving a region which occurs both in M and in the prism. 
	See Sec. 3. 
	The  involvement of the second homotopy group requires the introduction of extended presentations.
	
	\noindent Finally, we have to count $|Reg_2(M)|$. To this end it is enough to show that two bands cannot intersect more than once. (See Sec. 5). 
	Essentially we prove this simultaneously with the other main Proposition (See Sec. 5).\\
	The ideas developed here lead to further results for Artin groups A dealt with in the present work. Recall that a parabolic subgroup of an Artin group is a conjugate of a standard parabolic subgroup. 
	In \cite{10} we show that the intersection of parabolic subgroups is parabolic. 
	Also we describe fusion in A. In particular we show that in even Artin groups every standard parabolic subgroup P controls fusion in P(i.e. if two elements of P are conjugate in A then they are already conjugate in P).\\
	Together with further ideas we show in \cite{11} that locally reducible Artin groups (i.e. no (2.3.3),(2.3.4) and (2.3.5) type standard parabolics occur) have polynomial $(n^6)$ isoperimetric functions.

	\vspace {10pt}

	The work is organised as follows:
	\begin{description}
		\item [in Section 1] we recall and introduce basic results in Artin groups, presentations, the C(4)\&T(4) condition, 2 generated Artin groups, and bands.
		\item[in Section 2] we introduce the diagrams $\mathbb{M}$ and $\mathbb{M}_t$ and show that they are adequate. Also $M$ is adequate.
		\item[In Section 3] we consider relative extended presentations, Howie diagrams I-moves and banded diagrams.
		\item [In Section 4] we consider diagrams M with $Reg(M)=Reg_2(M)$. These are diagrams of right angled artin groups. We call them 
		{\it Abelian diagrams}.
		\item[In Section 5] we consider the connection between $\mathbb{M}_t$ and $\widetilde{\mathbb{M}^t}$ and prove the main Theorems.
	\end{description}
	\clearpage
	\addtocontents{toc}
	{\protect\contentsline{section} {Contents}
		{\thepage}{}
	}
	\tableofcontents
	
	\section*{Acknowledgments}
	I am very grateful to the Referee for his important remarks on a previous version of the work, in particular for pointing out several omissions.
	\vspace{10pt}
	
	%
	%
	
	\section{Preliminary Results}
	
	\tocexclude{
		\subsection{Basic results in Artin groups}
	}
	\ \par
	
	\addtocontents{toc}{\protect\contentsline{subsection}
		{\numberline{\thesubsection}Theorem \thesubsection\ (van der Lek)}
		{\thepage}{}
	}
	
	\noindent {\bf  Theorem 1.1 [van der Lek] \cite{13}}  
	{\it \begin{itemize} 
			\item[(a)] Let $A$ be an Artin group and let $B$ be a standard parabolic subgroup of $A$ with $B=\langle Y\rangle$, $Y \subseteq X$. Consider the Artin group $B'$ generated by $Y$ and related by all the relators $R_j$ such that $x_i, x_j \in Y$. Then the natural map $B' \rightarrow A$ is an embedding.
			\item[(b)] Let $B$ and $C$ be standard parabolic subgroups of $A$, $B=\langle Y\rangle$, $C=\langle Z\rangle$,
			$Y,Z \subseteq X$. Then $B \cap C = \langle Y \cap Z \rangle$.
		\end{itemize} 
	}
	
	Theorem 1.1 has several consequences which we collect in the Lemma below.
	For a word $W$ in $F(x_1, \ldots , x_n)$ let $Supp(W)$ be the subset of letters of $X$ which occur in $W$ or $W^{-1}$.
	
	\vspace{10pt}
	
		
		\begin{lemma}
		\addtocontents{toc}{\protect\contentsline{subsubsection}
			{Lemma \thesubsubsection\ (basic results on Artin groups)}
			{\thepage}{}
		}\ 
			\begin{itemize}
				\item [(a)]	Every standard generator has infinite order.
				\item[(b)] Let $U$, $V$ be non-empty words in
				$F(x_1, \ldots , x_n)$, $U \neq_A 1$ and $V \neq_A 1$. 
				If $U=_A V$ then $Supp(U) \cap Supp(V) \neq \emptyset$.
				\item[(c)] Let $1 \neq W \in F(x_1, \ldots , x_n )$ be a cyclically reduced word 
				which represents 1 in $A$. 
				If $x \in  Supp(W)$ then $x^{\pm 1}$ has at least two non-adjacent occurences in $W$.
				\item[(d)]  Let $W$ be a non-empty reduced word in $F(x_1, \ldots , x_n)$. If $W=_A1$ 
				then $||W|| \geq 4$.\\ 
				If $||W||=4$ then $W=x^{\alpha}y^{\beta}x^{-\alpha}y^{-\beta}, 
				[x,y]=1$, $ \alpha , \beta \in {\mathbb Z} \setminus \{0\}$.
			\end{itemize}
	\end{lemma}
	
	\vspace{10pt}
	
	\noindent {\bf Proof}
	\begin{itemize}
		\item [(a)] Apply Theorem 1.1 to $B=\langle x \rangle$, $x \in X$.
		
		\item [(b)] Since $U=_A V$ and $1 \neq_AU \in \langle Supp(U) \rangle \cap \langle Supp(V) \rangle = \langle Supp(U) \cap Supp (V) \rangle$, by part (b) of Theorem 1.1. Hence $Supp(U) \cap Supp(V) \neq \emptyset$.
		\item[] Here $\langle Supp (U) \rangle$ denotes the subgroup of $A$ generated by $Supp(U)$.
		
		\item [(c)] Let $W^*$ be a cyclic conjugate of $W$ such that 
		$W^*=x^{\alpha}W_1$, $\alpha \neq 0$, $h(W_1) \neq x^{\varepsilon}$,
		$\varepsilon \in \{1, -1 \}$ and $t(W_1) \neq x^{\varepsilon}$, 
		$\varepsilon \in \{1, -1 \}$, where $h(W)$ and $t(W)$ are the first and last letters of $W$, respectively.
		By part(a) $x^{\alpha}\neq_A 1$, 
		hence $W_1 \neq_A 1$ and $x^{-\alpha}=_A W_1$. Hence 
		$\{x\} = Supp(x^{-\alpha})\subseteq Supp(W_1)$ by part(b), i.e. $x \in Supp(W_1)$.	
		
		\item[(d)]  Assume first that $W$ is cyclically reduced and $||W|| \leq 3$.
		By part(a) $||W|| \in \{2,3\}$. If $||W||=2$ then $W=x^{\alpha}y^{\beta}$, $x \neq y$, $x,y \in X$, $\alpha , \beta \in {\mathbb Z}\setminus \{0\}$.
		If $\langle x,y \rangle$ is free then $W \neq_A1$, contrary to assumption, hence there is a single Artin relation $R$ between $x$ and $y$ and due to Theorem 1.1(a) $\langle x,y\rangle = \langle x,y \ | \ R \rangle$.
		Therefore, by \cite[Lemma 7]{1} $||W|| \geq ||R|| \geq 4$,
		i.e. $||W|| \geq 4$, a contradiction to the assumption $||W|| \leq 3$.
		
		Assume now that $||W|| =3$. It follows from part(c)  that 
		$3=||W|| \geq 2|Supp(W)|$, hence $|Supp(W)| =1$. But then $||W||=1$, a contradiction.

		If  $W$ is not cyclically reduced let $W'$ be its cyclically reduced form.
		Then $4 \leq ||W'|| \leq ||W||$ and the result follows. \\
		\ 	If $||W||=4$  then $W=x^{\alpha}y^{\beta}x^{-\alpha}y^{-\beta}$, because 
		$4 = ||W|| \geq 2 |Supp(W)| \geq 4$, i.e. $|Supp(W)| =2$,
		say $Supp(W) = \{x,y\}$.
		Then 
		$\langle x,y\rangle = \langle x,y \ | \ R \rangle$, where 
		$R$ is an Artin relator with $||R||=4$. There is only one such relator:
		$xyx^{-1}y^{-1}$. Hence, if $||W||=4$ and $W=_A1$ then $x$ and $y$ commute. \hfill $\Box$
		
	\end{itemize}

	\subsection{Presentations, Maps and Diagrams}
	We introduce some further notation and assumptions. For a vertex $v$ in a map $M$ let $N_M(v)$ be the set of all the regions of $M$ which contain $v$ on their boundaries. Similarly, for a path $\mu$ in $M$  we denote by $N_M(\mu)$ the set of all the regions $D$ with $\partial D \cap \mu \neq \emptyset$. 
	For a path $\mu$ of $M$ denote by $o(\mu)$ and $t(\mu)$ the origin and the terminus of $\mu$ in $M$, respectively. Thus,
	$\bar{\mu} = o(\mu) \cup \mu \cup t(\mu)$, where $\bar{\mu}$ is the closure of $\mu$ in $M$. 
	Let $D_1$ and $D_2$ be regions in $M$. Say that $D_1$ and $D_2$ are \emph{neighbours} or are \emph{adjacent} if they have an edge in common.
	\noindent Recall that every group $G$ has a free presentation
	${\cal P} = \langle X \ | \ {\cal R} \rangle$ by a set of 
	generators $X$ and a cyclically closed
	set of cyclically reduced defining relators ${\cal R} \subseteq F(X)$, where $F(X)$ is the free group, freely generated by $X$ and
	$G \cong F(X) / N$, where $N$ is the normal closure of $\cal R$ in $F(X)$.  (See \cite{14}). To every element $W \in N$ there corresponds a set of planar connected and 
	simply connected 2-complexes $M$ the edges of which are labelled by words in $F(X)$ the boundary of the 2-cells are labelled by $R \in \mathcal{R}$  and $W$ is a boundary label of $M$. 
	Such diagrams are called \emph{van Kampen diagrams}. 
	See \cite[Ch. V]{14}.
	We call a vertex of $M$ \emph{inner} if it is not on the boundary of $M$. We call a region $D$ \emph{inner} if all its boundary vertices are inner vertices of $M$. 
	In particular, $\partial D \cap \partial M = \emptyset$.
	Diagrams exist (under appropriate conditions) on every closed surface. They define a tessellation of the surface by 2-cells. Apart from van Kampen diagrams we shall need \emph{spherical diagrams} over $F(X)$. These are diagrams on the sphere $S^2$ which tessellate $S^2$ by 2-cells, whose boundaries are labelled by the defining relations. 
	
	We assume that 
	\begin{itemize} \item[(A1)] \emph{$M$ has no inner vertices with valency 2. }
		
		However, boundary vertices may have valency 2.
	\end{itemize}
	
	For a path $\mu$ in $M$ define the length $|\mu|_{M}$ of $\mu$ in $M$ as the number of vertices on $\mu$ in $M$ plus 1.
	When convenient, we shall write $||\mu||$ for $||\Phi(\mu)||$, where $\Phi$ is the labeling function $\Phi:M\longrightarrow F$.
	Since all our maps are underlying maps of diagrams, if $\theta$ is a boundary path with $\Phi(\theta)=W$ and $W$ decomposes by $W=W_1\cdots W_k,\quad k\geq 2$ reduced as written, then $\theta$ is decomposed into $\theta=\theta_1 v_1\theta_2\cdots v_{r-1}\theta_r$ with $\Phi(\theta_i)=W_i$ 
	and the vertices $v_i$ have valency 2 in $M,\ i=1,\ldots,r$. In particular, if the presentation satisfies the condition C(p),$p\geq 2$, 
	as defined above and $D$ is a boundary region of $M$ such that $\partial D$ is simple and $\theta=\partial D\cap \partial M$ is connected and the complement $\eta$ of $\theta$ on $\partial D$, decomposes by 
	$\eta=\eta_1\cdots\eta_q,\quad p-q\geq 2,\quad  \eta_i$ pieces, $1\leq i\leq q$, then $\theta$ contains at least $p-q-1\geq 1$ vertices with valency 2.
	\par\ \par
	\begin{lemma}
		Let $M$ be a connected, simply connected map which satisfies the condition C(p), $p\geq 2$ and let $N$ be a connected, simply connected submap with connected interior. 
		Let $D$ be a boundary region of $N$ such that $\partial D$ is simple and $\theta:=\partial D\cap \partial N$ is connected. 
		Assume $M\setminus N$ contains regions $E_1,\ldots,E_l$, such that $\theta\cap\partial E_i$ contains an edge. Then $i_N(D)\geq p-l$.
	\end{lemma}
	\addtocontents{toc}{\protect\contentsline{subsubsection}
		{Lemma \thesubsubsection\ ($i_N(D)\geq p-L$)}
		{\thepage}{}
	}
	Here $i_N(D)$ is the number of neighbours of $D$ in $N$.
	
	\noindent {\bf Proof:}\\
	Assume by way of contradiction that $i_N(D)< p-l$. The neighbours of $D$ in $M$ are those in $M\setminus N$ together with those in $N$.
	The neighbours of D in $M\setminus N$ are $E_1,\ldots,E_l$, hence, 
	$p\leq i_N(D)+l$. i.e. $p\leq i_N(D)+l<(p-l)+l=p$, a contradiction. Hence $i_N(D)\geq p-l$.
	\ \hfill $\Box$
	
	\par \ \\

	In certain cases it is convenient to ignore a part of the defining relators and to emphasize the others. This can be done when we can divide $X$ and $\cal R$ into subsets $X=X_1 \dot{\cup} X_2$ and ${\cal R}={\cal R}_1\dot{\cup} {\cal R}_2$ such that ${\cal R}_1 \subseteq F(X_1)$. Then we ignore ${\cal R}_1$ and emphasize ${\cal R}_2$ by considering the quotient of $H*F(X_2)$ by the normal closure of ${\cal R}_2$ in $H*F(X_2)$, where $H=F(X_1) / N_1$ and $N_1$ is the normal closure of ${\cal R}_1$ in $F(X_1)$. This gives a (relative)\emph{ presentation of $G$ relative to $H$}.
	We recall briefly the relevant definitions. For more details see \cite{7} and \cite{12}. 
	
	\vspace{10pt}
	\tocexclude{
		\subsubsection{Definitions}
	}
	\addtocontents{toc}{\protect\contentsline{subsubsection}
		{Definition \thesubsubsection\ (relative presentations, Howie diagrams)}
		{\thepage}{}
	}
	
	\begin{itemize}
		\item[(a)] \emph{A relative presentation $\mathbb P$} consists of a group $H$, a set $Y$, $Y \cap H = \emptyset$ and a set ${\cal R}'$ of cyclically reduced words in $H * F(Y)$ which are not in $H$. We denote
		\begin{equation}\tag{\rm I'}
			{\mathbb P} = \langle H, Y | {\cal R}' \rangle
		\end{equation}
		
		The group defined by $\mathbb P$ is
		$H * F(Y) / \langle \langle {\cal R}'\rangle \rangle$, where
		$\langle \langle {\cal R}' \rangle \rangle$ is the normal closure of 
		${\cal R}'$ in $H * F(Y)$. 
		\item[(b)] The \emph{Howie diagrams} which correspond to $\mathbb P$ are finite planar or spherical maps $M$ in which the 2-cells (regions) have interior homeomorphic to the open unit disc. 
		The edges of $M$ are labelled by elements of $\langle Y \rangle$, such that reading the corner labels of an inner vertex gives 1 in $H$ and reading the boundary label of a region clockwise including the corner labels gives one of the defining relation from 
		$\langle \langle {\cal R}'\rangle\rangle$. 
		Indeed, let $v$ be an inner vertex in $M$ and let $D_1,\ldots,D_m$ be the regions which carry $v$ on their boundaries. Let $h_i$ be the vertex label of $v$ at $D_i$. We have to show that $h_1h_2\cdots h_m=1$ in $H$. 
		We have $R_i=h_if_ik_ig_il_it_i,\quad h_i,k_i,l_i$ vertex labels, $g_i$ boundary labels. $f_i$ and $t_i$ are labels of common edges between $R_i$ and $R_{i-1}$, and $R_i$ and $R_{i+1}$ respectively. \\
		Let $x_i=k_ig_il_i,\quad i=1,\ldots, m$. Then 
		\begin{equation*}
			\tag{*}
			h_1\ldots h_m\cdot t_1^{-1}x_1\cdots x_m t_1=1_A,
		\end{equation*}
		since $h_i=t_i^{-1}x_if_i$ in A and $t_{i+1}=f_i$ in F. ($R_i=1$ in A). But $x_1\ldots x_m$ is a boundary label of the simply connected diagram consisting of $R_1,\ldots R_m$. Hence $x_1\cdots x_m=1$ in A and therefore $h_1\cdots h_m=1$ in A due to (*) Consequently, $h_1\cdots h_m=1$ in H by Theorem 1.1.

		We call the word obtained by the concatenation of the corner labels "vertex label". It is defined, up to cyclic permutation. See \cite[p.90]{7} and Remark 2.1.  Figure 4 below gives an example  where
		$H=\langle a,b \ | \ aba^{-1}b^{-1} \rangle$,
		$Y=\{t\}$ and 
		${\cal R}=\{R_1,R_2\}$.
		
	\end{itemize}
	
	
	\vspace{10pt}
	
	\begin{minipage}{1.5in}
		\begin{center}
			\begin{tikzpicture}[scale=.35]
				\draw[->] (2.5,5)--(0,0)--(5,0);
				\draw (5,0)--(10,0)--(7.5,5);
				\draw[->] (5,10)--(7.5,5);
				\draw[->] (5,10)--(2.5,5);
				\draw (1.5,0) arc(0:65:1.5);
				\draw (9.4,1.25) arc(85:179:1.25);
				\draw (4.5,9) arc(210:310:.75);
				\draw (1.1,.5)--(1.4,.5)--(1.4,.8);
				\draw (8.5,.5)--(8.2,.5)--(8.2,.9);
				\draw (4.7,8.5)--(5,8.7)--(4.7,9);
				\node[below] at (5,0) {$t$};
				\node[left] at (2.5,5) {$t$};
				\node[right] at (7.5,5) {$t$};
				\node[right] at (1,.75) {$a$};
				\node[left] at (8.5,.75) {$a$};
				\node[below] at (5,8.8) {$a$};
			\end{tikzpicture}
		\end{center}
		\[
		R_1=atat^{-1}a^{-1}t^{-1}
		\]
	\end{minipage}\hspace{.5cm}\begin{minipage}{1.5in}
		\begin{center}
			\begin{tikzpicture}[scale=.35]
				\draw[->] (2.5,5)--(0,0)--(5,0);
				\draw (5,0)--(10,0)--(7.5,5);
				\draw[->] (5,10)--(7.5,5);
				\draw[->] (5,10)--(2.5,5);
				\draw (1.5,0) arc(0:65:1.5);
				\draw (9.4,1.25) arc(85:179:1.25);
				\draw (4.5,9) arc(210:310:.75);
				\draw (1.1,.5)--(1.4,.5)--(1.4,.8);
				\draw (8.5,.5)--(8.2,.5)--(8.2,.9);
				\draw (4.7,8.5)--(5,8.7)--(4.7,9);
				\node[below] at (5,0) {$t$};
				\node[left] at (2.5,5) {$t$};
				\node[right] at (7.5,5) {$t$};
				\node[right] at (1,.75) {$b$};
				\node[left] at (8.5,.75) {$b$};
				\node[below] at (5,8.8) {$b$};
			\end{tikzpicture}
		\end{center}
		\[
		R_2=btbt^{-1}b^{-1}t^{-1}
		\]
	\end{minipage}\hspace{.5cm}
	\begin{minipage}{1.5in}
		\begin{center}
			\begin{tikzpicture}[scale=.35]
				\draw (0,0)--(8,0)--(10,10)--(2,10)--(0,0);
				\draw (0,0) --(10,10);
				\draw (8,0) --(2,10);
				\draw (5,5) circle (2);
				\draw (1.5,0) arc(0:79:1.5);
				\draw (8.25,1.25) arc(90:170:1.5);
				\draw (1.75,9) arc(240:370:1);
				\draw (9,10) arc(180:260:1);
				\draw (1,.8)--(1.3,.7)--(1.4,.9);
				\node[right] at (1.3,.7) {$a$};
				\draw (.85,1.1)--(.8,1.3)--(1,1.3);
				\node[above] at (.8,1.3) {$t$};
				\draw (2.8,5.2)--(3,5)--(3.2,5.2);
				\node[left] at (3,5) {$b$};
				\draw (6.8,5.2)--(7,5)--(7.2,5.2);
				\node[right] at (7,5) {$b$};
				\draw (4.8,2.8)--(5,3)--(4.8,3.2);
				\node[above] at (5,3) {$a$};
				\draw (4.8,6.8)--(5,7)--(4.8,7.2);
				\node[above] at (5,7) {$a$};
				\draw (7.8,1)--(8,1.2)--(7.8,1.3);
				\node[above] at (8,1.2) {$t$};
				\draw (6.9,.7)--(7,.5)--(7.3,.7);
				\node[left] at (7,.7) {$a$};
				\draw (9.3,9.1)--(9.5,9.1)--(9.5,9.3);
				\draw (9,9.4)--(9.1,9.6)--(9.3,9.5);
				\node[left] at (9.1,9.4) {$t$};
				\draw (2.4,8.7)--(2.2,8.9)--(2.4,9.1);
				\draw (3,9.6)--(3.2,9.6)--(3.3,9.3);
				\node[right] at (3.2,9.6) {$a$};
				
				\node[right] at (9,5) {$t$};
				\node[left] at (1,5) {$t$};
				\node[above] at (6,10) {$t$};
				\node[below] at (4,0) {$t$};

			\end{tikzpicture}
		\end{center}
		\[
		aba^{-1}b^{-1}=1 \ \mbox{in} \ H
		\]
	\end{minipage}
	
	\begin{center}
		Figure 4
	\end{center}
	
	\vspace{10pt}

	\begin{itemize}
		\item[(c)] 
		Let $M$ be a Howie diagram. $M$ is called \emph{reduced} if $M$ contains no adjacent regions $D_1$ and $D_2$ with a common edge and boundary cycles
		$v_1ev_2\omega_1$ and $v_1ev_2\omega_2$, respectively, such that 
		$\Phi(\omega_1) = \Phi(\omega_2)$ and
		$h_1 = h'_1$ and $h_2=h'_2$, where $h_1$ and $h'_1$ are the corner labels of $D_1$ and $D_2$ at $v_1$ and $h_2$ and $h'_2$ are the corner labels of $D_1$ and $D_2$ at $v_2$, respectively. See Figure 5.
	\end{itemize}

	Van Kampen diagrams are special cases of Howie diagrams. They are Howie diagrams with $H=1$. Hence, this definition is relevant for van Kampen diagrams, as well.

	\vspace{20pt}
	
	\begin{center}
		\begin{tikzpicture}[scale=.5]
			\draw (0,0) arc(30:330:2);
			\draw[rotate around={180:(0,0)}] (0,0) arc(330:30:2);
			\draw (0,0)--(0,-2);
			\draw[fill] (0,0) circle (.1);
			\draw (-.2,-.7)--(0,-1)--(.2,-.7);
			\draw (-3.9,-1.3)--(-3.7,-1)--(-3.5,-1.3);
			\draw (3.5,-1.3)--(3.7,-1)--(3.9,-1.3);
			\draw[fill] (0,-2) circle (.1);
			
			\node at (-4.5,-1.5) {$\omega_1$};
			\node at (4.5,-1.5) {$\omega_2$};
			\node[above=0.1cm] at (0,0) {$v_1$};
			\node[below=0.1cm] at (0,-2) {$v_2$};
			\node at (-.5,.2) {$h_1$};
			\node at (1,.2) {$h'_1$};
			\node at (-.5,-1.7) {$h_2$};
			\node at (.5,-1.7) {$h'_2$};
			\node at (-.5,-1) {$e$};
			\node at (-2,-1) {$D_1$};
			\node at (2,-1) {$D_2$};
		\end{tikzpicture}
		
		Figure 5
	\end{center}
	
	\vspace{20pt}
	
	\begin{remark}
		\addcontentsline{toc}{subsubsection}{Remark\quad\thesubsubsection}
		\ \\Howie diagrams were introduced by James Howie \cite{7} under the name "relative diagrams" (see \cite[p. 90]{7}) 
		for the purpose of studying solvability of equations over groups in one indeterminate $t$. Every such equation over a group $G$ is $r(t)=1$, 
		where $r(t)$ is an element  in the free product $G * \langle t \rangle$. 
		The equation $r(t)=1$ is solvable in an overgroup if and only if the natural map
		$G \rightarrow G * \langle t \rangle / \langle \langle r(t) \rangle \rangle$ is an embedding. Here $\langle \langle r(t) \rangle \rangle$ is the normal closure of $r(t)$ in $G * \langle t \rangle$. It follows that if $G$ is presented by 
		$\langle Z | S \rangle$ then $K:=G * \langle t \rangle / \langle \langle r(t) \rangle \rangle$ is isomorphic to $\langle Z \cup \{t\} \ | \ S \cup \{r(t) \} \rangle$.
		Let $Z_1 = Z \cup \{t\}$ and let $S_1 = S \cup \{r(t)\}$. In our work we start with
		$\langle Z_1 | S_1 \rangle$, the free presentation of $K$, and split $Z_1$ and $S_1$ as above to get the relative presentation of K:
		$\displaystyle Z / \langle \langle S \rangle \rangle * \langle t \rangle /\langle \langle r(t) \rangle \rangle $.
		Relative diagrams (Howie diagrams) are spherical diagrams with a distinguished vertex $v_0$. (See \cite[p. 90]{7}.)
		
		We use an extension of this construction to planar diagrams where the boundary of the diagram replaces the distinguished vertex $v_0$. 
		The boundary label of the diagram is a word in $G * \langle t \rangle$, rather than an element of $G$, as in the spherical case. 
		Also we allow any finite number of equations $r_1(t), \ldots , r_k(t) \in G * \langle t \rangle$. These are our relations.
	\end{remark}
	
	\subsection{C(4)\&T(4) maps}
	In this subsection we recall known results on C(4) \& T(4) maps in the form we need them. Recall that for a boundary region D of M, 
	$i(D)$ is the number of neighbours of $D$ in $M$.
	
	\vspace{10pt}
	
	\begin{definition}[Greendlinger regions]
		\addcontentsline{toc}{subsubsection}{Definition \thesubsubsection\ (Greendlinger regions)}
		Let $M$ be a map, let $k$ be a natural number, $k\in \{1,2,3\}$ and let $D$ be a boundary region of $M$. Say that $D$ is a \emph{$k$-Greendlinger region}
		if each of the following holds:
		
		\begin{enumerate}
			\item $\partial D$ is a simple closed curve;
			\item $\partial D \cap \partial M$ is connected;
			\item $i(D) = k$.
		\end{enumerate}
	\end{definition}
	Denote by $\mathcal D_k(M)$ the set of all the k-Greendlinger regions
	
	\vspace{10pt} 
	
	\begin{lemma}[Basic properties of C(4) \& T(4) Maps]
			\addcontentsline{toc}{subsubsection}{Lemma \thesubsubsection\ (Basic properties of C(4) \& T(4) Maps)}
			\ \\ Let $M$ be a connected, simply connected map with connected interior which contains at least two regions. Assume that $M$ satisfies the small cancellation condition C(4) $\&$ T(4). Then $M$ contains a set of Greendlinger regions 
			$\mathcal D = {\cal D}(M)={\cal D}_1(M)\cup {\cal D}_2(M)$, such that each of the following holds:
			\begin{itemize}
				\item[(a)] $i(D)\leq 2$ for $D\in{\cal D}$ and $\sum\limits_{D \in {\cal D}(M)}(3-i(D)) \geq 4$
				\item[(b)] $|{\cal D}| \geq 2$. If $|{\cal D}| = 2$ then $M$ is a one-layer map and $i(D_1) = i(D_2)=1$, where ${\cal D} = \{D_1, D_2 \}$ and  
				$M= \langle D_1, E_1, \ldots, E_k, D_2 \rangle$,\ 
				$E_i$ regions,\ $i=1,\ldots,k$.
				\item[(c)] If $|{\cal D}| =3$ then the dual $M^*$ of $M$ is a tripod. Moreover $i(D_j)=1$ for $j = 1,2,3$,  ${\cal D}=\{D_1, D_2, D_3\}$.
				\item[(d)] Assume that $\partial M$ decomposes into $\omega_1 v \omega_2 u \omega_3 z$, $v,u$ and $z$ vertices. Then there exists a $D \in {\cal D}(M)$ such that 
				$|\partial D \cap \omega_i| \geq 2$, for some $i,\ i=1,2,3$.
				\item[(e)] Let $D$ be a region in $M$. Then $\partial D$ is a simple closed curve.
				\item[(f)] Let $D_1$ and $D_2$ be regions in $M$. If $\partial D_1 \cap \partial D_2 \neq \emptyset$ then $\partial D_1 \cap \partial D_2$ is connected.
			\end{itemize}
	\end{lemma}
	\vspace{10pt}
	
	\begin{remark}
		\addcontentsline{toc}{subsubsection}{Remark \thesubsubsection\ (The Lemma follows from part (a))}
		As we shall see, the proofs of parts (b) - (f) of the Lemma rely only on part (a) and the  assumption that $M$ satisfies the condition C(4) \& T(4). Hence, if $M$ satisfies the condition C(4) \& T(4) and ${\cal D}_0(M)$ is any set of Greendlinger regions which satisfy part (a), then the Proposition will hold if we replace ${\cal D}(M)$ with ${\cal D}_0(M)$.
	\end{remark}
	\vspace{10pt}
	
	\noindent {\bf Proof} \begin{itemize} \item[(a)] This is \cite[p. 248]{14}.
		\item[] We prove parts (b)-(f) by induction on $|M|$, the case $|M|=2$ being clear.
		\item[(b)] Due  to part (a) $|{\cal D}| >1$. Assume $|{\cal D}| =2$. Then again due to part (a) $i(D_1)=i(D_2)=1$, ${\cal D}=\{D_1, D_2\}$. We prove by induction on $|M|$ that $M$ is a one layer map. If $|M|=2$ this is clear. Suppose $|M| \geq 3$ and the claim that $M$ is a one-layer map holds true for $M_1 = M \setminus \{D\}$, where $ D \in {\cal D}$.
		Since $i(D)=1$ and $i(K)\leq 2$ for every $K\in{\cal D}(M_1)$ the removal of $D$ from $M$ may change the status of being a Greendlinger region in $M_1$ and not being a Greendlinger region in M, for at most one region of $M_1$. Hence $|{\cal D}(M_1)|\leq\left( |{\cal D}(M)|-1\right)+1=|{\cal D}(M)|=2$. Hence $|{\cal D}(M_1)|=2$ and the induction hypothesis applies.
		
		Notice that due to the definition of Greendlinger regions, $M_1$ is connected, simply connected with connected interior. 
		Since $i(D)=1$, $D$ has one neighbour, $E$, in $M$. Let $M_1 = \langle E_1, \ldots , E_k \rangle$, $ k \geq 2$. Then either $E \in \{E_2, \ldots , E_{k-1}\}$ or $ E \in \{E_1, E_k\}$. See Figure 6 (a) and (b) respectively.
		
		\begin{center}
			
			\begin{minipage}{2.2in}
				\begin{center}
					\begin{tikzpicture}[scale=.7]
						
						\draw (0,0) ellipse (4 and 2);
						
						\draw (-3,-1.3)--(-3,1.3);
						\draw (-2,-1.75)--(-2,1.75);
						\draw (-1,-1.9)--(-1,1.9);
						\draw (0,-2)--(0,2);
						\draw (3,-1.3)--(3,1.3);
						\draw (2,-1.75)--(2,1.75);
						\draw (1,-1.9)--(1,1.9);
						
						\draw (1,-1.9)--(1,-3)--(2,-3)--(2,-1.7);
						
						\node at (-3.5,0) {$E_1$};
						\node at (3.5,0) {$E_k$};
						\node at (1.5,0) {$E$};
						\node at (1.5,-2.5) {$D$};

					\end{tikzpicture}
					
					(a)
			\end{center}\end{minipage}\hspace{10pt}\begin{minipage}{2in}
				\begin{center}
					\begin{tikzpicture}[scale=.7]
						
						\draw (0,0) ellipse (4 and 2);
						
						\draw (-3,-1.3)--(-3,1.3);
						\draw (-2,-1.75)--(-2,1.75);
						\draw (-1,-1.9)--(-1,1.9);
						\draw (0,-2)--(0,2);
						\draw (3,-1.3)--(3,1.3);
						\draw (2,-1.75)--(2,1.75);
						\draw (1,-1.9)--(1,1.9);
						
						\draw[rounded corners] (-3.75,.6)--(-4.5,1.5)--(-3.2,2.2)--(-3,1);
						
						\node at (-3.5,0) {$E_1$};
						\node at (3.5,0) {$E_k$};
						\node at (-3.7,1.3) {$D$};

					\end{tikzpicture}
					
					\vspace{.5cm}
					
					(b)
				\end{center}
			\end{minipage}
			
			\begin{center}
				\begin{tikzpicture}[scale=.7]
					
					\draw (0,0) ellipse (4 and 2);
					
					\draw (-3,-1.3)--(-3,1.3);
					\draw (-2,-1.75)--(-2,1.75);
					\draw (-1,-1.9)--(-1,1.9);
					\draw (0,-2)--(0,2);
					\draw (3,-1.3)--(3,1.3);
					\draw (2,-1.75)--(2,1.75);
					\draw (1,-1.9)--(1,1.9);
					
					\node at (-3.5,0) {$\Delta_1$};
					\node at (-2.5,0) {$\Delta_2$};
					\node at (3.5,0) {$\Delta_k$};
					\node at (-3,1.8) {$v_1$};
					\node at (-1.5,2.2) {$\omega_1$};
					
					\node at (-3,-1.9) {$u_1$};
					\node at (-1.5,-2.2) {$\omega_2$};
					
					\node[right] at (3,1.9) {$v_{k-1}$};
					\node[right] at (3,-1.9) {$u_{k-1}$};
					
					\draw[fill] (-3,1.3) circle(.1);
					\draw[fill] (3,1.3) circle(.1);
					\draw[fill] (-3,-1.3) circle(.1);
					\draw[fill] (3,-1.3) circle(.1);
					\draw[fill] (4,0) circle(.1);
					\draw[fill] (-4,0) circle(.1);
					
				\end{tikzpicture}
				
				(c)
			\end{center}
			
			Figure 6
		\end{center}
		
		If Case (b) of Figure 6 occurs then $M$ is one layer. If Case (a) of Figure 6 occurs then 
		$|{\cal D}| =3$\quad (${\cal D} = \{E_1, E_k, D \}$), violating the assumption that $|{\cal D}|=2$, hence in all cases the result follows.
		
		\item[(c)] $|{\cal D}| =3$ together with part (a) imply that there is a region $D \in {\cal D}$ with $i(D) =1$. Remove it and proceed as in part (b), 
		the notation of which we follow, thus $|{\cal D}(M_1)|\leq|{\cal D}(M)|=3$. Hence either $|{\cal D}(M_1)|=2$ or $|{\cal D}(M_1)|=3$. In the first case by part(b), 
		$M_1$ is a one layer map, 
		hence since $i(D)=1$ the dual of $M$ is a tripod. 
		In the second case the induction hypothesis applies.
		See Figures 7(a) and 7(b), resp.
	\end{itemize}
	
	\begin{minipage}[b]{1.6in}
		\begin{center}
			\begin{tikzpicture}[scale=.3]
				\draw[rounded corners] (0,1.3)--(0,4)--(4,6);
				\draw (4,6) .. controls (6.5,7) and (7.5,7) ..(10,6);
				\draw [rounded corners] (10,6)--(14,4)--(14,1)--(8,4)--(6,4)--(0,1.3);
				
				\draw (4,6)--(6,4);
				\draw (8,4)--(10,6);
				
				\draw (12,2)--(13,4.5);
				\draw (10,3)--(11,5.5);
				
				
				\draw (2,2)--(1,4.5);
				\draw (4,3)--(3,5.5);
				
				\draw[rounded corners] (13,4.5)--(13.3,6)--(11.5,7.5)--(11,5.5);
				
				\node at (.9,2.5) {$F_1$};
				\node at (13,2.3) {$F_3$};
				\node at (12,5.9) {$D$};
				
			\end{tikzpicture}
			
			(a)
	\end{center}\end{minipage}\hspace{10pt}\begin{minipage}[b]{1.6in}
		\begin{center}
			\begin{tikzpicture}[scale=.3]
				
				\draw[rounded corners] (0,1.3)--(0,4)--(4,6)--(6,8)--(6,12)--(8,12)--(8,8)--(10,6)--(14,4)--(14,1)--(8,4)--(6,4)--(0,1.3);
				
				\draw (4,6)--(6,4);
				\draw (6,8)--(8,8);
				\draw (8,4)--(10,6);
				
				\draw (12,2)--(12.5,4.7);
				\draw (10,3)--(11,5.5);
				
				\draw (6,10.5)--(8,10.5);
				\draw (6,9.5)--(8,9.5);
				
				\draw (2,2)--(1,4.5);
				\draw (4,3)--(3,5.5);
				
				\draw[rounded corners] (14,4)--(14.8,7)--(13,7.5)--(13,4.5);
				
				\node at (.9,2.5) {$F_1$};
				\node at (7,11.2) {$F_2$};
				\node at (13,2.4) {$F_3$};
				\node at (13.8,6) {$D$};

			\end{tikzpicture}

			(b)
		\end{center}
	\end{minipage}\hspace{10pt}\begin{minipage}[b]{1.5in}
		\begin{center}
			\begin{tikzpicture}[scale=.3]
				
				\draw[rounded corners] (0,1.3)--(0,4)--(4,6)--(6,8)--(6,12)--(8,12)--(8,8)--(10,6)--(14,4)--(14,1)--(8,4)--(6,4)--(0,1.3);
				
				\draw (4,6)--(6,4);
				\draw (6,8)--(8,8);
				\draw (8,4)--(10,6);
				
				\draw (12,2)--(12.5,4.7);
				\draw (10,3)--(11,5.5);
				
				\draw (6,10.5)--(8,10.5);
				\draw (6,9.5)--(8,9.5);
				
				\draw (2,2)--(1,4.5);
				\draw (4,3)--(3,5.5);

				\node at (.9,2.5) {$F_1$};
				\node at (7,11.2) {$F_2$};
				\node at (13,2.3) {$F_3$};

			\end{tikzpicture}
			
			(c)
		\end{center}
	\end{minipage}
	
	\begin{center}
		Figure 7
	\end{center}
	
	\begin{itemize}
		\item[] 
		
		the dual of $M^*$ is a tripod
		
		\item[(d)] If $|{\cal D}| \geq 4$ then by the pigeonhole principle (P.H.P) at least one region $D \in {\cal D}$ avoids $v$, $u$ and $z$. Hence, since $|\partial D\cap \partial M|\geq 2$ 
		due to part (a) and the C(4) \& T(4) condition,
		$D$ satisfies the requirement of part (d).
		If $|{\cal D}| =3$ then the result follows easily from part (c).
		Assume that $|{\cal D}|=2$. Then by part (b) of the Lemma $M$ is a one layer map $M = \langle D_1, \ldots , D_k \rangle$, $ k \geq 2$.
		In particular, $D_1$ and $D_k$ are in ${\cal D}(M)$ with $i(D_j)=1$, $ j \in \{1, k\}$. Hence, 
		$|\partial D_j \cap \partial M| = |\partial D_j| - i(D_j) \geq 4-1 =3$, due to the C(4) condition. 
		Hence $|\partial D_j \cap \omega_i| \geq 2$, for some $j \in \{1,k\}$
		and some $i$, $i=1,2$, due to the P.H.P.
		
		\item[(e)] Assume the claim of the Lemma is false. Then ${\mathbb E}^2 \setminus D = Q_1\cup \cdots \cup Q_m \cup Q_{\infty}$, where $Q_i$ are bounded simply connected connected components (with connected interior),
		$i=1, \ldots , m$, and $Q_{\infty}$ is the unbounded component. Observe that $Q_i \subseteq M$. There is a $Q_i$, say $Q_1$, such that $\partial Q_1$ contains a single double point $w$. See Figure 8. 
		By part (d), taking $\omega_1 = \emptyset$ and $u=v=z=w$ it follows that $|\partial E \cap \partial D | \geq 2$ for some $E \in {\cal D}(Q_1)$.
		Thus $\partial E\cap \partial Q_1=\partial E \cap \partial D$ as $\partial Q_1\subseteq\partial D$. Hence by Lemma 1.2.1 with $N=Q_1,p\geq 4$ and $l=1$ we get $i_{Q_1}(E)\geq 4-1=3$, violating that $i_{Q_1}(E)\leq 2$ by part (a) of the Lemma.
		Hence $\partial D$ is a simple curve.
		\item[(f)] Similar to part (e).
	\end{itemize}
	
	\ \hfill $\Box$
	
	\begin{center}
		\begin{tikzpicture} [scale =0.5]
			\filldraw (0,2) node[label={above:$w$}] (v2) {} circle [radius=0.15];
			\draw (0,0) node {$Q_1$} circle [radius=2];
			\node (v1) at (-1.5,3.5) {};
			\node (v3) at (1.5,3.5) {};
			\draw (v1) -- (v2.center) -- (v3);
		\end{tikzpicture}
		
		Figure 8
	\end{center}
	\vspace{10pt}
	
		\begin{proposition}
			\addcontentsline{toc}{subsubsection}{Proposition \thesubsubsection\ (The subset $\mathcal D_0(M)$ of $\mathcal D(M)$)}
			Let $M$ be a connected, simply connected map which contains at least two regions. Assume that $M$ satisfies the small cancellation condition C(4) $\&$ T(4). 
			Then ${\cal D}(M)$ contains a subset ${\cal D}_0(M)$ such that each of the following holds
			\begin{itemize}
				\item[(a)] If $D \in {\cal D}_0(M)$ and $i(D)=2$ then $\partial D \cap \partial M$ has an endpoint with valency 3.
				\item[(b)] $\sum\limits_{D \in {\cal D}_0(M)} (3-i(D)) \geq 4,\quad i(D)\leq 2$
				\item [(c)] If $|\mathcal D_0(M)|=2$ then $M$ is a one layer map and if $|\mathcal D_0(M)|=3$ then $M^*$ is a tripod.
			\end{itemize}
	\end{proposition}
	We need some preparatory results for the proof of Proposition 1.3.4.
	
	\vspace{10pt}
	
	\begin{definition}[Convex layer structure]
		\addcontentsline{toc}{subsubsection}{Definition \thesubsubsection\ (Convex layer structure)}
		Let $M$ be a connected and simply connected map with connected interior. For a connected submap $M_0$ denote by $N(M_0)$ the set of all the regions $ D$ in $M$ with $D \not\in Reg(M_0)$ and $\partial D \cap \partial M_0 \neq \emptyset$.
		
		Let $v_0$ be a vertex of $M$. 
		Define $S_0(v_0)= L_0 = \{ v_0\}$, $L_1 = N(v_0)$, $S_1 = S_0 \cup L_1$ and by induction 
		$L_i = N(S_{i-1})$ and $S_i = S_{i-1} \cup L_i$. 
		Since $M$ contains a finite number of regions and every region is contained in some $S_i$, it follows that there is a natural number $ p \geq 0$ such that $M = S_p$. 
		Denote $\Lambda (v_0) = (L_0, \ldots , L_p)$, $L_0 = \{v_0\}$. 
		We call $\Lambda (v_0)$ the \emph{layer decomposition of $M$ with layers $L_i$}. It is defined uniquely by $v_0$.\\
		We say that $\Lambda(v_0)$ is a \emph{convex layer structure} of $M$ (relative to $v_0$) if for every subset $S$ of regions in $L_i$, the submap $S_{i-1}\cup S$ is simply connected.
	\end{definition}
	We propose to show that every C(4) \& T(4) map has convex layer structure. See Proposition 1.3.7. We need the following Lemma.
	
	\vspace{10pt}
	
		\begin{lemma} 
			\addcontentsline{toc}{subsubsection}{Lemma \thesubsubsection\ (Basic properties of diagrams with convex layer structure)}
			Let $M$ be a connected, simply connected map with connected interior
			having at least two regions. Suppose that $M$ satisfies the small cancellation condition C(4) $\&$ T(4) and $M$ has a convex layer structure $\Lambda(v_0) = (L_0 , \ldots , L_p)$ relative to a vertex $v_0$. Then each of the following holds.
			\begin{itemize}
				\item[(a)] Let $D$ be a region in $L_i$ and let  $\omega_i = \partial S_i \cap \partial L_{i+1}$, $ i = 1, \ldots, p-1$,
				$\omega_p = \partial L_p \cap \partial M$. Then $\partial D \cap \omega_{i-1}$ is connected and $\beta (D) \leq 2$, where $\beta(D)$ is the number of neighbours of $D$ in $L_i$. 
				\item[(b)]  Let $ v \in \omega_i$. Then $d_{L_i}(v) \leq 3$. (Notice that $v\in\partial L_i$)
				\item[(c)]  Let $\alpha (D)$ be the number of neighbours of $D$ in $L_{i-1}$. Then $\alpha(D) \leq 1$.
			\end{itemize}
	\end{lemma}
	
	\vspace{10pt}
	
	\noindent {\bf Proof} By induction on $i$, for $i=1$ being clear due to Lemma 1.3.2 parts (e) and (f).  
	Assume $ i \geq 2$.
	\begin{itemize}
		\item[(a)] If $\partial D \cap \partial S_{i-1}$ is not connected then $S_{i-1} \cup \{D\}$ violates the convexity of the layer structure. Hence $\partial D \cap \partial S_{i-1}$ is connected. Suppose $\beta(D) \geq 3$. 
		Then there is a region $E$ in $L_i$ such that $\partial E \cap \partial D \neq \emptyset$, but $\partial E \cap \partial D \cap \partial S_{i-1} = \emptyset$.
		Let $P$ be a bounded connected component of ${\mathbb E}^2 \setminus (S_{i-1} \cup \{D \cup E\})$. Then $P \subseteq M$ and $P$ is simply
		connected and has connected interior. 
		Hence part (d) of Lemma 1.3.2 applies with $\omega_1 = \partial P \cap \partial S_{i-1}$, 
		$\omega_2 = \partial P \cap \partial E$ and $\omega_3 = \partial P \cap \partial D$.
		(Here, in the proof, $\omega_1$, $\omega_2$ and  $\omega_3$ are defined as in Lemma 1.3.2 and not as subpaths of  $\partial S_i$.)
		Hence there is a region $K$ in ${\cal D}(P)$ with $|\partial K \cap \omega_l| \geq 2$. $l\in\{1,2,3\} $. Say $l=1$.
		Hence $\partial K \cap \omega_1 = \gamma_1 v_1 \gamma_2 \cdots \gamma_k$,
		$ k \geq 2$, $\gamma_i$ edges,  $v_j$ boundary vertices of P with $d_P(v_j) =2$ and by the induction hypothesis for $S_{i-1}$,
		$d_{S_{i-1}}(v_j) \leq 3$. 
		Hence $d_M (v_j) \leq 2+3-2=3$. 
		But $v_j$ is an inner vertex of $M$. Hence the condition T(4) is violated, proving our claim.

		\item[(b)] Assume by way of contradiction that $d_{L_i}(v) \geq 4$. Then $L_i$ contains a region $E$ with $\partial E \cap \omega_i = \{v\}$. Since $\beta (E) \leq 2$ by part (a), it follows that $\alpha(E) \geq d(E) - \beta(E) \geq 4-2=2$. Let $K_1$ and $K_2$ be regions in $L_{i-1}$, adjacent to $E$, containing a common vertex $z$ on their boundary and let $\xi_j = \partial K_j \cap \partial E$, $j =1,2$. Then $\xi_j$ is connected and $\xi_1 z \xi_2 \subseteq \partial E \cap \omega_{i-1}$, where $z$ is the common endpoint of $\xi_1$ and $\xi_2$, $z \in \omega_{i-1}$, $z$ an inner vertex of $M$,
		$d_{L_i}(z) =2$ ($w_i$ as in part(a)).
		Hence, $d_M(z) = d_{L_i}(z) + d_{L_{i-1}}(z)-2$. 
		Therefore, $d_{L_{i-1}} (z) = d_M(z) - d_{L_i}(z) +2 = d_M(z)$. But $d_M(z) \geq 4$, as $z$ is an inner vertex of $M$. Hence $d_{L_{i-1}} (z) \geq 4$, violating the induction hypothesis  (b) for $S_{i-1}$. Hence $d_{L_i}(v) \leq 3$.
		\item[(c)] Assume by way of contradiction that $\alpha (D) \geq 2$. Let $\partial D \cap \omega_{i-1} = \xi_1 v_1 \cdots \xi_k$,
		$v_i$ vertices, $ k \geq 2$. Then $d_{L_{i-1}}(v_1) \geq 4$, by the argument at the end of the proof of part (b), violating part (b) for $L_{i-1}$. \hfill $\Box$
	\end{itemize}
	
	\vspace{10pt}
	
		\begin{proposition}
			\addcontentsline{toc}{subsubsection}{Proposition \thesubsubsection\ (C(4)\& T(4) implies convex layer structure)}
			Let $M$ be a connected, simply connected map with connected interior. If $M$ satisfies the condition C(4) $\&$ T(4), then $M$ has convex layer structure.
		\end{proposition}
	
	\vspace{10pt}
	
	\noindent {\bf Proof} Let $\Lambda = (L_0, \ldots , L_p)$ be a layer decomposition of $M$ relative to $L_0= \{v\}$. 
	We prove the Proposition by induction on $p$, the case $p=1$ being clear. Assume $ p \geq 2$ and the Proposition holds true for $S_{p-1}$. 
	Assume by way of contradiction that $L_p$ has a subset $S$ of regions such that $L_{p-1} \cup S$ is not simply connected. 
	Then there are regions $E_1$ and $E_2$ in $L_p$ such that $\partial E_1 \cap \partial E_2 \neq \emptyset$
	and $\partial E_1 \cap \partial E_2 \cap \omega_{p-1}= \emptyset$. 
	Hence, as in the proof of part (a) of the previous Lemma
	(with regions $E$ and $D$), it follows that the T(4) condition is violated. Hence $\Lambda$ is a convex layer structure. \hfill $\Box$
	
	\vspace{10pt}
	
	\begin{remark}
		\addcontentsline{toc}{subsubsection}{Remark \thesubsubsection\ ($L$ is one layer)}
		Each layer $L_i$ is a one-layer map.
	\end{remark}
	
	\vspace{10pt}
	
	\noindent {\bf Proof of Proposition 1.3.4} Let $\Lambda = (L_0, \ldots , L_p)$ be a convex layer structure of $M$, guaranteed by Proposition 1.3.7. 
	We prove the Proposition by induction on $p$, the result being clear for $p=1$. Assume $ p \geq 2$. We consider two cases according as $Int(L_p)$ is annular or simply connected.
	
	\underline{Case 1} \emph{$L_p$ has annular interior}. Then every layer has annular interior, since $M$ is simply connected. Let $ D \in {\cal D}_0 (S_{p-1})$. Then $i(D) = \alpha (D)  + \beta (D) \leq 2$. 
	Hence $(\alpha (D) , \beta (D)) \in \{(0,2), (1,1)\}$. Since $Int(L_p)$ is annular,
	$\beta (D) =2$, hence $\alpha(D) =0$ and $i(D)=2$. 
	Also, it follows from Lemma 1.3.6 part (b) that the endpoints of $\partial D \cap \omega_{p-1}$ have valency 3 in $L_{p-1}$. Moreover, let $\gamma(D)$ be the number of neighbours of $D$ in $L_p$. Then $\alpha (D) + \beta (D) + \gamma (D) \geq 4$, by the C(4) condition. Hence $\gamma (D) \geq 2$.
	Let $E_1$ and $E_2$ be regions in $L_p$ such that $\partial E_i \cap \partial D$ contains an edge, $i=1,2$ and such that $\partial E_1 \cap \partial E_2 \neq \emptyset$. 
	Then $\partial E_1 \cap \partial E_2 \cap \omega_{p-1}$ is an inner vertex $z$ and $d_M(z) = d_{L_{p-1}}(z) + d_{L_p}(z) - 2 = 2+ d_{L_p}(z) -2 \geq 4$, by the condition T(4). Thus $d_{L_p}(z) \geq 4$. Hence there is a region $K$ in $L_p$ with $\partial K \cap \omega_{p-1} = \{z\}$. 
	Consequently, $\alpha (K) = 0$, hence $i(K)=\alpha (K) + \beta (K) = 0+2 =2$ and hence, due to Lemma 1.3.6 part (b) $ K \in {\cal D}_0(L_p)$. 
	We say that \emph{$D$ in ${\cal D}_0(L_{p-1})$ induced $K$ in ${\cal D}_0(L_p)$}. 
	Notice that $\sigma (D) := 3-i(D) =1 = 3-i(K)$.
	Hence $\sigma(D) = \sigma(K)$ and the claim follows by induction.
	
	\vspace{10pt}
	\underline{Case 2} \emph{$Int(L_p)$ is the disjoint union of $k$ simply connected components} $C_i$, $i=1, \ldots , k$, $ k \geq 1$. 
	Then $C_i$ are one layer maps by Remark 1.3.8. 
	Consider $C=C_1$. $C = \langle D_1, \ldots , D_r \rangle$, $ r \geq 1$, $ \beta (D_j) =1$, $ \alpha (D_j) \leq 1$, $ j \in \{1,r\}$. 
	Hence, $ D_j \in {\cal D}_0(M)$, for $j \in \{1,r \}$ and if $ r \geq 2$ then $\sigma (D_j) \geq 1$, hence $\sigma (D_1) + \sigma (D_r) \geq 2$.  
	If $r=1$ then $\beta (D_j) =0$, hence $\sigma (D_j) = 2$. Here $\sigma$ has the same meaning as in case 1. 
	Therefore $\sigma (D_1) + \sigma (D_r) \geq 2$ in both cases and if $k \geq 2$ then 
	$ \sigma (M) = \sigma (S_p) \geq 2 + 2= 4$. Hence if $k\geq 2$ then part(b) of the Proposition follows. 
	Assume therefore that $Int(L_p)$ is connected and simply connected.
	Assume first that $Int(L_{p-1})$ is connected.
	Let $L_{p-1}= \langle E_1 , \ldots ,E_s \rangle$, $ s \geq 1$ as we saw, $E_1, E_s \in {\cal D}_0(S_{p-1})$. Consider $L_p\cup L_{p-1}$. 
	One of the following occurs:
	
	\vspace{10pt}
	
	\begin{minipage}{1.8in}

		\begin{tikzpicture}[scale=.35]
			\draw(0,0)--(10,0)--(10,6)--(6,6)--(6,4)--(2,4)--(2,0)--(0,0)--(0,2)--(2,2);
			\draw (8,6)--(8,0);
			\draw (6,4)--(6,0);
			\draw (4,4)--(4,0);
			\draw (2,2)--(10,2);
			\draw (6,4)--(10,4);

			\node at (5,5) {$L_p$};
			\node at (7,5) {$D_1$};
			\node at (0,3) {$L_{p-1}$};
			\node at (3,3) {$E_1$};
		\end{tikzpicture}
		
		\begin{center}
			(a)
			
		\end{center}
		
	\end{minipage}\hspace{3pt}\begin{minipage}{1.8in}

		\begin{tikzpicture}[scale=.35]
			\draw(0,0)--(10,0)--(10,6)--(6,6)--(6,4)--(2,4)--(2,0)--(0,0)--(0,2)--(2,2);
			\draw (8,6)--(8,0);
			\draw (6,4)--(6,0);
			\draw (4,4)--(4,0);
			\draw (2,2)--(10,2);
			\draw (6,4)--(10,4);
			
			\draw (4,4)--(4,6)--(6,6);
			
			\node at (11,5) {$L_p$};
			\node at (3,3) {$E_1$};
			\node at (5,5) {$D$};
			
		\end{tikzpicture}
		
		\begin{center}
			(b)
			
		\end{center}
		
	\end{minipage}\hspace{3pt}\hspace{3pt}\begin{minipage}{1in}

		\begin{tikzpicture}[scale=.35]
			\draw(0,0)--(6,0) --(6,2)--(4,2)--(4,4)--(2,4)--(2,2)--(0,2)--(0,0);
			\draw (2,0)--(2,2)--(4,2)--(4,0);
			
			\node at (5,3) {$L_p$};
			\node at (1,1) {$E_1$};
			\node at (5,1) {$E_s$};
			\node at (3,3) {$D$};
			
		\end{tikzpicture}
		
		\begin{center}
			(c)
			
		\end{center}
		
	\end{minipage}

	\begin{center}
		Figure 9
	\end{center}
	
	\begin{itemize}
		\item [(i)] $E_j\in \mathcal{D}_0(S_p)(=\mathcal{D}_0(M))$ such that $\sigma_{S_p}(E_j)=\sigma_{S_{p-1}}(E_{j})$,
		for $j=1$ or $j=s$.
		Hence, $\sigma_{S_p}(E_j)=3-i_{S_{p-1}}(E_j)$, where 
		$i_{S_{p-1}}(E_j)$ is the number of neighbours of $E_j$ in $S_{p-1}$. See Fig 9(a).
		\item [(ii)] Either $E_j$ induces a region $D$ in $S_p$ for $j=1$ such that $\sigma(E_j)\leq\sigma(D)$, or
		$\sigma(E_s)\leq\sigma(D_r)$, for $j=s$. See Fig. 9(b)
	\end{itemize}
	
	Taking into account this, 
	it follows that $ \sigma (S_p) \geq \sigma (S_{p-1})$.
	
	If $E_1$ and $E_s$ in $L_{p-1}$ induce the same region $D$ in $L_{p}$ then
	$\sigma(E_1) + \sigma(E_2) = (3-2) +(3-2) =2=3-1=\sigma(E)$, i.e. $\sigma (E_1) + \sigma (E_2) = \sigma (D)$. See Figure 9(c).
	
	\noindent Assume now that $Int(L_{p-1})$ is not connected.
	Then $Int (L_{p-1})=K_1\cup K_2\cup\cdots\cup K_{l},\ l\geq 2$.
	If $\partial K_m\cap\partial C_1=\emptyset$ for some $m, ,\in\{1,\ldots,l\}$, 
	then $K_1$ and $K_l$ belong to $\mathcal{D}(S_p)$
	and $\sigma(K_1)+\sigma(K_l)\geq 2$.
	By the same argument $\sigma(D_1)+\sigma(D_r)\geq 2$,
	hence $\sigma(M)=\sigma(S_p)\geq 2+2=4$.
	Hence, assume that $\partial C_1\cap \partial K_m\neq\emptyset$,
	for $m=1,\ldots,l$.
	Then it follows from the definitions of $L_{p-1}$ and $L_p$ that $l=1$, 
	i.e. $Int(L_{p-1})$ is connected, violating the assumption that
	$Int(L_{p-1})$ is not connected.
	
	\noindent Hence the result follows by the induction hypothesis. 
	
	{\bf (c)} Let $v_0$ be an inner vertex of $M$.
	Then $d(v_0)\geq 4$.
	Since $|\mathcal D_0(S_1)|\geq d(v_0)$, we get
	\begin{equation*}\tag {*}
		|\mathcal D_0(S_1)|\geq 4
	\end{equation*}
	Let $L_e$ be the layer of $\Lambda(v_0)$ with $Int(L_e)$ 
	connected and annular with $e$ maximal.
	Then due to (*), $e\geq 1$.
	If $S_e=M$ then using Induction of Greendlinger regions,
	it follows that $|\mathcal{D}_0(M)|\geq |\mathcal{D}_0(S_e)|
	\geq 4$.
	Hence $S_e\neq M$.
	Let $C_1,\ldots,C_m$ be the connected components of 
	$M\setminus S_e$. Since each component contains an element of $\mathcal{D}_0(M)$ (for example, in its layer with highest index), $1\leq m\leq 3$.
	Let $D_1,\ldots,D_n$ be the elements of $\mathcal D_0(S_e)$.
	Then $n\geq 4$.
	Checking all the relative positions of $D_i$ and $C_j$ 
	for Induction from $\mathcal{D}_0(S_0)$ to $\mathcal{D}_0(M)$, 
	it follows that each $D_j$ contributes at least 1 to $|\mathcal{D}_0(M)|$.
	Hence $|\mathcal{D}_0(M)|\geq 4$ and consequently, $M$
	cannot contain an inner vertex.
	Therefore, it follows, by an easy induction on $|M|$, that the 
	dual $M^*$ of $M$ is tree-like and every endpoint $D^*$ of $M^*$ defines an element $D$ of $\mathcal{D}_0(M)$ with
	$i(D)=1$.
	Hence, either $M$ is one-layer or with tripod dual.\hfill $\Box$
	\vspace{10pt}

	\begin{lemma}
			\addcontentsline{toc}{subsubsection}{Lemma \thesubsubsection\ (Improved version of Lemma 1.3.4)}
			Let $M$ be connected, simply connected diagram with
			connected interior which satisfies the condition C(4)\&T(4). Assume that $|{\cal D}_0(M)|\geq 4$. Assume also that $\partial M$ decomposes by $v\omega_1 w\omega_2$. 
			Then for $i = 1$ or $i = 2, N(\omega_i)$ contains a region $D\in {\cal D}_0(M)$ such that either $i(D) = 1$ or $i(D) = 2$ and $\partial D \cap \omega_i$ has an end-point $z$ with valency 3 , $z \in \omega_i$.
	\end{lemma}
	
	\begin{remark}\ 
		\addcontentsline{toc}{subsubsection}{Remark \thesubsubsection}
		\begin{itemize}
			\item [(a)] Notice that without the requirement $z\in\omega_i$ such a region D exists by Proposition 1.3.4.
			\item [(b)] In definition 1.3.5(b) we may replace $v_0$ by a region $D_0$ and even by  a path $\mu$ with certain properties 
			and define convex layer structure relative to $D_0$, or $\mu$ respectively, exactly as in definition 1.3.5,
			except that now $L_0=\{D\}$, or $L_0=\{\mu\}$, respectively. 
			Lemma 1.3.6 and Proposition 1.3.7 remain true with the same proofs, provided that we consider only regions from the same side of $\mu$ (left or right) and $L_1$ is a one-layer map. 
		\end{itemize}
	\end{remark} 
	\ \par
	\noindent {\bf Proof of Lemma 1.3.9}
	\begin{description}
		\item[Case 1] $M$ has no inner vertex. Then $M$ is "tree-like" i.e. the dual $M^*$ of $M$ is a tree. (See def. in \cite[p.245]{14}). 
		Then $i(D)=1$ for every $D\in {\cal D}_0(M)$. Hence the result follows by the Pigeon Hole Principle.
		\item [Case 2] M has an inner vertex $v_0$. Let $\Lambda = (L_0, L_1,... L_p),\quad p \geq 1$ be the layer structure of $M$ with $L_0=\{v_0\}$.
		\begin{description}
			\item[Subcase 1]  $Int (L_p)$ is annular. Then due to Lemma 1.3.6 and Proposition 1.3.7 both endpoints of $\partial D \cap \partial M$ for $D \in {\cal D}_0(M)$ have valency 3.
			Since $|{\cal D}_0(M)|\geq 4$, again the result follows from the PHP.
			
			\item [Subcase 2] $Int (L_p)$ is not annular. Let $i_0$ be the biggest number such that $Int(L_i)$ is annular. 
			Since $v_0$ is an inner vertex of $M, i_0 \geq 1$. 
			Define $V_i= \{\text{endpoints of\ } \partial D\cap\partial S_i\text{\ with valency }3, D \in {\cal D}_2(S_i)\}$. Then $|V_{i_0}|\geq 4$. 
			Observe that each $u\in V_i$ induces $u'$ in $V_{i+1}$ 
			such that if $u_1,u_2\in V_i,u_1\neq u_2$ then $u'_1,u'_2\in V_{i+1}\quad u'_1\neq u'_2$. See proof of Proposition 1.3.4.
			Hence $|V_p|\geq |V_{i_0}|\geq 4$. Now the result follows again by P.H.P. 
		\end{description}
	\end{description}	
	\hfill $\Box$\par 
	\noindent {\it We shall use this lemma without mentioning it.}
	
	\subsection{Diamond moves}
	\begin{definition}
		\addcontentsline{toc}{subsubsection}{Definition \thesubsubsection\ (Diamond moves)}
		\ \\ Let $u\alpha v\beta w$ be a path, where $u,v$ and $w$ are vertices
		and $\alpha$ and $\beta$ are paths. See Fig. 10.
		Recall that a \emph{diamond move \diamondmove{v}	at $v$}
		is a surgery in $M$ in which we split $\alpha$ and $\beta$ into $\alpha_1, \alpha_2$ and $\beta_1, \beta_2$ and split $v$ into $v_1$ and $v_2$, respectively, 
		then cut out a hole $K$ in $M$ with boundary $\alpha_1 v_1 \beta_1^{-1} w \beta_2 v_2 \alpha_2^{-1} u$ and sew back $K$ by identifying $u$ with $w$.

		\begin{center}
			\begin{tikzpicture}[scale=.5]
				\draw (-.5,0)--(0,.5)--(.5,0)--(0,-.5)--(-.5,0);
				\node at (0,0) {$v$};
				
				\draw (2,0)--(10,0);
				
				\draw (11,-1)--(13,-1)--(12.5,-.7);;
				\draw (13,-1)--(12.5,-1.3);
				
				\draw (15,0)--(18,3)--(21,0)--(18,-3)--(15,0);
				
				\node at (23,-1) {\huge $\rightsquigarrow$};
				
				\draw (26,-4)--(26,4);
				
				\draw[fill] (2,0) circle(.1);
				\draw[fill] (6,0) circle(.1);
				\draw[fill] (10,0) circle(.1);
				\draw[fill] (15,0) circle(.1);
				\draw[fill] (21,0) circle(.1);
				\draw[fill] (18,3) circle(.1);
				\draw[fill] (18,-3) circle(.1);
				\draw[fill] (26,0) circle(.1);
				\draw[fill] (26,4) circle(.1);
				\draw[fill] (26,-4) circle(.1);
				
				\draw (4.5,.5)--(5,0)--(4.5,-.5);
				\draw (7.5,.5)--(7,0)--(7.5,-.5);
				\draw (16.5,2)--(17,2)--(17,1.5);
				\draw (16.5,-2)--(17,-2)--(17,-1.5);
				\draw (19,-1.5)--(19,-2)--(19.5,-2);
				\draw (19,1.5)--(19,2)--(19.5,2);
				\draw (25.5,2)--(26,2.5)--(26.5,2);
				\draw (25.5,-2)--(26,-2.5)--(26.5,-2);
				
				\node[below] at (2,0) {$u$};
				\node at (4.5,1) {$\alpha$};
				\node[below] at (6,0) {$v$};
				\node at (7.5,1) {$\beta$};
				\node[below] at (10,0) {$w$};
				\node[left] at (15,0) {$u$};
				\node at (16,2) {$\alpha_1$};
				\node[above] at (18,3) {$v_1$};
				\node at (20,2) {$\beta_1$};
				\node[right] at (21,0) {$w$};
				\node at (20,-2) {$\beta_2$};
				\node[below] at (18,-3) {$v_2$};
				\node at (16,-2) {$\alpha_2$};
				\node[above] at (26,4) {$v_1$};
				\node at (25,2) {$\alpha_1$};
				\node at (27,2) {$\beta_1$};
				\node[right] at (26,0) {$u=w$};
				\node at (25,-2) {$\alpha_2$};
				\node at (27,-2) {$\beta_2$};
				\node[below] at (26,-4) {$v_2$};
				\node at (18,0) {$K$};

			\end{tikzpicture}
			
			\hspace{10pt}
			
			Figure 10
		\end{center}
	\end{definition}

	\subsection{2-generated Artin groups}
	
	We recall in Lemma 1.5.1 well known results from \cite{1} on 2 -generated Artin groups.
	
	\vspace{10pt}
	
		\begin{lemma}
			\addcontentsline{toc}{subsubsection}{Lemma \thesubsubsection\ (Basic results on 2 generated Artin groups)}
			Let $M$ be a reduced van Kampen $\cal R$-diagram with connected interior and cyclically reduced boundary label, 
			$\cal R$ the cyclic closure of a single Artin relation
			$R=(ab)^na \left(b^{-1}a^{-1}\right)^nb^{-1}$ or $R = (ab)^n\left(a^{-1}b^{-1}\right)^n$, $n \geq 1$. Then each of the following holds:
			\begin{itemize}
				\item[(a)] $\cal R$ satisfies the small cancellation condition C(4) \& T(4).
				\item[(b)] Let $\omega$ be a boundary cycle of $M$ with decomposition $\omega = u\mu v\nu^{-1}$, 
				$u$ and $v$ vertices. 
				If $||\Phi(\mu) || \leq \frac{1}{2}||R||$ then $|\mu| \leq |\nu|$ and if $||\Phi(\mu) || < \frac{1}{2}||R||$ then $|\mu| < |\nu|$.
				\item[(c)] $||\omega|| \geq |R|$
				\item[(d)] if $D_1$ and $D_2$ are adjacent regions then $|\partial D_1\cap \partial D_2|\leq \frac 1 2 ||R||-1$.
			\end{itemize}
		\end{lemma}
		
		\vspace{10pt}
		
		\noindent {\bf Proof} 
		
		\begin{itemize}
			\item[(a)] This is Lemma 3 in \cite{1}. 
			\item[(b)]  This is Lemma 5 in \cite{1}.
			\item[(c)]  This is Lemma 7 in \cite{1}.
			\item[(d)]  Follows from the proof of Lemma 3 in \cite{1} 	\hfill $\Box$
		\end{itemize}
		We remark that part (a) holds true for arbitrary connected diagrams.

		\begin{proposition}
			\addcontentsline{toc}{subsubsection}{Proposition \thesubsubsection\ ($\partial\Delta=\eta_1\eta_2\eta_3,\ \|\Phi(\eta_1)\|=1$ implies
				$\mathcal{D}^v(\eta_2)\neq\emptyset$ or $\mathcal{D}^v(\eta_3)\neq\emptyset$	)}
			Let $\Delta$ be a connected, simply connected reduced $\mathcal{R}_0$-diagram with
			connected interior and cyclically reduced boundary label $\mathcal{R}_0=R(a,b),\quad a,b\in X$. Assume that $\partial\Delta$ has a decomposition by 
			$\eta_1^{-1} v_1\eta_2 v_2\eta_3^{-1} v_3$ such that $\Phi(\eta_1)=x^\alpha,\quad\alpha\neq 0,\quad x\in\{a,b\}$,\ 
			$v_1,v_2$ and $v_3$ vertices.\\
			Then for some $i,\ i\in\{2,3\}$ the following holds:
			$N(\eta_i)$ contains a one-layer diagram L with $\theta:=\partial L\cap \eta_i$ connected, $||\theta||\geq n(\Delta), 
			|\xi|\leq |\theta|$, where $\xi$ is the complement of $\theta$ on $\partial L$. See Fig. 11.
		\end{proposition}
		
		\begin{center}
			\begin{tikzpicture} [scale =0.5]
				\draw (-0.4,3.5) 
				node [draw,fill=black,circle,minimum size=2mm,
				inner sep=0pt,label={left:$v_1$}] (v1){}
				.. controls (-0.4,6.5) and (0.1,8) .. (2.1,9.5) 
				node [draw,fill=black,circle,minimum size=2mm,
				inner sep=0pt,label={above:$v_2$}] (v2){}
				node[pos=.67,sloped,label={[label distance=-1mm]above:$\eta_2$} ]
				{\tikz{\draw (-3.2,7.65) -- (-3,7.5) -- (-3.2,7.35);}}
				node[pos=.37,sloped,label={above:$\theta$} ] 
				{\tikz{\draw (0,.15) -- (0.2,0) -- (0,-.15);}}
				node[pos=.1](p1){}
				node[pos=.5,draw,fill=black,circle,minimum size=2mm,
				inner sep=0pt](p2){}
				;
				\begin{scope}[xscale=-1, shift={(-7.6,0)}]
					\draw (3,3.5)
					node [draw,fill=black,circle,minimum size=2mm,
					inner sep=0pt,label={right:$v_3$}] (v3) {}
					.. controls (3,6.5) and (3.5,8) .. (5.5,9.5) 
					node [pos=.5,rotate =120,label={below:$\eta_3$} ]
					{\tikz{\draw (-3.2,7.65) -- (-3,7.5) -- (-3.2,7.35);}}
					
					;
				\end{scope}
				
				\draw (p1.center) .. controls ++(2,0.5) and ++(2,-1) .. (p2)
				node [pos=.7,sloped,label={[label distance=-1mm]above:$\xi$} ]
				{$<$}
				;
				\draw  (v1) -- (v3)
				node [pos=.5,label={below:$\eta_1$} ]
				{\tikz{\draw (-3.2,7.65) -- (-3,7.5) -- (-3.2,7.35);}}
				;
				\node at (2,5) {$\Delta$};
				\node at (0.5,5.5) {$L$};
			\end{tikzpicture}
			
			Figure 11
		\end{center}
		
		\noindent 
		We need the following Lemmas and definitions for the proof of the Proposition.

		\begin{definition}[terminal component]
			\addcontentsline{toc}{subsubsection}{Definition \numberline{\thesubsubsection} (terminal component)}
			\ \\ Let $\Lambda$ be a layer structure of $\Delta$. 
			Let $L_j$ be a layer of $\Lambda$ and assume 
			$Int(L_j)=Q_1\cup\cdots\cup Q_r,\quad r\geq 1,\quad Q_i$ the connected components.
			Say that $Q_i$ is a {\it terminal component of $L_j$} if 
			$\partial Q_i\cap\partial L_{j+1}$ is either empty or consists of at most two vertices.
			Thus, every connected component of $L_p$ is terminal.
			See Fig. 12. where $Q_1$ is terminal and $Q_2$ is not.
		\end{definition}
		\begin{center}
			
			\begin{tikzpicture} [scale =0.5]
				
				\node  (v1) at (-0.4,3.5) {};
				\node  (v3) at (7.5,3.5) {};

				\draw  (v1.center) -- (v3.center)
				node [pos=.5,label={below:$\eta_1$} ]
				{};
				
				\node at (-1,5) {$L_j$};
				\draw (v1.center) .. controls (-0.5,6.5) and (7.5,6.5) .. (v3.center);
				\draw (0.5,5.1) .. controls (0.5,8.1) and (3,8) .. (3.3,5.7);
				\node at (2.1,6.5) {$Q_1$};
				\draw (4.5,5.7) .. controls (4.5,8) and (7.5,7.5) .. (6.9,4.8);
				\draw  plot[smooth, tension=.7] coordinates {(5,7) (5.1,8.5) (6.5,8.5) (6.5,6.8)};
				\node at (5.5,6.5) {$Q_2$};
				\node at (7.5,8) {$L_{j+1}$};
			\end{tikzpicture}
			
			Figure 12
		\end{center}
		
		\begin{definition}[Maximality condition ]
			\addcontentsline{toc}{subsubsection}{Definition \numberline{\thesubsubsection}  (Maximality condition)}
			Let $D_1,\ldots,D_{k_i}$ be the regions of $L_i$ in $\Lambda$. Let 
			$\delta_i=\sum_{j=1}^{k_i-1}|\partial D_j\cap\partial D_{j+1}|$
			and let $\delta(\Lambda)=\sum_{i=1}^{p}\delta_i$. Say that $\Delta$ 
			{\em satisfies the maximality condition} if $\delta(\Lambda)$ is maximal among all $\delta(\Lambda')$, 
			where $\Lambda'$ is a convex layer structure of $\Delta'$ with $L'_0=\mu_1$ and is obtained from $\Lambda$ by diamond moves
		\end{definition}
		\begin{definition}[separating vertices]
			\addcontentsline{toc}{subsubsection}{Definition
				\numberline{\thesubsubsection}  (separating vertices)}
			Let $D$ be a region in $\Delta$. Then $\Delta$ has a boundary cycle
			$v\mu u\nu^{-1},\ \ v,u$ vertices such that $\Phi(\mu)$ and $\Phi(\nu)$ are positive words. 
			The vertices $u$ and $v$ are called {\em the separating vertices of D}. The vertices u and v are unique.
			It is easy to see that the separating vertices cannot be inner vertices of pieces.
			See Fig. 13
		\end{definition}
		
		\begin{center}
			\begin{tikzpicture}[scale = 0.5]
				
				\draw (-2.5,0) node (v1) {} .. controls (-3,2.5) and (1,5.5) .. (3.5,5) node (v2) {}
				node[pos=.67,sloped ]{$>$}
				;
				\draw[dashed]  (v1) edge (v2);
				\draw (v1.center) .. controls (0,-0.5) and (4,3) .. (v2.center)
				node[pos=.67,sloped ]{$>$}
				;
				\node [draw,fill=black,circle,minimum size=2mm,
				inner sep=0pt,label={left:$u$}] at (v1){};
				\node [draw,fill=black,circle,minimum size=2mm,
				inner sep=0pt,label={right:$v$}] at (v2){};
				\node at (-1,3.7) {$\mu$};
				\node at (1.6,1.1) {$\nu$};
				\node at (3.5,1.5) {$D$};
			\end{tikzpicture}
			
			Figure 13
		\end{center}
		
		\begin{lemma}
			\addcontentsline{toc}{subsubsection}{Lemma \numberline{\thesubsubsection} ($\Delta$ has convex layer structure based on $\eta$)}
			\ \\ Follow the notation of 1.3.5 and 1.3.6.\\
			$\Delta$ has convex layer structure $\Lambda=(L_0,L_1,\ldots,L_p),\quad
			L_0=\eta_1,\quad p\geq 1$ such that $\alpha(D)\leq 1,\beta(D)\leq 2$ and 
			$d_{L_i}(v)\leq 3$, for every region $D\in L_i,i=1,\ldots,p$ and vertex
			$v\in\partial L_i\cap \partial L_{i+1}$.
			See Remark 1.3.10(b)
		\end{lemma}	
		
		\ \par
		\noindent {\bf Proof of lemma 1.5.6}\\
		By induction on p.
		\begin{description}
			\item [$\mathbf {p=1}$] Consider the following condition
			\begin{equation*}\tag {*}
				\left.
				\begin{minipage}{4in}
					{\em If $D_1$ and $D_2$ are regions in $L_1$ such that $\partial D_1\cap \partial D_2\neq\emptyset$ then $\partial D_1\cap \partial D_2\cap \eta_1\neq\emptyset$}
				\end{minipage}
				\right\}
			\end{equation*}
			\noindent Observe that if condition (*) holds then for every subdiagram $S\subseteq L_1$, 
			$L_1\setminus S$ is simply connected. Hence it is enough to prove (*). 
			Suppose by way of contradiction that (*) does not hold true. 
			Then there are regions $D_1$ and $D_2$ in $L_1$ such that 
			$\partial D_1  \cap \partial D_2\neq \emptyset$ but 
			$\partial D_1  \cap \partial D_2\cap\eta_1= \emptyset$.
			See Fig 14(a).
			Let $\mathbb E^2\setminus(D_1\cup D_2\cup\eta_1)=
			Q_1\cup  Q_\infty , \ \ r\geq 1$, 
			where $Q_1$ is the bounded connected component and $Q_\infty$ is the unbounded component. Consider $Q_1$.
			If $E$ is a boundary region of $Q_1$ with $\partial E\cap 
			\partial Q_1$ connected,
			then $\partial E\cap\partial Q_1\subseteq
			\partial D_1\cup \partial D_2\cup \eta_1$.
			If $Q_1$ consists of a single region then since $\Phi(\eta_1)=$
			$a^l$, \quad $l\neq 0$, 
			we have $|\partial E\cap\eta_1|= l$.
			Also, $\partial E\cap\partial D_i$ is a piece $i=1,2$.
			Hence, due to Lemma 1.5.1(d), $|\partial E\cap\partial D_i|
			\leq n(\Delta)-1,\ i=1,2$, $E$ is a region of $\Delta$ with $\partial E$ the product of three pieces, 
			namely $\partial E\cap\partial D_1$,
			$\partial E\cap \partial D_2$ and $\partial E\cap \eta_1$
			such that $|\partial E|\leq 2 n(\Delta)-1$, violating $|\partial E|= 2 n(\Delta)$.
			Hence $|Q_1|\geq 2$.
			Therefore by Lemma 1.3.2(b) $|\mathcal D(Q_1)|\geq 2$.
			Assume $\mathcal{D}(Q_1)=\{E_1,E_2\}$.
			Then by Lemma 1.3.2(b) $i(E_1)=i(E_2)=1$, hence $\partial E_i\cap \partial Q_1$ is the
			product of at least three pieces (each) due to the C(4) condition.
			In particular $\partial E_i\cap \partial Q_1\not\subseteq\partial D_j,\quad
			j=1,2$ and $\partial E_i\cap\partial Q_1\not\subseteq\eta_1$.
			Consequently, each $\partial E_i\cap\partial Q_1$ contains at least two of 
			$\{u_1,u_2,u_3\}$ in Fig. 14(a).
			Since $Q_1$ is planar, this is not possible.
			The impossibility of the cases $|\mathcal{D} (Q_1)|\geq 3$ follows by similar arguments. We omit details.
			Hence (*) holds true.
			It also follows from (*) that $\beta(D)\leq 2$.
			It remains to show that $d_{L_1}(v)\leq 3$, for every vertex $v$ on
			$\theta:=\partial L_1\setminus\eta_1$.
			If not, then there is a region $F$ in $L_1$ such that 
			$\partial F\cap\theta=\{v\}$.
			Since $\beta(F)\leq 2$, the C(4) condition and Lemma 1.5.1(d)
			imply that $\|\partial F\cap\eta_1\|\geq 2$.
			But $\|\eta_1\|=1$, a contradiction.
			Hence $d_{L_1}(v)\leq 3$.
			\begin{center}
				\begin{tabular}{cc}
					\begin{tikzpicture}[scale=0.5]
						\node [draw,fill=black,circle,minimum size=2mm,inner sep=0pt
						,label={[label distance=-1mm]below right:$u_2$}]  (v1) at (-0.5,4) {};
						\draw (-4,2) coordinate  (v2) .. controls (0,-1.5) and (5.5,2) .. (3,6.5) coordinate  (v3)
						node[pos=0.5,sloped,label={below:$\eta_1$}] {$>$}
						;
						\draw  plot[smooth, tension=.7] coordinates {(v2) (-2.5,3.5) (v1)};
						\draw  plot[smooth, tension=.7] coordinates{(v1) (1,5.5) (v3) };
						\draw  plot[smooth, tension=.7] coordinates {(-2.45,1) (v1) (3.55,5)};
						
						\node at (1,2.5) {$Q_1$};
						\node at (-3,2) {$D_1$};
						\node at (2.5,5.5) {$D_2$};
						\node at (-1.7,1.3) {$u_1$};
						\node at (3,4.5) {$u_3$};
					\end{tikzpicture}
					&
					\begin{tikzpicture}[scale=0.4]
						\node [draw,fill=black,circle,minimum size=2mm,inner sep=0pt] (v1) at (-0.5,5) {};
						\draw (-6,2) coordinate  (v2) .. controls (-6.5,2.5) and (-5,10) .. (-3.5,10.5) coordinate  (v3)
						node[pos=0.5,sloped,label={above:$\eta_2$}] {$>$}
						node[pos=0.75] (v4){}
						node[pos=0.25] (v5){}
						;
						\draw  plot[smooth, tension=.7] coordinates {(v2) (-2.5,2.5) (v1)(-1,8.5) (v3)};

						\node at (-3.5,5.5) {$Q$};
						\node at (-5,2.5) {$E_2$};
						\node at (-3,9.5) {$E_1$};
						
						\draw (v4.center) .. controls (0.5,8.5) and (-2,1.5) .. (v5.center)
						node[pos=0.5,fill=black,circle,minimum size=2mm,inner sep=0pt] (v6){}
						;
						
						\draw (v1.center) -- (v6.center);
					\end{tikzpicture}\\
					(a) & (b)
				\end{tabular}
				
				Figure 14
			\end{center}
			\item [The case $\mathbf{p\geq 2}$] By the induction hypothesis $d_{S_{p-1}}(u)\leq 3$
			for every vertex $u$ on $\partial L_p\cap\partial L_{p-1}$.
			It follows that the arguments of the proof of 1.3.7 apply here, 
			showing that $\Lambda$ has convex layer structure.
			We show that $\alpha(D)\leq 1,\beta(D)\leq 2$ and
			$d_{L_p}(w)\leq 3$ for every  vertex $w\in \partial P\cap\partial\Delta$ and region $D$ in $L_p$.
			If $\beta(D)\geq 3$ then there is a region $D_1$ in $L_p$ such that 
			$\partial D\cap \partial D_1\neq \emptyset$ but 
			$\partial D_1\cap \partial D\cap \omega_{p-1}= \emptyset$,
			where $\omega_{p-1}=\partial L_p\cap\partial L_{p-1}$. Let 
			$\mathbb E^2\setminus(D\cup D_1\cup \omega_{p-1})=
			Q_1\cup Q_\infty$, where $Q_1$ is the bounded connected component and $Q_\infty$ is the unbounded component.
			At this point the arguments in the proof of Lemma 1.3.6 apply and show that 
			$Q_1=\emptyset$.
			Hence $\beta(D)\leq 2$.
			We show that $\alpha(D)\leq 1$.
			If $\alpha(D)\geq 2$ then $L_{p-1}$ contains regions $E_1$ and $E_2$ such that
			$\partial E_i\cap \partial D$ contains and edge, $i=1,2$ 
			and $\partial E_1\cap\partial E_2\cap\partial D$ is a vertex $w$ on $\omega_{p-1}$.
			By the induction hypothesis $d_{L_{p-1}}(w)\leq 3$.
			Since $w\in \partial D, d_{L_p}(w)=2$.
			Since $w\in L_{p-1}\cup L_p$, $w$ is an inner vertex of $\Delta$ and
			$d_\Delta(w)=d_{L_{p-1}}(w)+d_p(w)-2\leq 3+2-2=3$, violating    the T(4) condition.
			Hence $\alpha(D)\leq 1$. 
			Finally, assume that $d_{L_p}(w)\geq 4$, for a vertex $w$ in 
			$\partial L_p\cap \partial\Delta$.
			Then $L_p$ contains a region E with $\partial E\cap \partial\Delta=\{w\}$.
			Since $\beta(E)\leq 2$, it follows from C(4) condition that 
			$4\leq \beta(E)+\alpha(E)$.
			Hence $\alpha(E)\geq 2$, violating $\alpha(E)\leq 1$.
			This completes the proof of the Lemma.  \hfill $\Box$
		\end{description}
		\noindent {\bf Notation}\\
		\begin{itemize}
			\item [(a)] Denote by $\mathcal{D}^v(\eta_i),\ i=2,3$ the collection of all the regions D and all the one layer diagrams S in $N(\eta_i)$ which satisfy each of the following.
			Let $K\in\{S,D\}$. Then
			\begin{itemize}
				\item [i)] $\partial K\cap \eta_i$ is connected and if $K=S$ then
				$\partial K\cap \eta_i$ is a side of S
				\item[ii)] $\|\partial K\cap \eta_i\|\geq n(\Delta)$
				\item [iii)] $|\xi|\leq |\partial K\cap \eta_i|$ where $\xi$ is the complement of $\partial K\cap \eta_i$ on $\partial K$.
				Notice that if $\mathcal{D}^v(\eta_2)=\emptyset$
				and $\mathcal{D}^v(\eta_3)=\emptyset$
				then $|\mathcal{D}^v(\eta_2 v_2\eta_3)|\leq 1$.
			\end{itemize}
			\item [(b)] Let $\mu$ be a boundary path of $\Delta$ and
			consider $N(\mu)$. Denote the complement of $\mu$ on $\partial N(\mu)$ by $\mu^*$
			\item [(c)] 
			{\bf Definition}
			Let $L=\langle D_1,\ldots,D_k\rangle,\ k\geq 4$ be a one layer diagram with connected interior 
			and sides $\mu$ and $\nu$, See Fig 15.
			Say that {\it $L$ is of type K relative to $\mu$} if each of the following holds
			\begin{itemize}
				\item [i)] all vertices on $\mu$ have valency 3.
				\item [ii)] $\partial D_1\cap\partial D_2\cap\partial\nu$ and
				$\partial D_{k-1}\cap \partial D_k\cap\nu$ are vertices with valency 4 and all the other vertices on $\nu$ have valency 3.
			\end{itemize}
			Define {\it $L$ is of type K relative to $\nu$} accordingly.
			\begin{center}
				\begin{tabular}{cc}
					\begin{tikzpicture}[yscale=0.5,xscale =-0.5]
						
						\draw (-4,2.5) node (v1) {} -- (-2,2.5) node (v3) {} -- (0,2.5) node (v5) {} 
						node [pos=0.5] {$>$}
						--  (2,2.5) node (v7) {} -- (4,2.5) node (v9) {} -- (6,2.5) node (v11) {};
						\draw[yshift=-3cm] (-4,2.5) node (v2) {} -- (-2,2.5) node (v4) {} -- (0,2.5) node (v6) {} -- 
						(2,2.5) node (v8) {} -- (4,2.5) node (v10) {} -- (6,2.5) node (v12) {};

						\draw  (v1.center) edge (v2.center);
						\draw  (v3.center) edge (v4.center);
						\draw  (v5.center) edge (v6.center);
						\draw  (v7.center) edge (v8.center);
						\draw  (v9.center) edge (v10.center);
						\draw  (v11.center) edge (v12.center);
						\draw[dashed]  (v3.center) edge (v2.center);
						\draw[dashed]  (v5.center) edge (v4.center);
						\draw[dashed]  (v7.center) edge (v6.center);
						\draw[dashed]  (v9.center) edge (v8.center);
						\draw[dashed]  (v11.center) edge (v10.center);
						\node at (-1,3) {$\mu$};
						\node at (-1,-1) {$\nu$};
						\draw (v2.center) -- (-4,-3.5) node (v13) {} -- (-2,-3.5) -- (v4.center);
						\draw [dashed] (v13.center) edge (v4.center);
						\draw (v10.center) -- ++(0,-3) node (v15) {} -- ++(2,0) -- (v12.center);
						\draw [dashed] (v15.center) edge (v12.center);
						\node at (6.5,-4.5) {$L$};
					\end{tikzpicture}
					&
					\begin{tikzpicture}[scale=0.5]
						
						\draw (-4,2.5) node (v1) {} -- (-2,2.5) node (v3) {} -- (0,2.5) node (v5) {} -- (2,2.5) node (v7) {} -- (4,2.5) node (v9) {} -- 
						(6,2.5) node (v11) {};
						\draw[yshift=-3cm] (-4,2.5) node (v2) {} -- (-2,2.5) node (v4) {} -- (0,2.5) node (v6) {} -- 
						(2,2.5) node (v8) {} -- (4,2.5) node (v10) {} -- (6,2.5) node (v12) {};

						\draw  (v1.center) edge (v2.center);
						\draw  (v3.center) edge (v4.center);
						\draw  (v5.center) edge (v6.center);
						\draw  (v7.center) edge (v8.center);
						\draw  (v9.center) edge (v10.center);
						\draw  (v11.center) edge (v12.center);
						\draw[dashed]  (v3.center) edge (v2.center);
						\draw[dashed]  (v5.center) edge (v4.center);
						\draw[dashed]  (v7.center) edge (v6.center);
						\draw[dashed]  (v9.center) edge (v8.center);
						\draw[dashed]  (v11.center) edge (v10.center);
						\node at (1,3) {$\mu$};
						\node at (1,-1) {$\nu$};
						\draw (v2.center) -- (-4,-3.5) node (v13) {} -- (-2,-3.5) -- (v4.center);
						\draw [dashed] (v13.center) edge (v4.center);
						\draw (v10.center) -- ++(0,-3) node (v15) {} -- ++(2,0) -- (v12.center);
						\draw [dashed] (v15.center) edge (v12.center);
						\node at (6.5,-4.5) {$L$};
					\end{tikzpicture}\\
					(a) & (b)
				\end{tabular}
				
				Figure 15
			\end{center}
		\end{itemize}
		\ \par\ \par
		\begin{lemma}
			\addcontentsline{toc}{subsubsection}{Lemma \numberline{\thesubsubsection}(basic properties of separating vertices)}\ \\
			Let $\eta\in\{\eta_2,\eta_3\}$ and consider $N(\eta)$.
			Assume
			\begin{itemize}
				\item [(I)] $\Delta$ is reduced;
				\item [(II)] $\Phi(\partial\Delta)$ is cyclically reduced;
				\item [(III)] $\Delta$ satisfies the maximality condition;
				\item [(IV)] $\mathcal{D}^v(\eta)=\emptyset$.
			\end{itemize}
			Then each of the following holds
			\begin{itemize}
				\item [(a)]Let D and E  be adjacent regions in $L_i':=L_i\cap N(\eta), i=1,\ldots,p$.
				Assume 
				\begin{itemize}
					\item [(i)] $\alpha(D)=1$ and $\alpha(E)=1$
					\item [(ii)] $\partial D\cap\eta$ and $\partial E\cap\eta$
					contains an edge each.
					Let $u$ and $z$ be the endpoints of $\xi:=\partial D\cap\partial E,\ u\in\eta$ and $z\in\omega_{i-1}$.
					Then $z\in\eta^*$ and either z is a separating vertex of $D$ and $u$ is a separating vertex of $E$, or the other way around.
					
					\noindent We define $Or(D)=1$ if $u$ is a separating vertex of D and $z$ is a separating vertex of $E$ and
					define $Or(D)=-1$ if $u$ is a separating vertex of $E$ and $z$ is a separating vertex of $D$.
					Similarly for E. See Fig. 16.
				\end{itemize}
				\begin{center}
					\begin{tikzpicture}
						
						\draw (-4.5,3.5) node (v3) {} -- (-4.5,2) node (v4) {} -- 
						(-2.5,2) node[draw,fill=black,circle,minimum size=2mm,
						inner sep=0pt,label={below:$z$}] (v1) {} -- 
						(-2.5,3.5) node[draw,fill=black,circle,minimum size=2mm,
						inner sep=0pt,label={above:$u$}] (v2) {} -- 
						(-0.5,3.5) node (v5) {} -- (-0.5,2) node (v8) {} -- (v1.center) -- (v2.center) -- (v3.center);
						\draw[dashed]  (v2) edge (v4);
						\draw[dashed]  (v5) edge (v1);
						\node at (-3,2.5) {$D$};
						\node at (-1,2.5) {$E$};	
						\node (v7) at (0.2,3.5) {};
						\draw  (v5.center) edge (v7.center);
						\node (v9) at (0.2,2) {};
						\draw  (v8.center) edge (v9.center);
						\node at (-0.15,3.75) {$\eta$};
						\node at (0.2,2.25) {$\omega_{i-1}$};
						
						\begin{scope}[xshift=6cm]
							\draw (-4.5,3.5) node (v3) {} -- (-4.5,2) node (v4) {} -- 
							(-2.5,2) node[draw,fill=black,circle,minimum size=2mm,
							inner sep=0pt,label={below:$z$}] (v1) {} -- 
							(-2.5,3.5) node[draw,fill=black,circle,minimum size=2mm,
							inner sep=0pt,label={above:$u$}] (v2) {} -- 
							(-0.5,3.5) node (v5) {} -- (-0.5,2) node (v8) {} -- (v1.center) -- (v2.center) -- (v3.center);
							\node[yshift=0.5cm] at (-3,2.5) {$D$};
							\node[yshift=0.5cm] at (-1,2.5) {$E$};
							\node (v7) at (0.2,3.5) {};
							\draw  (v5.center) edge (v7.center);
							\node (v9) at (0.2,2) {};
							\draw  (v8.center) edge (v9.center);
							\node at (-0.15,3.75) {$\eta$};
							\node at (0.2,2.25) {$\omega_{i-1}$};
						\end{scope}
						\draw[dashed]  (v3) edge (v1);

						\draw[dashed]  (v2) edge (v8);

					\end{tikzpicture}
					
					Figure 16
				\end{center}
				\item[(b)] Let $C$ be a terminal connected  component of $L_i$ in $L_i'$.
				Then i) and ii) of part (a) hold for every region D of C.
				\item [(c)] $C$ in part (b) contains a region from $\mathcal{D}^v(\eta)$.
				In particular, due to (IV), $C\not\subseteq N(\eta)$.
				\item [(d)] $Int (L_i)$ is connected.
				\item [(e)] $\alpha(D)=1$ for every region of $N_1(\eta):=
				N(\eta)\setminus D_0$, where $D_0$ is the region which contains $v_2$.
				\item [(f)] If $D$ and $E$ are adjacent regions in $N(\eta)$ with $\xi=\partial D\cap \partial E$ having endpoints $u$ and $w$,\ $u\in\eta, w\in\eta^*$, then either $u$ is a separating vertex of $D$
				and $w$ is a separating vertex of $E$ or the other way around, 
				such that $Or(E)=Or(D)$.
				\item [(g)] Every vertex on $\eta$ has valency at most 4.
				\item [(h)] $N_1(\eta)$ is a one-layer diagram without subdiagram of type K.
			\end{itemize}
		\end{lemma}
		\noindent {\bf Proof} \ Assume $\eta=\eta_2$.
		The case $\eta=\eta_3$ is dealt with similarly.
		\begin{itemize}
			\item [(a)] Due to 1.5.5 and Lemma 1.5.1(d), if $i\geq 2$ then $\eta$ contains a separating vertex of $D$ and a separating vertex of $E$.
			If $i=1$ then this is true due to 1.5.5 and Lemma 1.5.1 and the fact that $|\partial D\cap \eta_1|=|\partial E\cap \eta_1|=1<n$.
			But then there is a cancellation in $\Phi(\partial \Delta)$ at $u$, violating (II).
			Hence $u$ is either a separating vertex of $D$ or of $E$, 
			but due to (I) not for both.
			Suppose $u$ is a separating vertex of $D$.
			Then the other separating vertex of $D$ cannot occur on $\eta$,
			since then (V)  is violated.
			Hence $D$ has a separating vertex on $\eta^*$.
			Due to 1.5.5, it is not $z$.
			If $z$ is not a separating vertex of $E$ then (III) is violated.
			Hence $w$ is a separating vertex of $E$.
			The other separating vertex of $E$ is on $\eta$.
			If $u$ is a separating vertex of $E$ then the above arguments apply to $E$ in place of $D$ and the result follows.
			
			\item [(b) and (c)] Let $C=\langle D_1,\ldots, D_k\rangle,\ k\geq 1$.
			By Lemma 1.3.6 $z_j:=\partial D_j\cap\partial D_{j+1}\cap \eta$ have valency 3, $j=1,\ldots, k-1$.
			Hence $\partial D_j\cap \eta$ contains an edge, for every $D_j\in C$.
			Assume $\alpha(D_j)=0$, for some $j$.
			Let $j_0$ be the smallest index for which $\alpha(D_j)=0$.
			$k\geq 2$, since $\Delta$ has connected interior.
			If $j_0=1$ then $\beta(D_1)=1$, hence $i(D_1)=\alpha(D_1)
			+\beta(D_1)\leq 1$.
			Consequently, $D_1\in \mathcal{D}^v(\eta)$.
			Hence $j_0\geq 2$ and $\alpha(D_j)=1$, for $j=1,\ldots,j_0-1$.
			Consider $D_{j_0}$.
			Either it has both of its separating vertices on $\eta$ or one
			on $\eta$ and the other on $\eta^*$.
			In the first case $D_j\in\mathcal{D}^v(\eta)$, violating (IV),
			hence we consider the second case.
			Let $z'_1:=\partial D_{j_0-1}\cap\partial D_{j_0}\cap \eta$.
			$z'_1$ is not a separating vertex of $D_{j_0}$, due to 1.5.5.
			If $z'_1$ is not a separating vertex of $D_{j_0-1}$ then 
			$\Phi(\partial\Delta)$ is not reduced at $z'_1$.
			Hence $z'_1$ is a separating vertex of $D_{j_0-1}$.
			If $j_0-1=1$ then $\beta(D_{j_0-1})=1$ and $D_{j_0-1}\in \mathcal{D}^v(\eta)$, violating (IV).
			Hence $j_0-1\geq 2$.
			Let $\theta_l=\partial D_l\cap \partial D_{l+1}$.
			Then $\theta_l$ is connected by Lemma 1.3.2(f).
			Let $u_l=\theta_l\cap\eta$ and let $w_l=\theta_l\cap\eta^*,\ 
			i\leq l\leq j_0-1$.
			Then i) and ii) of part (a) are satisfied by $D_l$,\ 
			$1\leq l\leq j_0-1$.
			Hence part (a) applies.
			
			It follows by successive applications of part a) that either $D_l$ has $u_l$ and $w_{l-1}$ as separating vertices, or $w_l$ and $u_{l+1}$ as separating vertices.
			Due to (I) $u_l$ and $w_{l-1}$ are the separating vertices of $D_l$.
			It follows that $D_1\in \mathcal{D}(\eta)$, violating (IV).
			See Fig. 17.
			\begin{center}
				\begin{tikzpicture} [scale=0.5]
					
					\draw (-4,2.5) node (v1) {} -- (-2,2.5) node (v3) {} -- (0,2.5) node (v5) {} -- (2,2.5) node (v7) {} -- (4,2.5) node (v9) {} -- (6,2.5) node (v11) {} -- (8.5,4.5) node (v13) {};
					\draw[yshift=-3cm] (-4,2.5) node (v2) {} -- (-2,2.5) node (v4) {} -- (0,2.5) node (v6) {} -- (2,2.5) node (v8) {} -- (4,2.5) node (v10) {} -- (6,2.5) node (v12) {} -- (8.5,4.5) node (v14) {};

					\draw  (v1.center) edge (v2.center);
					\draw  (v3.center) edge (v4.center);
					\draw  (v5.center) edge (v6.center);
					\draw  (v7.center) edge (v8.center);
					\draw  (v9.center) edge (v10.center);
					\draw  (v11.center) edge (v12.center);
					\draw  (v13.center) edge (v14.center);
					\draw[dashed]  (v3.center) edge (v2.center);
					\draw[dashed]  (v5.center) edge (v4.center);
					\draw[dashed]  (v7.center) edge (v6.center);
					\draw[dashed]  (v9.center) edge (v8.center);
					\draw[dashed]  (v11.center) edge (v10.center);
					\draw[dashed]  (v13.center) edge (v12.center);
					\node at (-3,-1) {$D_1$};
					\node at (7.5,0) {$D_{j_0}$};
					\draw[line width=2pt] (-4,-2) -- (v1.center) -- (v11.center) -- (10,5.7);
					\node at (-4.5,-1.25) {$\eta$};
				\end{tikzpicture}
				
				Figure 17
			\end{center}
			Hence, $\alpha(D)=1$ for every region D of C and $D_1\in\mathcal{D}^v(\eta)$
			\item [(d)] Due to part (c) and assumption (IV)
			\begin{equation*}\tag {**}
				\text{\it $L_i'$ does not contain a connected component of $L_i$}
			\end{equation*}
			Hence either $L'_i$ is contained in a connected component C of $L_i$ or in two (adjacent) connected components $C_1$ and $C_2$.
			In the first case $N(\eta)$ contains a terminal connected component $C_3$ of $L_j$ for some $j>i$, violating (**). See Fig 18(a).\\
			Assume the second case.
			See Fig 18(b).
			\begin{center}
				
				\begin{tabular}{cc}
					\begin{tikzpicture}[scale=0.5]
						\draw (-3,-1) -- (3,-1);
						\coordinate (c1) at (-1.5,0);
						\coordinate (c2) at (1.5,0);
						\draw  plot[smooth, tension=.7] coordinates 
						{(-2,-1)  (c1) (c2) (2,-1)};
						\coordinate (d1) at (-1,1);
						\coordinate (d2) at (1,1);
						\draw  plot[smooth, tension=.7] coordinates {(c1) (d1) (d2) (c2)};
						\draw  plot[smooth, tension=.7] coordinates {(d1) (-0.5,2) (0.5,2) (d2)};
						\node at (0,1.5) {$C$};
						\node at (1.25,1.5) {$L_j$};
						\node at (2.4,-0.5) {$L_i$};
						\node[rotate=80] at (-1.95,-0.5) {$>$};
						\node[rotate=70] at (-1.42,0.5) {$>$};
						\node[rotate=60] at (-0.7,1.7) {$>$};
					\end{tikzpicture}
					&	
					\begin{tikzpicture}[scale=0.5]
						
						\draw (-4,1) -- (-2,1) 
						.. controls (-1,3) and (2,3) .. (3,1) node [pos=0.7,sloped,
						label={left:$\eta_2$}] 
						{\tikz {\draw (0,.15) -- (0.2,0) -- (0,-.15);}} .. controls (4,3) and (5.5,2.5) .. (7.5,2.5);
						\node at (0.5,1.5) {$C_1$};
						
						\node at (5.5,1.5) {$C_2$};
						\node at (-3,1.7) {$L_i$};
					\end{tikzpicture}
					\\
					(a) & (b)
				\end{tabular}
				
				Figure 18
			\end{center}
			Let $H(L_i)$ be the left most region of $L_i$.
			Now, either  $N(\eta)$ contains a connected component of 
			$L_j, j>i$, or $C_1$ is the connected component which contains
			$H(L_i)$.
			(Observe that $H(L_i)\in N(\eta)$).
			Hence $C_1\subseteq L'_i$, again violating (**).
			Consequently, $L'_i$ is connected.
			
			\item [(e)] Let $D\in L'_i$.
			Then either $\partial D\cap \partial L_{i+1}=\emptyset$ or 
			$\partial D\cap \partial L_{i+1}\neq\emptyset$.
			Consider the first case.
			Then $\partial D\cap\eta$ contains an edge, due to the C(4)
			condition and Lemma 1.3.2.
			The arguments of the proof in part (b) apply here 
			and imply that $\alpha(D)=1$.
			Finally, assume $\partial D\cap \partial L_{i+1}\neq\emptyset$
			See Fig. 19.
			\begin{center}
				\begin{tabular}{cp{2cm}c}
					\begin{tikzpicture} [scale=0.5]
						
						\draw (-2.5,1.5) node (v1) {} -- (-2.5,-0.5) -- (1.5,-0.5) -- (1.5,3.5) -- (-0.5,3.5) -- (-0.5,1.5) node (v2) {} -- (v1.center);
						\node (v4) at (1.5,1.5) {};
						\node (v3) at (-0.5,-0.5) {};
						\draw  (v2.center) edge (v3.center);
						\draw  (v2.center) edge (v4.center);
						\node at (-1.5,0.5) {$E$};
						\node at (0.5,0.5) {$D$};
						\node at (0.5,2.5) {$F$};
					\end{tikzpicture}
					&  &
					\begin{tikzpicture} [scale=0.5]
						
						\draw (-3,-1.5) -- (-1,-1.5) node (v2) {} -- (1,-1.5) node (v4) {} -- (3,-1.5) node (v6) {} -- (5,-1.5)
						--   (5,0.5) node (v7) {} -- (5,2.5) -- (3,2.5) node (v5) {} -- (1,2.5) -- (1,0.5) node (v3) {} -- (-1,0.5) node (v1) {} -- (-3,0.5)-- cycle;
						\draw  (v1.center) edge (v2.center);
						\draw  (v3.center) edge (v4.center);
						\draw  (v5.center) edge (v6.center);
						\draw  (v3.center) edge (v7.center);
						\node at (-2,-0.5) {$D$};
					\end{tikzpicture}\\
					(a) & & (b)
				\end{tabular}
				
				Figure 19
			\end{center}
			If $\alpha(D)=0$ then $D$ induces a region $F$ in $L_{i+1}$ with 
			$F\in \mathcal{D}^v(\eta)$.
			Hence $\alpha(D)=1$ in this case as well.
			\item [(f)] If $D$ and $E$ are both in $L_i$ then this follows by part (a).
			Assume $D\in L_i$ and $E\in L_{i+1}$.
			If $u$ is a separating vertex neither of $D$ nor of $E$ then $\Phi(\partial\Delta)$ is not reduced, violating II.
			Hence $u$ is a separating vertex of $E$ or of $D$ (but not both).
			Consider $w$.
			Either it is an inner vertex or it is a boundary vertex.
			In the first case $E$ has a neighbour $E_1$ in $L_{i+1}$ and $D$ has a neighbour $D_1$ in $L_i$.
			If $w$ is neither a separating vertex of $D$ nor of $E$ then by part (a), $w$ is a separating vertex for both $E_1$ and $D_1$,
			implying that $E_1$ and $D_1$ cancel each other, violating (I).
			Hence $w$ is a separating vertex for D or for $E$.
			Finally, assume $w$ is a boundary vertex of $D$.
			If neither $D$ nor $E$ have neighbours in $L_i$ and $L_{i+1}$ respectively, then $\Phi(\partial\Delta)$ is not reduced, 
			violating (II),
			Hence assume $E$ has a neighbour $E_1$ in $L_{i+1}$ and
			$w$ is a separating vertex of $E_1$.
			Therefore, by part (a) for $E$ and $E_1$ (in $L_{i+1}$), $z$ is a separating vertex for $E$, where $z=\partial E\cap \partial E_1\cap \eta$.
			But then   $z$ and $v$ are the separating vertices of $E$, where $v\in\eta$.
			Hence, $E\in\mathcal{D}^v(\eta)$, violating (IV).
			Hence $w$ is a separating vertex for $E$ or for $D$.
			If $D$ has a neighbour $D_1$ in $L_i$ then the same arguments hold.
			\item [(g)] if $v\in L_1$ then by Lemma 1.3.6, $d_{L_i}(v)=3$.
			Hence, if $v\in\partial L_i$ but $v\not\in \partial L_{i+1}$ then $d_\Delta(v)=3$.
			If $v\in\partial L_i\cap\partial L_{i+1}$ then there are regions $D,D_1$ in $L_i$ and region $E$ in $L_{i+1}$, such that $v=\partial D\cap\partial D_1\cap\partial E$.
			Hence $d_\Delta(v)=d_D(v)+d_{D_1}(v)+d_E(v)-2=2+2+2-2=4$.
			\item [(h)] 
			Assume by way of contradiction that $N(\eta_2)$ is not a one-layer subdiagram. 
			Then it contains a region which has more than two neighbours in $N(\eta_2)$.
			This implies that there are regions $E_1$ and $E_2$ in $N(\eta_2)$ such that $\partial E_1 \cap \partial E_2$
			contains an edge in common and $\partial E_1 \cap \partial E_2\cap \eta_2=\emptyset$. See Fig. 14(b).  
			Let Q be the bounded connected component of $\mathbb E^2\setminus
			(\eta_2\cup E_1\cup E_2)$.
			Then we may assume that $Int(Q)$ is connected, due to Lemma 1.3.2(f).
			Consider the subdiagram $P:=Q\cup E_1\cup E_2$.
			It is connected and simply connected with connected interior.
			Hence Proposition 1.3.7 applies to it relative to any vertex v. 
			We choose  
			$v=\partial Q\cap \partial E_1\cap \partial E_2$.
			Then $E_1$ and $E_2$ are in the first layer $K_1$
			of $\Lambda(v):=(K_0,K_1,\ldots,K_q),\ q\geq 1,K_0=\{v\}$. 
			
			\noindent If  $q\geq 2$ then $K_q\subseteq N(\eta)$, violating part (c) of the Lemma.
			If $q=1$ then $K_1$ contains at least 4 regions due to the T(4) condition (Notice that $v$ is an inner vertex).
			Hence, it contains at least 2 regions apart of $E_1$ and $E_2$.
			By the definition of $K_1$, every region in $K_1$ contains $v$ and
			it follows from part (f) of the Lemma that one of these regions is in $\mathcal{D}^v(\eta)$, violating (IV).
			Hence, $N(\eta)$ is a one-layer diagram.
			Finally, if $L$ is of type $K$ relative to $\mu$ and $Or(L)=1$ then $D_2\in\mathcal{D}^v(\mu)$ and if
			$Or(L)=-1$ then $D_{k-1}\in\mathcal{D}^v(\mu)$.
			Similarly, 
			if $L$ is of type $K$ relative to $\nu$ and $Or(L)=1$ then $D_1\in\mathcal{D}^v(\nu)$ and if
			$Or(L)=-1$ then $D_{k}\in\mathcal{D}^v(\nu)$.
			In particular, $L\not\subseteq N(\eta_2)$ and
			$L\not\subseteq N(\eta_3)$, due to assumption (IV) in Lemma 1.5.7. See Fig. 15. (Replace $L$ with $N(\mu)$).)
			
		\end{itemize}
		\ \hfill $\Box$
		

		\vspace{10pt}
		
		\noindent {\bf Proof of the Proposition }\\
		Assume (I),(II),(III) and (IV) of Lemma 1.5.7
		We claim that $v_2\in\partial L_p$.
		If not, then either $L_p\subseteq N(\eta_2)$ or
		$L_p\subseteq N(\eta_3)$, violating (**) in part (d) of Lemma 1.5.7.
		Hence, $v_2\in L_p$.
		Also, it follows from parts (e) and (f) of Lemma 1.5.7 that one of the following holds.
		\begin{itemize}
			\item [i)] $Or(L)= -1, v_2 \in \partial T_p, \eta_2$ contains a separating vertex of $T_p$ and $\eta_3$ contains a separating vertex of $T_p$.
			Here $T_p$ is the right most region in $L_p$. 
			\item [ii)] $Or(L)=1, v_2 \in \partial H_p$ and $\eta_3$ contains a separating vertex of $H_p$ and $\eta_2$ contains a separating vertex of $H_p$.
			Here $H_p$ is the left most region of $L_p$. 
		\end{itemize}
		\noindent Here $L=N(\eta_2)$ or  $L=N(\eta_3)$
		%
		\noindent Let $\eta'^*_2$ be the complement of $\eta'_2:=\eta_2\setminus \partial T_p$, where $v_2\in T_p$, on $N(\eta_2)\setminus D_0$
		and define $\eta'^*_3$ accordingly.
		We claim that:\\ 
		in Case (i),
		$|\eta'_2|-|\eta'^*_2|=2(x-1)$, where $x=|\xi|,\ x\geq 1,\ \xi$
		is the tail of $\eta_2$ which start at $z_1$.
		Moreover, $\|\eta_2\|\geq n(\Delta)$. 
		Hence $N(\eta_2)\in\mathcal{D}^v(\eta_2)$, violating (IV).\\
		in Case (ii), 
		$|\eta'_3|-|\eta'^*_3|=2(y-1)$, where $y=|\theta|,\ y\geq 1,\  \theta$ is the tail of $\eta_3$ which starts at $v_2$. 
		Moreover, $\|\eta_3\|\geq n(\Delta)$. 
		Hence $N(\eta_3)\in\mathcal{D}^v(\eta_3)$, violating (IV).

		\noindent Consider the first. 
		It follows that $x\geq 1$, due to (IV). See Fig. 20.
		\begin{center}
			\begin{tabular}{cc}
				\begin{tikzpicture}[yscale=0.5,xscale=0.3]
					
					\draw (-2,-2.5) node (v1) {} .. controls (-2,2) and (-2,2) .. (3,2) 
					node[label={above:$v_2$}]{$\bullet$}
					node [pos=0.3,sloped,label={above:$\eta_2$}] {$>$}
					node [pos=0.5,label={[label distance=-4mm]above left:$z_1$}] (vl){$\bullet$}
					.. controls (8,2) and (8,2) .. (8,-2.5) 
					node [label={[label distance=-3mm]below right:$z_2$}](v2) {}
					node [pos=0.7,sloped,label={above:$\eta_3$}] {$<$}
					;
					\draw  (v1.center) edge (v2.center);
					\draw [dashed] (vl.center) edge (v2.center);
					\node at (-0.5,3) {$\xi$};
					\node at (0.5,-1.5) {$T_p$};
				\end{tikzpicture}
				&
				\begin{tikzpicture}[yscale=0.5,xscale=0.3]
					
					\draw (-2,-2.5) node [label={[label distance=-4mm]below 	left:$w_2$}] (v1) {} .. controls (-2,2) and (-2,2) .. (3,2) 
					node[label={above:$v_2$}]{$\bullet$}
					node [pos=0.3,sloped,label={above:$\eta_2$}] {$>$}
					.. controls (8,2) and (8,2) .. (8,-2.5) node (v2) {}
					node [pos=0.7,sloped,label={above:$\eta_3$}] {$<$}
					node [pos=0.5,label={[label distance=-4mm]above 	right:$w_1$}] (vr){$\bullet$}
					;
					\draw  (v1.center) edge (v2.center);
					\draw [dashed] (vr.center) edge (v1.center);
					\node at (7,2.8) {$\theta$};
					\node at (6,-0.5) {$w_p$};
					\node at (3.5,-1.8) {$H_p$};
				\end{tikzpicture}\\
				(a) & (b)
			\end{tabular}
			
			Figure 20
		\end{center}
		
		\noindent Let $N(\eta_2)=\langle D_1,\ldots, D_l\rangle,\ D_l=T_p$, 
		let $w_1$ and $w_2$ be the separating vertices of $H_p$. 
		Then $w_1\in \eta_2$ and $w_2\in \eta^*_2$. 
		Let $z_1$ and $z_2$ be the separating vertices of $T_p, z_1\in\eta_2$ and $z_2\in\eta^*_2$.
		Let $\eta'_2$ be the subpath of $\eta_2$ which starts at $w_1$ and ends at $z_1$ 
		and let $\eta'^*_2$ be the subpath of $\eta^*_2$ which starts at $w_2$ and ends at $z_2$. 
		Then it follows by an easy induction on $l$ using Lemma 1.5.7 that $|\eta'_2|=|\eta'^*_2|$. 
		Now $|\eta_2|=n(\Delta)-1+|\eta'_2|+x$ and 
		$|\eta^*_2|=1+|\eta'^*_2|+n(\Delta)-x$. 
		Hence $|\eta_2|-|\eta^*_2|=(n(\Delta)-1+x)-
		(1+n(\Delta)-x)+|\eta'_2|-|\eta'^*_2|=2(x-1)$
		as required. 
		Finally $\|\eta_2\|\geq n(\Delta)-1+x\geq n(\Delta)$. 
		The second case is dealt with similarly. \\
		To complete the proof of the Proposition we notice that we may choose $\Delta$ such that (I),(II) and (III) hold true.
		Hence, the contradictions obtained contradict (IV) as required.

		\subsection{Bands}
		
		We extend the definition of \simBt to $\underset{B_a}{\sim}$ for $a \in X$ (instead of $t \in T$) and denote $\kappa(D) = Int(\cup [D]_{B_a})$. 
		
		We call $\kappa(D)$ an \emph{a-band} if $\partial \kappa (D)$ is a simple closed curve with cyclically reduced boundary label and call $\kappa (D)$ a 
		\emph{closed  a-band} if $\kappa(D)$ is a one layer annular diagram with connected interior, simple closed boundary components, with cyclically reduced boundary labels. See Figure 21.

		\begin{minipage}{2.2in}
			\begin{center}
				\begin{tikzpicture}[scale=.3]
					
					\draw[rounded corners] (0,0)--(3,2)--(6,2)--(9,0)--(12,0)--(14,2);
					\draw[rounded corners] (0,3)--(3,5)--(6,5)--(9,3)--(12,3)--(14,5);
					
					\draw (0,0)--(0,3);
					\draw (2,1.3)--(2,4.3);
					\draw (4,2)--(4,5);
					\draw (6,1.9)--(6,4.9);
					\draw (8,.7)--(8,3.7);
					\draw (10,0)--(10,3);
					\draw (12,0.2)--(12,3.2);
					\draw (14,2)--(14,5);

					\draw (-.2,.7)--(0,1)--(.2,.7);
					\draw (1.8,2.7)--(2,3)--(2.2,2.7);
					\draw (13.8,2.7)--(14,3)--(14.2,2.7);
					
					\node[above] at (5,5) {$\nu$};
					
					\node at (5,1) {$\mu$};
					
					\node[left] at (0,1) {$\alpha_1$};
					\node[right] at (14,3) {$\alpha_2$};
					\draw (5.4,5.3)--(5.7,5)--(5.4,4.7);
					\draw (5.4,2.3)--(5.7,2)--(5.4,1.7);
				\end{tikzpicture}
				
				(a)
			\end{center}
		\end{minipage}\hspace{10pt}\begin{minipage}{1.1in}
			\begin{center}
				\begin{tikzpicture}[scale=.8]
					
					\draw (0,0) circle(2);
					\draw (0,0)--(2,0);
					\draw (0,0)--(-2,0);
					\draw (0,0)--(0,2);
					\draw (0,0)--(0,-2);
					\draw (0,0)--(1.73,1);
					\draw (0,0)--(1.73,-1);
					\draw (0,0)--(-1.73,1);
					\draw (0,0)--(-1.73,-1);
					
					\draw (0,0)--(1,1.73);
					\draw (0,0)--(1,-1.73);
					\draw (0,0)--(-1,1.73);
					\draw (0,0)--(-1,-1.73);
					
					\draw[thick] (0,0) circle(1.01);
					
					\draw[fill,white] (0,0) circle(1);
					
					\draw (-.2,1.5)--(0,1.7)--(.2,1.5);
					
					\node[above] at (0,2) {$\nu$};
					\node[left] at (0,1.5) {$a$};

				\end{tikzpicture}

				(b)
			\end{center}
		\end{minipage}\hspace{10pt}\begin{minipage}{1.5in}
			\begin{center}
				\begin{tikzpicture}[scale=1]
					\draw (0,0)--(2,0)--(2,2)--(0,2)--(0,0);
					\draw (0,1)--(2,1);
					
					\draw (-.2,.8)--(0,.5)--(.2,.8);
					\draw (-.2,1.2)--(0,1.5)--(.2,1.2);
					\draw (1.8,1.2)--(2,1.5)--(2.2,1.2);
					\draw (1.8,.8)--(2,.5)--(2.2,.8);
					\draw (.7,.8)--(1,1)--(.7,1.2);
					\draw (.7,1.8)--(1,2)--(.7,2.2);
					\draw (.7,-.2)--(1,0)--(.7,.2);
					
					\node at (-.2,1.75) {$x$};
					\node at (-.2,.25) {$x$};
					\node at (1,2.2) {$a$};
					\node at (1,-.2) {$a$};
					\node at (1,1.5) {$D_1$};
					\node at (1,.5) {$D_2$};
					\node at (2.2,1.5) {$x$};
					\node at (2.2,.5) {$x$};
					\node at (1,1.2) {$a$};

				\end{tikzpicture}

				(c)
			\end{center}
		\end{minipage}
		
		\begin{center}
			Figure 21
		\end{center}
		
		Notice that an a-band cannot intersect itself, because $aaa^{-1}a^{-1}$ is not a defining relator. Notice also that if $M$ is reduced then $\Phi(\partial B)$ is cyclically reduced, $B$ an a-band, since if $D_1$ and $D_2$ are adjacent regions in $B$ such that $\Phi(\omega_1) = a x a^{-1}x^{-1}$ and $\Phi(\omega_2) = a x^{-1}a^{-1}x$, $\omega_1$ and $\omega_2$ are boundary cycles of $D_1$ and $D_2$ then $(D_1, D_2)$ is a cancelling pair, contradicting the assumption that $M$ is reduced. See Figure 21(c).
		
		Let $B$ be an a-band, $B = \langle D_1, \ldots , D_k \rangle$, $ k \geq 1$. Then $B$ has boundary cycle $\alpha_1 \mu \alpha_2^{-1}\nu^{-1}$ with $\Phi(\alpha_1) = a^{\varepsilon}$, $\Phi(\alpha_2) = a^{-\varepsilon}$, $ \alpha_1 \subseteq \partial D_1$, $\alpha_2 \subseteq \partial D_k$ $ \varepsilon \in \{1,-1\}$, $\mu = \mu_1 \cdots \mu_k$, $\nu = \nu_1 \cdots \nu_k$, $ \mu_i = \mu \cap \partial D_i$, $\nu_i = \nu \cap D_i$. Call $\alpha_1$ and $\alpha_2$ the \emph{poles of $B$} and call $\mu$ and $\nu$ the \emph{sides of $B$}. See figure 7(a).
		
		\vspace{10pt}
		
		\tocexclude{\subsubsection {Definitions}}
		\addcontentsline{toc}{subsubsection}{\numberline{\thesubsubsection}  Definitions\quad  (band-bundles)}
		\begin{itemize}
			\item[(a)] Let $B_1$ and $B_2$ be bands in an $\cal R$-diagram $M$. Say that $B_1$ and $B_2$ are \emph{adjacent} if $B_1$ and $B_2$ have sides $\mu_1$ and $\nu_2$ respectively such that $o(\mu_1) = o (\nu_2)$ and either $\mu_1 \subseteq \nu_2$ or $\nu_2 \subseteq \mu_1$. 
			See Figure 22(a). (Recall from 1.1 that $o(\mu))$ denotes the origin of $\mu$ and $t(\mu)$ denotes its terminus.
		\end{itemize}
		\vbox{
			\begin {center}
			\begin{tabular}{ccccc}
				\begin{tikzpicture}[scale=.3]
					\draw (5,0)--(0,0)--(0,15)--(5,15)--(5,0)--
					(10,0)--(10,10)--(5,10);
					\node at (2.3,10) {$B_1$};
					\node at (7,7) {$B_2$};
					\node at (6,4) {$\nu_2$};
					\node[left] at (0,4) {$\nu_1$};
					\node[left] at (5,4) {$\mu_1$};
				\end{tikzpicture}
				& &
				\begin{tikzpicture}[scale=.3]
					\draw (0,0)--(0,15)--(10,15)--(10,0)--(0,0);
					\draw (5,0)--(5,15);
					\node at (2.3,7) {$B_1$};
					\node at (7,7) {$B_2$};
					\node at (6.3,4.5) {$\nu_2$};
					\node at (4,4) {$\mu_1$};
				\end{tikzpicture}
				& &
				\begin{tikzpicture}[scale=.5]
					\draw (4.5,4.5) arc (50:130:7);
					\draw (4.5,0) arc (50:130:7);
					\draw (-3,1) -- (-3,5.5);
					\draw (3,5.5) -- (3,1);
					\draw[fill] (-1,6) circle (0.15);
					\draw (-3,3) -- (3,3);
					\draw (-3,4.5) -- (3,4.5);
					\draw (-1,6) -- (-1,1.5);
					\draw (1,6) -- (1,1.5);
					\node at (-2,2) {$B_1$};
					\node at (0,2.1) {$B_2$};
					\node at (2,2) {$B_3$};
					\node at (-1,6.5) {$v_1$};
				\end{tikzpicture}\\
				(a)& &(b)& &(c)
			\end{tabular}
			\ \\ \ 
			
			Figure 22
		\end{center}
	}
		%
	%
		%
		%
		%
		%
		%
		%
		%
	%
	
	\begin{itemize}
		\item[(b)] Let $B_1$ and $B_2$ be adjacent bands. Say that $B_1 \cup B_2$ is \emph{homogeneous} if $t(\mu_1) = t(\nu_2)$ (hence $\mu_1 = \nu_2$). See Figure 22(b).
		\item[(c)] Let $B_1, \ldots , B_k$ be bands in $M$. Say that $(B_1, \ldots , B_k)$ form a \emph{band-bundle} $\mathbb B$ if $(B_i, B_{i+1})$ are homogeneous adjacent bands. See Figure 22(c).
		We call $k$ the \emph{width of $\mathbb B$} and we call the length of a band in $\mathbb B$ the \emph{height of 
			${\mathbb B}$}, and we call configuration 8(a) a \emph{stair} if $t(\mu_1) \neq t(\nu_2)$.
	\end{itemize}
	
	\begin{definition}[B-connected edges]
		\addcontentsline{toc}{subsubsection}{\numberline{\thesubsubsection} Definition \thesubsubsection\quad (B-connected edges)}\ \\
		Let $\Delta \in Reg_{4^+}({\mathbb M})$ and let $e$ be an edge of $\partial \Delta$. Say that $e$ is \emph{B-connected to $f$, $f$ an edge of $\partial M$} 
		if there is a band that starts on $e$ and terminates on $ f \in \partial M$. Let $\Delta_1 \in Reg_{4^+}(M)$. 
		Thus, one pole of the connecting band is $e$ and the other one is $f$.
		Say that $\Delta $ is \emph{B-connected to $\Delta_1$} if there is a band that starts on $\Delta$ and terminates on $\Delta_1$. We say that $\Delta$ and $\Delta_1$ are \emph{$B$-neighbours}.
	\end{definition}
	
	\vspace{10pt}
	
	\noindent {\bf Notation} $\partial \Delta \cap_B \partial M := \bigcup \{e \subseteq \partial \Delta \ , \ e \  \mbox{is B-connected to} \ f \subseteq \partial M \}$.
	We denote $f$ by $\beta(e)$ and denote $\cup \beta (e)$, $ e \in \partial \Delta$, $e$ B-connected to $\partial M$ by $\beta(\partial \Delta \cap_B \partial M$).
	
	\vspace{10pt}
	
	\begin{definition}[B-interior edge]
		\addcontentsline{toc}{subsubsection}{\numberline{\thesubsubsection} Definition \thesubsubsection\quad (B-interior edges)}
		\ \\Let notation be as in the previous definition and let $e$ be an edge of $\partial \Delta$. Say that $e$ is \emph{B-interior} if $ \{e\}\cap (\partial \Delta \cap_B \partial M) = \emptyset$. Denote by $i_B(\Delta)$ the number of $B$-interior neighbours of $\Delta$.
		
		Every edge of $\partial \Delta$ is B-connected either to an inner region or to the boundary of $M$.
	\end{definition}
	
	
	\section {The diagrams $\mathbb{M}$ and $\mathbb{M}_t, t \in T(M)$}
	\tocexclude{\subsection{The diagrams $\mathbb{M}$ and $\mathbb{M}_t, t \in T(M)$}}
	\addcontentsline{toc}{subsection}{\numberline{\thesubsection} $\mathbb{M}$ and $\mathbb{M}_t$}
	Let $M$ be a  simply connected $\cal R$-diagram over $F$, $\cal R$ given by ({\rm III}) and let $t\in T(M)$. 
	Denote by $Reg_{H_t}(M)$ the set of all the regions of $M$ the boundaries of which are labelled only by letters from $H_t$ ($H_t = \langle X-\{t\} \rangle$). 
	Denote by $Reg_t(M)$ the set of all the regions of $M$ the boundaries of which have labels from $\langle t \rangle$ and also from $H_t$. 
	Thus, $Reg(M) = Reg_t(M) \cup Reg_{H_t}(M)$. We subdivide $Reg_t(M)$ further by 
	\begin{equation*}
		Reg_t(M) = Reg_{4^+,t}(M) \cup Reg_{2,t}(M)
	\end{equation*}
	where $Reg_{4^+,t}(M)$ is the set of regions in $Reg_t(M)$ with boundary length at least 8 and $Reg_{2,t}(M)$ is the set of regions in $Reg_t(M)$ with boundary length 4. We have for $t \in T$
	\[ \tag{1}
	Reg(M) = Reg_{4^+,t}(M) \cup Reg_{H_t}(M) \cup Reg_{2,t}(M)
	\]
	For the definition of the diagrams $\mathbb M$ and $\mathbb M_t$ we introduce equivalence relations
	"$\underset{t}{\sim}$", "$\underset{H_t}{\sim}$" and "$\underset{B_t}{\sim}$" on $Reg_{4^+,t}(M)$, $Reg_{H_t}(M)$ and $Reg_{2,t}(M)$, respectively.

	\vspace{10pt}
	
	\begin{sssection}[Definition of $\underset{t}{\sim}$]
		\addcontentsline{toc}{subsubsection}{\numberline{\thesubsubsection}  
			Definition of $\protect\underset{t}{\sim}$}
		\ \\ Let $D_1$ and $D_2$ be regions in $Reg_{4^+,t}(M)$. Say that they are \emph{$t$-friends} if the following hold
		\begin{itemize}
			\item[($i$)] $\partial D_1 \cap \partial D_2$ contains an edge.
			\item[($ii$)] $Supp (\partial D_1) = Supp (\partial D_2)$.
		\end{itemize}
		Let \simt be the transitive closure of $t$-friendness. It is easy to see that \simt  is an equivalence relation on $Reg_{4^+,t}(M)$. Denote by $[D]_t$ the \simt equivalence class of $D$.
	\end{sssection} 
	\vspace{10pt}
	
	\begin{sssection}[Definition of \simHt]
		\addcontentsline{toc}{subsubsection}{\numberline{\thesubsubsection}  
			Definition of $\protect\underset{H_t}{\sim}$}
		\ \\ Let $D_1$ and $D_2$ be regions in $Reg_{H_t}(M)$. Say that they are \emph{$H_t$-friends} if
		$\partial D_1 \cap \partial D_2$ contains an edge. Let \simHt be the transitive closure of $H_t$-friendness. Then \simHt is an equivalence relation on $Reg_{H_t}(M)$. Denote by $[D]_{H_t}$ the equivalence class of $D \in Reg_{H_t}(M)$.
	\end{sssection} 
	\vspace{10pt}
	
	\begin{sssection}[Definition of \simBt]
		\addcontentsline{toc}{subsubsection}{\numberline{\thesubsubsection}  
			Definition of $\protect\underset{B_t}{\sim}$}
		\ \\ Let $D_1$ and $D_2$ be regions in $Reg_{2,t}(M)$. Say that $D_1$ and $D_2$ are \emph{$B_t$-friends} if
		$\partial D_1 \cap \partial D_2$ consists of an edge with label $t^{\pm 1}$. Let \simBt be the transitive closure of $B_t$-friendness. 
		Then \simBt is an equivalence relation. For $D \in Reg_{2,t}(M)$ denote by $[D]_{B_t}$ the equivalence class of $D$ in $Reg_{2,t}(M)$.
	\end{sssection}
	\vspace{10pt}
	\begin{remark}
		\addcontentsline{toc}{subsubsection}{\numberline{\thesubsubsection}  
			Remark: $M = \cup[D]_Z$ is a partition of $Reg(M)$}
		\text{ }
		\begin{itemize}
			\item [(a)]
			Let $D_1$ and $D_2$ be regions of $M$. Then $[D_j]_{Z} \neq \emptyset$ if and only if $D_j \in Reg_{Z}(M)$ where $[D_j]_Z$ is either $[D_j]_t$ or $[D_j]_{B_t}$ or $[D_j]_{H_t}$. 
			We shall denote these three possibilities by $Z=t^{4^+}$ or $Z=B_t$ or $Z=H_t$, respectively, i.e. $Z\in\{t^{4^+},B_t,H_t\}$. 
			Also, $[D_1]_{Z_1}= [D_2]_{Z_2} \neq \emptyset$ if and only if $Z_1 = Z_2$ and $D_1 \in [D_2]_{Z_2}$. If $Z_1 \neq Z_2$ or $D_1 \not\in [D_2]_{Z_2}$ then $[D_1]_{Z_1} \cap [D_2]_{Z_2} = \emptyset$.
			Clearly $\cup[D]_Z$, where $D \in Reg_Z(M)$ and $Z \in \{t^{4^+}, H_t, B_t\}$, is the whole of $M$. Hence $M = \cup[D]_Z$ given in (1) is a partition of $Reg(M)$.
			\item[(b)]
			Definition 2.1.1 makes sense without the requirement that $D_1$ and $D_2$ are in $Reg_{4^+,t}(M)$ namely, $D_1\sim D_2$ if $Supp(\partial D_1)=Supp(\partial D_2)$
			and $\partial D_1\cap\partial D_2$ contains an edge.
			We shall call the corresponding equivalence relation "Equivalence". 
			For $D \in Reg_Z(M)$ define $\kappa_0(D) = Int (\cup \{\overline{E} \ | \ E \in [D]_Z \}$. Then $\kappa_0(D)$ is connected.
			If $Z = t^{4^+}$ then define $n(\kappa_0(D)) = n(D)$. Notice that $n(\kappa_0(D))$ is well defined.
			Also, denote $\kappa_0(D)$ by $\Delta(D)$, for this case.
		\end{itemize}
	\end{remark}
	
	\vspace{10pt}
	
	\tocexclude{\subsubsection{Definition (adequate)} }
	\addcontentsline{toc}{subsubsection}{\numberline{\thesubsubsection}  
		Definition (adequate diagrams)}
	Let $M$ be a connected, simply connected $\cal R$-diagram over $F$, $M\in \mathcal{M}_3(W)$. 
	Say that $M$ is \emph{adequate}, if for every $D \in Reg(M),\quad D$ is homeomorphic to the open unit disc, $\partial D$ is simple,
	such that the boundary label is not $1$, and is cyclically reduced.
	
	\begin{definition}[${\mathbb M}_t$] 
		\addcontentsline{toc}{subsubsection}{\numberline{\thesubsubsection}  
			Definition (${\mathbb M}_t$)}	
		\ \\ Let $M$ be a van Kampen $\cal R$-diagram over $F$ and let $ t \in T(M)$. 
		If $\kappa_0(D)$ is homeomorphic to the open unit disc, for every $D\in Reg(M)$, 
		then we may consider $\{ \kappa_0(D) \ | \ D \in Reg(M) \}$ as a set of regions such that their totality generates a diagram which we denote by ${\mathbb M}_t$. 
		Thus $Reg({\mathbb M}_t) = \{ \kappa_0 (D) \ | \ D \in Reg(M) \}$. 
		We have  $\partial {\mathbb M}_t = \partial M$. Also,
		\begin{equation*}\tag{2}
			Reg(\mathbb M_t)=Reg_{4^+,t}(\mathbb M_t)\cup Reg_{B_t}(\mathbb M_t)\cup Reg_{H_t}(\mathbb M_t)\text {\ for every\ }
			t\in T(M).
		\end{equation*}
	\end{definition}
	The main results of this section, are Propositions 2.1.7 and 2.1.8 
	
	\begin{proposition}
		\addcontentsline{toc}{subsubsection}{\numberline{\thesubsubsection}  
			Proposition (${\mathbb M}_t$ is adequate)}		
		\ \\ Let $M \in \mathcal{M}_3(W)$. Then $\mathbb{M}_t$ is adequate, for every $t \in T(M)$. Moreover, if $K_1$ and $K_2$ are regions of $\mathbb{M}_t$, not both bands, then $\partial K_1 \cap \partial K_2 \neq \emptyset$ 
		implies that $\partial K_1 \cap \partial K_2$ is connected.
	\end{proposition}
	
	\begin{proposition}
		\addcontentsline{toc}{subsubsection}{\numberline{\thesubsubsection}  
			Proposition ($M$ is adequate)}		
		\ \\ 	Let $M\in \mathcal M_3(W)$. Then M  is adequate.  Moreover, if $\partial D_1\cap \partial D_2 \neq\emptyset$ then $\partial D_1\cap \partial D_2$ is connected.
	\end{proposition}
	\vspace{10pt}
	
	\begin{remark}
		\addcontentsline{toc}{subsubsection}{\numberline{\thesubsubsection}  
			Remark (The passage from ${\mathbb M}_t$ to 
			${\widetilde{\mathbb M}}^t$)}			
		\ \\ Let $K \in Reg_{4^+,t}({\mathbb M}_t)$. When passing from ${\mathbb M}_t$ to $\widetilde{{\mathbb M}^t}$ by shrinking edges labelled by elements of $H_t$ to a point and shrinking $t$-bands to edges, it may happen that the interior of the obtained region $\widetilde{K^t}$  in $\widetilde{{\mathbb M}^t}$ becomes disconnected. 
		Since $M$ is adequate, we can replace  $\widetilde{K^t}$ by a disc (and we shall do so) the boundary of which is labelled by $\Phi(\partial (\widetilde{K^t}))$.
	\end{remark}
	
	We prove Propositions 2.1.7 and 2.1.8 by simultaneous induction.
	
	\subsubsection{The induction hypothesis $\mathcal H$}
	{\it \linespread{1.5}\selectfont Let $M$ be a connected, simply connected ${\cal R}$-diagram over $F$, ${\cal R}$ given by (III), with cyclically reduced boundary label $W$ such that $M \in {\cal M}_3(W)$. 
		Let $N$ be a connected, simply connected ${\cal R}$-diagram over $F$ with cyclically reduced boundary label $U$ such that $N \in {\cal M}_3(U)$. If $|N|<|M|$ then the results of Theorem A, B, C and
		Propositions 2.1.7 and 2.1.8  hold for N.\par\bigskip	
	}
	
	In the rest of the work we show that under assumption ${\cal H}$ theorems A, B, C and Propositions 2.1.7 and 2.1.8 hold true.
	
	\vspace{10pt}
	
	\begin{remark}
		\addcontentsline{toc}{subsubsection}{\numberline{\thesubsubsection}
			Remark (replacements do not cause reduction)}
		\ \\ In this section we are going to show that if  $M\in \mathcal{M}_3(W)$ 
		then $M$ and  $\mathbb{M}_t$ are adequate. $t \in T(M)$.
		In the course of the proof we replace a subdiagram P with cyclically reduced boundary label $U$ 
		by a diagram $P'\in \mathcal{M}_3(U)$ and we want to show that the new diagram 
		$M':=(M\setminus P)\cup P'$ satisfies $M'\in \mathcal{M}_3(U)$. 
		One of the problems with such replacements in general, is that while $M$ may be reduced, $M'$ need not be.
		In the replacement we carry out here, P is surrounded by regions in 
		$Reg_{4^+}(M)$, hence since $M\in \mathcal{M}_2(W)$, $M$ is $4^+$-reduced, hence the replacement does not cause reduction. 
		Thus, in order to show that $M'\in\mathcal{M}_3(W)$ we have to show
		\begin{itemize}
			\item [(a)] for every band $B'$ in $M'$,\,$\partial B'$ is a simple closed curve, and $\Phi(\partial B')$ is cyclically reduced.
			\item [(b)] $M'\in\mathcal{M}_2(W)$.
		\end{itemize}
		Now (a) follows since P is surrounded by regions from $Reg_{4^+}(M)$ and therefore the replacement does not effect $\partial B'$.
		Hence we shall concentrate on showing (b).
	\end{remark}
	
	\tocexclude{\subsection{Preparatory Results}}
	\addcontentsline{toc}{subsection}{\numberline{\thesubsection}  
		Preparatory Results for the proof of \\ Propositions 2.1.7 and 2.1.8}
	\begin{lemma}
		\addcontentsline{toc}{subsubsection}{\numberline{\thesubsubsection}  
			Lemma (replacement of $Q$ by $Q'$ leaves $M$ in $\mathcal{M}_3(W)$)}
		Let $M\in \mathcal{M}_3(W)$ and let Q be a connected, $4^+$-reduced simply connected subdiagram with
		$U:=\Phi ( \partial Q)$ cyclically reduced. Denote by $|Inn(Q)|_{4^+}$ 
		the set of all the Equivalence classes of Q which are not fragments of Equivalence classes of M and denote by $|Q|_{4^+}$ the number of Equivalence classes in Q.
		\begin{itemize}
			\item [(a)] if $|Inn(Q)|_{4^+}=|Q|_{4^+}$ then $Q\in \mathcal{M}_3(U)$.
			\item [(b)] if Q is replaced by $Q'\in \mathcal{M}_3(U)$ then 
			$M'\in \mathcal{M}_3(W),\text{\ where\ } M'=(M\setminus Q)\cup Q'$	
		\end{itemize} 
	\end{lemma}
	
	\noindent{\bf Proof of Lemma 2.2.1} 
	\begin{itemize}
		\item [a)]  Replace Q with $Q' \in \mathcal{M}_3(W)$ and let $M'=(M\setminus Q)\cup Q'$. We have
		\begin{itemize}
			\item [1)] $|M|_{4^+}=|M\setminus Q|_{4^+}+|Inn(Q)|_{4^+}$ because all the regions $\Delta$ in $Reg_{4^+}(Q)$ outside $Inn(Q)$ are fragments of regions from $Reg_{4^+}(M\setminus Q)$, hence are already counted in $|M\setminus Q|_{4^+}$.\\ For the same reason: 
			\item [2)] $|M'|_{4^+}=|M\setminus Q|_{4}+|Inn(Q')|_{4^+}$.  Also, 
			\item [3)] $|M'|_{4^+}\geq |M|_{4^+}$ since $M\in \mathcal{M}_3(W)$ and
			\item [4)] $|Q'|_{4^+}\leq |Q|_{4^+}$\\ \ \\
			from (1) - (4) we get \\
			$0\leq |M'|_{4^+}-|M|_{4^+}=|Inn (Q')|_4-|Inn (Q)|_4$, i.e. 
			\item [5)] $|Inn(Q)|_{4^+}\leq |Inn(Q')|_{4^+}$.\\
			Hence we have, due to the assumptions $|Inn Q|_{4^+}=|Q|_{4^+}$\\
			
			$|Q|_{4^+}=|Inn(Q)|_{4^+}\leq |Inn(Q')|_{4^+}\underset{5)}\leq |Q'|_{4^+}$, i.e. 
			$|Q|_{4^+}\leq |Q'|_{4^+}$. Hence, by (4)
			\item [6)] $|Q|_{4^+}=|Q'|_{4^+}$.\\ But $Q'\in\mathcal{M}_3(U)$, by assumption on $Q'$. Hence $Q\in\mathcal{M}_3(U)$.
		\end{itemize}
		\item [b)] By 4) and 6) $|M'|_{4^+}=|M\setminus Q|_{4^+}+|Inn(Q')|_{4^+}=|M\setminus Q|_{4^+}+|Inn(Q)|_{4^+}=|M|_{4^+}$ 
		, i.e. $|M'|_{4^+}=|M|_{4^+}$\\
		Hence $M'\in \mathcal{M}_3(W)$.
	\end{itemize}
	\hfill $\Box$

	\begin{lemma}
		\addcontentsline{toc}{subsubsection}{\numberline{\thesubsubsection}  
			Lemma (endpoints of a path $\mu$ in $\mathbb M_t$ labelled by $M_t$ on $\Delta$ are adjacent)}
		\ \\ Let $M\in \mathcal{M}_3(W)$ and consider $\mathbb{M}_t,\, t\in T(M)$. Assume $\mathbb{M}_t$ is adequate and has connected interior. Let $\Delta\in Reg_{4^+,t}(\mathbb M_t)$ 
		and let $\mu$ be a simple path in $\mathbb{M}_t$ with endpoints $w_1$ and $w_2$ on $\partial\Delta$. Assume that $\Phi(\mu)$ is reduced. 
		Let $\omega_2$ be the boundary path of $\Delta$ with endpoints $w_1$ and $w_2$ in the closed bounded component of $\mathbb E^2\setminus (\mu\cup\Delta)$. 
		If $Supp(\Phi(\mu))\subseteq H_t$ then either $\mu\subseteq\partial\Delta$ or $\Phi(\omega_2)\in H_t$.
	\end{lemma} 
	
	\noindent{\bf Proof of the Lemma}\\
	Let $\mathcal{N}$ be the collection of all the diagrams $M_1\in\mathcal{M}_3(W)$ which contain a pair $(\Delta,\mu)$ as in the Lemma, 
	such that $\mu\not\subseteq\partial\Delta$. 
	If $\mathcal{N}=\emptyset$ then $\mu\subseteq\partial\Delta$.
	Hence we have to show that if $\mathcal{N}\neq\emptyset$ then
	$\Phi(\omega_2)\in H_t$.
	Assume $\mathcal{N}\neq\emptyset$.
	Let P be the bounded connected component of $\mathbb{E}^2\setminus (\Delta\cup\mu)$. Denote $P=P(M,\Delta,\mu)$ and choose $M_0\in\mathcal{N},
	\Delta\in Reg_{4^+,t}(\mathbb{M}_t)$ and $\mu$ in $\mathbb{M}_t$ such that $|P(M_0,\Delta,\mu)|$ is minimal, 
	over all $M\in\mathcal{N},\Delta\in Reg_{4^+,t}(\mathbb{M}_t)$ and $\mu$ in $\mathbb{M}_t$. 
	
	Let $U = \Phi(\partial P)$. By carrying out diamond moves at the endpoints of $\Phi(\mu)$, if needed, we may assume that $\mu$ is cyclically reduced. We claim that $P\in\mathcal{M}_3(U)$. Replace P by $P' \in \mathcal{M}_3(U)$  and let $M'=(M\setminus P)\cup P'$. Since $Int(\mathbb{M}_t)$ is connected, 
	every edge in $\mu$ is either on the boundary of $\mathbb M_t$ or is common to neighbouring regions in $\mathbb{M}_t$. 
	Hence every edge of $\mu$ either separates two adjacent regions in $\mathbb M_t$, or is on the boundary.
	In particular, $\mu$ separates regions of $\mathbb M_t$ which are in P and those which are outside P, hence P cannot contain fragments of regions from $M\setminus P$.
	Hence for every $\Delta \in Reg(\mathbb{M}_t)$ either $\Delta \in Reg(P)$  or $\Delta \cap P = \emptyset$. Therefore
	\begin{itemize}
		\item [(1)] $|M|_{4^+}=|M\setminus P|_{4^+}+|P|_{4^+}$
		\item [(2)] $|M'|_{4^+}=|M\setminus P|_{4^+}+|P'|_{4^+}$\\
		By the assumptions on $M'$ and $P'$ we have:
		\item [(3)] $|M'|_{4^+}\geq |M|_{4^+}$
		\item [(4)] $|P'|_{4^+}\leq |P|_{4^+}$
	\end{itemize}
	It follows from (1)-(4) that $|P|_{4^+}= |P'|_{4^+}$ and $|M|_{4^+}= |M'|_{4^+}$. Consequently $P\in \mathcal{M}_3(U)$, since $P'\in \mathcal{M}_3(U)$. 
	Hence theorem C applies to P with $\omega_1=\mu$ and $\omega_2=\partial P \cap \partial\Delta$. (as defined above). 
	Hence there is a $P''\in \mathcal{M}_3(U)$ with the regions $K_1$ and $K_2$ mentioned in Theorem C. 
	If $K_i\in N(\omega_2)$ then either $n(K_i)=2$ in which case $n(\Delta)=2$ or $\Delta$ is not a "full" Equivalence class in M, $i=1,2$. 
	By assumption both cases are false.
	Let $\omega_1=\mu$. 
	Then $K_i\in N(\omega_1)$. Since $||\partial K \cap \omega_1|| \geq 2$ by Theorem C, $Supp(K)\subseteq Supp(\omega_1) \subseteq H_t$.
	Consequently if $\mu$ decomposes by $\mu=\mu_1 \mu_2 \mu_3,\quad \mu_2=\omega_1\cap \partial K$ and $\mu'_2$ is the complement of $\mu_2$ on $\partial K$,  then $Supp (\Phi(\mu'_2))\subseteq H_t$. 
	If $|P(M,\mu,\Delta)|\geq 2$ then $|P(M,\mu',\Delta)|<|P(M,\mu,\Delta)|$ where $\mu'=\mu_1 \mu_2'\mu_3$ since 
	$|P(M,\mu',\Delta)|=|P(M_t,\mu,\Delta)|-|\{K\}|$. 
	This violates the minimality of $|P|$.
	Hence $|P(M,\mu,\Delta)|=1$ i.e.  $P(M,\mu,\Delta)=\{K\}$ and $Supp(K)\subseteq H_t$
	Consequently, $Supp(\Phi(\partial\Delta\cap\partial K))\subseteq H_t$ and hence $\Phi(\omega_2)=\Phi(\partial\Delta\cap\partial K)=a^m, a\in X_t$ for some $m\neq 0$. 
	Hence $\Phi(\omega_2)\in H_t$, as required.
	\hfill $\Box$
	
	\vspace{10pt}

	\subsection{Proof of Proposition 2.1.7}
	\noindent{\bf Proof } 
	\ We have to show that if $\Delta\in Reg (\mathbb{M}_t)$ then \begin {itemize}
	\item [a)] $Int(\Delta)$ is homeomorphic to the open unit disc, and $\partial\Delta$ is a simple closed curve.
	\item [b)] $\Phi(\partial \Delta)$ is cyclically reduced.
\end{itemize}
\begin{sssection}[(a)]
	Consider first the case $D\in Reg_{4^+,t}(M)$ and let 
	$\Delta=\kappa_0(D)$
\end{sssection}
\noindent Assume by way of contradiction that $\partial \Delta$ is not a simple closed curve.
Consider the set $\mathcal{N}(\Delta)$ of all the diagrams $N\in \mathcal{M}_3(W)$ 
which contain $\Delta$ and in which $\partial \Delta$ is not simple closed. Then for 
every $N$ in $\mathcal{N}(\Delta),\  \mathbb{E}^2\setminus \Delta= Q_1(N)\cup Q_2(N)\ldots \cup Q_{k(N)}(N)\cup Q_\infty$, where $Q_i(N),i=1,\ldots,k(N),\quad k(N)\geq 1$ are the bounded connected components of $\mathbb{E}^2\setminus \Delta$ in $N,\quad Q_i(N)\subseteq N$. 
Let $\rho_\Delta(N)=\sum_{i=1}^{k(N)}|\partial Q_i(N)|$. Choose $N_0$ in $\mathcal N(\Delta)$ with minimal $\rho_\Delta(N)$  among all $N\in\mathcal N(\Delta)$ and denote $Q_i=Q_i(N_0), i = 1,\ldots,k(N_0)$. Let $U_i=\Phi(\partial Q_i)$. We claim that
\begin{itemize}
	\item [(*)] {\it $U_i$ is cyclically reduced} 
\end{itemize}  
Suppose not. 
Then there is a $Q_i$ such that $\partial Q_i$ has at most one double point $w$ 
in which cancellation occurs. 
Carrying out a diamond move \diamondmove{w} 
at $w$ reduces $|U_i|$. 
\noindent Since diamond moves do not alter the number of regions in $N_0$ and since $\partial Q_i \subseteq \partial \Delta$, it follows that the diamond move \diamondmove{w} 
does not alter the number of Equivalence classes, 
hence $N'_0 \in \mathcal{M}_3(W)$, where $N'_0$ is the diagram obtained from $N_0$ by \hbox{\diamondmove{w}.} 
Thus $N'_0 \in \mathcal N(\Delta)$. 
But $\rho_\Delta(N'_0)<\rho_\Delta(N'_\Delta)$ because $\diamondmove{w}$ reduces $|{U_i}|$ and does not alter 
$|U_j|,\ j\neq i$, violating the minimality of $\rho_\Delta(N_\Delta)$. Hence $U_i$ is cyclically reduced, as claimed.

Now, it follows from Lemma 2.2.1 that ${Q_i} \in \mathcal{M}_3(U_i)$. 
Hence Theorem B is applicable to $Q_i$ and therefore $Supp(E)=Supp(\Delta)$, for every 
$E\in Reg({Q_i})$, as $||U_i|| \geq 4$ by Lemma 1.1.1(c). 
But this contradicts the definition of $\Delta$, showing that $\partial\Delta$ is a simple closed curve. \\ \ \\
{\bf (b)} We show that using diamond moves on $\partial\Delta$ we may obtain from $\Delta$ a subdiagram $\Delta'$
of $\mathbb M_t$ with boundary label V, V cyclically reduced, such that V is the cyclically reduced form of $\Phi(\partial\Delta)$, and such that if K is a neighbour of $\Delta$ in $\mathbb M_t$ 
and $\Phi(\partial K)$ is cyclically reduced then it remains cyclically reduced after the diamond moves.


\noindent Suppose $\Delta \in Reg_{4^+,t}({\mathbb M}_t)$ and $\Phi(\partial \Delta)$ is not cyclically reduced.
Since $M\in\mathcal{M}_3(W)$, $\Phi(\partial \Delta) \neq 1$. 
Hence, $\partial \Delta$ contains a subpath 
$\mu_0 u \alpha v \beta w \nu_0$, $\mu_0$, $ \nu_0$, $\alpha$ and $\beta$ non-empty paths, $u$, $v$ and $w$ vertices, such that $\Phi(\alpha) = \Phi(\beta)^{-1}$ and $\Phi(\mu_0)\Phi(\nu_0)$ is reduced as written. 
Consider first the effect of a diamond move
\diamondmove{v} at $v$
on $\Delta$. See Figure 23.

\begin{minipage}{1.7in}
	\begin{center}
		
		\begin{tikzpicture}[scale=.1]
			\draw[rounded corners] (-3,1)--(3,1)--(4,2);
			
			\draw (4,2) arc(140:35:20);
			
			\draw (36,1)--(41,1);
			\draw (40,2)--(41,1)--(40,0);
			\draw (0,2)--(1,1)--(0,0);
			
			\draw[fill] (4,2) circle(.8); 
			\draw[fill] (36,1) circle(.8);
			\draw[fill] (20,9) circle(.8);
			
			\draw (9,7)--(10,7)--(9,5);
			\draw (29,8)--(30,6)--(29,5);
			
			\node at (0,4) {$\mu_0$};
			\node[above] at (4,2) {$u$};
			\node[above] at (9,8) {$\alpha$};
			\node[above] at (20,9) {$v$};
			\node at (30,10) {$\beta$};
			\node[above] at (36,1) {$w$};
			\node at (41,3) {$\nu_0$};

			\node at (20,1) {$\Delta$};
			
		\end{tikzpicture}
	\end{center}
	
\end{minipage}\hspace{.3cm}\begin{minipage}{3in}
	\begin{tikzpicture}[scale=.3]
		
		\draw (0,0)--(3,0)--(6,5)--(6,-5)--(3,0);
		\draw (6,0)--(9,0);
		
		\draw[fill] (3,0) circle(.25);
		\draw[fill] (6,0) circle(.25);
		\draw[fill] (6,5) circle(.25);
		\draw[fill] (6,-5) circle(.25);
		
		\draw (3.5,-1.5)--(4.1,-2)--(4.2,-1.2);
		\draw (3.5,1.5)--(4.1,2)--(4.3,1.5);
		\draw (5.5,-1.5)--(6,-2)--(6.5,-1.5);
		\draw (5.5,1.5)--(6,2)--(6.5,1.5);
		\draw (8.5,.5)--(9,0)--(8.5,-.5);
		
		\node at (1,.5) {$\mu_0$};
		\node at (2.3,1) {$u$};
		\node at (3,3) {$\alpha_1$};
		\node[above] at (6,5) {$v_1$};
		\node at (7,3) {$\beta_1$};
		\node at (6.5,1) {$w$};
		\node at (9,1) {$\nu_0$};
		\node at (7,-3) {$\beta_2$};
		\node[below] at (6,-5) {$v_2$};
		\node at (2,-3) {$\alpha_2$};
		\node at (6,-8) {$\Delta$};
		
		\node at (11,-2) {\huge $\rightsquigarrow$};
		\node at (-2,-2) {\huge $\rightsquigarrow$};

		\draw (14,0)--(22,0);
		\draw (18,-4)--(18,4);
		\draw[fill] (18,0) circle(.2);
		\draw[fill] (18,4) circle(.2);
		\draw[fill] (18,-4) circle(.2);
		
		\draw (15.5,.5)--(16,0)--(15.5,-.5);
		\draw (19.5,.5)--(20,0)--(19.5,-.5);
		
		\draw (17.5,1.5)--(18,2)--(18.5,1.5);
		\draw (17.5,-1.5)--(18,-2)--(18.5,-1.5);
		
		\node at (16,1) {$\mu_0$};
		\node[above] at (18,4) {$v_1$};
		\node at (20,1) {$\nu_0$};
		\node at (20,-1) {$u=w$};
		\node[below] at (18,-4) {$v_2$};
		\node at (18,-7) {$\Delta$};

	\end{tikzpicture}
\end{minipage}\hspace{10pt}\begin{minipage}{.5cm}
	\begin{tikzpicture}[scale=.25]
		\draw (1,0)--(2,1)--(1,2)--(0,1)--(1,0);
		\node at (1,1) {$v$};
	\end{tikzpicture}
\end{minipage}

\vspace{10pt}

\begin{center}
	Figure 23
\end{center}

If $\Delta''$ is the diagram obtained from $\Delta$ then $Int(\Delta'')$ does not contain $\alpha_1$ ($=\beta_1$) and $\mu_0$ is followed by $\nu_0$, hence no cancellation occurs at $u$. If the boundary of $\overline{Int(\Delta'')}$ is not cyclically reduced then carry out diamond moves at the appropriate vertices, as in 
\begin{tikzpicture}[scale=.25]
	\draw (1,0)--(2,1)--(1,2)--(0,1)--(1,0);
	\node at (1,1) {$v$};
\end{tikzpicture}
above, until we get a cyclically reduced boundary label. Notice, that in this process there can be no loops, since
at each step the boundary  is shortened by $|\alpha| \neq 0$. Hence, the procedure ends after at most $\frac{1}{2}|\partial \Delta|$ steps. By construction, clearly $\Phi (\partial \Delta')$ is the cyclic reduction of $\Phi (\partial \Delta)$.

Finally, consider the effect of a diamond move
\begin{tikzpicture}[scale=.25]
	\draw (1,0)--(2,1)--(1,2)--(0,1)--(1,0);
	\node at (1,1) {$v$};
\end{tikzpicture}
carried out on $\partial \Delta$,
on regions $\Delta_j$, distinct from $\Delta$. Clearly, if $\partial \Delta_j \cap u \alpha v \beta w = \emptyset$ then 
\begin{tikzpicture}[scale=.25]
	\draw (1,0)--(2,1)--(1,2)--(0,1)--(1,0);
	\node at (1,1) {$v$};
\end{tikzpicture}
does not change $\overline{\Delta_j}$. Hence, we consider $\Delta_j$ which are neighbours of $\Delta$ and have nontrivial intersection with $u\alpha v \beta w$. 

\vspace{40pt}

\begin{minipage}{2in}
	
	\begin{center}
		
		\begin{tikzpicture}[scale=.15]
			
			\draw (4,2) arc(140:35:20);
			
			\draw (5,3)--(4,12)--(10,13)--(11,7.2);
			\draw (10,13)--(16,14)--(17,9);
			\draw (16,14)--(18,18)--(20,14);
			\draw (17,9)--(20,14)--(24,14)--(23,8.7);
			\draw (24,14)--(28,14)--(28.5,7);
			
			\draw[fill] (8,5.5) circle(0.2);
			\draw[fill] (17,9) circle(0.2);
			\draw[fill] (28.5,7) circle(0.2);
			
			\draw (10,6.2)--(11,7.2)--(10,7.5);
			\draw (20,8.7)--(19,9.2)--(20,9.7);
			
			\node[below] at (8,5.5) {$u$};
			\node at (11,5.8) {$\alpha$};
			\node[below] at (17,9) {$v$};
			\node at (20,6.6) {$\beta$};
			\node[below] at (28.5,7) {$w$};
			\node at (17,4) {$\Delta$};
			\node at (8,9) {$\Delta_1$};
			\node at (13,10) {$\Delta_2$};
			\node at (18,14) {$\Delta_3$};
			\node at (21.5,11) {$\Delta_4$};
			\node at (26,10) {$\Delta_5$};

		\end{tikzpicture}
		
		\vspace{10pt}
		
		(a) 
		
	\end{center}
	
\end{minipage}\hspace{.5cm}\begin{minipage}{2in}
	
	\begin{center}
		\begin{tikzpicture}[scale=.15]
			
			\draw (4,2) arc(140:35:20);
			
			\draw (10,7)--(9,12)--(13,13)--(15,8.6)--(17,13)--(20,14)--(19,9);
			\draw (20,14)--(28,13)--(27.5,7.3);
			\draw (28,13)--(33,12)--(32,4.7);
			
			\draw (13,13)--(15,15)--(17,13);
			
			\draw[fill] (23.2,8.7) circle(.15);
			\draw[fill] (30,6.1) circle(.15);
			
			\draw (21,9.4)--(22,9)--(21,8.7);
			\draw (25.7,8.5)--(25,8.4)--(25.7,7.8);
			
			\node at (20,3) {$\Delta$};
			\node at (15,7) {$\alpha$};
			\node[below] at (23.2,8.7) {$v$};
			\node at (25,6) {$\beta$};
			\node[below] at (30,6.1) {$w$};
			\node at (8,9) {$u$};
			\node at (12,10) {$ \Delta_1$};
			\node at (15,12.4) {$\Delta_2$};
			\node at (18,11) {$\Delta_3$};
			\node at (25,11) {$\Delta_4$};
			\node at (30,9) {$\Delta_5$};

		\end{tikzpicture}
		
		\vspace{10pt}
		
		(b)
	\end{center}
\end{minipage}

\begin{center}
	Figure 24
\end{center}

We consider two cases according to $v$ being an endpoint or an inner vertex of the common edge of $\Delta$ and $\Delta_i$. See Figure 24 (a) and (b), respectively. In both cases
\begin{tikzpicture}[scale=.25]
	\draw (1,0)--(2,1)--(1,2)--(0,1)--(1,0);
	\node at (1,1) {$v$};
\end{tikzpicture}
does not change the boundary labels of the $\Delta_j$. Also, in case of figure 24(b) $\Phi(\omega_4)$ is not cyclically reduced, where $\omega_4 = \partial \Delta_4$, 
hence if by assumption $\Delta_4$ has cyclically reduced boundary label then only case (a) may occur and in this case 
\begin{tikzpicture}[scale=.25]
	\draw (1,0)--(2,1)--(1,2)--(0,1)--(1,0);
	\node at (1,1) {$v$};
\end{tikzpicture}
leaves cyclically reduced boundary labels cyclically reduced. See Figures 25(a),  25(b) and 25(c) for the result of the diamond moves.

\begin{minipage}{1.5in}
	\begin{center}
		
		\begin{tikzpicture}[scale=.7]
			\draw (0,0) ellipse(2 and 1);
			
			\draw (0,0) --(0,3);
			\draw (-1,.85)--(-.9,2)--(-1.1,3.5)--(-.5,4)--(0,3)--(1,4)--(1.1,2.5)--(1,.89);
			\draw (-.5,4)--(1,4);
			\draw (-.9,2)--(0,2)--(1.1,2.5);
			
			\draw[fill] (0,0) circle(.1);
			
			\draw (-.2,.5)--(0,.2)--(.2,.5);
			\draw (-.2,2.5)--(0,2.8)--(.2,2.5);
			
			\node at (0,-.5) {$\Delta$};
			\node at (-.5,1.5) {$\Delta_1$};
			\node at (-.5,3) {$\Delta_2$};
			\node at (0.1,3.7) {$\Delta_3$};
			\node at (.5,3) {$\Delta_4$};
			\node at (.5,1.5) {$\Delta_5$};
			\node at (1,.3) {$u=v$};
			
		\end{tikzpicture}
		
		\vspace{10pt}
		
		(a)
	\end{center}
\end{minipage} \hspace{.2cm} \begin{minipage}{1.5in}
	\begin{center}
		\begin{tikzpicture}[scale=.3]
			\draw (0,0) ellipse(4 and 2);
			
			\draw (0,0) -- (0,6);
			
			\draw (0,6) circle(2);
			
			\draw (-2,1.7)--(-2,2.5)--(-3,3.25)--(-2,4)--(-1.8,5);
			\draw (-2,2.5)--(-1,3.25)--(-2,4);
			
			\draw (1.9,5.2)--(2.3,3)--(1.8,1.7);
			
			\draw (-1,3.25)--(0,3.25);
			
			\draw[fill] (0,0) circle(.1);
			\draw[fill] (0,6) circle(.1);
			
			\node at (-1,-1) {$\Delta$};
			\node at (-.7,2.5) {$\scriptstyle \Delta_1$};
			\node at (-2,3.3) { $\scriptstyle \Delta_2$};
			\node at (-.7,3.8) {$\scriptstyle \Delta_3$};
			\node at (1,6) {$\scriptstyle \Delta_4$};
			\node at (1,3) {$\scriptstyle \Delta_5$};

		\end{tikzpicture}
		
		\vspace{10pt}
		
		(b)
	\end{center}
\end{minipage}\hspace{.2cm} \begin{minipage}{1.8in}
	
	\begin{center}
		
		\begin{tikzpicture}[scale=.3]
			
			\draw (0,0) circle(4);
			\draw (0,0)--(0,8);
			
			\draw[fill] (0,0) circle(.15);
			\draw[fill] (0,4) circle(.15);
			\draw[fill] (0,8) circle(.15);
			
			\draw (-.3,1.3)--(0,1)--(.3,1.3);
			\draw (-.3,6.7)--(0,7)--(.3,6.7);
			\draw (2.6,-2.7)--(2.6,-3)--(2.9,-3);
			
			\node[above] at (0,8) {$v_1$};
			\node at (-1,7) {$\alpha_1$};
			\node at (1,7) {$\beta_1$};
			\node at (1,5.7) {$\Delta_2$};
			\node at (-.5,3) {$u=w$};
			\node at (-1,-1.5) {$\Delta_1$};
			\node at (2,-1) {$B_1$};
			\node at (3.8,-3) {$w$};

		\end{tikzpicture}
		
		\vspace{10pt}
		
		(c)
	\end{center}

\end{minipage}

\begin{center}
	Figure 25
\end{center}

\hfill $\Box$

\begin{sssection}
	For the case where $\Delta\in Reg_{H_t}(M)$  
	we follow 2.3.1 and in the same way derive the same contradiction to the definition of $[\ ]_{H_t}$. 
	We omit the details. 
\end{sssection}
\begin{sssection}
	Consider the case $\Delta\in Reg_{2,t}(M)$. 
	$\partial \Delta$ is simple closed due to the assumption that $M\in\mathcal{M}_3(W)$.
	As in 2.3.1 we apply diamond moves to get cyclically reduced boundary label. 
\end{sssection}


\begin{sssection}
	Now we show that if $K_1, K_2 \in Reg (\mathbb M_t), K_1$ and $K_2$ are not both bands, then $\partial K_1 \cap \partial K_2 \neq \emptyset$ implies that $\partial K_1 \cap \partial K_2$ is connected.\ \\ 
	We have the following cases to check
	\begin{itemize}
		\item [1)] $K_1, K_2 \in Reg_{4^+,t} (\mathbb M_t)$
		\item [2)] $K_1, K_2 \in Reg_{H_t} (\mathbb M_t)$
		\item [3)] $K_1 \in Reg_{4^+,t} (\mathbb M_t)$ and $K_2 \in Reg_{H_t} (\mathbb M_t)$
		\item [4)] $K_1 \in Reg_{B_t} (\mathbb M_t)$ and $K_2 \in Reg_{4^+,t} (\mathbb M_t)$
		\item [5)] $K_1 \in Reg_{B_t} (\mathbb M_t)$ and $K_2 \in Reg_{H_t} (\mathbb M_t)$
	\end{itemize} 
\end{sssection}

\begin{sssection}[Case 1]
	Write $\Delta_i$ for $K_i,\ i=1,2$. 
	Assume $\partial\Delta_1\cap\partial\Delta_2\neq\emptyset$ and $\partial\Delta_1\cap\partial\Delta_2$ is not connected. Consider the set $\mathcal N(\Delta_1,\Delta_2)$ of all diagrams $N\in \mathcal{M}_3(W)$ which contain $\Delta_1\cup \Delta_2$
	with $\partial\Delta_1\cap\partial\Delta_2\neq\emptyset$ such that $\partial\Delta_1\cap\partial\Delta_2$ is not connected. 
	Then for every $N\in \mathcal N(\Delta_1,\Delta_2)$ 
	we have $\mathbb E^2\setminus ( \partial{\Delta_1}\cup\partial{\Delta_2})=Q_1(N)\cup Q_2(N)\ldots\cup Q_{k(N)}(N)\cup Q_\infty$,
	where $k(N)\geq 1,\quad \Delta_i(N)=Int(\cup[D_i]_N), [D_i]_N$ stands for the Equivalence class of $D_i$ in $N$ and $Q_j(N)$ are the bounded connected components of $\mathbb{E}^2\setminus (\Delta_1(N)\cup \Delta_2(N))$. 
	Notice that $M\in \mathcal N(\Delta_1,\Delta_2)$ and $\Delta_i(M)=\Delta_i,\quad i=1,2$. Let $\rho_{\Delta_1,\Delta_2}(N)=\sum_{i=1}^{k(N)}|\partial Q_i(N)|$ and denote by $\mathcal N_0(\Delta_1,\Delta_2)$ the set of all $N$ in 
	$\mathcal N(\Delta_1,\Delta_2)$ with minimal $\rho_{\Delta_1,\Delta_2}(N)$. Choose $N_0\in \mathcal N_0(\Delta_1,\Delta_2)$ such that $\rho^v(N_0):=\sum_{i=1}^{k(N_0)}|Q_i(N_0)|$ is minimal. 
	Denote $Q_i:=Q_i(N_0)$ and assume that $\partial\Delta_1\cap\partial\Delta_2\cap Q_1$ consists of two  vertices, $w_1$ and $w_2$.
	Let $U_1$ be the label of the boundary path $\mu_1$ of $Q_1$ which starts at $w_1$ and terminates at $w_2$ and let $U_2$ be the label of the boundary path $\mu_2$ of $Q_1$ which starts at $w_2$ and terminates at $w_1$, $\mu_1\neq \mu_2^{-1}$. Let $U=U_1 U_2$ reduced as written. Then $U$ is a boundary label of $Q_1$.
	We claim that $U$ is cyclically reduced. Suppose not. Since $U_1$ and $U_2$  are reduced, there is a cancellation at $w_1$  or at $w_2$. 
	Hence, carrying out a diamond move \diamondMove{w_1} or \diamondMove{w_2}, reduces $|U|$. 
	Since diamond moves do not alter the number of regions in $N_0$ and since $\partial Q_1 \subseteq \partial \Delta_1 \cup \partial \Delta_2$, it follows that \diamondMove{w_1}  or \diamondMove{w_2} does not alter $|N_0|_{4^+}$. 
	Let $N’_0$  be the result of applying \diamondMove{w_1} or \diamondMove{w_2}  on $N_0$.  
	Then $N'_0 \in \mathcal{N}( \Delta_1, \Delta_2)$. But $\rho_{\Delta_1, \Delta_2}(N’_0)<\rho_{\Delta_1, \Delta_2}(N_0)$, violating the minimality of $\rho _{\Delta_1, \Delta_2}(N_0)$,  hence $U$ is cyclically reduced.\\
	It follows from Lemma 2.2.1 that $Q_1\in \mathcal M_3(U)$.  Hence Theorem C  applies to $Q_1$  with $\mu_1$ and $\mu_2$  in place of $\omega_1$ and $\omega_2$, respectively. 
	Consequently, there is a $Q''_1 \in \mathcal M_3(U)$  with $N''_0:=\left(N_0\setminus Q_1\right)\cup Q_1''\in \mathcal M_3(W)$ due to Lemma 2.2.1 such that $Q''_1$  has a boundary region or
	boundary modified region K such that $\partial K\cap \partial N''_0$ is connected and $\partial K\cap \partial N''_0\subseteq \mu_i, i=1$ or $i=2$ 
	and $||\partial K\cap \partial N''_0||\geq 2$. Suppose $\partial K\cap \partial N''_0\subseteq \mu_1$ and consider $\Delta_1 \cup K$.
	Let $D_1$ be a region in $N''_0$.
	Denote  the Equivalence class of $D_1$ in $N''_0$ by $[D_1]_{N''_0}$. Then $[D_1]_{N''_0}\supseteq \Delta_1\cup K$. 
	Hence if we denote by $[\Delta_1]_{N''_0}$ the Equivalence class of $\Delta_1$ in $N''_0$ then $[\Delta_1]_{N''_0}=[D_1]_{N''_0}$. 
	Consequently $[\Delta_1]_{N''_0}\supseteq \Delta_1\cup K$.
	Consequently $\left|Q_1^{\prime \prime}\left(N_0''\right)\right| \leqslant\left|Q_1\left(N_0\right)\right|-|K|<\left|Q_1\left(N_0\right)\right|$. 
	Hence $|Q_1''(N_0)|\leq|Q_1(N_0)$.
	Now,  for $i\geq 2$ we have $\left|Q_i^{\prime \prime}\left(N_0^{\prime \prime}\right)\right|=\left|Q_i\left(N_0\right)\right|$,
	because for $i\geq 2$, $Q_i \subseteq N \backslash Q_1$, hence the replacement of $Q_1$ by $Q''_1$  has no effect on $|Q_i|$. 
	Hence,
	$\rho^v\left(\Delta_1\left(N_0^{\prime \prime}\right), \Delta_2\left(N_0^{\prime \prime}\right)\right)
	=
	\sum_{i=1}^{k\left(N_0\right)}\left|Q_i^{\prime \prime}\left(N_0^{\prime \prime}\right)\right|=\left|Q_1^{\prime \prime}\left(N_0^{\prime \prime}\right)\right|+\sum_{i=2}^{k\left(N_0\right)} \left|Q_i\left(N_0\right)\right|<|Q_1\left(N_0\right)| +\sum_{i=2}^{k\left(N_0\right)}\left|Q_i\left(N_0\right)\right|=\sum_{l=1}^{k\left(N_0\right)}\left|Q_i\left(N_0\right)\right|=\rho^v\left(\Delta_1\left(N_0\right), \Delta_2\left(N_0\right)\right)$,
	Contradicting the minimality assumption on $\rho^v\left(\Delta_1, \Delta_2\right)$. 
	Hence $\partial\Delta_1\cap \partial\Delta_2$  is connected. \\ 
	For the remaining cases we consider $\mathbb{E}^2\setminus(K_1\cup K_2)$. We have, as before $\mathbb{E}^2\setminus(K_1\cup K_2)=Q_1\cup Q_2\cdots\cup Q_k\cup Q_\infty,k\geq 1$, 
	where $Q_i$ are the bounded connected components, $1\leq i\leq k$, and $Q_\infty$ the unbounded connected component.
\end{sssection}
\begin{sssection}[Case 2]
	Due to Lemma 2.2.1 the proof of case 1 in 2.3.5 goes through, for case 2 in  2.3.4
\end{sssection}
\begin{sssection}[Case 3]
	Let $\mu_1=\partial K_1\cap \partial Q_1$ and let $\mu_2=\partial K_2\cap \partial Q_1$. 
	As in the previous cases, we may assume that $U:=\Phi(\partial Q_1)$ is cyclically reduced. 
	Also $|Inn(Q_1)|_{4^+}=|Q_1|_{4^+}$. 
	Hence by Lemma 2.2.1 $Q_1\in \mathcal{M}_3(U)$.\\ 
	Consequently, Theorem C applies to $Q_1$. 
	Hence we may replace $Q_1$ with $Q'_1\in \mathcal{M}_3(U)$ such that $Q'_1$ contains a boundary region $ K$ with $\partial K\cap\partial Q_1$ connected. 
	By Theorem C either $\partial K\cap\partial Q_1\subseteq\mu_1$ or 
	$\partial K\cap\partial Q_1\subseteq\mu_2$.
	Since $\Phi(\mu_2)\subseteq H_t$, by Lemma 2.2.2 $\Phi(\mu)\subseteq H_t$.
	But $\|\partial K\cap\mu_1\|\geq n\geq 2$. Hence
	$\partial K\cap\partial Q_1\subseteq \mu_2$.
	Now, $M':=(M\setminus Q_1)\cup Q_1'$ is in $\mathcal{M}_3(W)$, due to Lemma 2.2.1(b). 
	Hence we may repeat the argument of 2.3.5 with $\rho^v$, leading to a contradiction which shows that $\partial K_1\cap \partial K_2$ is connected.
\end{sssection}
\begin{sssection}[Case 4]
	Recall that $C_\Gamma(t)=\{x\in X, n(x,t)=2\}$.
	Assume $\partial K_1\cap \partial K_2\neq\emptyset$ and $\partial K_1\cap \partial K_2$ is not connected. 
	Then $\mathbb{M}_t\setminus( K_1\cup K_2)$ contains a simply connected component Q with $\partial Q=\mu_1\mu_2$, where $\mu_1=\partial Q\cap K_1$ and $\mu_2=\partial Q\cap K_2$. 
	Now $\Phi(\mu_1)=a^\alpha$ due to Lemma 2.2.2, where $Supp(K_1)=\{a,t\}$,
	while $\Phi(\mu_2)\in C_\Gamma(t)$. 
	Consequently $\|\Phi(\partial Q)\|_a \leqslant\left\|\Phi\left(\mu_1\right)\right\|_a+\left\|\Phi\left(\mu_2\right)\right\|_a=1+0=1$, violating Lemma 1.1.1(c). 
	Hence $\partial K_1\cap\partial K_2$ is connected.
\end{sssection}
\begin{sssection}[Case 5]
	Assume $\partial K_1\cap\partial K_2\neq \emptyset$ and $\partial K_1\cap\partial K_2$ is not connected. 
	Let $Q,\mu_1$ and $\mu_2$ be as in 2.3.8 and among all such Q choose one with minimal $|Q|$. As in previous cases we may assume that $U:=\Phi(\partial Q)$ is cyclically reduced. 
	Notice that $Q\subseteq\mathbb{M}_t$. 
	Hence it follows from Lemma 2.2.1 that $Q\in \mathcal{M}_3(U)$. 
	Consequently Theorem C applies to Q. 
	It follows from Theorem C that Q contains a region K with $\partial K\cap\mu_1$ connected, 
	$K\in Reg(\mathbb M)\cup Reg_2(M)$ and $||\partial K\cap \mu_i||\geq n(K)$ for $i=1$ or for $i=2$. 
	If $i=1$ then $Q_1:=Q\setminus K$, violates the minimality of $|Q|$ and if $i=2$ then remove $K$ from $Q$ via $K_1$ ($K_1$ is a t-band) and $Q_2:=Q\setminus K_2$ violates the minimality of $|Q|$. 
	Hence $\partial K_1\cap \partial K_2$ is connected.
\end{sssection}
\hfill$\Box$

\subsection{Proof of Proposition 2.1.8}
We have to prove:
\begin{itemize}
	\item [(a)] $\partial D$ is a simple closed curve, for every $D\in Reg(M)$
	\item [(b)] if $\partial D_1\cap \partial D_2\neq \emptyset$ then $\partial D_1\cap \partial D_2$ is connected.
\end{itemize}
\ \\
\begin{itemize}
	\item [(a)] The proof of part (a) follows the proof of part (a) of Proposition 2.1.7, with two exceptions:
	\begin{itemize}
		\item [i)] Write everywhere $D$ in place of $\Delta$.
		\item [ii)] Replace the last sentence by the following:\\
		"In particular, $Q_1$ is a diagram of a 2-generated Artin group. 
		Hence by Lemma 1.5.1(c) $||U||\geq 2 n(D)$. But $U$ is a proper subword of $\Phi(\partial D)$ and $|\Phi(\partial D)|=\| \Phi(\partial D)\|=2 n(D)$. A contradiction, showing that $\partial D$ is a simple closed curve."
	\end{itemize}
	\item[(b)] {\bf 1)\quad} Suppose not. 
	Then $M\setminus({D_1}\cup {D_2})=Q_1\cup Q_2\cup\cdots\cup Q_k\cup Q_\infty, k\geq 1$, where $Q_i$ are the bounded connected components and $Q_\infty$ is the unbounded component. 
	As in part 2.3.1, we may choose $Q_1$ such that $\partial D_1\cap \partial D_2\cap \partial Q_1$  consists of two vertices, $w_1$ and $w_2$ and $U:=\Phi(\partial Q_1)$ is cyclically reduced.
	Let $\mu_1=\partial Q_1\cap\partial D_1$ and let  $\mu_2=\partial Q_1\cap\partial D_2$. Let $U_1=\Phi(\mu_1)$ and $U_2=\Phi(\mu_2)$. 
	Then $U=U_1U_2$ is reduced as written. 
	Let $\Delta_1=\kappa_M(D_1)$ and let $\Delta_2=\kappa_M(D_2)$,
	where $\kappa_M(D_i)=\kappa_0(D_i)$ in $M$.
	(See 2.1.4 for definition). Then $\partial \Delta_1\cap \partial\Delta_2\neq\emptyset$  since $w_1$ and $w_2$ are in the intersection. Hence $\partial \Delta_1\cap\partial\Delta_2$ is connected by Proposition 2.1.7. 
	This implies that $Q=\Delta'_1\cup\Delta'_2$ where $\Delta'_i=\Delta_i\cap Q_1, i=1,2$.
	Assume 
	$||\Phi(\partial\Delta'_1\cap\partial\Delta'_2)||\geq 2$.
	Then $Supp(D_1)=Supp(D_2)=\{a,b\}$. Hence $P:=D_1\cup D_2\cup Q_1$ is a connected diagram (not necessarily simply connected) over F(a,b). 
	Consequently, by Lemma 1.5.1(a) it satisfies the condition    C(4)\&T(4).  (See Remark after Lemma 1.5.1.)
	
	\noindent {\bf 2)\quad} Consider the subdiagram $Q_1$ of P. 
	Due to Proposition 1.3.4, if $|Q_1|\geq 2$ then $Q_1$ has a boundary region E with $\theta:=\partial E\cap\partial D_i$, connected, $i=1$ or $i=2$, and either $\theta=\partial E\cap\partial Q_1$ 
	or $\theta$ contains $w_1$ or $w_2$ such that $i_{Q_1}(E)\leq 2$ or $i_{Q_1}(E)=1$, respectively. But then $\partial E=\theta\cup\eta_i, i=1$ or 
	$i=2$, where $\eta_i$ is the product of at most 2 pieces 
	($\eta_i$ are the complements of $\theta$ on $\partial E$). Since $\theta=\partial E\cap \partial D_i$, $\theta_i$ is a piece, hence E is an inner region of P with at most 3 neighbours, contradicting the condition C(4). 
	Hence $\|\Phi(\partial\Delta'_1\cap\partial\Delta'_2)\|=1$ (if $|Q|=1$ then $\Phi(\partial Q)$ is the product of two pieces, violating the C(4) condition).
	If $Supp(D_1)=Supp(D_2)$ then the arguments for the case $\|\Phi(\partial\Delta'_1\cap\partial\Delta'_2)\|\geq 2$ apply. 
	Hence assume $Supp(D_1)\neq Supp(D_2)$.
	Then we can repeat the above argument with $\Delta'_1$ in place of $Q_1$. 
	We have $\Delta'_1\cup D_1$ and $\Delta'_2\cup D_2$  are diagrams over two generated Artin groups. Since $\|\Phi(\partial\Delta'_1\cap\partial\Delta'_2)\|= 1$ 
	hence $\Delta'_i$ has a boundary region $E_1$ with $\partial E_1\cap \partial Q=\partial E_1\cap\partial D_1$ as above, in 2). 
	$i=1$ or $i=2$ leading to the same contradiction.
	Hence $\partial D_1\cap \partial D_2$ is connected.	
\end{itemize}
\ \hfill $\Box$

\vspace{10pt}

\section {Relative extended presentations, I-moves and banded diagrams}
\subsection{Identities among relations}
We recall identities among relations from \cite[p. 157]{4}. Let $C$ be a 2-complex. Let $c_0$ be a basepoint for $C$ (some designated vertex of $C$). A \emph{sequence $\sigma$ over $C$ at $v_0$} is $\sigma = (q_1, \ldots , q_n)$, 
where $q_i$ is a path in $C$ consisting of a concatenation $\gamma_i \partial \Delta_i \gamma_i^{-1}$, $\Delta_i$ being a region (face) of $C$ and $\gamma_i$ a path from $v_0$ to $v(\Delta_i)$, $1 \leq i \leq n$, 
where $v(\Delta_i)$ is the basepoint of $\Delta_i$. We call $\sigma$ an \emph{identity sequence} if the concatenation $q_1\cdots q_n$ is freely equivalent to the trivial path at $v_0$.

Given a presentation ${\cal P}= \langle X | {\cal R} \rangle$, following \cite{14} we construct a boquet of tailed circles $C$ with basepoint $v_0$. Let $\sigma = (q_1, \ldots , q_n)$ be a sequence over $C$ at $v_0$, $q_i=\gamma_i\Delta_i\gamma_i^{-1}$, $\gamma_i$ tails. 
Label $\partial \Delta_i$ by $R_i \in {\cal R}$ and label $\gamma_i$ by a word in $F(X)$, via a labeling function $\Phi_0: C \rightarrow F$. After carrying out all free cancellations in $\Phi_0(\sigma) = \Phi_0(q_1) \cdots \Phi_0(q_n)$ we get an $\cal R$-diagram $M$ over $F$ the boundary label of which is the cyclically reduced word for $\Phi_0(q_1\cdots q_n)$ in $F$. 

If $\Phi_0(q_1 \cdots q_n) \neq 1_F$ then $M$ is a connected and simply connected planar diagram (van Kampen diagram) and if $\Phi_0(q_1 \cdots q_n)= 1$ in $F$ then $M$ is the union of a finite number of spines and tessellated spheres ($\cong S^2 \subseteq {\mathbb R}^3$) by the $\Delta_i$. In this case we call $M$ a \emph{spherical diagram}. We call $M$ a \emph{singular sphere} or just a sphere, if it consists of a single sphere. When the $q_i$ are labelled with $R_i$ and $\phi(q_1 \cdots q_n) = R_1 \cdots R_n = 1$ in $F$ then we call the sequence $S=(R_1, \ldots, R_n)$ an \emph{identity sequence over $\cal R$} and call it a \emph{simple identity sequence},  if the corresponding spherical diagram is a simple sphere.

\vspace{10pt}

\begin{definition}[Extended presentations]
	\addcontentsline{toc}{subsubsection}{\numberline{\thesubsubsection}  
		Definition (Extended presentations) }
	\ \\ Let ${\cal P} = \langle X \ | \ {\cal R} \rangle$ be a finite presentation of a group $G$ and let $\cal I$ be a set of simple identity sequences over $\cal R$. 
	The corresponding \emph{extended presentation} is the triple $\langle X \ | \ {\cal R} \ | \ {\cal I} \rangle$. 
	This notion was introduced in \cite{5} by Roger Fenn. Now let $A= A(\Gamma)$ be an Artin group given by $\langle X \ | \ {\cal R} \rangle$ in ({\rm III}). The identitites among relation $\cal I$ we choose are obtained from the 3-generated 
	standard
	parabolic subgroups $A(i,j,k)$ with defining subgraphs $\Gamma_{i,j,k}$ with the vertex set $V(\Gamma_{i,j,k}) = \{ x_i, \ x_j , \ x_k \}$ and the edge set $E(\Gamma_{i,j,k}) = \{e_{ij}, \ e_{i,k}, \ e_{j,k} \}$, $ 1 \leq i < j< k \leq n$, with labels $\lambda(e_{ij}) = \lambda(e_{j,k}) = 2 $ and $\lambda(e_{i,k})=m$, $ m \geq 4$. 
	The corresponding tessellated simple spheres $\Sigma_{i,j,k}$ are prisms with two regions of degree $2m$ (upper and lower) labelled with $R_{ik}$, which we call the \emph{big regions}, and a band of length $2n$ in which the regions are labelled by $R_{ij}$ and $R_{jk}$ in an alternating manner. See Figure 26 for $n=4$.
	
	\begin{center}
		\begin{tikzpicture}[scale =.4]
			\draw (3,0)--(6,0)--(9,3)--(9,6)--(6,9)--(3,9)--(0,6)--(0,3)--(3,0);
			\draw (3,1.5)--(6,1.5)--(7.5,3)--(7.5,6)--(6,7.5)--(3,7.5)--(1.5,6)--(1.5,3)--(3,1.5);
			\draw (0,3)--(1.5,3);
			\draw (0,6)--(1.5,6);
			\draw (3,9)--(3,7.5);
			\draw (6,9)--(6,7.5);
			\draw (7.5,6)--(9,6);
			\draw (7.5,3)--(9,3);
			\draw (6,0)--(6,1.5);
			\draw (3,0)--(3,1.5);
			\node at (4,4) {$R_{ik}$};
			\node at (4.5,8.2) {$R_{ij}$};
			\node at (7.1,7) {$R_{jk}$};
			\node at (8.3,4) {$R_{ij}$};
		\end{tikzpicture}
		
		Figure 26
	\end{center}
	
	We introduce one of our main tools, \emph{$\cal I$-moves}.
	
	Let $\Sigma = \Sigma_{i,j,k}$ and let $\omega$ be a simple closed curve on $\Sigma$. Then  $\omega$ subdivides $\Sigma$ into two simply connected subdiagrams, $\omega_R$ and $\omega_L$, with connected interiors, where $\omega_R$ is the submap of all regions to the right of $\omega$ and $\omega_L$ is the submap of all the regions to the left of $\omega$. Clearly $\overline{\omega_R} \cap \overline{\omega_L} = \omega$ and $\omega_R \cup \omega \cup \omega_L = \Sigma$.
	
	Let $S(\Sigma)$ be the set of all the connected, simply connected submaps of $\Sigma$ with connected interior and let $\beta_{\omega} : S(\Sigma) \rightarrow S(\Sigma)$ be the function which for every closed curve $\omega$ on $\Sigma$ sends $\omega_R$ to $\omega_L$ and $\omega_L$ to $\omega_R$ and leaves $\omega$ unaltered. 
	($\beta_{\omega}$ has nothing to do with $\beta$ in section 1.3)
	
	Let $M$ be a van Kampen $\cal R$-diagram over $F$ and let $M_0$ be a simply connected subdiagram of $M$ 
	with connected interior and simple boundary cycle $\mu$. 
	Suppose that $\overline{M_0}$ can be embedded into $\Sigma=\Sigma_{i,j,k}$, 
	for some (uniquely defined) $i$, $j$ and $k$ as an $\cal R$-diagram with image $N$ in $\Sigma$ 
	with $\partial N$ simple closed such that $\mu$ is mapped onto $\partial N$. Then $ N \in S(\Sigma)$.
\end{definition}
\vspace{10pt}

\begin{definition}[$\cal I$-moves]
	\addcontentsline{toc}{subsubsection}{\numberline{\thesubsubsection}  
		Definition ($\cal I$-moves) }	
	\ \\ Let $\cal I$ be a simple identity among relation realised by a uniquely defined sphere $\Sigma = \Sigma_{i,j,k}$. 
	Let $M$, $M_0$, $N$, $S(\Sigma)$ and $\beta_{\omega}$ be as defined above, 
	for some $\omega$. 
	An \emph{$\cal I$-move at $M_0$} is a surgery in $M$ by which we cut out $M_0$ and fill in instead by $\beta_{\omega}(N)$. Since $\partial N = \partial (\beta_{\omega} N)$, this is well defined. 
\end{definition}
\noindent Below we consider some examples which illustrate the way I-moves are used.
\ \\

\noindent {\bf Example \sssectionnum}
\addcontentsline{toc}{subsubsection}{\numberline{\thesubsubsection}  
	Example}
\ \\ $X=\{a,b,c\}$, ${\cal R}=\{R_1,R_2,R_3\}$, where $R_1=aba^{-1}b^{-1}$,
$R_2=aca^{-1}c^{-1}$, $R_3 =bcb^{-1}c^{-1}$.
Here $G\cong {\mathbb Z} \oplus {\mathbb Z} \oplus {\mathbb Z}$.
Let $W=abca^{-1}b^{-1}c^{-1}$. Then $W$ represents 1 in $G$. Figure 27 shows corresponding van Kampen diagrams.

\begin{minipage}{2in}
	\begin{center}
		\begin{tikzpicture}[scale=.45]
			\draw (4,0)--(7,2)--(7,6)--(4,8)--(1,6)--(1,2)--(4,0)--(4,4)--(1,6);
			\draw (4,4)--(7,6);
			\draw (2.5,.75)--(2.5,1)--(2.8,1);
			\node[below] at (2.5,.75) {$a$};
			\draw (5.5,.75)--(5.5,1)--(5.2,1);
			\node[below] at (5.5,.75) {$c$};
			\draw (2.5,4.75)--(2.5,5)--(2.8,5);
			\draw (5.5,4.75)--(5.5,5)--(5.2,5);
			\draw (2.5,6.75)--(2.5,7)--(2.2,7);
			\node[left] at (2.2,7) {$c$};
			
			\draw (.8,3.75)--(1,4)--(1.2,3.75);
			\node[left] at (1,4) {$b$};
			
			\draw (6.8,3.75)--(7,4)--(7.2,3.75);
			\node[right] at (7,4) {$b$};
			
			\draw (3.8,1.75)--(4,2)--(4.2,1.75);
			
			\draw (6.3,6.75)--(6,6.75)--(6,6.45);
			\node[right] at (6.3,7) {$a$};
			
			\node at (4,6) {$D_1$};
			\node at (2.5,3) {$D_2$};
			\node at (5.5,3) {$D_3$};
			\draw[fill] (4,0) circle(3pt);

		\end{tikzpicture}
		
		$M_0$
		
		(a)
	\end{center}
\end{minipage}\hspace{.2cm} \begin{minipage}{2in}
	\begin{center}
		\begin{tikzpicture}[scale=.45]
			\draw (3,0)--(6,2)--(6,6)--(3,8)--(0,6)--(0,2)--(3,0);
			\draw (0,2)--(3,4)--(3,8);
			\draw (3,4)--(6,2);
			\draw (1.4,.7)--(1.4,1)--(1.8,1);
			\node[left] at (1.4,.7) {$a$};
			\draw (3.9,.8)--(4.2,.8)--(4.2,.5);
			\draw (-.3,3.7)--(0,4)--(.3,3.7);
			\draw (5.7,3.7)--(6,4)--(6.3,3.7);
			\draw (1.2,3)--(1.5,3)--(1.5,2.7);
			\draw (5.1,2.8)--(4.8,2.8)--(4.8,2.5);
			\draw (1.2,7)--(1.5,7)--(1.5,6.7);
			\draw (4.5,6.7)--(4.5,7)--(4.8,7);
			\node at (5,3.5) {$a$};
			\node at (5,8,) {$a$};
			\node at (-1,4) {$b$};
			\node at (7,4) {$b$};
			\node at (5,1) {$c$};
			\node at (1.5,2.5) {$c$};
			\node at (1,8) {$c$};
			\node at (3,2) {$D'_1$};
			\node at (4.5,5) {$D'_2$};
			\node at (1.5,5) {$D'_3$};
			\draw[fill] (3,0) circle(3pt);
		\end{tikzpicture}
		
		$M_1$
		
		(b)
	\end{center}
\end{minipage}

\begin{center}
	Figure 27
\end{center}

\noindent {\bf Example \sssectionnum}
\addcontentsline{toc}{subsubsection}{\numberline{\thesubsubsection}  
	Example}
\ \\ Now let ${\cal E} = \langle a,b,c \ | \ {\cal  R} \ | \ I \rangle$, $ {\cal R} = \{ R_1,R_2,R_3 \}$, $I=R_1$, $R_2^b$, $R_3$,
$(R_1^{-1})^c$, $R_2^{-1}$, $(R_3^{-1})^a$ be an extended presentation,
where $R^x = xRx^{-1}$ and $R_1$, $R_2$, $R_3$ are as in Example 3.1.3.
The effect of $I$ on the diagrams of Figure 27 is the replacement of the diagram in Figure 27(a) by the diagram in
Figure 27(b).

\vspace{20pt}

\begin{minipage}{1.2in}
	\begin{center}
		\begin{tikzpicture}[scale=.3]
			\draw (4,0)--(7,2)--(7,6)--(4,8)--(1,6)--(1,2)--(4,0)--(4,4)--(1,6);
			\draw (4,4)--(7,6);
			\draw (2.5,.75)--(2.5,1)--(2.8,1);
			\node[below] at (2.5,.75) {$a$};
			\draw (5.5,.75)--(5.5,1)--(5.2,1);
			\node[below] at (5.5,.75) {$c$};
			\draw (2.5,4.75)--(2.5,5)--(2.8,5);
			\draw (5.5,4.75)--(5.5,5)--(5.2,5);
			\draw (2.5,6.75)--(2.5,7)--(2.2,7);
			\node[left] at (2.2,7) {$c$};
			
			\draw (.8,3.75)--(1,4)--(1.2,3.75);
			\node[left] at (1,4) {$b$};
			
			\draw (6.8,3.75)--(7,4)--(7.2,3.75);
			\node[right] at (7,4) {$b$};
			
			\draw (3.8,1.75)--(4,2)--(4.2,1.75);
			
			\draw (6.3,6.75)--(6,6.75)--(6,6.45);
			\node[right] at (5,6.5) {$a$};
			
			\node at (4,6) {$D_1$};
			\node at (2.5,3) {$D_2$};
			\node at (5.5,3) {$D_3$};
			
			\draw (7,2)--(11,6)--(4,8);
			\draw (9,3.7)--(9,4)--(8.7,4);
			\draw(8,6.6)--(8,6.9)--(8.3,7);
			\node at (9,3) {$a$};
			\node at (8,8) {$b$};
			\node at (8,5) {$D_0$};
			
		\end{tikzpicture}
		
		$abcb^{-1}a^{-1}c^{-1}$
		
		$M$
	\end{center}
\end{minipage}$\bf \begin{array}{c}I\\ \longrightarrow\end{array}$\begin{minipage}{1.2in}
	
	\begin{center}
		\begin{tikzpicture}[scale=.3]
			\draw (3,0)--(6,2)--(6,6)--(3,8)--(0,6)--(0,2)--(3,0);
			\draw (0,2)--(3,4)--(3,8);
			\draw (3,4)--(6,2);
			\draw (1.4,.7)--(1.4,1)--(1.8,1);
			\node[left] at (1.4,.7) {$a$};
			\draw (3.9,.8)--(4.2,.8)--(4.2,.5);
			\draw (-.3,3.7)--(0,4)--(.3,3.7);
			\draw (5.7,3.7)--(6,4)--(6.3,3.7);
			\draw (1.2,3)--(1.5,3)--(1.5,2.7);
			\draw (5.1,2.8)--(4.8,2.8)--(4.8,2.5);
			\draw (1.2,7)--(1.5,7)--(1.5,6.7);
			\draw (4.5,6.7)--(4.5,7)--(4.8,7);
			\node at (5,3.5) {$a$};
			\node at (5,8,) {$a$};
			\node at (-1,4) {$b$};
			\node at (6.3,4.1) {$b$};
			\node at (5,1) {$c$};
			\node at (1.5,2.5) {$c$};
			\node at (1,8) {$c$};
			\node at (3,2) {$D'_1$};
			\node at (4.5,5) {$D'_2$};
			\node at (1.5,5) {$D'_3$};
			\draw (6,2)--(9,8)--(3,8);
			\node at (7,7) {$D_0$};
			\draw (7,4.7)--(7.3,4.7)--(7.4,4.3);
			\node at (8,4.7) {$a$};
			\draw (6.3,8.2)--(6,8)--(6.3,7.8);
			\node at (6,8.5) {$b$};
			\draw[red,thick] (3,4)--(6,2)--(9,8)--(3,8)--(3,4);

		\end{tikzpicture}
		
		$M^I$
		
	\end{center}

\end{minipage} \ $\bf \begin{array}{c}\longrightarrow\\ \mbox{free}\\ \mbox{reduction}\end{array}$\begin{minipage}{1.2in}
	\begin{center}
		\begin{tikzpicture}[scale=.4]
			
			\draw (0,0)--(3,0)--(3,6)--(0,6)--(0,0);
			\draw (0,3)--(3,3);
			\draw (1.2,.3)--(1.5,0)--(1.2,-.3);
			\draw (1.2,3.3)--(1.5,3)--(1.2,2.7);
			\draw (1.2,6.3)--(1.5,6)--(1.2,5.7);
			\draw (-.3,1.7)--(0,2)--(.3,1.7);
			\draw (2.7,1.7)--(3,2)--(3.3,1.7);
			\draw (-.3,4.7)--(0,5)--(.3,4.7);
			\draw (2.7,4.7)--(3,5)--(3.3,4.7);
			
			\node at (-.5,2) {$a$};
			\node at (3.5,2) {$a$};
			\node at (-.5,5) {$b$};
			\node at (3.5,5) {$b$};
			\node at (1.5,-.5) {$c$};
			\node at (1.5,6.5) {$c$};
			\node at (1.5,2) {$D'_1$};
			\node at (1.5,5) {$D'_3$};

		\end{tikzpicture}

	\end{center}
	
	$M^I$ after free reduction
\end{minipage}

\begin{center}
	Figure 28
\end{center}

\vspace{20pt}

This changes the positions of the regions realising $R_1$, $R_2$ and $R_3$ and may reduce the number of regions. See Figure 28. 

\vspace{10pt}
\noindent {\bf Example \sssectionnum}
\addcontentsline{toc}{subsubsection}{\numberline{\thesubsubsection}  
	Example}

\noindent Let ${\cal E} = \langle a,b,c \ | \ (ab)^2(a^{-1}b^{-1})^2 \ , \ cac^{-1}a^{-1} \ , \ cbc^{-1}b^{-1} \ | \ I \rangle$ be an extended presentation, where

$I=R_0, R_1, R_2^{h_1}, R_1^{h_2}, R_2^{h_3}, R_1^{-1 h_5}, R_2^{-1 h_6}, R_1^{-1 h_7},
R_2^{-1 h_8}, cR_0^{-1}c^{-1}$.

and where $R_0=(ab)^2(a^{-1}b^{-1})^2$ , $R_1=cac^{-1}a^{-1}$, $R_2=cbc^{-1}b^{-1}$ and $h_i$ is the prefix of $R_0$ of length $i$.

\vspace{20pt}

\begin{minipage}{1.6in}
	\begin{center}
		\begin{tikzpicture}[scale=.25]
			\draw (2,0)--(5,0)--(7,3)--(7,6)--(5,9)--(2,9)--(0,6)--(0,3)--(2,0);
			\draw (7,3)--(10,3)--(10,7)--(9,10)--(5,9);
			\draw (7,6)--(10,7);
			\draw (3.7,-.3)--(4,0)--(3.7,.3);
			\node at (4,.7) {$a$};
			
			\draw (3.7,8.7)--(4,9)--(3.7,9.3);
			\node at (4,8.4) {$a$};
			
			\draw (-.3,4.7)--(0,5)--(.3,4.7);
			\node at (.6,5) {$a$};
			
			\draw (6.7,4.7)--(7,5)--(7.3,4.7);
			\node at (6.4,5) {$a$};
			
			\draw (.9,1.2)--(.9,1.6)--(1.3,1.6);
			\node at (2,2) {$b$};
			
			\draw (5.7,1.4)--(6,1.4)--(6,1.1);
			\node at (5,2) {$b$};
			
			\draw (.9,7.8)--(1.2,7.8)--(1.2,7.5);
			\node at (2,7) {$b$};
			
			\draw (5.8,7.4)--(5.9,7.8)--(6.3,7.7);
			\node at (5,7) {$b$};
			
			\draw (7.1,9.7)--(7.5,9.6)--(7.5,9.2);
			\node at (8,9) {$c$};
			
			\draw (9.2,8.2)--(9.5,8.5)--(9.8,8.4);
			\node at (8.5,7.5) {$b$};
			
			\draw (9.7,4.7)--(10,5)--(10.3,4.7);
			\node at (9,5) {$a$};
			
			\draw (8.2,3.3)--(8.5,3)--(8.2,2.7);
			\node at (8.5,4) {$c$};

		\end{tikzpicture}
		(a)
	\end{center}
\end{minipage}$\bf \longrightarrow$ \hspace{.2cm}\begin{minipage}{2in}
	\begin{center}
		\begin{tikzpicture}[scale=.2]

			\draw (4,3)--(7,3)--(9,6)--(9,9)--(7,12)--(4,12)--(2,9)--(2,6)--(4,3);
			\draw (7,12)--(7,14)--(3,14)--(0,9)--(0,6)--(3,1)--(8,1)--(11,4)--(9,6);
			\draw (8,1)--(7,3);
			\draw (3,1)--(4,3);
			\draw (0,6)--(2,6);
			\draw (0,9)--(2,9);
			\draw (4,12)--(3,14);
			
			\draw (5.7,14.3)--(6,14)--(5.7,13.7);
			\node at (6,15) {$a$};
			
			\draw (4.7,11.7)--(5,12)--(4.7,12.3);

			\draw (6.7,13.3)--(7,13)--(7.3,13.3);
			\node at (8,13) {$c$};
			
			\draw (8,10.2)--(8,10.5)--(8.3,10.5);
			\node at (7,10) {$b$};
			
			\draw (8.7,7.7)--(9,8)--(9.3,7.7);
			\node at (8,8) {$a$};
			
			\draw (10,4.7)--(10,5)--(10.3,5);
			\node at (10,6) {$c$};
			
			\draw (9.7,3)--(10,3)--(10,2.7);
			\node at (10,2) {$b$};

		\end{tikzpicture}
		
		(b)
		
	\end{center}
\end{minipage}

\begin{center}
	Figure 29
\end{center}

\vspace{20pt}

In Figure 29 an $I$-move changes the diagram in Figure 29(a) to the diagram in Figure 29(b). As in Figure 28 this may reduce the number of regions in the diagram 29.

But even if it does not reduce the number of 
regions, 
an $I$-move may reduce the number of equivalence classes $\{\Delta(D) \}$. See Figure 30.

\vspace{20pt}
\vbox{
	\begin{minipage}{1.6in}
		\begin{center}
			\begin{tikzpicture}[scale=.2]
				\node at (3,5) {$D_1$};
				\node at (12,10) {$D_2$};
				
				\draw (2,0)--(5,0)--(7,3)--(7,6)--(5,9)--(2,9)--(0,6)--(0,3)--(2,0);
				\draw (7,3)--(10,3)--(10,7)--(9,10)--(5,9);
				\draw (7,6)--(10,7);
				\draw (3.7,-.3)--(4,0)--(3.7,.3);
				\node at (4,.7) {$a$};
				
				\draw (3.7,8.7)--(4,9)--(3.7,9.3);
				\node at (4,8.4) {$a$};
				
				\draw (-.3,4.7)--(0,5)--(.3,4.7);
				\node at (.6,5) {$a$};
				
				\draw (6.7,4.7)--(7,5)--(7.3,4.7);
				\node at (6.4,5) {$a$};
				
				\draw (.9,1.2)--(.9,1.6)--(1.3,1.6);
				\node at (2,2) {$b$};
				
				\draw (5.7,1.4)--(6,1.4)--(6,1.1);
				\node at (5,2) {$b$};
				
				\draw (.9,7.8)--(1.2,7.8)--(1.2,7.5);
				\node at (2,7) {$b$};
				
				\draw (5.8,7.4)--(5.9,7.8)--(6.3,7.7);
				\node at (5,7) {$b$};
				
				\draw (7.1,9.7)--(7.5,9.6)--(7.5,9.2);
				\node at (8,9) {$c$};
				
				\draw (9.2,8.2)--(9.5,8.5)--(9.8,8.4);
				\node at (8.5,7.5) {$b$};
				
				\draw (9.7,4.7)--(10,5)--(10.3,4.7);
				\node at (9,5) {$a$};
				
				\draw (8.2,3.3)--(8.5,3)--(8.2,2.7);
				\node at (8.5,4) {$c$};
				
				\draw (10,3)--(12,3)--(15,6)--(16,10)--(14,14)--(11,13)--(9,10);
				
			\end{tikzpicture}
			
			\vspace{1cm} 
			
			(a)

			\[
			[D_1]_t \neq [D_2]_t
			\]
		\end{center}
	\end{minipage}$\bf \longrightarrow$ \hspace{.2cm}\begin{minipage}{2in}
		\begin{center}
			\begin{tikzpicture}[scale=.25]
				\node at (5,8) {$D'_1$};
				\node at (11,10) {$D'_2$};
				
				\draw (4,3)--(7,3)--(9,6)--(9,9)--(7,12)--(4,12)--(2,9)--(2,6)--(4,3);
				\draw (7,12)--(7,14)--(3,14)--(0,9)--(0,6)--(3,1)--(8,1)--(11,4)--(9,6);
				\draw (8,1)--(7,3);
				\draw (3,1)--(4,3);
				\draw (0,6)--(2,6);
				\draw (0,9)--(2,9);
				\draw (4,12)--(3,14);
				
				\draw (5.7,14.3)--(6,14)--(5.7,13.7);
				\node at (6,15) {$a$};
				
				\draw (4.7,11.7)--(5,12)--(4.7,12.3);

				\draw (6.7,13.3)--(7,13)--(7.3,13.3);
				\node at (8,13) {$c$};
				
				\draw (8,10.2)--(8,10.5)--(8.3,10.5);
				\node at (7,10) {$b$};
				
				\draw (8.7,7.7)--(9,8)--(9.3,7.7);
				\node at (8,8) {$a$};
				
				\draw (10,4.7)--(10,5)--(10.3,5);
				\node at (11,5) {$c$};
				
				\draw (9.7,3)--(10,3)--(10,2.7);
				\node at (10,2) {$b$};
				
				\draw[green,thick] (9,6)--(14,5)--(16,9)--(15,12)--(13,15)--(10,14)--(7,12)--(4,12)--(2,9)--(2,6)--(4,3)--(7,3)--(9,6);
				
				\node at (14,7) {$a$};
				\draw (14.7,6.9)--(15,7.2)--(15.3,6.9);
				\node at (11.5,6.3) {$b$};
				\draw (11.5,5.3)--(11.8,5.6)--(11.5,5.9);
				
				\node at (15,10) {$b$};
				\draw (14.7,11.7)--(15,12)--(15.3,11.7);
				
				\node at (13,12.5) {$a$};
				\draw (13.7,13.5)--(14,13.8)--(14.3,13.5);
				
				\node at (12,13.7) {$b$};
				\draw (11.4,14.8)--(11.7,14.5)--(11.4,14.1);
				
				\node at (11,12.5) {$a$};
				\draw (8.7,13.2)--(9,13.2)--(9,12.9);

			\end{tikzpicture}
			
			(b)
			
			\[
			[D'_1]_t = [D'_2]_t
			\]
		\end{center}
		
	\end{minipage}
	
	\begin{center}
		Figure 30
	\end{center}
}
\vspace{10pt}

\noindent {\bf Remark \sssectionnum}
\addcontentsline{toc}{subsubsection}{\numberline{\thesubsubsection}  
	Remark}
\begin{itemize}
	\item[(a)]  Let $\Delta \in Reg_{4^+}({\mathbb M})$ and let $B$ be a band such that $||\partial \Delta \cap \theta || \geq 2$, for one of the sides $\theta$ of $B$. 
	Considering $\Delta$ as a subdiagram of $M$, it follows by induction on $|\Delta|$ using Lemma 1.5.1  that there exists a band $B'$ such that $B \cup B'$ is a closed (annular) band and there exists a copy $\Delta'$ of $\Delta$ 
	such that $\theta_1 = \partial \Delta'$, where $\theta_1$ is a side of $B \cup B'$ such that $\Delta \cup(B \cup B') \cup \Delta'$ is the concatenation of spheres ${ S}_D$, for every $D \in Reg(\Delta)$. We call the replacement of $B \cup \Delta$ with $B' \cup \Delta'$ the \emph{extended I-move}. 
	It is a sequence of I-moves. Thus, the extended I-move $J$ sends $\Delta \in Reg_{4^+}({\mathbb M})$ to $\Delta^J \in Reg_{4^+}({\mathbb M}^J)$. 
	Here $\Delta^J$ is the diagram obtained from $\Delta$ by the action of $J$ on it. $\Delta^J=\Delta'$.
	\item[(b)]  Let $\mathbb B$ be a band-bundle emanating from $\partial \Delta$ and ending on $\partial M$. Considering $\Delta \cup {\mathbb B}$ we can replace it with $\Delta' \cup {\mathbb B}'$ via a sequence $J$ of I-moves, where
	${\mathbb B} \cup {\mathbb B}'$ is an annular band-bundle and $\Delta \cup \Delta' \cup B \cup B'$ is the concatanation of the corresponding prisms.
	Then $\Delta'$ is a boundary region of ${\mathbb M}^J$ with $\partial \Delta' \cap \partial M^J$ the top of $\mathbb B$. 
	We call the replacing of $\Delta \cup {\mathbb B}$ by $\Delta' \cup {\mathbb B}'$ a  \emph{pushing up of $\Delta$ to $\theta$}. If $\theta = \partial M$ we call the replacing of $\Delta \cup {\mathbb B}$ by $\Delta' \cup {\mathbb B}'$, a 
	{\it pushing up to $\partial M$}.

\end{itemize}

\vspace{10pt}


\subsection{Banded Subdiagrams and transfer of regions}
In this subsection we assume $M\in \mathcal M_3(W)$ and assumption $\mathcal{H}$

\vspace{10pt}

\begin{lemma}
	\addcontentsline{toc}{subsubsection}{\numberline{\thesubsubsection}  
		Lemma ($I$-moves leave the diagrams in 
		$\mathcal{M}_2(W)$  )}
	Let $M \in {\cal M}_2(W)$ and let I be an I-move or extended I-move. Then $M^I \in {\cal M}_2(W)$.
\end{lemma}

\vspace{10pt}

\noindent {\bf Proof} \quad We have to show that $I$  increases  neither $|Reg_{4^+}(M)|$ nor $|Reg_{4^+}({\mathbb M})|$. 
Suppose that $I$ replaces the simply connected subdiagram $N$ with connected interior by a diagram $N'$. 

We use notation of Definition 3.1.1. 
If $n=2$ then clearly I cannot alter $|Reg_{4^+}(M)|$ and $|Reg_{4^+}({\mathbb M})|$. If $ n \geq 4$ then $N$ contains exactly one big region $\Delta$  and $N'$ also contains exactly one big region $\Delta'$, 
which is a copy of $\Delta$. Since $M \in {\cal M}_2(W)$ the result follows.
\hfill $\Box$

\vspace{10pt}


\subsubsection{Definition (Banded Subdiagrams)}
Let $ B = \langle E_1, \ldots , E_k \rangle$ be a band in $M$ with sides $\sigma_1$ and $\sigma_2$ and let $\theta$ be a boundary path of a simply connected subdiagram $M_0$ with connected interior.

\begin{itemize}
	\item[(a)] Say that \emph{$M_0$ is banded at $\theta$ by $B$} if $||\theta|| \geq 1$ and $\theta \subseteq \sigma_1 \cup \sigma_2$. See Figure 31(a).
\end{itemize}

\hspace{-20pt}\begin{minipage}{1.6in}
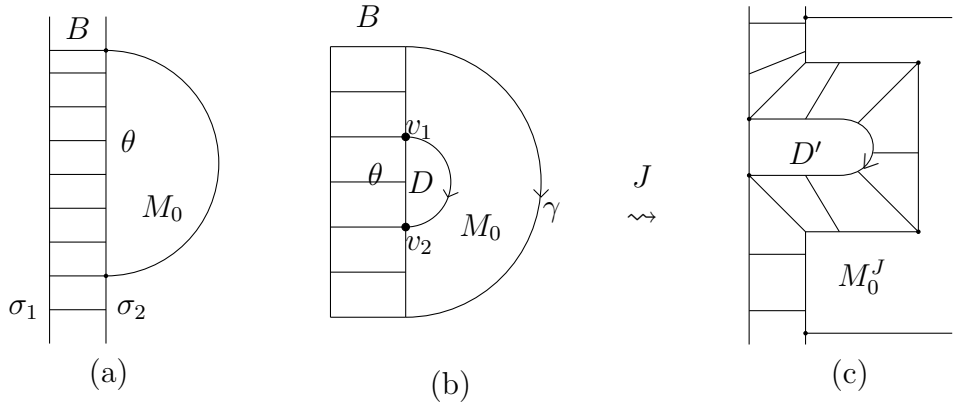
	\begin{center}
		\begin{tikzpicture}[scale=.15]
			
			\draw[rotate around={180:(5,6)}] (5,6) arc(90:270:10);
			\draw (0,0)--(0,29);
			\draw (5,0)--(5,29);
			\draw (0,3)--(5,3);
			\draw (0,6)--(5,6);
			\draw (0,9)--(5,9);
			\draw (0,12)--(5,12);
			\draw (0,15)--(5,15);
			\draw (0,18)--(5,18);
			\draw (0,21)--(5,21);
			\draw (0,24)--(5,24);
			\draw (0,26)--(5,26);
			
			\draw[fill] (5,6) circle (.15);
			\draw[fill] (5,26) circle (.15);
			
			\node[left] at (0,3) {$\sigma_1$};
			\node[right] at (5,3) {$\sigma_2$};
			\node at (10,12) {$M_0$};
			\node at (7,18) {$\theta$};
			\node at (2.5,28) {$B$};

		\end{tikzpicture}
		
		(a)
	\end{center}
\end{minipage}\hspace{10pt}\begin{minipage}{1.7in}
	\begin{center}
		\begin{tikzpicture}[scale=.2]
			\draw[rotate around={180:(5,6)}] (5,6) arc(90:270:9);
			\draw[rotate around={180:(5,12)}] (5,12) arc(90:270:3);
			\draw (0,6)--(0,24);
			\draw (5,6)--(5,24);
			
			\draw (0,6)--(5,6);
			\draw (0,9)--(5,9);
			\draw (0,12)--(5,12);
			\draw (0,15)--(5,15);
			\draw (0,18)--(5,18);
			\draw (0,21)--(5,21);
			\draw (0,24)--(5,24);
			
			\draw (7.5,14.5)--(7.8,14)--(8.5,14.3);
			\draw (13.5,14.5)--(14,14)--(14.5,14.5);

			\draw[fill]  (5,12) circle(.25);
			\draw[fill]  (5,18) circle(.25);
			\node at (5.9,10.8) {$v_2$};
			\node at (5.9,18.5) {$v_1$};
			\node at (3,15.5) {$\theta$};
			\node at (14.7,13) {$\gamma$};
			
			\node at (10,12) {$M_0$};
			\node at (6,15) {$D$};
			\node at (2.5,26) {$B$};
		\end{tikzpicture}
		
		\vspace{.5cm}
		
		(b)
	\end{center}
\end{minipage}$\begin{array}{c} J\\ \bf \leadsto \end{array}$ \begin{minipage}{1.8in}
	\begin{center}
		\begin{tikzpicture}[scale=.15]
			\draw (0,0)--(0,15)--(8,15) .. controls (12,15) and (12,20) .. (8,20)
			node [xscale=0.7,pos=0.3,sloped] {$<$}
			--(0,20)--(0,30);
			\draw (5,0)--(5,10)--(15,10)--(15,25)--(5,25)--(5,30);
			
			\draw (5,29)--(18,29);
			\draw[fill] (5,29) circle (.15);
			\draw (5,1)--(18,1);
			\draw[fill] (5,1) circle (.15);
			
			\draw[fill] (0,15) circle (.15);
			\draw[fill] (0,20) circle (.15);
			\draw[fill] (15,25) circle (.15);
			\draw[fill] (15,10) circle (.15);
			
			\draw (0,15)--(0,20);
			
			\node at (10,6) {$M_0^J$};
			\node at (5,17) {$D'$};
			
			\draw (0,15)--(5,10);
			\draw (0,20)--(5,25);
			\draw (0,3)--(5,3);
			\draw (0,8)--(5,8);
			\draw (0,28.5)--(5,28.5);
			\draw (0,24)--(5,26);
			
			\draw (11.,17)--(15,17);
			\draw (9.6667,15.3333)--(15,10);
			\draw (9.6667,19.6667)--(15,25);
			
			\draw (5,20)--(8,25);
			\draw (5,15)--(8,10);
			
		\end{tikzpicture}
		
		(c)
	\end{center}
\end{minipage}

\begin{minipage}{1.6in}

	\
	
\end{minipage}\hspace{10pt}\begin{minipage}{1.7in}
	\begin{center}
		\begin{tikzpicture}[scale=.2]
			\draw[rotate around={180:(5,6)}] (5,6) arc(90:270:9);
			\draw[rotate around={180:(5,12)}] (5,12) arc(90:270:3);
			\draw (0,6)--(0,24);
			\draw (5,6)--(5,24);
			
			\draw (0,6)--(5,6);
			\draw (0,9)--(5,9);
			\draw (0,12)--(5,12);
			\draw (0,15)--(5,15);
			\draw (0,18)--(5,18);
			\draw (0,21)--(5,21);
			\draw (0,24)--(5,24);
			
			\draw (7.5,14.5)--(7.8,14)--(8.5,14.3);
			\draw (13.5,14.5)--(14,14)--(14.5,14.5);

			\draw[fill]  (5,12) circle(.25);
			\draw[fill]  (5,18) circle(.25);
			\node at (5.9,10.8) {$v_2$};
			\node at (5.9,18.5) {$v_1$};
			\node at (3,15.5) {$\theta$};
			\node at (14.7,13) {$\gamma$};
			
			\node at (10,12) {$M_0$};
			\node at (6,15) {$\Delta$};
			\node at (2.5,26) {$B$};
		\end{tikzpicture}
		
		\vspace{.5cm}
		
		(b')
	\end{center}
\end{minipage}$\begin{array}{c} J\\ \bf \leadsto \end{array}$ \begin{minipage}{1.8in}
	\begin{center}
		\begin{tikzpicture}[scale=.15]
			\draw (0,0)--(0,15)--(8,15) .. controls (12,15) and (12,20) .. (8,20)
			node [xscale=0.7,pos=0.3,sloped] {$<$}
			--(0,20)--(0,30);
			\draw (5,0)--(5,10)--(15,10)--(15,25)--(5,25)--(5,30);
			
			\draw (5,29)--(18,29);
			\draw[fill] (5,29) circle (.15);
			\draw (5,1)--(18,1);
			\draw[fill] (5,1) circle (.15);
			
			\draw[fill] (0,15) circle (.15);
			\draw[fill] (0,20) circle (.15);
			\draw[fill] (15,25) circle (.15);
			\draw[fill] (15,10) circle (.15);
			
			\draw (0,15)--(0,20);
			
			\node at (10,6) {$M_0^J$};
			\node at (5,17) {$\Delta'$};
			
			\draw (0,15)--(5,10);
			\draw (0,20)--(5,25);
			\draw (0,3)--(5,3);
			\draw (0,8)--(5,8);
			\draw (0,28.5)--(5,28.5);
			\draw (0,24)--(5,26);
			
			\draw (11.,17)--(15,17);
			\draw (9.6667,15.3333)--(15,10);
			\draw (9.6667,19.6667)--(15,25);
			
			\draw (5,20)--(8,25);
			\draw (5,15)--(8,10);
			
		\end{tikzpicture}
		
		(c')
	\end{center}
\end{minipage}

\begin{center}
	Figure 31
\end{center}

\begin{itemize}
	\item[(b)] Let notation be as in part (a) and let $\theta = v_0 \theta_1v_1 \theta_2 v_2 \cdots \theta_rv_r$ be a decomposition of $\theta$, $v_i$ vertices, $r \geq 1$. Say that $M_0$ is \emph{banded at $(\theta_1, \ldots , \theta_r)$ by bands $B_i$}, respectively, if $M_0$ is banded at $\theta_i$ by $B_i$ for every $i$, $ i = 1, \ldots , r$.
\end{itemize}

\vspace{10pt}

\begin{definition}	[transfer of regions from $M_0$]
	\addcontentsline{toc}{subsubsection}{\numberline{\thesubsubsection}  
		Definition (transfer of regions from $M_0$)}
	\  \\  
	Let $M$ be a connected, simply connected $\cal R$-diagram over $F$ with connected interior and let $M_0$ be a connected, simply connected $\cal R$-sub-diagram over $F$ with connected interior. 
	Let $D$ be a boundary region of $M_0$ with $\mu := \partial D \cap \partial M_0$ connected and let $B= B_1 \cup B_2 \cup B_3$ be such that $\partial D \cap \partial B = \partial D \cap \partial B_2$. Observe that $D \cup B_2$
	embeds to the prism $P_m$, $m = 2n(D)$, $P_m = D \cup B_2 \cup D' \cup B'_2$, where $B_2 \cup B'_2$  is an annulus with $m$ regions from $Reg_2(M)$ (a closed band) and $D$ is a region with $n(D') = n(D)$ and $Supp (D) = Supp(D')$. 
	(We identify $D$ and $B_2$ with their images in $P_m$).  
	Notice that $B_2$ and $B'_2$ have the same poles $\alpha$ and $\beta$, hence $B' = B_1 \cup B'_2 \cup B_3$ is a band. 
	Replace $B \cup D$ with $B' \cup D'$ (which have the same boundaries) and denote the resulting diagram by $M^J_0$, where $J$ is the $I$-move which replaced $B \cup D$ with $B' \cup D'$. Thus $M_0 \cup B$ is replaced by $\left((M_0 \cup B) \setminus (D \cup B)\right) \cup (D' \cup B')$.
	Hence
	
	\vspace{10pt}
	
	$M_0^J = M_0 \setminus D$, $B^J = B'$ and $M_0^J$ is banded by $B'$. \hfill ({\rm V})

	\vspace{10pt}
	
	In this situation we say that \emph{$J$ transferred $D$ from $M_0$}. 
	See Figure 31(b) and 31(c). 
	If every region can be transferred from $M_0$, we say that $M_0$ is \emph{transferable}.
	The removal of $D \in Reg(M_0)$ naturally extends to transfer of regions $\Delta$ from $\mathbb M$.
	(See figures 31(b') and 31(c').)
\end{definition}

Notice that if instead of  $B$ we have  a band-bundle then $D$ can be transferred by a sequence of $I$-moves beyond the band-bundle, without destroying the bundle structure. 
See Figure 32(a) and 32(b) for $D$, and 32(a') and 32(b') for $\Delta$.

\begin{minipage}{1.7in}
	\begin{center}
		\begin{tikzpicture}[scale=.2]
			\draw[rotate around={180:(5,6)}] (5,6) arc(90:270:9);
			\draw[rotate around={180:(5,12)}] (5,12) arc(90:270:3);
			\draw (0,6)--(0,24);
			\draw (5,6)--(5,24);
			\draw (2.5,6)--(2.5,24);

			\draw (0,6)--(5,6);
			\draw (0,9)--(5,9);
			\draw (0,12)--(5,12);
			\draw (0,15)--(5,15);
			\draw (0,18)--(5,18);
			\draw (0,21)--(5,21);
			\draw (0,24)--(5,24);

			\node[left] at (0,15) {$\theta$};
			\node at (10,12) {$M_0$};
			\node at (6,15) {$D$};
			\node at (1,26) {$B_1$};
			\node at (4,26) {$B_2$};
		\end{tikzpicture}
		
		(a)
		
	\end{center}
\end{minipage}$\begin{array}{c} J\\ \bf \leadsto \end{array}$ \begin{minipage}{1.8in}
	\begin{center}
		\begin{tikzpicture}[scale=.15]
			\draw (0,0)--(0,15)--(8,15) .. controls (12,15) and (12,20) .. (8,20)
			--(0,20)--(0,30);
			\draw (5,0)--(5,10)--(15,10)--(15,25)--(5,25)--(5,30);
			\draw (2.5,0)--(2.5,12.5)--(12.5,12.5)--(12.5,22.5)--(2.5,22.5)--(2.5,30);
			
			\draw (5,29)--(16,29);
			\draw (5,1)--(16,1);
			\draw[rotate around={180:(16,1)}] (16,1) arc(90:270:14);

			\draw (0,15)--(0,20);
			
			\node at (10,6) {$M_0$};
			\node at (5,17) {$D'$};
			\node at (1,32) {$B'_1$};
			\node at (4,32) {$B'_2$};
			
			\node[left] at (0,17) {$\theta$};
			
			\draw (0,3)--(5,3);
			\draw (0,8)--(5,8);
			\draw (0,28.5)--(5,28.5);
			\draw (0,26)--(5,26);
			
			\draw (0,15)--(5,10);
			\draw (0,20)--(5,25);
			
			\draw (11,17)--(15,17);
			\draw (9.6,15.4)--(15,10);
			\draw (9.6,19.6)--(15,25);
			
			\draw (5,20)--(8,25);
			\draw (5,15)--(8,10);
		\end{tikzpicture}
		
		(b)
		
	\end{center}
\end{minipage}

\begin{minipage}{1.7in}
	\begin{center}
		\begin{tikzpicture}[scale=.2]
			\draw[rotate around={180:(5,6)}] (5,6) arc(90:270:9);
			\draw[rotate around={180:(5,12)}] (5,12) arc(90:270:3);
			\draw (0,6)--(0,24);
			\draw (5,6)--(5,24);
			\draw (2.5,6)--(2.5,24);

			\draw (0,6)--(5,6);
			\draw (0,9)--(5,9);
			\draw (0,12)--(5,12);
			\draw (0,15)--(5,15);
			\draw (0,18)--(5,18);
			\draw (0,21)--(5,21);
			\draw (0,24)--(5,24);
			
			\node[left] at (0,15) {$\theta$};
			
			\node at (10,12) {$M_0$};
			\node at (6,15) {$\Delta$};
			\node at (1,26) {$B_1$};
			\node at (4,26) {$B_2$};
		\end{tikzpicture}
		
		(a')
		
	\end{center}
\end{minipage}$\begin{array}{c} J\\ \bf \leadsto \end{array}$ \begin{minipage}{1.8in}
	\begin{center}
		\begin{tikzpicture}[scale=.15]
			\draw (0,0)--(0,15)--(8,15) .. controls (12,15) and (12,20) .. (8,20)
			--(0,20)--(0,30);
			\draw (5,0)--(5,10)--(15,10)--(15,25)--(5,25)--(5,30);
			\draw (2.5,0)--(2.5,12.5)--(12.5,12.5)--(12.5,22.5)--(2.5,22.5)--(2.5,30);
			
			\draw (5,29)--(16,29);
			\draw (5,1)--(16,1);
			\draw[rotate around={180:(16,1)}] (16,1) arc(90:270:14);

			\draw (0,15)--(0,20);
			
			\node at (10,6) {$M_0$};
			\node at (5,17) {$\Delta'$};
			\node at (1,32) {$B'_1$};
			\node at (4,32) {$B'_2$};
			
			\node[left] at (0,17) {$\theta$};
			
			\draw (0,3)--(5,3);
			\draw (0,8)--(5,8);
			\draw (0,28.5)--(5,28.5);
			\draw (0,26)--(5,26);
			
			\draw (0,15)--(5,10);
			\draw (0,20)--(5,25);
			
			\draw (11,17)--(15,17);
			\draw (9.6,15.4)--(15,10);
			\draw (9.6,19.6)--(15,25);
			
			\draw (5,20)--(8,25);
			\draw (5,15)--(8,10);
		\end{tikzpicture}
		
		(b')
		
	\end{center}
\end{minipage}

\begin{center}
	Figure 32
\end{center}

We are interested in a special type of banded subdiagrams (2-banded subdiagrams) which we consider in the next subsection.

\vspace{10pt}

\begin{definition}[2-banded subdiagram] (See Figure 33(c).)
	\addcontentsline{toc}{subsubsection}{\numberline{\thesubsubsection}  
		Definition (2-banded subdiagram)}
	\ \\ Let $ W \neq 1$ be a cyclically reduced word in $F$ which represents 1 in $A$. Let $M \in {\cal M}_3(W)$  and let $Q$ be a connected, 
	simply connected subdiagram of $M$ with connected interior. 
	Say that $Q$ is \emph{2-banded} if
	\begin{itemize}
		\item[(1)] there is an a-band $B_l$ and a b-band $B_r$ in $Q$, $ a,b \in X$. ($a=b$ is not excluded.)
		\item[(2)] there are connected and simply connected subdiagrams $Q_r$, $Q_l$ and $P$ in $Q$ and paths $\theta_r$ and $\theta_l$ such that
	\end{itemize}
	\[\left.\begin{array}{cl}
		(i) & Q=Q_l \cup B_{l }\cup_{\theta_l} P \cup_{\theta_r}B_r \cup Q_r \ \mbox{ and} \ 
		\partial Q = \omega_1 \alpha_l^{-1} \gamma \alpha_r^{-1} \omega_3 \beta_r\beta_l,\\
		& \mbox{where} \ 
		\theta_l = \partial B_l \cap \partial P, \  \theta_r = \partial B_r \cap \partial P, \  \theta_l \ \mbox{ and} \ \theta_r
		\ \mbox{ are sides of}  \\
		&  B_l \ \mbox{ and} \ B_r \ \mbox{ respectively,} \ \omega_1 = \partial Q_l \cap \partial Q, \  \gamma = \partial P \cap \partial Q, \
		\\
		& \omega_3 = \partial Q_r \cap \partial Q \ \mbox{ are connected}, \  \alpha_l \ \mbox{ and} \  \beta_l \ \mbox{ are poles of} \\  
		& \ B_l \ \mbox{ and} \ \alpha_r \ \mbox{and} \  \beta_r \ 
		\mbox{ are poles of} \  B_r \\
		(ii)& o(\theta_l) = o(\theta_r)
	\end{array}\right\} ({\rm VI})
	\]
	\begin{itemize}
		\item[(3)] If $Int(P) \neq \emptyset$ then $\Phi(\partial P)$ is cyclically reduced. 
	\end{itemize} 
	
	Call $Q$ \emph{closed 2-banded} if $t(\theta_l) = t(\theta_r)$. Let $U$ be a cyclically reduced word in $F$ which represents 1 in A.
	Denote by ${\cal Q}_M(U)$ the collection of all the 2-banded subdiagrams of $M$ with boundary label $U$.
\end{definition}

We shall show that under certain conditions there is a boundary region of $M_0$ which we can transfer from $M_0$. We wish to repeat the process of transferring until either no regions remain in $M_0$, 
or $M_0$ has a special structure. The next Proposition shows that $I$-moves preserve these "certain conditions" hence we can repeat transferring regions  
as long as there are regions in $M_0^J$, or $M_0^J$ does not have the special structure.
Here $M_0^J$ is the result of applying $J$ on $M_0$, $J$ a sequence of $I$-moves.

\vspace{10pt}

\begin{definition} [Standard I-moves]\ 
	\addcontentsline{toc}{subsubsection}{\numberline{\thesubsubsection} Definition (Standard I-moves)}
	Say that an $I$-move $J$ is \emph{standard} if $\theta$ is a maximal piece in the sense that 
	$\theta$ on Figure 32(b) cannot be extended as a common boundary of $D$ and $B$.
\end{definition}

\vspace{10pt}

\begin{proposition}
	\addcontentsline{toc}{subsubsection}{\numberline{\thesubsubsection} Proposition (I-moves leave banded diagrams in $\mathcal M_2(W)$)}
	Let $W$ be a cyclically reduced word in $F$, $ W \neq 1$, which represents 1 in $A$ and let $M \in {\cal M}_3(W)$. Let $Q$ be a subdiagram of $M$ with cyclically reduced boundary label $U$ and suppose that $Q \in {\cal Q}_M(U)$. 
	(Thus $Q$ is 2-banded.) Let notation be as in Figure 33(c).
	
	Assume that $ P$ 
	contains a boundary region $\Delta \in Reg_{4^+,t}({ P})$ or $D \in Reg_2({ P})$ such that $\partial \Delta \cap \partial P$ or $\partial D \cap \partial P$, respectively, are connected and 
	$\partial \Delta \cap \partial P = \partial \Delta \cap \theta$, or $\partial D \cap \partial P = \partial D \cap \theta$, $ \theta = \theta_r$ or $\theta = \theta_l$ such that $||\partial \Delta \cap \theta|| \geq2$, or 
	$||\partial \theta \cap \partial D|| \geq 2$. Let $I$ be the standard $I$-move on $B \cup \Delta$, or $B \cup D$, respectively, where $B=B_r$ or $B=B_l$ according as $\theta = \theta_r$ or $\theta = \theta_l$, respectively, which transfers $D$ or $\Delta$, respectively, from $Q$. 
	Then $Q^I \in {\cal Q}_{M^I}(U)$.
\end{proposition}

\vspace{10pt}

\noindent {\bf Proof} \\
By the definition of ${\cal Q}_M(U)$ we have to prove each of the following
\begin{itemize}
	\item[(1)] $M^I \in {\cal M}_2(W)$ and $B^I$ are bands, $B \in \{B_r, B_l\}$.
	\item[(2)] $B^I$, $\theta^I$, $Q^I$, $P^I$ satisfy ($i$) and ($ii$) of ({\rm VI}), where $\theta \in \{ \theta_l, \theta_r\}$.
	\item[(3)] If $Int(P^I) \neq \emptyset$ then $\Phi(\partial (P^I))$ is cyclically reduced.
\end{itemize}
We prove these  in turn. Assume that $P$ contains a boundary region $K \in Reg_2(M) \cup Reg_{4^+,t}({ P})$.

\begin{itemize}
	\item[(1)] $M^I \in {\cal M}_2(W)$ by Lemma 3.2.1.
	If $B^I$ is not an adequate band then $\partial B$ is not simple. This may happen only if $\partial \Delta \cap \partial P$ is not connected, violating the assumption on $\Delta$.
	
	\item[(2)] This follows from (1).
	
	We follow the notation of Definition 3.2.3 with $P$ in place of $M_0$.
	\item[(3)] Let $B=B_1 \cup B_2 \cup B_3$, $B_i$  subbands, $ i =1,2,3$. Let $\theta = \partial B_2 \cap \partial P$ and let $v_1 \theta \gamma$ be a boundary cycle of $P$. Then by ({\rm V}) \  $B^I = (B \setminus B_2) \cup B_2^I$. Let $\eta$ be the complement of $\theta$ in $\partial \Delta$.
	Then $B_2^I \cap \partial \Delta = \eta$ and $\partial P^I = v_1 \eta v_2 \gamma$, $v_1, v_2$ vertices. Now suppose by way of contradiction that $\Phi (\eta v_2 \gamma)$ is not cyclically reduced. Since $\partial \Delta$ is cyclically reduced and $\partial P$ is cyclically reduced.
	Then there is a cancellation around $v_1$ and/or $v_2$, which extends $\theta$ as a common piece. 
	But by the definition of standard $I$-move (3.2.5), $\theta$ is maximal in the sense that it cannot be extended along $\partial B$. Hence, this is not possible. 
	Hence $\Phi(\eta v_2\gamma)$ is cyclically reduced.
	\hfill $\Box$
\end{itemize}

\section{Abelian diagrams}
\setcounter{subsection}{1}
\begin{definition}
	\addcontentsline{toc}{subsubsection}{\numberline{\thesubsubsection}  
		Definition (a-abelian diagrams, permutational diagrams)}
	Assume $M\in\mathcal{M}_3(W)$. Let  $M$ be an $\cal R$-diagram over $F$, $\cal R$ given by ({\rm III}) and let $ a \in X$.
	\begin{itemize}
		\item[(a)] Say that $M$ is \emph{a-abelian} if for every region $D$ in $M$ with $ a \in Supp(D)$ $D \in Reg_2(M)$.
		\item[(b)] Say that $M$ is \emph{abelian} if it is a-abelian for every $ a \in \cup Supp (D)$ $D \in Reg(M)$. Thus $M$ is abelian if and only if $Reg(M) = Reg_2(M)$. 
		Hence for every $D\in Reg(M)$ with $Supp(D)=\{a,b\}\subseteq X$ the subdiagrams generated by $[D]_s,\ s\in\{a,b\}$ are bands.
		\item[(c)] Let $M$ be an abelian diagram with boundary cycle $\omega$. Let $ \omega = \omega_1 v_1 \omega_2^{-1} v_2$ be a decomposition of $\omega$, $v_1, v_2$ boundary vertices, 
		$\omega_1$ and $\omega_2$ boundary paths. Say that $M$ is \emph{$(\omega_1, \omega_2)$-permutational} if every band that starts on $\omega_i$, ends on $\omega_{i+1}$ (indices taken modulo 2).
	\end{itemize}
\end{definition}

\vspace{10pt}

\begin{lemma}
	\addcontentsline{toc}{subsubsection}{\numberline{\thesubsubsection}  
		Lemma (Properties of permutational and banded diagrams)} \ 
	\begin{itemize}
		\item[(a)] Let $M$ be a connected, simply connected abelian $\cal R$-diagram with connected interior and boundary cycle $\omega$ with decomposition $ \omega = \omega_1 v_1 \omega_2^{-1} v_2$, $v_1, v_2$ vertices, $|M| \geq 2$.
		If $M$ is a $(\omega_1, \omega_2)$-permutational diagram, then there are sequences $J_1$ and $J_2$ of I-moves on $M$ such that $N_{M_1}(\omega_1) \cap {\cal D}_2(M_1) \neq \emptyset$ and $N_{M_2}(\omega_2) \cap {\cal D}_2(M_2) \neq \emptyset$ where $M_1 = M^{J_1}$ and $M_2= M^{J_2}$.
		\item[(b)] If $P$ is a closed 2-banded subdiagram of a simply connected abelian diagram, then $P$ is transferable. 
		Moreover, if $P$ is banded by $B_l$ and $B_r$ and $|Reg_2(M)|$ is minimal among all the diagrams with boundary label $\Phi(\partial M)$, 
		then $B_r$ and $B_l$ cannot be extended to bands $\widehat{B_r}$ and $\widehat{B_l}$ respectively, such that $|\widehat{B}_l \cap \widehat{B}_r | \geq 2$ where $\widehat{B}_r = B_r \cup D_1 \cup D_2$ \ \ and $\widehat{B}_l = B_l \cup D_1 \cup D_2$. 
		See Figure 33 (a).
	\end{itemize}
\end{lemma}

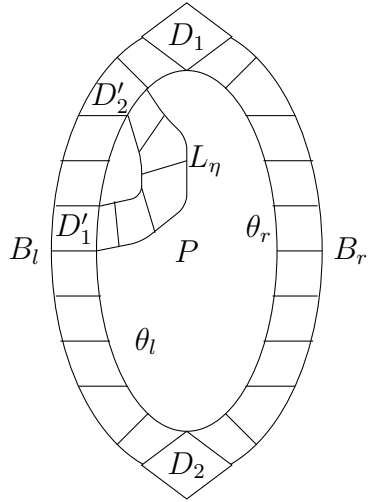
\begin{center}
	\begin{tikzpicture}[scale=.6]
		\draw (0,0) ellipse (2 and 4);
		\draw (0,0) ellipse (3 and 5);
		
		\draw[fill=white] (0,4)--(-1,4.75)--(0,5.5)--(1,4.75)--(0,4);
		\draw[fill=white] (0,-4)--(-1,-4.75)--(0,-5.5)--(1,-4.75)--(0,-4);
		
		\draw (-3,0)--(-2,0);
		\draw (-2.9,-1)--(-1.9,-1);
		\draw (-2.8,-2)--(-1.7,-2);
		\draw (-2.4,-3)--(-1.3,-3);
		\draw (-2.9,1)--(-1.9,1);
		\draw (-2.8,2)--(-1.7,2);
		\draw (-2.4,3)--(-1.3,3);
		
		\draw (-.85,3.6)--(-1.55,4.3);
		\draw (-.85,-3.6)--(-1.55,-4.3);
		\draw (.85,3.6)--(1.55,4.3);
		\draw (.85,-3.6)--(1.55,-4.3);
		
		\draw[rounded corners] (-2,0)--(-1,.2)--(0,1)--(0,2.5)--(-.5,3)--(-.9,3.6);
		\draw[rounded corners] (-1.9,1)--(-1,1.2)--(-1,2)--(-1.3,3);
		\draw (-1.5,0.1)--(-1.6,1.1);
		\draw (-.7,.4)--(-1,1.4);
		\draw (-1,1.7)--(0,2);
		\draw (-1.05,2.2)--(-.5,3);
		
		\node at (-1.7,3.4) {$D'_2$};
		\node at (-2.5,.5) {$D'_1$};
		\node at (0.4,2) {$L_{\eta}$};
		\node at (1.6,0.5) {$\theta_r$};
		
		\draw (3,0)--(2,0);
		\draw (2.9,-1)--(1.9,-1);
		\draw (2.8,-2)--(1.7,-2);
		\draw (2.4,-3)--(1.3,-3);
		\draw (2.9,1)--(1.9,1);
		\draw (2.8,2)--(1.7,2);
		\draw (2.4,3)--(1.3,3);
		
		\node at (0,0) {$P$};
		\node at (0,-4.75) {$D_2$};
		\node at (0,4.75) {$D_1$};
		\node[left] at (-3,0) {$B_l$};
		\node[right] at (3,0) {$B_r$};
		\node at (-.9,-2) {$\theta_l$};
	\end{tikzpicture}
	
	Figure 33 (a)
	
	$D_1$ and $D_2$ extend $B_l$ and $B_r$
\end{center}

{\it
	\begin{itemize}
		\item[(c)] Let $M$ be a 2-banded abelian diagram $M=(Q_lB_l \theta_lP\theta_r B_rQ_r)$ given by ({\rm VI}). Assume that $\gamma(P)$ is a boundary path of a band $B_m$ which crosses $B_l$ and $B_r$ in regions $D_1$ and $D_2$, respectively. If 
		$ Int(P) \neq \emptyset$ then $P$ is transferable. See Figure 33(b).
\end{itemize}}

\begin{center}
	\begin{tikzpicture}[scale=.65]
		\draw (0,0)--(2,4)--(-2,4)--(0,0);
		\draw (0,0)--(1,-1)--(3,5)--(-3,5)--(-1,-1)--(0,0);
		
		\draw (-1,5)--(-1,4);
		\draw (0,5)--(0,4);
		\draw (1,5)--(1,4);
		\draw (-2.55,3.6)--(-2,4)--(-2,5);
		\draw (2.55,3.6)--(2,4)--(2,5);
		\draw (-2.2,2.6)--(-1.5,3);
		\draw (-1.9,1.6)--(-1,2);
		\draw (-1.5,.6)--(-.5,1);
		\draw (2.2,2.6)--(1.5,3);
		\draw (1.9,1.6)--(1,2);
		\draw (1.5,.6)--(.5,1);
		
		\draw[rotate around={163:(1,-1)}] (1,-1) arc(90:268:3.1);
		\draw (-1,-1) arc(270:127:3.3);

		\node[below] at (1,4) {$\gamma'$};
		\node[above] at (1,5) {$\gamma$};
		\node[below] at (0,4) {$P$};
		\node at (-3,3) {$Q_l$};
		\node at (3,3) {$Q_r$};
		\node at (-1.6,2.4) {$B_l$};
		\node at (1.6,2.4) {$B_r$};
		\node at (-.4,1.5) {$\theta_l$};
		\node at (.45,1.5) {$\theta_r$};
		\node at (-2.4,4.5) {$D_1$};
		\node at (2.4,4.5) {$D_2$};
		\node at (-0.5,4.4) {$B_m$};
	\end{tikzpicture}
	
	Figure 33 (b)
\end{center}

\begin{center}
	\begin{tikzpicture}[scale=.65]
		\draw (0,0)--(2,4);
		\draw (-2,4)--(0,0);
		\draw (0,0)--(1,-1)--(3,5)--(2,5);
		\draw (-2,5)--(-3,5)--(-1,-1)--(0,0);
		
		\draw (-2.55,3.6)--(-2,4)--(-2,5);
		\draw (2.55,3.6)--(2,4)--(2,5);
		\draw (-2.2,2.6)--(-1.5,3);
		\draw (-1.9,1.6)--(-1,2);
		\draw (-1.5,.6)--(-.5,1);
		\draw (2.2,2.6)--(1.5,3);
		\draw (1.9,1.6)--(1,2);
		\draw (1.5,.6)--(.5,1);
		
		\draw[rotate around={163:(1,-1)}] (1,-1) arc(90:268:3.1);
		\draw (-1,-1) arc(270:127:3.3);
		
		\draw (-2,5)--(2,5);
		\draw (-.3,5.2)--(0,5)--(-.3,4.8);
		
		\node[above] at (0,5.2) {$\gamma$};
		\node[below] at (0,4) {$P$};
		\node at (-3,3) {$Q_l$};
		\node at (3,3) {$Q_r$};
		\node at (-1.6,2.4) {$B_l$};
		\node at (1.6,2.4) {$B_r$};
		\node at (-.4,1.5) {$\theta_l$};
		\node at (.45,1.5) {$\theta_r$};
		\node at (-2.4,4.5) {$D_1$};
		\node at (2.4,4.5) {$D_2$};
		\node at (-2.4,5.3) {$\alpha_l$};
		\node at (2.4,5.3) {$\alpha_r$};
		\node at (-.4,-1) {$\beta_l$};
		\node at (.45,-1) {$\beta_r$};
		\node at (-4.5,1) {$\omega_1$};
		\node at (5.4,2.5) {$\omega_3$};
		\draw (-.2,-.3)--(-.5,-.5)--(-.4,-.2);
		\draw (.5,-.7)--(.5,-.5)--(.7,-.5);
		
		\draw (-4.2,.7)--(-4,1)--(-3.7,.7);
		\draw (4.8,2.3)--(5,2)--(5.2,2.3);
	\end{tikzpicture}
	
	Figure 33 (c)
\end{center}

\vspace{10pt}

\begin{itemize}
	\item[(d)] 
	{\it Let $M$ be a 2-banded abelian diagram given by ({\rm VI}), $M=(Q_lB_l \theta_lP\theta_rB_rQ_r)$. Then there is a sequence $J$ of $I$-moves such that either $P^J= \emptyset$ or $P^J$ is $(\gamma, \mu)$-permutational, where 
		$ \mu = \theta_l\theta_r^{-1}$ and where $\gamma = \gamma(P)$ is the connected component of $\partial P \cap \partial M$ which contains $t(\theta_r)$ and $t(\theta_l)$.} See Figure 33(c).
\end{itemize}

\vspace{10pt}

\noindent {\bf Proof} \\
By induction on $|M|$, the case $|M|=2$ being clear. Assume first that $Int(M)$ is connected.

\vspace{10pt}
\begin{itemize}
	\item[(a)] Since $M$ is abelian, every edge on $\partial M$ is connected by a band to another edge on $\partial M$. Since $M$ is
	$(\omega_1, \omega_2)$-permutational, every band that starts on $\omega_1$, ends on $\omega_2$ and every band that starts on $\omega_2$, 
	ends on $\omega_1$. Let $\omega_1 = \eta_1 \cdots \eta_n$, $\omega_2 = \nu_1 \cdots \nu_n$ $\eta_i$ and $\nu_i$ edges. Every band $L$ in $M$ 
	which connects $\omega_1$ with $\omega_2$
	defines a pair of numbers $(i,j)$, where if $L$ starts on $\eta_i$, and it ends at $\nu_j$. This pair is uniquely defined by $L$ since $M$ is planar. 
	Hence the totality of the connecting bands define a permutation $\sigma$ on $\{1, \ldots , n \}$. We claim that there is an $i$, $ 1 \leq i \leq n-1$ such that $\sigma(i) > \sigma (i+1)$. Suppose not. Then for every $i$, $i = 1, \ldots , n-1$, $\sigma(i) \leq \sigma (i+1)$. Since $\sigma$ is a permutation, hence $\sigma(i) < \sigma (i+1)$ $i = 1, \ldots , n-1$, i.e. $\sigma(1) < \sigma(2) \cdots < \sigma(n-1) < \sigma(n)$.  
	Consequently $\sigma(1)=1$. Hence, there is a band $L$ in $M$ such that $L$ starts on $\eta_1$ and ends on $\nu_1$. But then $u$ is a double point of $\partial L$, a contradiction to $M\in\mathcal{M}_3(W)$. See Figure 34.
\end{itemize}

\begin{center}
	
	\begin{tikzpicture}[scale=1.2]
		
		\draw (-1.5,0) ellipse (3.5 and 1.5);
		\draw (-3,0) ellipse (2 and 1);
		
		\draw[rounded corners] (-5,0)--(-4,.5)--(-2,.5)--(-1.5,0)--(-2,-.5)--(-4,-.5)--(-5,0);
		
		\draw[fill] (-5,0) circle(.05);
		\draw[fill] (2,0) circle(.05);
		
		\draw (-4,.5)--(-4,.85);
		\draw (-3,.5)--(-3,1);
		\draw (-2.3,.5)--(-2.3,.95);
		\draw (-1.75,.25)--(-1.4,.6);
		
		\draw (-4,-.5)--(-4,-.85);
		\draw (-3,-.5)--(-3,-1);
		\draw (-2.3,-.5)--(-2.3,-.95);
		\draw (-1.75,-.25)--(-1.4,-.6);
		
		\draw (-1,0)--(-1.55,0);
		
		\draw (-4.4,.3)--(-4.65,.55);
		\draw (-4.4,-.3)--(-4.65,-.55);
		\draw (-1.8,-1.3)--(-1.5,-1.5)--(-1.8,-1.7);
		
		\node at (-2,-.65) {$L$};
		\node[above] at (-1.5,1.5) {$\mu_1$};
		\node[below] at (-1.2,-1.5) {$\nu_1$};
		\node[left] at (-5,0) {$u$};
	\end{tikzpicture}
	
	Figure 34
	
\end{center}

Hence, there is an $i$, $ 1 \leq i \leq n-1$, such that $\sigma(i) > \sigma(i+1)$. Then if we denote by $L_i$ and $L_{i+1}$ the bands that emanate from $\eta_i$ and $\eta_{i+1}$ respectively, then they cross. See Figure 35(a).

\vspace{10pt}

\begin{minipage}{1.9in}
	\begin{center}
		\begin{tikzpicture}[scale=.2]
			
			\draw[rounded corners] (0,0)--(5,5)--(5,10)--(0,15)--(-8,20);
			\draw[rounded corners] (0,0)--(-5,5)--(-5,10)--(0,15)--(8,20);
			
			\draw[rounded corners] (0,0)--(10,0)--(10,10)--(4,18)--(-5,24);
			\draw[rounded corners] (0,0)--(-10,0)--(-10,10)--(-4,18)--(5,24);
			
			\draw[red, thick] (0,15.4)--(3.6,18)--(0,20.3)--(-3.6,18)--(0,15.4);
			
			\draw[fill] (0,0) circle (.2);
			
			\node[below] at (-3,0) {$L_i$};
			\node[below] at (3,0) {$L_{i+1}$};
			\node at (0,8) {$P$};
			\node at (0,18) {$D$};
			
			\draw (-3,3)--(-10,3);
			\draw (-5,6)--(-10,6);
			\draw (-5,9)--(-10,9);
			\draw (-3,12)--(-7.4,13.5);
			\draw (-1,14)--(-5.5,16);
			
			\draw (3,3)--(10,3);
			\draw (5,6)--(10,6);
			\draw (5,9)--(10,9);
			\draw (3,12)--(7.4,13.5);
			\draw (1,14)--(5.5,16);
			\draw (-13,0)--(13,0);
			\draw (-12.3,-.2)--(-12,0)--(-12.3,.2);
			\node[below] at (-12,0) {$\omega_1$};
			
		\end{tikzpicture}
		
		\vspace{1cm}
		
		(a)
	\end{center}
\end{minipage}$\begin{array}{c} \mbox{\textcolor{red} J}\\ \leadsto \end{array}$\begin{minipage}{1.4in}
	\begin{center}
		\begin{tikzpicture}[scale=.25]
			\draw (0,0)--(0,18)--(5,21)--(0,24)--(-5,21)--(0,18);
			\draw (0,0)--(5,0)--(5,21);
			\draw (0,3)--(5,5);
			\draw (0,6)--(5,8);
			\draw (0,9)--(5,11);
			\draw (0,12)--(5,14);
			\draw (0,15)--(5,17);
			
			\draw (0,0)--(-5,0)--(-5,21);
			\draw (0,3)--(-5,5);
			\draw (0,6)--(-5,8);
			\draw (0,9)--(-5,11);
			\draw (0,12)--(-5,14);
			\draw (0,15)--(-5,17);
			
			\draw (5,0)--(6,0);
			\draw (-5,0)--(-6,0);
			\draw[red, thick] (0,18.3)--(4.7,21)--(0,23.7)--(-4.7,21)--(0,18.3);
			
			\node[below] at (-3,0) {$L_i^J$};
			\node[below] at (3,0) {$L_{i+1}^J$};
			\node at (0,21) {$D$};
			
			\draw (-7,0)--(7,0);
			\draw (-6.3,-.2)--(-6,0)--(-6.3,.2);
			\node[below] at (-6,0) {$\omega_1$};
			\draw (5.7,-.2)--(6,0)--(5.7,.2);
			
		\end{tikzpicture}

		(b)
	\end{center}
\end{minipage}$ \leadsto $ \begin{minipage}{1.4in}
	\begin{center}
		\begin{tikzpicture}[scale=.25]
			
			\draw (0,0)--(0,-18)--(5,-21)--(0,-24)--(-5,-21)--(0,-18);
			\draw (0,0)--(5,-2)--(5,-21);
			\draw (0,-3)--(5,-5);
			\draw (0,-6)--(5,-8);
			\draw (0,-9)--(5,-11);
			\draw (0,-12)--(5,-14);
			\draw (0,-15)--(5,-17);
			
			\draw (0,0)--(-5,-2)--(-5,-21);
			\draw (0,-3)--(-5,-5);
			\draw (0,-6)--(-5,-8);
			\draw (0,-9)--(-5,-11);
			\draw (0,-12)--(-5,-14);
			\draw (0,-15)--(-5,-17);

			\draw[red, thick] (0,-18.3)--(4.7,-21)--(0,-23.7)--(-4.7,-21)--(0,-18.3);
			
			\node[above] at (-3,0) {$L_i^{J_1}$};
			\node[above] at (3,0) {$L_{i+1}^{J_1}$};
			\node at (0,-21) {$E$};
			
			\draw (-7,-21)--(-5,-21);
			\draw (-6.3,-21.2)--(-6,-21)--(-6.3,-20.8);
			\node[below] at (-6,-21) {$\omega_1$};
			\draw (7,-21)--(5,-21);
			\draw (5.7,-21.2)--(6,-21)--(5.7,-20.8);
		\end{tikzpicture}
		
		(c)
	\end{center}
\end{minipage}

\begin{center}
	Figure 35
\end{center}

Since $P$ is a closed 2-banded diagram,  by the induction assumption part (b)  we may transfer every region from $P$ in Figure 35(a) by a sequence $J$ of $I$-moves and by a further sequence $J_1$, 
of $I$-moves we can push the region $D$ in Figure 35(b)  to the region $E$ on Figure 35(c).

\noindent $E \in N_{M_1}(\omega_1) \cap {\cal D}_2(M_1)$ where $M_1 = M^{J_0J_1}$. By a similiar argument, 
by replacing $\omega_1$ with $\omega_2$, we get a sequence $J_2$ of $I$-moves, such that if $M_2 = M^{J_2}$ then  $N_{M_2}(\omega_2) \cap {\cal D}_2(M_2) \neq \emptyset$. 
If $int(M)$ is not connected then we may apply the above arguments to each connected component.

\noindent (b) \quad Assume first that $Int(M)$ is connected.
From every edge $\eta$ of $\theta_l$ there emanates a band $L_{\eta}$. 
We claim that  $L_{\eta}$ ends on $\theta_r$. Otherwise, there are regions $D'_1$ and $D'_2$ in $B_l$ which extend $L_{\eta}$ to $\widehat{L}_{\eta}= L_{\eta} \cup \{D'_1, D'_2 \}$ such that
$|\widehat{L}_{\eta} \cap \theta_1 | \geq 2$, violating the induction  assumption of part (b) of the Lemma. 
See Figure 36(a).
Hence $L_{\eta}$ ends on $\theta_r$ and it follows that $P$ is $(\theta_r, \theta_l)$ permutational. Consequently, by part (a), there is a sequence $J$ of $I$-moves (which do not change $|Reg_2(P)|$) such that $N(\theta_1) \cap {\cal D}_2(P_1) \neq \emptyset$, where $P_1 = P^J$
and $\theta_1 = \theta_l^J$. But then if $E \in N(\theta_1) \cap {\cal D}_2(P_1)$ then an $I$-move $J_1$ transfers $E$ from $P_1$. See Figure 36.

\begin{center}
	\begin{tikzpicture}[scale=.4]
		\draw (0,7)--(12,7)--(12,4)--(0,4)--(0,7);
		\draw (3,7)-- (3,4);
		\draw (6,7)-- (6,4);
		\draw (9,7)-- (9,4);
		\draw (3,4)--(6,0)--(9,4);
		
		\draw[fill] (6,4) circle (.15);
		
		\node at (14,5) {$\rightsquigarrow$};
		
		\node[left] at (0,5.5) {$B_l$};
		\node[below] at (6,4) {$E$};
		\node[right] at (8,2) {$P_1$};
		\node[above] at (14,5) {$J_1$};
		\node at (24,6) {$E^1$};
		\node at (26,2) {$P_1^{J_1}$};

		\draw (19,7)--(29,7)--(29,4)--(26,4)--(24,0)--(22,4)--(19,4)--(19,7);
		\draw (22,7)-- (22,4);
		\draw (26,7)-- (26,4);
		\draw (24,3)-- (24,0);
		\draw (22,7)--(24,3)--(26,7);
	\end{tikzpicture}
	
	Figure 36
\end{center}

Doing this repeatedly, using Proposition 3.2.6, we get $|P_1^{J_3}|=0$. Hence $P$ is transferable, showing part (b). We prove the ''Moreover" part. Now by part (b) we may assume that $|P|=0$, hence $B_l$ and $B_r$ are adjacent. If there are regions $D_1$ and $D_2$ in $Reg_2(M)$ which extend $B_l \cup B_r$ such that $|\widehat{B}_l \cap \widehat{B}_r| \geq 2$, $\widehat{B}_l = B_l \cup \{D_1, D_2\}$, $\widehat{B}_r = B_r \cup \{D_1, D_2\}$, then by a sequence $J_5$ of $I$-moves we can push $D_1$ to $D_2$ to be adjacent.  But then  $\left(D_1^{J_5}, D_2\right)$ is a cancelling pair violating the minimality assumption. See Figure 37.

\vspace{10pt}

\begin{minipage}{1.5in}
	
	\begin{tikzpicture}[scale=.25]
		
		\draw[red, thick] (0,0)--(2,-2)--(4,0)--(2,2)--(0,0);
		\draw[red, thick] (12,0)--(14,-2)--(16,0)--(14,2)--(12,0);
		\draw (4,0)--(12,0);
		\draw (5,-2)--(6,0)--(5,2);
		\draw (8,-2)--(8,2);
		\draw (11,-2)--(10,0)--(11,2);
		\draw (2,2)--(14,2);
		\draw (2,-2)--(14,-2);
		
		\node at (2,0) {$D_1$};
		\node at (14,0) {$D_2$};
		
	\end{tikzpicture}
	
\end{minipage}$\begin{array}{c}  J_5\\ \leadsto \end{array}$\begin{minipage}{1.9in}
	\begin{center}
		\begin{tikzpicture}[scale=.35]
			
			\draw (2,0)--(10,0)--(12,-2)--(16,0)--(12,2)--(10,0);
			\draw (12,-2)--(13,0)--(12,2);
			\draw (5,-2)--(4,0)--(5,2);
			\draw (7,-2)--(6,0)--(7,2);
			\draw (9,-2)--(8,0)--(9,2);
			
			\draw (12,-2)--(2,-2)--(2,2)--(12,2);

			\node at (11.5,0) {$D_1^{J_5}$};
			\node at (13.8,0) {$D_2^{J_5}$};
			
		\end{tikzpicture}
	\end{center}
\end{minipage}$ \begin{array}{c}f.c.\\ \leadsto \end{array}$ \begin{minipage}{1in}
	
	\begin{tikzpicture}[scale=.25]
		
		\draw (0,0)--(8,0)--(8,4)--(0,4)--(0,0);
		\draw (0,2)--(8,2);
		\draw (2,0)--(2,4);
		\draw (4,0)--(4,4);
		\draw (6,0)--(6,4);

	\end{tikzpicture}
	
\end{minipage}

\vspace{10pt}

\hfill $f.c.$ is free cancellation
\begin{center}
	Figure 37
\end{center}

If $Int(M)$ is not connected then apply the arguments on each component.
\vspace{10pt}

\begin{itemize}
	\item[(c)] 
	Assume first that $Int(P)$ is connected. 
\end{itemize}

Due to Part (b) of the Proposition we may assume that

\vspace{10pt}

{\it No band in $P$ starts and ends on $B$, where $B \in \{B_l, B_r, B_m\}$} \hfill (**)

(We follow the notation of Figure 33(b).)

\vspace{10pt}

Next, assume there is a band $L$, $L \neq B_m$ which starts at $B_l$ and ends on $B_r$. 

Let $P_0$ be the subdiagram of $P$ which is bounded by $B_r$, $B_l$ and $L$. Then by part (c) of the induction hypothesis 
there is a sequence $J$ of $I$-moves such that the diagram $P_0^J$ contains no regions. Hence $L$ is adjacent to $B_r^J$ and $B_l^J$.
But then the vertex $z$ in Figure 38(a) has valency 3 and we may carry out a standard $I^2$-move $I_1$ around $z$. See Figure 38 (b) by which we transfer $E'$ and we get
$|P_0^{I_1}|=|P_0|-1$ and the result follows. Hence we may assume that every band that starts on $B_l \cup B_r$ ends on $B_m$ and every band that starts on $B_m$ ends on $B_l \cup B_r$. Consequently, $P_0$ is permutational
and by part (a) there is a region 
$D \in {\cal D}_2(M) \cap N(\gamma)$. Then we can transfer $D$ via $B_m$. Hence it follows from the induction hypothesis that $P$ is transferable. 
If $Int(M)$ is not connected, then it follows from part (b) that there is an $I$-move $J_4$ such that $Int(P_0^{J_4})$ is connected ($\theta_l$ and $\theta_r$ are connected).

\vspace{10pt}

\begin{minipage}{2.3in}
	\begin{center}
		\begin{tikzpicture}[scale=.6]
			\draw (2,0)--(4,0)--(6,3)--(6,7)--(0,7)--(0,3)--(2,0);
			\draw (2,3)--(2,7);
			\draw (2,0)--(4,3)--(4,7);
			\draw (0,3)--(6,3);
			
			\draw[fill] (2,3) circle(.1);
			
			\node at (2,1.3) {$E$};
			\node at (4,1) {$B_l$};
			\node at (2.2,3.3) {$z$};
			\node at (3,5) {$B_m$};
		\end{tikzpicture}

		(a)
\end{center}\end{minipage}$\begin{array}{c} I_1\\ \leadsto\end{array}$\begin{minipage}{2.3in}
	\begin{center}
		\begin{tikzpicture}[scale=.65]
			
			\draw (2,0)--(4,2)--(6,2)--(6,5)--(0,5)--(0,2)--(2,0);
			\draw (0,5)--(2,2)--(4,5)--(4,2);
			\draw (2,0)--(2,2);
			
			\node at (2,4) {$E'$};
			\node at (5,3) {$B_l$};

		\end{tikzpicture}

		(b)
	\end{center}
\end{minipage}

\begin{center}
	Figure 38
\end{center}

\noindent (d) \quad Assume first that $|P| \neq 0$ and $Int(P)$ is connected. See Figure 33(c).

If every band which starts on $\mu := \theta_l \theta_r^{-1}$ ends on $\gamma$ and every band that starts on $\gamma$ ends on $\mu$ then $P$ is already permutational. Hence we may assume that there is a band $L$ in $P$ which either starts on $\gamma$ and ends on $\gamma$ or starts on $\mu$ and ends on $\mu$. In the first case
assume that no band in $P$ starts and ends inside the sub-diagram $P_0$ bounded by $L$ and $\gamma$. Then $P_0$ is $(\mu_1, \nu_1)$-permutational, where $\mu_1 = \partial P_0 \cap \gamma$ and $\nu_1$ is the complement of $\mu_1$ on $\partial P_0$. 
Hence, by part (a) of the Lemma  $P$ contains a region $ D \in N(\gamma) \cap {\cal D}_2(M)$. 
Let $P_1= P \setminus \{D\}$. Then $P_1$ is a 2-banded abelian diagram with $|P_1| < |P|$, hence by the induction hypothesis there exist a sequence $J$ of $I$-moves such that either $P_1^J = \emptyset$ or $P_1^J$ is permutational. From this it easily follows that $M^J$ is $(\gamma, \mu)$-permutational.
Hence after carrying out a sequence of  $I$-moves $J_0$ we may assume that no band in $P^{J_0}$ starts on $\gamma$ and ends on $\gamma$. Notice that $|P^{J_0}| < |P|$. Assume therefore that $L$ starts on $\mu$ and ends on $\mu$. If $L$ starts on $\theta_l$ then due to 
the ``Moreover" part of Part (b) 
we may assume that $L$ ends on $\theta_r$. Then $P_L$  is removable by part (c) of the Lemma, 
where $P_L$ is the subdiagram bounded by $\theta_r$, $\theta_l$ and $L$. 
Hence as long as there are such bands we can remove regions from $P$. In other words, there is a sequence $J$ of $I$-moves on $P$ such that if $P^{J} \neq \emptyset$ then every band starting on $\mu$ ends on $\gamma$ and every band that starts on $\gamma$ ends on $\mu$, i.e. $P^J$ is $(\gamma, \mu)$-permutational.

\ \hfill $\Box$
\tocexclude
	{
		\section{Lifting $\mathcal{D}(\widetilde{\mathbb M}^t)$ to $\mathbb M_t,\ t\in T(M)$ and the proofs of the Main Results}
	}
	\addcontentsline{toc}{section}{\numberline{\thesection}  
		Lifting $\mathcal{D}(\widetilde{\mathbb M}^t)$}
	The idea behind showing that $\widetilde{\mathbb M}^t$ satisfies the condition C(4)\&T(4)
	is to use it in M. Obviously, we cannot lift the condition C(4)\&T(4) from $\widetilde{\mathbb M}^t$ to M,
	however we can lift a version of Greendlinger's Lemma. This is the content of theorem C. 
	The first step towards such a lifting is to establish the basic connection between $\widetilde{\mathbb M}^t$ and $\mathbb M_t$. We consider $\mathbb M_t$ as a map and consider $(\mathbb M_t)^1$, its 1-skeleton. 
	Following the definition of $\widetilde{\mathbb M}^t$ we consider $\left(\widetilde{\mathbb M}^t\right)^1$ as a  quotient of $\left({\mathbb M}_t\right)^1$. Then we extend this from the 1-skeleton to the whole map and show that
	$|Reg_{4^+,t}(\mathbb M_t)|=|Reg(\widetilde{\mathbb M}^t)|$.
	Since $\widetilde{\mathbb M}^t$ satisfies the condition C(4)\&T(4), this allows to estimate $|Reg_{4^+,t}(\mathbb M)|$
	\subsection {Preliminary Results}
	\tocexclude{
		\subsubsection {The mappings $\Psi_t$ and $\Psi_t^{-1},\quad t\in T(M)$}
	}
	\addcontentsline{toc}{subsubsection}{\numberline{\thesubsubsection}
		The mappings $\Psi_t$ and $\Psi_t^{-1}$}
	
	Let $M$ be a connected, simply connected R-diagram, $M\in \mathcal M_3(W)$ 
	and $t\in T(M)$.
	Consider the diagram $\mathbb M_t$ (Defined in 2.1).  Assume first that it has connected interior. 
	Denote by $(\mathbb M_t)^1$  the 1-skeleton  of $\mathbb M_t$.  
	Define an equivalence relation "$\apx 1$" on the vertices of $(\mathbb M_t)^1$ to be the transitive  closure of the relation "$\Sim 1$" defined as follows: vertices $v_1$ and $v_2$ satisfy $v_1 \Sim1 v_2$ if either $v_1=v_2$, or there is an edge e labelled  by a letter from $H_t$ ($H_t$-edge) with $\iota(e)=v_1$ and $\tau (e)=v_2$. 
	Here $\iota (e)$ denotes the initial vertex of $e $ and 
	$\tau(e )$ denotes the terminal vertex of $e $.
	Edges $e_1$ and $e_2$ satisfy $e_1\Sim 1 e_2$ if $e_1=e_2$.
	Then by the definition of $\mathbb M^t$, $(\mathbb M_t)^1\Big/ \apx1$ is isomorphic to $\mathbb (\mathbb M^t)^1$.
	
	\noindent Let $E_t$ be the collection of all the edges of 
	${\mathbb M}_t$ which are labelled by $t^{\pm 1}$ and define $E_{H_t}$ accordingly.
	
	\noindent Denote the projection mapping of $(\mathbb M_t)^1$ 
	which realises $\apx 1$ by \\ 
	$\phi_t:(\mathbb M_t)^1 \longrightarrow (\mathbb M ^t)^1$. Then $\phi_t$ sends a vertex $v$ in $(\mathbb M_t)^1$ to the vertex $\phi_t(v)$ which we denote by $v^t$.
	If $e\in E_t\left( (\mathbb M_t)^1\right)$ then $\phi_t(e)$ is an edge which we denote by $e^t$. 
	If $e\in E_{H_t}\left( (\mathbb M_t)^1\right)$ then
	$e^t:=\phi_t(e)$ is a loop with $\iota(e)=\tau(e)$, which we remove.
	We can do this because there is no interrelation between a loop and an edge (or loop).
	So we define
	$\phi_t(e)=\emptyset$ and $\phi_t(\bar e)=v_1^t=v_2^t$ where $\iota(e)=v_1$ and $\tau (e)=v_2$. 
	Notice that this does not invalidate Lemma 5.1.2 below.
	We define $\phi_t^{-1}:\left( \mathbb M^t \right)^1\longrightarrow \left( \mathbb M_t \right)^1$ by its set theoretical definition as follows:
	
	\noindent Let $v^t$ be a vertex in $\mathbb M^t, v^t=\phi_t(v), v$ a vertex of $\mathbb M_t$. 
	\begin{itemize}
		\item [a)]  Define  
		$\phi_t^{-1}({v^t})=\{w| w\,\text{is}\,H_t\text{-connected to v}\}$.
		\item[b)] if $e^t$ is an edge in \mt1 then $\phi_t^{-1}(e^t)$ is the unique edge $e\in E_t\left(({ \mathbb M}^t )^1\right)$ such that $\phi_t(e)=e^t$.
		(Notice that $e$ is necessarily in $E_t\left(\left({ \mathbb M}_t \right)^1\right)$.)
	\end{itemize}
	If $Int(\mathbb M_t)$ is not connected then consider each connected component, individually.
	
	\ \\ \ \\
	\noindent Consider now $\widetilde{\mathbb M}^t$. Assume first that it has connected interior. Define an equivalence relation "$\apx2$" on $(\mathbb M^t)^1$ to be the transitive closure of the relation "$\Sim2$" defined as follows:
	\begin{itemize}
		\item [1)] vertices $v_1^t$ and $v_2^t$ of \mt1 satisfy $v_1 \Sim2 v_2$ if and only if $v_1=v_2$.
		\item[2)] edges $e_1^t$ and $e_2^t$ satisfy $e_1^t\Sim2 e_2^t$ if and only if $\tau(e_1^t)=\tau(e_2^t)$ and $\iota(e_1^t)=\iota (e_2^t)$.  
		
	\end{itemize}
	
	\noindent It follows from the definition of $\widetilde{\mathbb M}^t$ in the introduction that \mt1 $\Big/ "\apx2"$ is isomorphic to $\left(\widetilde{ \mathbb M}^t \right)^1$.
	
	Denote by $\xi_t$ the projection from $\left({ \mathbb M}^t \right)^1$ onto $\left(\widetilde{ \mathbb M}^t \right)^1$ which realises $\apx 2$. Then we have
	\begin{itemize}
		\item [1.] if $v^t \in V\left(\left({ \mathbb M}^t \right)^1\right)$ then $\xi_t(v^t)=\widetilde{v},\  \widetilde{v}$ a uniquely defined vertex in $\widetilde{ \mathbb M^t}^1$.
		\item [2.] if $e^t \in E\left(\left({ \mathbb M}^t \right)^1\right)$ then $\xi_t(e^t)=\widetilde{e},\  \widetilde{e}$ a uniquely defined edge in $\widetilde{ \mathbb M^t}^1$.
	\end{itemize}
	\noindent We define $\xi_t^{-1}:\widetilde{ \mathbb M}^t \longrightarrow { \mathbb M}^t$ as follows:
	\begin{itemize}
		\item [1.] if $\widetilde{v} \in V\left(\left(\widetilde{ \mathbb M^t} \right)^1\right)$ then $\xi_t^{-1}(\widetilde{v})=v^t$, where $v^t$ is the unique 
		vertex in $\left(\mathbb M^t\right)^1$ with $\xi_t(v^t)=\widetilde{v}$.
		\item [2.] if $\widetilde{e} \in E\left(\left(\widetilde{ \mathbb M^t} \right)^1\right)$ then $\xi_t^{-1}(\widetilde{e})
		=\{e^t\in\mathbb{M}^t|\iota(e^t)=\iota(e_0^t)\text{\ and\ }\tau(e^t)=\tau(e_0^t) \}$, where $\xi_t(e_0^t)=\widetilde{e}$.
	\end{itemize}
	Now, by definition, $\xi_t$ does not move vertices.
	In particular, it does not identify vertices. Hence,
	\begin{equation*}\tag*{\it(1)}
		\hbox{
			If $\gamma^t$ is a simple closed curve then 
			$\xi_t(\gamma^t)$ is a simple closed curve
		}
	\end{equation*}
	Let $\mu^t$ be a simple closed curve  in $({ \mathbb M}^t)^1$.
	
	\noindent Define: 
	\begin{align*}
		RN(\mu^t) & =\{e\in N(\mu)\big| e\text{ an edge}, e\text{ is to the right of }\mu^t\},\\
		LN(\mu^t) & =\{e\in N(\mu)\big| e\text{ an edge}, e\text{ is to the left of }\mu^t\}
	\end{align*}
	Then since $\xi_t$ does not move vertices, 
	\begin{equation*}\tag*{\it(2)}
		\xi(RN(\mu^t))=RN(\xi_t(\mu^t)) \text{ and } 
		\xi(LN(\mu^t))=LN(\xi_t(\mu^t))
	\end{equation*}
	Since $\widetilde{ \mathbb M}^t$ is planar, it follows from (2) that
	\begin{equation*}\tag*{\it(3)}
		\left.
		\begin{minipage}{9cm}
			If $\widetilde{\gamma}^t$ is simple closed and minimal then
			$\xi_t^{-1}(\widetilde{\gamma}^t)$ contains a simple closed minimal curve $\gamma^t$.
		\end{minipage}
		\right\}
	\end{equation*}
	In other words, all the simple closed and minimal curves in $\left(\widetilde{ \mathbb M}^t \right)^1$ come from
	$\left({ \mathbb M}^t \right)^1$.
	Moreover, 
	\begin{equation*}\tag*{\it(4)}
		\left.
		\begin{minipage}{10cm}
			if $\gamma_1^t$ and $\gamma_2^t$ are distinct simple 	minimal closed curves with $\|\gamma_1^t\|\geq 4$
			and $\|\gamma_2^t\|\geq 4$, then
			$\xi_t(\gamma_1^t)\neq  \xi_t(\gamma_2^t)$.
		\end{minipage}
		\right\}
	\end{equation*}

	\noindent If $\mathbb{\widetilde{M}}^t$ has non-connected interior, then consider each component separately.
	
	\noindent Finally, define $\Psi_t: (\mathbb M_t)^1 \rightarrow \left(\widetilde{ \mathbb M}^t \right)^1$ by $\Psi_t=\xi_t\circ\phi_t$. 
	Then  $\Psi^{-1}_t=\phi^{-1}_t\circ\xi_t^{-1}$.
	It follows by checking Euler characteristics that 
	$\widetilde{\mathbb M}^t$ is planar.
	We have the following
	\par\ \\ \ \\
	
	\begin{lemma}\ 
		\addcontentsline{toc}{subsubsection}{\numberline{\thesubsubsection}  
			Lemma: Basic Properties of $\left(\widetilde{M}^t \right)^1$, $\left(\widetilde{ \mathbb M}^t \right)^1$ and $\left({ \mathbb M}_t \right)^1$
		}
		\begin{itemize}
			\item [1)] Let $e\in E_t\left(\left({ \mathbb M}_t \right)^1\right)$
			and let $B^1_e$ be the 1-skeleton of the t-band which contains $e$. Then $\Psi_t^{-1}\circ\Psi_t(\bar e)=B^1_e$. In particular $\bar e \subseteq \Psi_t^{-1}\circ\Psi_t(\bar e)$
			\item [2)] Let ${e}\in E_t\left(\left({ \mathbb M}_t \right)^1\right)$. Then $\Psi_t^{-1}\circ\Psi_t(e)={e}$
			\item [3)] Let $v\in V\left(\left({ \mathbb M}_t \right)^1\right)$. Then $\Psi_t^{-1}\circ\Psi_t(v)=[v]_{H_t}$. In particular $v\in \Psi_t^{-1}\circ\Psi_t(v)$
			\item [4)] Let $\widetilde{v}\in V\left(\left(\widetilde{\mathbb M}^t\right)^1\right)$. Then $\Psi_t\circ\Psi_t^{-1}(\widetilde{v})=\widetilde{v}$. 
			\item [5)] $\Psi_t(\iota(e))=\iota(\Psi_t(e))$ and $\Psi_t(\tau(e))=\tau(\Psi_t(e))$ for every edge $e\in V_t\left(\left({ \mathbb M}_t \right)^1\right)$
		\end{itemize}
	\end{lemma}
	
	\noindent {\bf Proof\quad} Immediate by the definitions of $\Psi_t$ and  $\Psi^{-1}_t$ \hfill $\Box$ 
	
	\begin{proposition}
		\addcontentsline{toc}{subsubsection}{\numberline{\thesubsubsection} Proposition: $|Reg_{4^+,t}(\mathbb M_t)|=|Reg(\widetilde{\mathbb M}_t)|$
		}
		\ \\ Let notation and assumptions be as above. Then
		$|Reg_{4^+,t}(\mathbb M_t)|= |Reg(\widetilde{\mathbb M}^t)|$.
	\end{proposition}
	
	\noindent{\bf Proof}\\
	We extend $\Psi_t$ from $\left({ \mathbb M}_t \right)^1$ to $\mathbb M_t$ and extend $\Psi^{-1}_t$ from $\left(\widetilde{ \mathbb M}^t \right)^1$ to $\widetilde{\mathbb M^t}$. 
	First extend $\phi_t$ and $\phi^{-1}_t$.
	
	\noindent If $\gamma$ is a closed curve in $({ \mathbb M}_t)^1$ then $\phi_t(\gamma)$ is a closed curve in $\left({ \mathbb M}^t \right)^1$. 
	If $\gamma$ is a minimal simple closed curve in $({ \mathbb M}_t)^1$ then 
	$\gamma$ is the boundary of a region $K$ of ${ \mathbb M}_t$, with a disc attached to $\gamma$. 
	Also, if $\gamma$ is the boundary of a region $K$ of $\mathbb{M}_t$ then $\gamma$ is simple due to Proposition 2.1.8 and it is minimal. 
	Let $\gamma$ be a simple closed minimal curve in $\left({ \mathbb M}_t \right)^1,\quad \gamma=\partial K, K\in Reg(\mathbb M_t)$ 
	with $\Phi(\gamma)\notin H_t$.
	We claim that $\phi_t(\gamma)$ is simple closed.
	Assume $\phi_t(\gamma)$ 
	is not  a non-trivial simple closed curve. 
	Then there is a path $\mu$ in $\mathbb M_t$ which connects distinct vertices on $\gamma,\Phi(\mu)\in H_t$.
	If $K\in Reg_{4^+,t}(\mathbb{M}_t)$ then due to Lemma 2.2.2, $\mu$ connects adjacent vertices $v_1$ and $v_2$ such that $v_1ev_2$ is the closure of an edge of $\partial K$, 
	labelled either by $t^{\pm 1}$ or by $h\in X_t$. 
	The first option violates Lemma 1.1.1(c), since $t\notin H_t$. 
	In the second option $v_1$ and $v_2$ are adjacent and become the same vertex under $\phi_t$, hence $\phi_t(\gamma)$ remains simple closed. 
	If $K\in Reg_{B_t}(\mathbb M_t)$ then $\mu$ connects either two vertices on the same side of $K$ or it connects two vertices on different sides of $K$. In the first case the two endpoints of $\mu$ become the same under $\phi_t$ hence connecting them cannot make $\phi_t(\partial K)$ non-simple. 
	In the second case $||\Phi(\mu\delta)||_t=1$, violating Lemma 1.1.1(c). 
	Here $\delta$ is a simple boundary path of $B_t$ which connects the endpoints of $\mu$ on $\partial K$.\\
	Hence, If $\gamma$ is a minimal simple closed curve then $\phi_t(\gamma)$ is a simple closed curve in $\mathbb M^t$. 
	We claim that $\phi_t(\gamma)$ is minimal.
	Since $M$ and $M^t$ are diagrams of group presentations over $F$ and over $H_t*<t>$, respectively,
	$\gamma$ is the boundary of a region $D_\gamma$ and similarly, if $\gamma^t$ is a simple minimal closed curve, then
	$\gamma^t=\partial D_{\gamma^t}$.
	It follows from the definition of $\phi_t$, the way $\Delta=
	\Delta(D)$ is constructed and the construction of $M^t$ from $M$ that if $\Delta^t$ is the region of $\mathbb M^t$
	obtained from $\Delta\in Reg_2(\mathbb M_t)\cup Reg_{4^+,t}(\mathbb M_t)$ then $\phi_t(\partial\Delta)=\partial\Delta^t$.
	Hence if $\gamma=\partial\Delta$ then $\phi_t(\gamma)=\gamma^t,
	\; \gamma^t=\partial\Delta^t$.
	Hence,
	\begin{equation} \tag{$*$}
		\text {if } \gamma \text{ is simple, closed and minimal then so is } \phi_t(\gamma)
	\end{equation}
	
	\noindent Hence if $\gamma=\partial\Delta,\quad\Delta$ a region in $Reg_{4^+,t}\left({ \mathbb M}_t \right)\cup Reg_{B_t}(\mathbb M_t)$ 
	then due to ($*$) we may attach to $\phi_t(\gamma)$ a disc $\Delta^t$ such  that $\partial\Delta^t=\phi_t(\gamma)$. More precisely, 
	let $\gamma=v_1\alpha_1 w_1\beta_1 v_2\alpha_2 w_2\beta_2\ldots v_k\alpha_k w_k\beta_k$, 
	where $\Phi(\alpha_i)=t^{f_i}, \Phi(\beta_i)=a^{g_i},\quad f_i,g_i\neq 0$, $v_i=\iota(\alpha_i)=\tau(\beta_i),\ i=1,\ldots, k$. 
	Then
	$\phi_t(\gamma)=\widetilde{v_1}\widetilde{\alpha_1}\widetilde{v_2}\widetilde{\alpha_2}\ldots \widetilde{v_k}\widetilde{\alpha_k}$, where $\widetilde v_i=\phi_t(v_i)=\phi_t(w_i), \widetilde{\alpha_i}=\phi_t(\alpha_i)$.
	Consequently, $\Delta^t$ is a k-gon (while $\Delta$ is a 2k-gon)  which we attach to $\phi_t(\gamma)$. 
	This way we may extend the definition of $\phi_t$ to $\mathbb M_t$ by defining $\phi_t(\Delta)=\Delta^t$  and get that
	\begin{equation*} \tag{5}
		|Reg_{4^+,t}(\mathbb M_t)|=|Reg_{4^+}(\mathbb M^t)|
	\end{equation*}
	Now we extend $\xi_t:(\mathbb{M}^t)^1\longrightarrow(\widetilde{\mathbb{M}}^t)^1$ to $Reg(\mathbb M^t)$ and extend $\xi_t^{-1}$ from $Reg (\widetilde{\mathbb M}^t)$ to $Reg_{4^+,t} ({\mathbb M}^t)$.
	Let $\widetilde{\gamma}^t$ be a simple minimal closed curve in
	$\left(\widetilde{ \mathbb M}^t \right)^1$.
	By (3), $\widetilde{\gamma}^t=\xi_t(\gamma^t)$.
	But  $\gamma^t$ is a simple minimal closed curve,
	$\gamma^t=\phi_t(\gamma)$, hence we may attach $K_{\gamma^t}$ to
	$\gamma^t$ and since $\xi_t$ does not alter t-edges, we may 
	attach $K_{\gamma^t}$ to $\tilde\gamma^t$. 
	We define $\xi_t(K_{\gamma^t})=K_{\xi_t}(\gamma^t)$.
	Then $\xi_t$ is defined on the whole of $\mathbb M^t$.
	Hence, as for $\phi_t$, 
	we may extend $\xi_t$ on $\mathbb M^t$ by first attaching to $\widetilde{\gamma}^t$ a disc $\widetilde{\Delta}^t$ 
	such that $\widetilde{\gamma}^t=\partial\widetilde{\Delta}^t$ 
	and then define $\xi_t(\Delta^t)=\widetilde{\Delta^t}$ and $\xi_t^{-1} (\widetilde{\Delta}^t)=\Delta^t$. 
	As for $\phi_t$, we get, due to (4), that
	$\xi_t$ is one to one on $Reg_{4^+,t}(\mathbb M^t)$, onto $Reg (\widetilde{\mathbb M}^t)$. Thus 
	\begin{equation*} \tag{6}
		|Reg(\widetilde{\mathbb M}^t)|=|Reg_{4^+}(\mathbb M_t)|
	\end{equation*}
	It follows from (5) and (6) that since $\Psi_t=\xi_t\circ\phi_t$ 
	\begin{equation*} \tag{7}
		|Reg_{4^+,t}(\mathbb M_t)|=|Reg(\widetilde{\mathbb M}^t)|
	\end{equation*}
	\ \hfill $\Box$
	
	\tocexclude{
		\subsection{The subdiagrams $Opt(\Delta),Opt_r(\Delta)$ and $Opt_l(\Delta),\Delta=\Psi_t^{-1}(\widetilde{\Delta}),\widetilde{\Delta}\in \mathcal D(\widetilde{\omega}_i),\ i=1,2.$}
	}
	\addtocontents{toc}{\protect\contentsline{subsection}{\numberline{\thesubsection}The subdiagrams $Opt(\Delta),Opt_r(\Delta)$,$Opt_l(\Delta)$ and $\widehat{Opt}(\Delta)$}{\thepage}{}}
	
	Let $M\in \mathcal M_3(W)$ and assume $Int(M)$ is connected. In general this does not imply that $\widetilde{\mathbb M}^t$ has connected interior. In the rest of the work we aim to prove Theorems A,B and C. 
	It is enough to prove  these Theorems for every connected component of $Int(\widetilde {\mathbb M}^t)$. 
	Hence, without loss of generality we may work under the assumption that $\widetilde{\mathbb M}^t$ has connected interior, and we shall do so. Thus,\\ 
	{\em for the rest of the work we shall assume that $\widetilde{\mathbb M}^t$ has connected interior, unless said otherwise.}
	
	\begin{definition}[$Opt(\Delta),Opt_r(\Delta),Opt_l(\Delta),\widehat{Opt}(\Delta)$]
		\addcontentsline{toc}{subsubsection}{\numberline{\thesubsubsection} 
			Definition ($Opt(\Delta),Opt_r(\Delta)$,$Opt_l(\Delta)$ and $\widehat{Opt}(\Delta)$)}
		\ \\ Let $\widetilde{\Delta}\in \mathcal D(\widetilde{\mathbb{M}}^t)$ and let $\Delta=\Psi_t^{-1}(\widetilde{\Delta})$. 
		Let $\partial\widetilde{\Delta}\cap\partial\widetilde{\mathbb{M}}^t=\widetilde{v}_1\widetilde{e}_1\widetilde{v}_2\widetilde{e}_2\widetilde{v}_3\widetilde{e}_3\ldots \widetilde{v}_k,\quad k\geq 2$, $\widetilde{v}_i$ vertices, $\widetilde{e}_i$ edges. 
		Define $Opt(\Delta)=\Psi_t^{-1}(\widetilde{e}_1\widetilde{v}_2\ldots \widetilde{e}_{k-1})$, define $Opt_r(\Delta)=\Psi_t^{-1}(\widetilde{e}_1\widetilde{v}_2\ldots \widetilde{e}_{k-1}\widetilde{v}_k)$ and define $Opt_l(\Delta)=\Psi_t^{-1}(\widetilde{v}_1\widetilde{e}_1\widetilde{v}_2\ldots \widetilde{e}_{k-1})$. 
		Thus $Opt_r(\Delta)=Opt(\Delta)\cup\Psi^{-1}_t(\widetilde{v}_k)$ and $Opt_l(\Delta)=Opt(\Delta)\cup\Psi^{-1}_t(\widetilde{v}_1)$.
		Define $\widehat{Opt}(\Delta)= \Psi^{-1}_t(\widetilde{v}_1\widetilde{e}_1\ldots\widetilde{v}_k)$.
		It follows from the definition of $\psi_t^{-1}$ and Lemma 5.1.2 that 
		$Opt(\Delta)=
		\Psi^{-1}_t\left(\tilde{e}_1\right)\cup
		\Psi^{-1}_t\left(\tilde{v}_2\right)\cup
		\Psi^{-1}_t\left(\tilde{e}_2\right)\cup \ldots \cup 
		\Psi^{-1}_t\left(\tilde{e}_{k-1}\right)=
		B_1\cup Z_2\cup B_2\cup Z_3\ldots \cup B_{k-1}$ where 
		$B_i$ are t-bands emanating from $\partial\Delta\cap \partial Opt(\Delta)$ and $Z_i$ are connected unions of regions from $Reg_{H_t}(\mathbb M_t)$.
		
		\noindent Thus, $\Delta\cup Opt(\Delta)$,  $\Delta\cup Opt_r(\Delta), \Delta\cup Opt_l(\Delta)$ and  $\widehat{Opt}(\Delta)$
		are connected subdiagrams of $\mathbb M_t$. 
		Due to Proposition 2.1.7 $\partial M\cap\partial B_i,\ \partial Z_i\cap \partial\Delta$ and $\partial B_i\cap\partial\Delta$ 
		are connected and it is easy to see that 
		$\partial\Delta\cap \partial Opt(\Delta),\ \partial\Delta\cap \partial Opt_r(\Delta),\ \partial\Delta\cap \partial Opt_l(\Delta)$ 
		and $\partial\Delta\cap\partial\widehat{Opt}(\Delta)$
		are connected.
		
		\noindent We are interested in $Opt_r(\Delta) (Opt_l(\Delta))$ for the special case when $\widetilde{v}_k\ (\widetilde{v}_1)$ has valency 3 in $\widetilde{\mathbb M}^t$.
	\end{definition} 
	
	\begin{lemma}
		\addcontentsline{toc}{subsubsection}{\numberline{\thesubsubsection} 
			Lemma ($P$ and $Q$ are connected and simply connected with cyclically reduced boundary labels)}
		Let  $M\in\mathcal M_3(W)$. Assume that $\partial M $  decomposes by $\partial M =u\omega_1 w\omega_2^{-1}$, $u,w$ vertices and $t\in T(M)$.
		Let $ \widetilde\Delta_1^t$ and  $ \widetilde\Delta_2^t$ be boundary regions of $\widetilde{\mathbb{M}}^t$, $ \widetilde\Delta_1^t$ and  $ \widetilde\Delta_2^t$ in  
		$N(\tilde \omega_1), \Delta_i=\Psi^{-1}_t(\tilde \Delta^t_i),
		Supp(\Delta_i)=\{a,t\}, \widetilde{\omega}_1=\Psi_t(\omega_1)$.
		Assume that $\partial\widetilde{\Delta}_1\cap\partial\widetilde{\Delta}_2\cap\partial\widetilde{\omega}_1$ is a vertex $\widetilde{v}^t$ 
		with valency 3 in $\widetilde{\mathbb M}^t$. Let $Q=\Psi_t^{-1}(\widetilde{v})$ and let $P=\Psi_t^{-1}(N(\widetilde{v}))$. 
		Thus $P=\Psi_t^{-1}\left(\tilde{\Delta^t}_1 \cup \tilde{\Delta^t_2}\cup
		(\partial\tilde{\Delta^t}_1\cap \tilde{\partial\Delta^t}_2)\cup
		\widetilde{e}_k\cup\widetilde{f}_1\cup\widetilde{v}\right)$ 
		where $\widetilde e_k$ and $\widetilde{f}_1$ are the boundary paths of $\widetilde{\Delta}_1^t$ and $\widetilde{\Delta}_2^t$ with common endpoint $\widetilde{v}$. Then each of the following holds:
		\begin{itemize}
			\item [(a)] $Q$ is connected and simply connected with cyclically reduced boundary label such that each connected component $Q_i, i=1,\ldots,k$ of Q has cyclically reduced boundary labels $U_i$. Moreover $Q_i\in \mathcal{M}_3(U_i), i=1,\ldots, k$.
			
			\item[(b)] P is connected and simply connected with connected interior and cyclically reduced boundary label V.
			$P\in \mathcal{M}_3(V)$. 
			\item [(c)] $\partial Q=\theta_1 p_1 \alpha_1 p_2 \theta_2 p_3 \alpha_2 p_4\theta_3 p_5\gamma^{-1}$, 
			where $p_i$ are vertices, $\theta_i$ are sides of t-bands $i=1,2,3$ 
			and $\alpha_j=\partial\Delta_j\cap\partial Q, j=1,2, \Phi(\alpha_j)=a^{c_j},\ c_j\neq 0$ 
			and $\gamma=\partial Q\cap \omega_1$.
		\end{itemize} 
	\end{lemma}

	\noindent {\bf Proof of Lemma 5.2.2}
	\begin{itemize}
		\item [(a)] Follows from the definition of $\Psi_t^{-1}(\widetilde{v})$ and Proposition 2.1.7.
		We turn to the "moreover" part. As shown above, $U_i$ is cyclically reduced, $i=1,\ldots ,k$. 
		Since $Supp(\Delta_i)=\{a,t\}, i=1,2$ while $Supp(Q_i)\subseteq H_t$, the conditions of Lemma 2.2.1 are satisfied, 
		hence $Q_i\in\mathcal{M}_3(U_i), i=1,\ldots,k$.
		
		\item[(b)] By the definition of $\Psi_t^{-1}$ we have $P=\Delta_1\cup\Delta_2\cup Q\cup L_1\cup L_2\cup L_3$, 
		where $\Delta_i=\Psi_t^{-1}(\widetilde{\Delta}_i), i=1,2,\  L_1=\Psi_t^{-1}(\widetilde{e}_k), L_2=\Psi_t^{-1}(\partial\widetilde{\Delta}_1 \cap \partial\widetilde{\Delta}_2)$ and $L_3=\Psi_t^{-1}(\widetilde{f}_1)$. 
		$L_1,L_2$ and $L_3$  are t-bands emanating from $\partial\Delta_1$ or from $\partial\Delta_2$.
		It also follows from the definition of $\Psi_t^{-1}$ that each of the following holds:
		
		\begin{equation*}
			\tag{Int}
			\left.
			\begin{minipage}{10cm}
				$\partial\Delta_1\cap\partial L_1\neq\emptyset,\partial \Delta_1\cap\partial L_2\neq\emptyset, \partial\Delta_2\cap\partial L_2\neq\emptyset$,\\ $\partial \Delta_2\cap\partial L_3\neq\emptyset, \partial L_i\cap \partial Q\neq\emptyset,
				i=1,2,3$,\\ 
				$\partial\Delta_j\cap\partial Q\neq\emptyset, j=1,2$.
			\end{minipage}
			\right\}
		\end{equation*}
		
		\noindent Since $\Delta_i,L_j$ and $Q$ are regions in $\mathbb M_t$, Proposition 2.1.7 applies and hence all the intersections in (Int) are connected. Also, $\partial L_1\cap \omega_1$ and 	$\partial L_3\cap \omega_1$ are connected. We have 
		\begin{equation*}\tag{1}
			\left.
			\begin{minipage}{10cm}
				$\partial P=\mu\cup\gamma^{-1}$, where $\mu=\partial P\cap
				\partial (L_1\cup \Delta_1\cup\Delta_2\cup L_2\cup L_3)$ and $\gamma=\partial Q\cap \omega_1$ 
			\end{minipage}
			\right\}
		\end{equation*}
		See Fig. 39.
		\begin{center}
			\begin{tikzpicture}[scale=.4]
				
				\draw (0,0) circle (2);
				\draw (9,0) circle (2);
				
				\draw (2,0)--(7,0);
				\draw (1.7,-1)--(7.2,-1);
				\draw (5.5,0)--(5.5,-1);
				\draw (3.5,0)--(3.5,-1);
				\draw (4.5,0)--(4.5,-1);

				\draw (0,2)--(0,7)--(1,7)--(1,1.7);
				\draw (9,2)--(9,7)--(10,7)--(10,1.7);
				
				\draw (0,3)--(1,3);
				\draw (0,4.5)--(1,4.5);
				\draw (0,6)--(1,6);
				
				\draw (9,3)--(10,3);
				\draw (9,4.5)--(10,4.5);
				\draw (9,6)--(10,6);

				\draw (1,7)--(10,7);
				
				\node at (0,0) {$\Delta_1$};
				\node at (9,0) {$\Delta_2$};

				\node  at (4.5,-2) {$L_3$};
				\node[right] at (1.5,1.5) {$\alpha_1$};
				\node[left] at (7.5,1.5) {$\alpha_2$};
				
				\node[above] at (5,3.5) {$Q$};
				\node[above] at (0.5,7.5) {$L_1$};
				\node[above] at (9.5,7.5) {$L_2$};
				
				\node at (11,-3) {$P$};
				
			\end{tikzpicture}
			
			Figure 39
		\end{center}
		\noindent First we show that P is connected and simply connected.\\
		P is clearly connected. Also, $L_i\cap L_j=\Delta_i\cap \Delta_j=\Delta_i\cap L_j=\emptyset,\ i\neq j$. Hence if we define $L_i'=L_i\cup Q,\ \Delta_j'=\Delta_j\cup Q$, then
		\begin{equation*}\tag{2}
			L'_i\cap L'_j=Q,\ L'_i\cap\Delta'_j=Q\text{ and } \Delta'_i\cap \Delta'_j=Q
		\end{equation*}
		Since $\partial L_i\cap\partial Q$ and $\partial\Delta_i\cap\partial Q$ are simple and connected by Proposition 2.1.7 and $Q,L_i$ and $\Delta_j$ are simply connected, hence
		\begin{equation*}\tag {3}
			Q,L'_i \text{ and } \Delta_j'\text{ are simply connected.} \ i=1,2,3,\ j=1,2
		\end{equation*}
		Hence it follows from Seifert van Kampen's Theorem that
		$P=Q\cup L_1'\cup L_2'\cup L_3'\cup \Delta_1\cup \Delta_2$ is simply connected.
		Next, we show that $\partial P$ is a simple closed curve. 
		$\gamma$ is a simple curve, as $\gamma\subseteq\omega_1\subset\partial M$. 
		We claim that $\mu$ is a simple curve. 
		Since $\partial L_i$ and $\partial \Delta_j$ are simple
		and 
		$\partial L_1\cap\partial\Delta_1,\ 
		\partial L_3\cap\partial\Delta_1,\ 
		\partial L_2\cap\partial\Delta_2$
		are connected by  Proposition 2.1.7,
		if $\mu$ is not simple then one of the following holds:
		\begin{itemize}
			\item [i)] $\partial L_1\cap \partial L_3\neq \emptyset$
			\item [ii)] $\partial L_1\cap \partial \Delta_2\neq \emptyset$
			\item [iii)] $\partial L_1\cap \partial L_2\neq \emptyset$
			\item [iv)] $\partial \Delta_1\cap \partial \Delta_2\neq \emptyset$
		\end{itemize}
		
		Consider these cases in turn. See Fig 40.
	\end{itemize}
	\begin{center}
		\begin{tabular}{cc}
			\begin{tikzpicture} [scale=0.35]
				
				\begin{scope}[shift={(1,2.5)}]
					
					\draw (0.5,-0.5) node (circo) {} circle (1);
					\draw (9.5,-0.5) circle (1);
					
					\draw (1.35,0)--(8.55,0);
					\draw (1.35,-1) node (arcs) {}--(8.55,-1);   
					\draw (6.45,0)--(6.45,-1);
					\draw (3.5,0)--(3.5,-1);
					\draw (5.1,0)--(5.1,-1);

					\draw (9,0.35)--(9,7)--(10,7)--(10,0.35);

					\draw (9,2)--(10,2);
					\draw (9,3.5)--(10,3.5);
					\draw (9,5)--(10,5);

					\draw (-2,7)--(10,7);
					
					\node at (0.5,-0.5) {$\Delta_1$};
					\node at (9.5,-0.5) {$\Delta_2$};

					\node [above right] at (1.21,0.2) {$\alpha_1$};
					\node [above left] at (8.73,0.2) {$\alpha_2$};
					
					\node[above] at (5,3.5) {$Q$};
					\node[above] at (-2.5,7.5) {$L_1$};
					\node[above] at (9.5,7.5) {$L_2$};

				\end{scope}

				\draw (-1,9.5) node  (us) {} .. 
				controls (-3.5,6) and (-5.5,1.5) .. (-4,-1.5) node (ue) {}
				node[pos=.15 ] (u1) {} 
				node[pos=.3 ] (u2) {}
				node[pos=.45 ] (u3) {}
				node[pos=.6 ] (u4) {}
				
				node[pos=.75 ] (u5) {}
				node[pos=.87 ] (u6) {}
				;
				\draw (-2.5,9.5)  node (ls) {} .. 
				controls (-5,7) and (-7,1) .. (-5,-2.5) node (le) {}
				node[pos=.15 ] (l1) {} 
				node[pos=.3 ] (l2) {}
				node[pos=.45 ] (l3) {}
				node[pos=.6 ] (l4) {}
				
				node[pos=.75 ] (l5) {}
				node[pos=.87 ] (l6) {}
				;
				\draw  (us.center) edge (ls.center);	
				\draw (u1.center) -- (l1.center);
				\draw (u2.center) -- (l2.center);
				\draw (u3.center) -- (l3.center);
				\draw (u4.center) -- (l4.center);
				\draw (u5.center) -- (l5.center);
				\draw (u6.center) -- (l6.center);
				\draw  (ue.center) edge (le.center);
				\draw (5.5,-2.5) node  (us) {} .. 
				controls (3,-4) and (-2.5,-3.5) .. (ue.center) {}
				node[pos=.15 ] (u1) {} 
				node[pos=.3 ] (u2) {}
				node[pos=.45 ] (u3) {}
				node[pos=.6 ] (u4) {}
				
				node[pos=.75 ] (u5) {}
				node[pos=.87 ] (u6) {}
				;
				\draw (6.5,-3)  node (ls) {} .. 
				controls (4,-5.5) and (-3.5,-5) .. (le.center) {}
				node[pos=.15 ] (l1) {} 
				node[pos=.3 ] (l2) {}
				node[pos=.45 ] (l3) {}
				node[pos=.6 ] (l4) {}
				
				node[pos=.75 ] (l5) {}
				node[pos=.87 ] (l6) {}
				;
				\draw  (us.center) edge (ls.center);	
				\draw (u1.center) -- (l1.center);
				\draw (u2.center) -- (l2.center);
				\draw (u3.center) -- (l3.center);
				\draw (u4.center) -- (l4.center);
				\draw (u5.center) -- (l5.center);
				\draw (u6.center) -- (l6.center);
				
				\draw   (us.center)  .. 
				controls (8.5,0.5) and (6.5,2) .. (0.5,-2) node (ue) {}
				node[pos=.15 ] (u1) {} 
				node[pos=.3 ] (u2) {}
				node[pos=.45 ] (u3) {}
				node[pos=.6 ] (u4) {}
				
				node[pos=.75 ] (u5) {}
				node[pos=.87 ] (u6) {}
				;
				\draw  (ls.center) .. 
				controls (10.5,1) and (6.5,3.5) .. (0.5,-0.5) node (le) {}
				node[pos=.15 ] (l1) {} 
				node[pos=.3 ] (l2) {}
				node[pos=.45 ] (l3) {}
				node[pos=.6 ] (l4) {}
				
				node[pos=.75 ] (l5) {}
				node[pos=.87 ] (l6) {}
				;
				\draw  (ue.center) edge (le.center);	
				\draw (u1.center) -- (l1.center);
				\draw (u2.center) -- (l2.center);
				\draw (u3.center) -- (l3.center);
				\draw (u4.center) -- (l4.center);
				\draw (u5.center) -- (l5.center);
				\draw (u6.center) -- (l6.center);
				
				
				\draw (-1.5,5.5) node  (us) {} .. 
				controls (-3.5,4) and (-4,-4) .. (ue.center) 
				node[pos=.15 ] (u1) {} 
				node[pos=.3 ] (u2) {}
				node[pos=.45 ] (u3) {}
				node[pos=.6 ] (u4) {}
				
				node[pos=.75 ] (u5) {}
				node[pos=.87 ] (u6) {}
				;
				\draw (-1,4.5)  node (ls) {} .. 
				controls (-2.5,3.5) and (-2.5,-3) .. (le.center) 
				node[pos=.15 ] (l1) {} 
				node[pos=.3 ] (l2) {}
				node[pos=.45 ] (l3) {}
				node[pos=.6 ] (l4) {}
				
				node[pos=.75 ] (l5) {}
				node[pos=.87 ] (l6) {}
				;
				\draw  (us.center) edge (ls.center);	
				\draw (u1.center) -- (l1.center);
				\draw (u2.center) -- (l2.center);
				\draw (u3.center) -- (l3.center);
				\draw (u4.center) -- (l4.center);
				\draw (u5.center) -- (l5.center);
				\draw (u6.center) -- (l6.center);
				
				
				\draw (us.center) .. controls (0.5,6.5) and (2,4.5) .. (2,2.85)
				node[pos=.2 ] (u1) {} 
				node[pos=.4 ] (u2) {}
				node[pos=.6 ] (u3) {}
				node[pos=.8 ] (u4) {}
				;
				
				\draw (ls.center) .. controls (0,5) and (1,4.5) .. (1,2.85)
				node[pos=.2 ] (l1) {} 
				node[pos=.4 ] (l2) {}
				node[pos=.6 ] (l3) {}
				node[pos=.8 ] (l4) {}
				;
				\draw (u1.center) -- (l1.center);
				\draw (u2.center) -- (l2.center);
				\draw (u3.center) -- (l3.center);
				\draw (u4.center) -- (l4.center);
				\begin{scope}[green,scale=1,line width=5, opacity=0.3,shift={(0,0)}]
					\draw (-1,4.5) .. controls (0,5) and (1,4.5) .. (1,2.85) node (arce) {}; 
					\draw (-1,4.5)   .. 
					controls (-2.5,3.5) and (-2.5,-3) .. (0.5,-0.5) ;
					\draw  (0.5,2) arc (180:120:1);
					\draw  (0.5,2) arc (180:330:1)--++(3,0).. controls (3,1) and (2,0.5) .. (0.5,-0.5);
				\end{scope}
				
				
				\node at (-1,2) {$Q_1$};
				\node at (-7,4.5) {};
			\end{tikzpicture}
			&
			\begin{tikzpicture}[scale=0.35]
				
				\begin{scope}[shift={(1,2.5)}]
					
					\draw (0.5,-0.5) node (circo) {} circle (1);
					\draw (9.5,-0.5) circle (1);
					
					\draw (1.35,0)--(8.55,0);
					\draw (1.35,-1) node (arcs) {}--(8.55,-1);   
					\draw (6.45,0)--(6.45,-1);
					\draw (3.5,0)--(3.5,-1);
					\draw (5.1,0)--(5.1,-1);

					\draw (9,0.35)--(9,7)--(10,7)--(10,0.35);

					\draw (9,2)--(10,2);
					\draw (9,3.5)--(10,3.5);
					\draw (9,5)--(10,5);

					\draw (-2,7)--(10,7);
					
					\node at (0.5,-0.5) {$\Delta_1$};
					\node at (9.5,-0.5) {$\Delta_2$};

					\node[above right] at (1.21,0.2) {$\alpha_1$};
					\node[above left] at (8.73,0.2) {$\alpha_2$};
					
					\node[above] at (5,3.5) {$Q$};
					\node[above] at (-2.5,7.5) {$L_1$};
					\node[above] at (9.5,7.5) {$L_2$};

				\end{scope}

				\draw (-1,9.5) node  (us) {} .. 
				controls (-3,6.5) and (-4,0.5) .. (-1.5,-1.5) node (ue) {}
				node[pos=.15 ] (u1) {} 
				node[pos=.3 ] (u2) {}
				node[pos=.45 ] (u3) {}
				node[pos=.6 ] (u4) {}
				
				node[pos=.75 ] (u5) {}
				node[pos=.87 ] (u6) {}
				;
				\draw (-2.5,9.5)  node (ls) {} .. 
				controls (-4,7) and (-5.5,0.5) .. (-2.5,-2) node (le) {}
				node[pos=.15 ] (l1) {} 
				node[pos=.3 ] (l2) {}
				node[pos=.45 ] (l3) {}
				node[pos=.6 ] (l4) {}
				
				node[pos=.75 ] (l5) {}
				node[pos=.87 ] (l6) {}
				;
				\draw  (us.center) edge (ls.center);	
				\draw (u1.center) -- (l1.center);
				\draw (u2.center) -- (l2.center);
				\draw (u3.center) -- (l3.center);
				\draw (u4.center) -- (l4.center);
				\draw (u5.center) -- (l5.center);
				\draw (u6.center) -- (l6.center);
				\draw  (ue.center) edge (le.center);
				\draw (5.5,-2.5) node  (us) {} .. 
				controls (3,-4) and (0.5,-3.5) .. (ue.center) {}
				node[pos=.15 ] (u1) {} 
				node[pos=.3 ] (u2) {}
				node[pos=.45 ] (u3) {}
				node[pos=.6 ] (u4) {}
				
				node[pos=.75 ] (u5) {}
				node[pos=.87 ] (u6) {}
				;
				\draw (6.5,-3)  node (ls) {} .. 
				controls (3,-5.5) and (-0.5,-4) .. (le.center) {}
				node[pos=.15 ] (l1) {} 
				node[pos=.3 ] (l2) {}
				node[pos=.45 ] (l3) {}
				node[pos=.6 ] (l4) {}
				
				node[pos=.75 ] (l5) {}
				node[pos=.87 ] (l6) {}
				;
				\draw  (us.center) edge (ls.center);	
				\draw (u1.center) -- (l1.center);
				\draw (u2.center) -- (l2.center);
				\draw (u3.center) -- (l3.center);
				\draw (u4.center) -- (l4.center);
				\draw (u5.center) -- (l5.center);
				\draw (u6.center) -- (l6.center);
				
				\draw   (us.center)  .. 
				controls (8.5,-0.5) and (12,1) .. (6,-1) node (ue) {}
				node[pos=.15 ] (u1) {} 
				node[pos=.3 ] (u2) {}
				node[pos=.45 ] (u3) {}
				node[pos=.6 ] (u4) {}
				
				node[pos=.75 ] (u5) {}
				node[pos=.87 ] (u6) {}
				;
				\draw  (ls.center) .. 
				controls (12.5,0.3) and (12,2.5) .. (6,0) node (le) {}
				node[pos=.15 ] (l1) {} 
				node[pos=.3 ] (l2) {}
				node[pos=.45 ] (l3) {}
				node[pos=.6 ] (l4) {}
				
				node[pos=.75 ] (l5) {}
				node[pos=.87 ] (l6) {}
				;
				\draw  (ue.center) edge (le.center);	
				\draw (u1.center) -- (l1.center);
				\draw (u2.center) -- (l2.center);
				\draw (u3.center) -- (l3.center);
				\draw (u4.center) -- (l4.center);
				\draw (u5.center) -- (l5.center);
				\draw (u6.center) -- (l6.center);
				
				
				\draw (-1.5,5) node  (us) {} .. 
				controls (-3,3.5) and (0,-3.5) .. (ue.center) 
				node[pos=.15 ] (u1) {} 
				node[pos=.3 ] (u2) {}
				node[pos=.45 ] (u3) {}
				node[pos=.6 ] (u4) {}
				
				node[pos=.75 ] (u5) {}
				node[pos=.87 ] (u6) {}
				;
				\draw (-1,4)  node (ls) {} .. 
				controls (-1.5,2) and (1.5,-2) .. (le.center) 
				node[pos=.15 ] (l1) {} 
				node[pos=.3 ] (l2) {}
				node[pos=.45 ] (l3) {}
				node[pos=.6 ] (l4) {}
				
				node[pos=.75 ] (l5) {}
				node[pos=.87 ] (l6) {}
				;
				\draw  (us.center) edge (ls.center);	
				\draw (u1.center) -- (l1.center);
				\draw (u2.center) -- (l2.center);
				\draw (u3.center) -- (l3.center);
				\draw (u4.center) -- (l4.center);
				\draw (u5.center) -- (l5.center);
				\draw (u6.center) -- (l6.center);
				
				
				\draw (us.center) .. controls (0.5,6.5) and (2,4.5) .. (2,2.85)
				node[pos=.2 ] (u1) {} 
				node[pos=.4 ] (u2) {}
				node[pos=.6 ] (u3) {}
				node[pos=.8 ] (u4) {}
				;
				
				\draw (ls.center) .. controls (0,5) and (1,4.5) .. (1,2.85)  
				node[pos=.2 ] (l1) {} 
				node[pos=.4 ] (l2) {}
				node[pos=.6 ] (l3) {}
				node[pos=.8 ] (l4) {}
				;
				\draw (u1.center) -- (l1.center);
				\draw (u2.center) -- (l2.center);
				\draw (u3.center) -- (l3.center);
				\draw (u4.center) -- (l4.center);
				
				\begin{scope}[green,scale=1,line width=5, opacity=0.3,shift={(0,0)}]
					\draw  (0.5,2) arc (180:120:1);
					\draw  (0.5,2) arc (180:330:1)--++(7.25,0) arc (210:250:1) .. 
					controls (8.5,1) and (9,1) .. (6,0);
					\draw (-1,4)  .. controls (-1.5,2) and (1.5,-2) .. (6,0);
					\draw (-1,4) .. controls (0,5) and (1,4.5) .. (1,2.85);
				\end{scope}

				\node at (3.5,0.5) {$Q_2$};
			\end{tikzpicture}\\
			
			(i) & (ii)\\

			\begin{tikzpicture}[scale=0.35]

				
				\draw (-1,9.5) node  (us) {} .. 
				controls (-3.5,6) and (-5,2) .. (-3.5,-1) node (ue) {}
				node[pos=.15 ] (u1) {} 
				node[pos=.3 ] (u2) {}
				node[pos=.45 ] (u3) {}
				node[pos=.6 ] (u4) {}
				
				node[pos=.75 ] (u5) {}
				node[pos=.87 ] (u6) {}
				;
				\draw (-2,9.5)  node (ls) {} .. 
				controls (-4.5,7) and (-6,1.5) .. (-4,-2) node (le) {}
				node[pos=.15 ] (l1) {} 
				node[pos=.3 ] (l2) {}
				node[pos=.45 ] (l3) {}
				node[pos=.6 ] (l4) {}
				
				node[pos=.75 ] (l5) {}
				node[pos=.87 ] (l6) {}
				;
				\draw  (us.center) edge (ls.center);	
				\draw (u1.center) -- (l1.center);
				\draw (u2.center) -- (l2.center);
				\draw (u3.center) -- (l3.center);
				\draw (u4.center) -- (l4.center);
				\draw (u5.center) -- (l5.center);
				\draw (u6.center) -- (l6.center);
				\draw  (ue.center) edge (le.center);
				\draw (2,-3.5) node  (us) {} .. 
				controls (-0.5,-3.5) and (-1.5,-3.5) .. (ue.center) {}
				node[pos=.15 ] (u1) {} 
				node[pos=.3 ] (u2) {}
				node[pos=.45 ] (u3) {}
				node[pos=.6 ] (u4) {}
				
				node[pos=.75 ] (u5) {}
				node[pos=.87 ] (u6) {}
				;
				\draw (2,-4.5)  node (ls) {} .. 
				controls (-1,-4.5) and (-1.5,-4.5) .. (le.center) {}
				node[pos=.15 ] (l1) {} 
				node[pos=.3 ] (l2) {}
				node[pos=.45 ] (l3) {}
				node[pos=.6 ] (l4) {}
				
				node[pos=.75 ] (l5) {}
				node[pos=.87 ] (l6) {}
				;
				\draw  (us.center) edge (ls.center);	
				\draw (u1.center) -- (l1.center);
				\draw (u2.center) -- (l2.center);
				\draw (u3.center) -- (l3.center);
				\draw (u4.center) -- (l4.center);
				\draw (u5.center) -- (l5.center);
				\draw (u6.center) -- (l6.center);
				
				\draw   (us.center)  .. 
				controls (4,-3) and (5.5,-2) .. (2,-2) node (ue) {}
				node[pos=.15 ] (u1) {} 
				node[pos=.3 ] (u2) {}
				node[pos=.45 ] (u3) {}
				node[pos=.6 ] (u4) {}
				
				node[pos=.75 ] (u5) {}
				node[pos=.87 ] (u6) {}
				;
				\draw  (ls.center) .. 
				controls (5,-4) and (7,-1) .. (2,-1) node (le) {}
				node[pos=.15 ] (l1) {} 
				node[pos=.3 ] (l2) {}
				node[pos=.45 ] (l3) {}
				node[pos=.6 ] (l4) {}
				
				node[pos=.75 ] (l5) {}
				node[pos=.87 ] (l6) {}
				;
				\draw  (ue.center) edge (le.center);	
				\draw (u1.center) -- (l1.center);
				\draw (u2.center) -- (l2.center);
				\draw (u3.center) -- (l3.center);
				\draw (u4.center) -- (l4.center);
				\draw (u5.center) -- (l5.center);
				\draw (u6.center) -- (l6.center);
				
				
				\draw (-2.2,5.3) node  (us) {} .. 
				controls (-4,4) and (-2.5,-2) .. (ue.center) 
				node[pos=.15 ] (u1) {} 
				node[pos=.3 ] (u2) {}
				node[pos=.45 ] (u3) {}
				node[pos=.6 ] (u4) {}
				
				node[pos=.75 ] (u5) {}
				node[pos=.87 ] (u6) {}
				;
				\draw (-1.7,4.3)  node (ls) {} .. 
				controls (-2.5,3.5) and (-2,-1) .. (le.center) 
				node[pos=.15 ] (l1) {} 
				node[pos=.3 ] (l2) {}
				node[pos=.45 ] (l3) {}
				node[pos=.6 ] (l4) {}
				
				node[pos=.75 ] (l5) {}
				node[pos=.87 ] (l6) {}
				;
				\draw  (us.center) edge (ls.center);	
				\draw (u1.center) -- (l1.center);
				\draw (u2.center) -- (l2.center);
				\draw (u3.center) -- (l3.center);
				\draw (u4.center) -- (l4.center);
				\draw (u5.center) -- (l5.center);
				\draw (u6.center) -- (l6.center);
				
				
				\draw (us.center) .. controls (0.5,6.5) and (1,4.5) .. (1,2.9)
				node[pos=.2 ] (u1) {} 
				node[pos=.4 ] (u2) {}
				node[pos=.6 ] (u3) {}
				node[pos=.8 ] (u4) {}
				;
				
				\draw (ls.center) .. controls (0,5) and (0,4.5) .. (0,2.9)
				node[pos=.2 ] (l1) {} 
				node[pos=.4 ] (l2) {}
				node[pos=.6 ] (l3) {}
				node[pos=.8 ] (l4) {}
				;
				\draw (u1.center) -- (l1.center);
				\draw (u2.center) -- (l2.center);
				\draw (u3.center) -- (l3.center);
				\draw (u4.center) -- (l4.center);

				\begin{scope}[green,scale=1,line width=5, opacity=0.3,shift={(0,0)}]
					\coordinate (ls) at (5,-2.25)  {} ;
					\draw  (ls.center) .. 
					controls (5,-1.5) and (4.5,-1) .. (2,-1) node (le) {};
					\draw (-1.7,4.3)  node (ls) {} .. 
					controls (-2.5,3.5) and (-2,-1) .. (le.center) ;
					\draw (ls.center) .. controls (0,5) and (0,4.5) .. (0,2.85) arc (120:330:1) 
					coordinate (t1) --(5.1,1.5);
					
				\end{scope}

				\draw  (0.5,2) node {$\Delta_1$} circle ( 1);
				\node at (-1.5,10.5) {$L_1$};
				\node[above right] at (1.3,2.65) {$\alpha_1$};

				\begin{scope}[xscale=-1, shift={(-10.1,0)}]
					
					\draw (-1,9.5) node  (us) {} .. 
					controls (-3.5,6) and (-5,2) .. (-3.5,-1) node (ue) {}
					node[pos=.15 ] (u1) {} 
					node[pos=.3 ] (u2) {}
					node[pos=.45 ] (u3) {}
					node[pos=.6 ] (u4) {}
					
					node[pos=.75 ] (u5) {}
					node[pos=.87 ] (u6) {}
					;
					\draw (-2,9.5)  node (ls) {} .. 
					controls (-4.5,7) and (-6,1.5) .. (-4,-2) node (le) {}
					node[pos=.15 ] (l1) {} 
					node[pos=.3 ] (l2) {}
					node[pos=.45 ] (l3) {}
					node[pos=.6 ] (l4) {}
					
					node[pos=.75 ] (l5) {}
					node[pos=.87 ] (l6) {}
					;
					\draw  (us.center) edge (ls.center);	
					\draw (u1.center) -- (l1.center);
					\draw (u2.center) -- (l2.center);
					\draw (u3.center) -- (l3.center);
					\draw (u4.center) -- (l4.center);
					\draw (u5.center) -- (l5.center);
					\draw (u6.center) -- (l6.center);
					\draw  (ue.center) edge (le.center);
					\draw (2,-3.5) node  (us) {} .. 
					controls (-0.5,-3.5) and (-1.5,-3.5) .. (ue.center) {}
					node[pos=.15 ] (u1) {} 
					node[pos=.3 ] (u2) {}
					node[pos=.45 ] (u3) {}
					node[pos=.6 ] (u4) {}
					
					node[pos=.75 ] (u5) {}
					node[pos=.87 ] (u6) {}
					;
					\draw (2,-4.5)  node (ls) {} .. 
					controls (-1,-4.5) and (-1.5,-4.5) .. (le.center) {}
					node[pos=.15 ] (l1) {} 
					node[pos=.3 ] (l2) {}
					node[pos=.45 ] (l3) {}
					node[pos=.6 ] (l4) {}
					
					node[pos=.75 ] (l5) {}
					node[pos=.87 ] (l6) {}
					;
					\draw  (us.center) edge (ls.center);	
					\draw (u1.center) -- (l1.center);
					\draw (u2.center) -- (l2.center);
					\draw (u3.center) -- (l3.center);
					\draw (u4.center) -- (l4.center);
					\draw (u5.center) -- (l5.center);
					\draw (u6.center) -- (l6.center);
					
					\draw   (us.center)  .. 
					controls (4,-3) and (5.5,-2) .. (2,-2) node (ue) {}
					node[pos=.15 ] (u1) {} 
					node[pos=.3 ] (u2) {}
					node[pos=.45 ] (u3) {}
					node[pos=.6 ] (u4) {}
					
					node[pos=.75 ] (u5) {}
					node[pos=.87 ] (u6) {}
					;
					\draw  (ls.center) .. 
					controls (5,-4) and (7,-1) .. (2,-1) node (le) {}
					node[pos=.15 ] (l1) {} 
					node[pos=.3 ] (l2) {}
					node[pos=.45 ] (l3) {}
					node[pos=.6 ] (l4) {}
					
					node[pos=.75 ] (l5) {}
					node[pos=.87 ] (l6) {}
					;
					\draw  (ue.center) edge (le.center);	
					\draw (u1.center) -- (l1.center);
					\draw (u2.center) -- (l2.center);
					\draw (u3.center) -- (l3.center);
					\draw (u4.center) -- (l4.center);
					\draw (u5.center) -- (l5.center);
					\draw (u6.center) -- (l6.center);
					
					
					\draw (-2.2,5.3) node  (us) {} .. 
					controls (-4,4) and (-2.5,-2) .. (ue.center) 
					node[pos=.15 ] (u1) {} 
					node[pos=.3 ] (u2) {}
					node[pos=.45 ] (u3) {}
					node[pos=.6 ] (u4) {}
					
					node[pos=.75 ] (u5) {}
					node[pos=.87 ] (u6) {}
					;
					\draw (-1.7,4.3)  node (ls) {} .. 
					controls (-2.5,3.5) and (-2,-1) .. (le.center) 
					node[pos=.15 ] (l1) {} 
					node[pos=.3 ] (l2) {}
					node[pos=.45 ] (l3) {}
					node[pos=.6 ] (l4) {}
					
					node[pos=.75 ] (l5) {}
					node[pos=.87 ] (l6) {}
					;
					\draw  (us.center) edge (ls.center);	
					\draw (u1.center) -- (l1.center);
					\draw (u2.center) -- (l2.center);
					\draw (u3.center) -- (l3.center);
					\draw (u4.center) -- (l4.center);
					\draw (u5.center) -- (l5.center);
					\draw (u6.center) -- (l6.center);
					
					
					\draw (us.center) .. controls (0.5,6.5) and (1,4.5) .. (1,2.9)
					node[pos=.2 ] (u1) {} 
					node[pos=.4 ] (u2) {}
					node[pos=.6 ] (u3) {}
					node[pos=.8 ] (u4) {}
					;
					
					\draw (ls.center) .. controls (0,5) and (0,4.5) .. (0,2.9)
					node[pos=.2 ] (l1) {} 
					node[pos=.4 ] (l2) {}
					node[pos=.6 ] (l3) {}
					node[pos=.8 ] (l4) {}
					;
					\draw (u1.center) -- (l1.center);
					\draw (u2.center) -- (l2.center);
					\draw (u3.center) -- (l3.center);
					\draw (u4.center) -- (l4.center);

					\begin{scope}[green,scale=1,line width=5, opacity=0.3,shift={(0,0)}]
						\coordinate (ls) at (5,-2.25)  {} ;
						\draw  (ls.center) .. 
						controls (5,-1.5) and (4.5,-1) .. (2,-1) node (le) {};
						\draw (-1.7,4.3)  node (ls) {} .. 
						controls (-2.5,3.5) and (-2,-1) .. (le.center) ;
						\draw (ls.center) .. controls (0,5) and (0,4.5) .. (0,2.85) arc (120:330:1) 
						coordinate (t1) --(5.1,1.5);
						
					\end{scope}

					\draw  (0.5,2) node {$\Delta_2$} circle ( 1);
					\node at (-1.5,10.5) {$L_2$};
					\node [above left] at (1.25,2.7) {$\alpha_2$};

				\end{scope} 

				\begin{scope}
					\path (t1) coordinate (ls) --++(0,1) coordinate (us);
					\draw (us.center)  -- ++(-7.4,0)
					node[pos=.3 ] (u1) {} 
					node[pos=.5 ] (u2) {}
					node[pos=.7 ] (u3) {}
					;
					
					\draw (ls.center)  -- ++(-7.4,0)
					node[pos=.3 ] (l1) {} 
					node[pos=.5 ] (l2) {}
					node[pos=.7 ] (l3) {}
					
					;
					\draw (u1.center) -- (l1.center);
					\draw (u2.center) -- (l2.center);
					\draw (u3.center) -- (l3.center);
					\draw (u4.center) -- (l4.center);
				\end{scope}
				\node at (5,0.5) {$Q_3$};
				\draw (-1,9.5) -- (11.1,9.5); 
				\node at (5,3.2) {$L_3$};
				\node at (5,6.5) {$Q$};
			\end{tikzpicture} &
			\begin{tikzpicture}[scale=.35]
				
				\draw (0,2)--(0,7)--(1,7)--(1,2.25);
				\draw (9,2.25)--(9,7)--(10,7)--(10,2);
				
				\draw (0,3)--(1,3);
				\draw (0,4.5)--(1,4.5);
				\draw (0,6)--(1,6);
				
				\draw (9,3)--(10,3);
				\draw (9,4.5)--(10,4.5);
				\draw (9,6)--(10,6);

				\draw (1,7)--(10,7);

				\node[above ] at (1.5,2.2) {$\alpha_1$};
				\node[above ] at (8.6,2.2) {$\alpha_2$};
				
				\node[above] at (5,5) {$Q$};
				\node[above] at (0.5,7.5) {$L_1$};
				\node[above] at (9.5,7.5) {$L_2$};

				\draw (5,0.5) coordinate (v1) .. controls (-3,6.5) and (-3,-5.5) .. (v1);
				\draw[xscale=-1] (-5,0.5) coordinate (v1) .. controls (-13,6.5) and (-13,-5.5) .. (v1);
				\node at (1.5,0.5) {$\Delta_1$};
				\node at (8.5,0.5) {$\Delta_2$};
				
				\draw (1.75,2.15) node  (us) {} .. 
				controls (3,3) and (7,3) .. (8.25,2.15) node (ue) {}
				node[pos=.15 ] (u1) {} 
				node[pos=.3 ] (u2) {}
				node[pos=.45 ] (u3) {}
				node[pos=.6 ] (u4) {}
				
				node[pos=.75 ] (u5) {}
				node[pos=.87 ] (u6) {}
				;
				\draw (2.55,1.9)  node (ls) {} .. 
				controls (4,2.5) and (6,2.5) .. (7.45,1.9) node (le) {}
				node[pos=.15 ] (l1) {} 
				node[pos=.3 ] (l2) {}
				node[pos=.45 ] (l3) {}
				node[pos=.6 ] (l4) {}
				
				node[pos=.75 ] (l5) {}
				node[pos=.87 ] (l6) {}
				;
				\draw (u1.center) -- (l1.center);
				\draw (u2.center) -- (l2.center);
				\draw (u3.center) -- (l3.center);
				\draw (u4.center) -- (l4.center);
				\draw (u5.center) -- (l5.center);
				\draw (u6.center) -- (l6.center);
				\begin{scope}[green,scale=1,line width=5, opacity=0.3,shift={(0,0)}]
					\draw (2.55,1.9)  node (ls) {} .. 
					controls (4,2.5) and (6,2.5) .. (7.45,1.9) node (le) {};
					\draw (ls.center) .. controls (3.9,1.3) and (3.9,1.3) .. (v1.center);
					\begin {scope} [xscale=-1, shift={(-10,0)}]
					\coordinate   (ls) at (2.55,1.9) {};
					\coordinate   (v1) at (5,0.5) {};
					\draw (ls.center) .. controls (3.9,1.3) and (3.9,1.3) .. (v1.center);
					\end {scope}
				\end{scope}

				\node at (5,1.5) {$Q_4$};
				\node at (5,3.35) {$L_3$};
			\end{tikzpicture}\\
			
			(iii) & (iv)
		\end{tabular}
		
		Figure 40
	\end{center}
	
	\begin{itemize}\item []
		\begin{itemize}
			\item [i)] Let $Q_1$ be the bounded component of $\mathbb E^2\setminus(\Delta_1\cup L_2\cup L_3)$. 
			Then\\
			$\partial Q_1=\rho_1\rho_2\rho_3$, 
			$\rho_1=\partial Q_1\cap\partial L_1$, 
			$\rho_2=\partial Q_1\cap\partial\Delta_1$, 
			$\rho_3=\partial Q_1\cap\partial L_3$.
			Then $\rho_i$ is connected and
			$\Phi(\rho_1)=c_1\in C_\Gamma(t)$,
			$\Phi(\rho_3)=c_3\in C_\Gamma(t)$,
			$\Phi(\rho_2)\in \langle a,t\rangle$,
			$\|\Phi(\rho_2)\|_t=\|\Phi(\partial\Delta_1)\|_t-2\geq 2$,
			$\|\Phi(\rho_2)\|_a=\|\Phi(\partial\Delta_1)\|_a-1\geq 1$.

			%
			Consider $\widetilde{Q}_1^t$. 
			Since $\|\Phi(\rho_2)\|_t\geq 2, |\widetilde{Q}_1^t|\neq 0$,
			we may assume that $V_1=\Phi(\partial Q_1)$ is cyclically reduced by carrying out diamond moves if needed. 
			Also, it is easy to see that the conditions of Lemma 2.2.1 are satisfied, hence may assume $Q_1\in \mathcal{M}_3(V)$. 
			Consequently, Theorem A applies to $Q_1\,(\subsetneq M)$. 
			Hence $\widetilde{Q}_1^t$ satisfies the condition C(4)\&T(4).
			This leads to a contradiction as in the proof of Lemma 1.3.2(e).
			(Notice that $L_1$ and $L_3$ become edges in $\widetilde{Q}_1^t$.) 
			
			
			\item [ii)] See Fig. 40(ii). Let $Q_2$ be the bounded component of $\mathbb E^2\setminus(\Delta_1\cup L_3\cup \Delta_2)$.
			Then $\partial Q_2=\rho_1\rho_2\rho_3$, where $\rho_1=\partial Q_2\cap\partial \Delta_1$,
			$\rho_2=\partial Q_2\cap\partial L_3$ and
			$\rho_3=\partial Q_2\cap\partial \Delta_2$
			and where
			$\Phi(\rho_1)\in \langle a,t\rangle$, 
			$\Phi(\rho_2)\in C_\Gamma(t)$ and
			$\Phi(\rho_3)\in \langle a,t\rangle$.
			Then $\rho_i$ are connected.
			Consider $\tilde{Q}_2^t$.
			Since $\Phi(\rho_1)\neq 1,\ |\tilde{Q}_2^t|\geq 1$.
			As in part (i), we may assume that $V_2:=\partial Q_2$ is cyclically reduced and $Q_2\in \mathcal{M}_3(V_2)$.
			Also, as in part (i), $\tilde{Q}_2^t$ satisfies the condition C(4)\& T(4) 
			and this leads to a contradiction, as in Lemma 1.3.2(f).
			(See remark at the end of part (i))
			\item [iii)] See Fig. 40(iii).
			Let $Q_3$ be the bounded connected component of
			$\mathbb{E}^2\setminus(L_1\cap\Delta_1\cap L_3\cap \Delta_2)$.
			Let $\partial Q_3=\rho_1\rho_2\rho_3\rho_4$,
			where 
			$\rho_1=\partial Q_3\cap\partial L_1$,
			$\rho_2=\partial Q_3\cap\partial \Delta_1$,
			$\rho_3=\partial Q_3\cap\partial L_3$ and
			$\rho_4=\partial Q_3\cap\partial \Delta_2$.
			From this point on the proof is similar to the proof of part (ii).
			\item [iv)]  See Fig. 40(iv).
			Let $Q_4$ be the bounded connected component of
			$\mathbb{E}^2\setminus(L_1\cup\Delta_1\cup L_3\cup \Delta_2\cup L_2)$.
			Let $\partial Q_4=\rho_1\rho_2\rho_3\rho_4\rho_5$,
			where 
			$\rho_1=\partial Q_4\cap\partial L_1$,
			$\rho_2=\partial Q_4\cap\partial \Delta_1$,
			$\rho_3=\partial Q_4\cap\partial L_3$,
			$\rho_4=\partial Q_4\cap\partial \Delta_2$ and
			$\rho_5=\partial Q_4\cap\partial L_2$
			From this point on the proof is similar to the proof of part (ii).
			We omit details. Thus, none of the cases (i)-(iv) is possible.
			
		\end{itemize}

		Hence $\mu$ is a simple curve.  
		Since $N_P(\mu)=L_1\cup\Delta_1\cup L_2\cup\Delta_2\cup L_3$ and since $L_i$ and $\Delta_j$ are regions in $\mathbb M_t$, they have connected interior.
		Hence $\mu\cap\gamma=\emptyset$, as $\gamma\subseteq\partial Q$.
		By the definition of $\mu$ and $\gamma$, they have the same endpoints. Consequently, since $M$ is planar, $\mu w\gamma^{-1}$ is a simple closed curve, $w$ a point.
		Checking $\Phi(\mu)$ shows that $\Phi(\mu)$ is reduced. Clearly, $\Phi(\gamma)$ is reduced. 
		Let $z_1$ and $z_2$ be the first and last letters of $\Phi(\mu)$, respectively and let $s_1$ and $s_2$ 
		be the first and last letters of $\Phi(\gamma)$, respectively. 
		If $z_1=t^{\pm 1}$ then $z_1s_1^{-1}$ is reduced  since $s_1\subseteq H_t$. 
		If $z_1\neq t^{\pm 1}$ then $z_1s_1^{-1}$ is a subword of $\Phi(\omega_1)$, hence is reduced. Similarly, for $z_2$ and $s_2$. Hence $\partial P$ is cyclically reduced. 
		Hence
		\begin{equation*}
			\left.
			\begin{minipage}{10cm}
				$P$ is connected, simply connected with simple boundary cycle and with cyclically reduced boundary label.
			\end{minipage}
			\right\}
		\end{equation*}
		Finally, it easily follows from Lemma 2.2.1 that $P\in\mathcal{M}_3(V)$ where $V:=\Phi(\partial P)$
		\item [(c)] Follows from part (a)
	\end{itemize}
	\ \hfill $\Box$

	\subsection {Definitions and assumptions}
	\subsubsection* {Generalised Greendlinger regions}
	
	In Section 1.3 we introduced (classical) Greendlinger regions. In the rest of the work we shall need more general notions of Greendlinger regions. Recall from Section 1.3 $\mathcal D_i(M),i=1,2$ and $\mathcal D(M)$, 
	where M is a connected, simply connected map with connected interior which contains at least two regions.\\
	If  $\partial M=u\mu v\nu^{-1},\mu$ and $\nu$ boundary paths, we shall call M {\em bigonal with sides $\mu$ and $\nu$}. We define $\mathcal D_i(\mu),\mathcal D_i(\nu),i=1,2$ and $\mathcal D(\mu)$ and $\mathcal D(\nu)$ as follows:
	
	\vspace{10pt}
	\tocexclude{
		\subsubsection{Definitions}
	}
	\addcontentsline{toc}{subsubsection}{\numberline{\thesubsubsection} 
		Definition (Generalised Greendlinger regions)}
	Let M be a bigonal map with sides $\mu$ and $\nu$. Define 
	\begin{itemize}
		\item[(i)] $\mathcal D_j(\mu)=\{ D\in Reg(M),\partial D\cap \mu$ is connected and $i(D)=j\}\quad j=1,2$. $\mathcal D_j(\nu)$ is defined accordingly. 
		Also, define $\mathcal D(\mu)=\mathcal D_1(\mu)\cup \mathcal D_2(\mu)\cup \mathcal D_0(\mu)$ where $\mathcal D_0(\mu)$ is defined only if $M=\{D\}$. 
		In this case $\mathcal D_0(\mu)=\{D\}$ if $|\mu|\geq |\nu|$ and $\mathcal D_0(\nu)=\{D\}$ if $|\nu|\geq |\mu|$. In case $|\mu|=|\nu|$ we choose $\mathcal D_0(\mu)=\{D\}$.
		\item [(ii)] Let $M$ be an $\mathcal R$-diagram, $\mathcal R$ given by (I). Assume $M\in \mathcal M_3(W)$ and $\partial M=\mu\nu^{-1}$. 
		
		\noindent Define 
		\begin{equation*}
			\mathbb D(\mu)=\left\{\Delta \in Reg_{4^+}(\mathbb M) \middle|\; 
			\begin{minipage}{6cm}
				$\partial \Delta \cap\mu$ is connected, 
				$||\partial \Delta \cap \mu||\geq n(\Delta)$  and 	$|\Phi(\xi)|\leq |\Phi(\partial\Delta\cap \mu)|$
			\end{minipage}
			\right\}
		\end{equation*} 
		Here $\xi$ is the complement of $\partial\Delta\cap\mu$ on $\partial\Delta$. Define $\mathbb D(\nu)$ accordingly.
		\item[(iii)] Define 
		\begin{equation*}
			\mathcal D^v(\mu)=\left\{\Delta \in Reg_{4^+}(\mathbb M) \middle|\; 
			\begin{minipage}{6cm}
				$\Delta$ contains a subdiagram S with 
				$S \in O L\left(\theta_1, \theta_2\right)$,
				$S\subseteq N(\mu)$,
				$\partial S \cap \mu=\theta_1$,
				$\|\theta_1\| \geq n\left(\Delta\right)$ $\geq 4$ and 
				$\left|\Phi\left(\theta_2\right)\right| \leq\left|\Phi\left(\theta_1\right)\right|$. 
			\end{minipage}
			\right\}
		\end{equation*}
		Define $\mathcal D^v(\nu)$  accordingly.
		\item [(iv)] Define\\ $\mathcal{D}^1_2(\mu)=\{D\in \mathcal{D}_1(\mu),n(D)=2\}$,
		$\mathcal{D}^2_2(\mu)=\{D\in \mathcal{D}_2(\mu),n(D)=2\}$.
		\item [(v)] Define $\sigma:\text{boundary paths\ }\rightarrow\{0,1\}$ 
		as follows:\\
		$\sigma(\mu)=1$ if one the following holds:
		\begin{itemize}
			\item [1)] $\mathcal{D}_2^1(\mu)\neq\emptyset$
			\item [2)] $|\mathcal{D}_2^2(\mu)|\geq 2$
			\item [3)] $|\mathcal{D}_2^2(\mu)|=1$ and there is a region $\Delta$ in 
			$\mathbb M_{4^+}$ and a sequence $J$ of I-moves such that $\partial\Delta^J\cap\mu$ is connected and 
			$|\partial\Delta^J\cap\mu|\geq 2$
			\item [4)] $\mathcal{D}_2^2(\mu)=\emptyset$ and 
			$\mathbb D(\mu)\cup \mathcal{D}^v(\mu)\neq\emptyset$
		\end{itemize}
		\noindent In all other cases $\sigma(\mu)=0$.
		
		\item [(vi)] Let $\tilde{\mu}=\Psi_t(\mu)$. Define
		\begin{equation*}
			\mathcal{D}'(\tilde{\mu})=
			\{\tilde{\Delta}\in Reg(\widetilde{\mathbb M}^t)|\Delta\in\mathcal{D}^v(\mu),
			\Delta=\Phi^{-1}_t(\tilde\Delta)\}
		\end{equation*}
	\end{itemize}

	\subsubsection{Assumptions}
	The proofs of Theorems A,B and C rely on other results as well(below) which  are proved by simultaneous induction. They are
	\begin{itemize}
		\item [(I)] Theorems A,B and C
		\item [(II)] Theorem S, Theorem L and Theorem T.
	\end{itemize}
	\begin{sssection}[Theorem S]
		\addcontentsline{toc}{subsubsection}{\numberline{\thesubsubsection} 
			Theorem S (Structure of $Opt(\Delta)$)}
		{\it \quad Let $M\in\mathcal{M}_3(W)$ and assume that $\partial M$  decomposes by 
			$\omega_1 v_1\omega_2^{-1} v_2$, where 
			$\omega_1$ and $\omega_2$ are boundary paths and $v_1$ and $v_2$ are boundary vertices. 
			Let $t\in T(M)$ and let $\widetilde{\Delta}\in\mathcal{D}(\widetilde{\omega}_1)$, 
			$\widetilde{\omega}_1=\Psi_t(\omega_1)$,
			$\Delta=\Psi_t^{-1}(\widetilde{\Delta}), Supp(\Delta)=\{a,t\}$
			Assume $\sigma(\omega_1)=0$ in parts (a),(b) and (c) . 
			Then there is a sequence J of I-moves such that 
			$\Delta^J=\Delta, M^J\in\mathcal{M}_3(W)$ and each of the following holds 
			\begin{itemize}
				\item [(a)] If $\beta$ is an edge of 
				$\left(\partial\Delta\cap\partial \widehat{Opt}(\Delta)\right)^J$ 
				then there is a band which has a pole on $\beta$ and the other pole on 
				$\omega_1$. 
				In particular, $(\partial\Delta\cap_B \omega_1)^J$ is connected.
				\item [(b)] If $B_1$ and $B_2$ are bands emanating from $(\partial\Delta\cap_B \omega_1)^J$
				with adjacent poles on $\partial\Delta$ then 
				$\partial B_1\cap \partial B_2$ is connected. Moreover, 
				let $L_1,\ldots, L_k$ be the bands that emanate from $(\partial\Delta\cap_B \omega_1)^J$. 
				Then $\partial L_i\cap \partial L_{i + 1}$ is connected, $i=1,\ldots,k-1$, where $L_i$ and $L_{i + 1}$ have adjacent poles on $\partial\Delta$.
				\item [(c)] 
				Let $Z:=\widehat{Opt}(\Delta)^J\cup \omega_1\setminus\{L_1,\ldots,L_k\}=Q_1\dot\bigcup\cdots\dot\bigcup Q_l$, 
				where $Q_i$ are the connected components.
				\begin{itemize}
					\item [(i)] $\partial Q_i$ are simple closed curves and 
					$\partial Q_i$ decomposes by $\partial Q_i=\rho_i\xi_i\gamma_i^{-1}$, 
					where $\rho_i=\partial Q_i\cap\partial L_j$,
					$\xi_i=\partial Q_i\cap \partial L_{j+1}$
					and $\gamma_i=\partial Q_i\cap \omega_1$, $j\in \{1,\ldots k\}, i=1,\ldots, l$
					\item [(ii)] Up to diamond moves $U_i:=\Phi(\partial Q_i)$ 
					is cyclically reduced and $Q_i\in\mathcal{M}_3(U_i),i=1,\ldots, l$
					and either $Q_i=\{D_i\}, n(D_i)=2$ or
					$Q_i$ is transferable via $L_j\cup L_{j+1}$. See Fig. 41.
				\end{itemize}
			\end{itemize}
			\vbox{\rm
				\begin{center}
					\begin{minipage}{1in}
						\begin{center}
							\begin{tikzpicture}[scale=.4]
								\draw (2,0)--(4,0)--(4,4)--(6,7)--(4.5,7)--(3,4)--(1.5,7)--(0,7)--(2,4)--(2,0);
								\draw (2,1)--(4,1);
								\draw (2,2)--(4,2);
								\draw (2,3)--(4,3);
								\draw (2,4)--(4,4);
								\draw (1.6,4.5)--(2.35,5.3);
								\draw (1.3,5)--(2,6);
								\draw (1,5.5)--(1.7,6.5)  node[right] {$\rho_i$};
								\draw (4.35,4.5)--(3.5,5);
								\draw (4.8,5.2)--(3.9,5.8);
								\draw (5.2,5.8)--(4.2,6.4)  node[left=-1mm] {$\xi_i$};

								\draw[->] (1.5,7)--(3,7);
								\draw (3,7)--(4.5,7);

								\node[above] at (1,7) {$L_j$};
								\node[above] at (5,7) {$L_{j+1}$};
								\node at (3,7.5) {$\gamma_j$};
								\node at (3,5.5) {$Q_i$};
								\draw (3,0)--(3,4);
								
							\end{tikzpicture}
							
							(0)
						\end{center}
					\end{minipage}
					\begin{minipage}{1in}
						\begin{center}
							\begin{tikzpicture}[scale=.4]
								
								\draw (2,0)--(4,0)--(4,4)--(6,7)--(4.5,7)--(3,4)--(1.5,7)--(0,7)--(2,4)--(2,0);
								\draw (2,1)--(4,1);
								\draw (2,2)--(4,2);
								\draw (2,3)--(4,3);
								\draw (2,4)--(4,4);
								\draw (1.7,4.5)--(2.35,5.3);
								\draw (1.3,5)--(2,6);
								\draw (1,5.5)--(1.7,6.5);
								\draw (4.35,4.5)--(3.5,5);
								\draw (4.8,5.2)--(3.9,5.8);
								\draw (5.2,5.8)--(4.2,6.4);

								\draw[thick, ->] (1.5,7)--(2.3,5.5);
								\draw[thick] (2.3,5.5)--(3,4)--(4.5,7);

								\node[above] at (1,7) {$B_l$};
								\node[above] at (5,7) {$B_{r}$};
								\node at (3,5.9) {$\gamma$};
								
								\draw (3,0)--(3,4);
								\node at (1,6.5) {$E_1$};
								\node at (5,6.5) {$E_2$};
								
							\end{tikzpicture}
							
							(1)
						\end{center}
					\end{minipage}
					\hspace{.1cm}
					\begin{minipage}{1in}
						
						\begin{center}
							\begin{tikzpicture}[scale=.4]
								\draw (2,0)--(4,0)--(4,4)--(6,7)--(4.5,7)--(3,4)--(1.5,7)--(0,7)--(2,4)--(2,0);
								\draw (2,1)--(4,1);
								\draw (2,2)--(4,2);
								\draw (2,3)--(4,3);
								\draw (2,4)--(4,4);
								\draw (1.6,4.5)--(2.35,5.3);
								\draw (1.3,5)--(2,6);
								\draw (1,5.5)--(1.7,6.5);
								\draw (4.35,4.5)--(3.5,5);
								\draw (4.8,5.2)--(3.9,5.8);
								\draw (5.2,5.8)--(4.2,6.4);

								\draw (1.5,7)--(3,8);
								\draw (3,8)--(4.5,7);

								\node[above] at (1,7) {$B_l$};
								\node[above] at (5,7) {$B_{r}$};
								\node at (3,8.5) {$\gamma$};
								\node at (3,5.7) {$Q_i$};
								\draw (3,0)--(3,4);
								
							\end{tikzpicture}

							(2)
						\end{center}
					\end{minipage}
					\hspace{.1cm}
					\begin{minipage}{1in}
						\begin{center}
							\begin{tikzpicture}[scale=.15]
								\draw (5,0)--(0,0)--(0,15)--(5,15)--(5,0)--(10,0)--(10,10)--(5,10);
								\node at (2.3,16.5) {$B_l$};
								\node at (7.6,11.5) {$B_{r}$};

								\node at (2.3,12.5) {$E$};
								\draw (0,3)--(10,3);
								\draw (0,7)--(10,7);
								\draw (0,10)--(5,10);
							\end{tikzpicture}
							
							\vspace{.5cm}
							
							(3)
						\end{center}
					\end{minipage}\hspace{.1cm}\begin{minipage}{1in}
						
						\begin{center}
							\begin{tikzpicture}[scale=.15]
								\draw (0,0)--(0,15)--(10,15)--(10,0)--(0,0);
								\draw (5,0)--(5,15);
								\node at (2.3,16.5) {$B_l$};
								\node at (7.3,16.5) {$B_{r}$};
								
							\end{tikzpicture}

							\vspace{.5cm}
							
							(4)
						\end{center}
					\end{minipage}
					
					Figure 41
				\end{center}
			}	
			\begin{itemize} \item []
				\begin{itemize}
					\item[(iii)] For every $D\in Reg(Q_i),i=1,\ldots,l,
					\ Supp(\Phi(\partial D))\subseteq X\setminus\{a,t\}$		
				\end{itemize}
				\item [(d)] Let $\eta=(\partial\Delta\cap_B\omega_1)^J$.
				Then the following hold:
				\begin{itemize}
					\item [(i)] If $\Delta\in\mathbb{D}(\eta)$ then $\sigma(\omega_1)=1$.
					\item [(ii)] If $\widetilde{\mathbb M}^t\in O.L,
					\widetilde{\mathbb M}^t=\langle\widetilde{\Delta}_1,\ldots, \widetilde{\Delta}_m\rangle$, $m\geq 2$ and 
					$\Delta_j\in\mathcal{D}^v(\eta)$, for $j=1$ and $j=m$,
					then $\sigma(\omega_1)=1$.
				\end{itemize}

			\end{itemize} 
		}
	\end{sssection}
	\noindent {\bf Remarks}
	\begin{itemize}
		\item [1)] It follows from part (c) that if $\mathcal{D}_2^2(\omega_1)=\emptyset$ then $Opt(\Delta)$ is a homogeneous band-bundle, 
		hence we can push up $\Delta$ by a sequence $J_1$ of I-moves to $\omega_1$, 
		such that $\Delta^j\in\mathbb D(\mu)$, obtaining the classical result on Greendlinger regions.
		\item [2)] The same results hold if we use $\nu$ in place of $\mu$.
	\end{itemize}
	
	\begin{sssection}[Theorem L]
		\addcontentsline{toc}{subsubsection}{\numberline{\thesubsubsection} 
			Theorem L (Lifting $\mathcal{D}(\widetilde{\mathbb M}^t)$ to
			$\mathcal{D}^v({\mathbb M}_t)$)}
		{\it
			\ \\ Let $M\in \mathcal{M}_3(W)$, $\partial M=u\omega_1 v\omega^{-1}_2,\ u$ 
			and $v$ vertices, $\omega_1$ and $\omega_2$ boundary paths of $M$.
			Assume $T(M) \neq\emptyset$ and let $t\in T(\mathbb{M})$.
			Let $\widetilde{\omega}_1=\Psi_t(\omega_1)$ and let
			$\widetilde{\Delta}\in Reg (\widetilde{\mathbb M}^t)$.
			Let $\Delta=\Psi^{-1}_t(\widetilde{\Delta})$ and let 
			$\eta=\partial\Delta\cap_B\omega_1$.
			Let $Supp (\Delta)=\{a,t\}$.
			\begin{itemize}
				\item [a)] If $\widetilde{\Delta}\in\mathcal{D}(\widetilde{\omega}_1)$
				then $\sigma(\omega_1)=1,\|\eta\|_a\geq 2$ and $\|\eta\|_t\geq 2$
				\item [b)] If $\widetilde{\Delta}\in\mathcal{D}'(\widetilde{\omega}_1)$
				then $\sigma(\omega_1)=1,\|\eta\|_a\geq 2$ and $\|\eta\|_t\geq 2$
			\end{itemize}
		}
	\end{sssection}
	\begin{sssection}[Theorem T]
		\addcontentsline{toc}{subsubsection}{\numberline{\thesubsubsection} 
			Theorem T (Transfer of transferable regions)}	
		{\it 
			\ \\ Let $M\in \mathcal{M}_3(W)$ and let P be a simply connected, connected subdiagram of M with at least one region. Assume that P is  closed 2-banded  in M by bands  $B_r$ and $B_l$. 
			If $\partial P$ has cyclically reduced label U and $P\in\mathcal{M}_3(U)$ then P is transferable.
			In particular, there is a sequence J of I-moves such that $B_r^J\cap B_l^J$ is connected and 
			$B_r^J\cup B_l^J$ is a band-bundle. Moreover 
			$|B_r^J|=|B_r|$ and $|B_l^J|=|B_l|$
		}
	\end{sssection}

	In each of these Theorems $M\in \mathcal{M}_3(W)$ and $\partial M=u\mu v\nu^{-1}$, $u$ and $v$ vertices. 
	The induction hypothesis  claims that if N is a proper subdiagram of M with cyclically reduced boundary 
	label $U_N$ such that $N\in\mathcal{M}_3(U_N)$ then the corresponding result holds for N. We shall denote these hypotheses by $\mathcal H(A),\mathcal{H}(B),
	\mathcal{H}(C)$, $\mathcal{H}(S), \mathcal{H}(L)$ 
	and $\mathcal{H}(T)$, respectively. 
	So we assume in the sequel that
	\begin{equation*}
		\tag {$\mathcal H_2$}
		\mathcal H(A),
		\mathcal{H}(B),
		\mathcal{H}(C),
		\mathcal{H}(S),
		\mathcal{H}(L) \text{ and } \mathcal{H}(T) \text { hold}.
	\end{equation*}
	
	\subsection {The proof of Theorem A}
	\begin{lemma}\ \\
		\addcontentsline{toc}{subsubsection}
		{\numberline{\thesubsubsection} 
			Lemma \quad (If $M\in \mathcal{M}_3(W)$ then poles of bands are not adjacent)}
		Let W be a cyclically reduced word in F, $W=1$ in A and $W\neq 1$ in $A(\Gamma_2)$. 
		Let M be a van Kampen $\mathcal R$-diagram over F, $\mathcal R$ given by (I). 
		Assume that $M\in\mathcal{M}_1(W)$ and all bands are adequate. 
		Let $\Delta_1$ and $\Delta_2$ be regions in $Reg_{4+}(\mathbb M)$, with 
		$Supp(\Delta_1)=Supp(\Delta_2)=\{a,t\},\Delta_1\neq \Delta_2$. 
		Assume that there are bands $B_1$ and $B_2$ connecting adjacent edges $e_1v_1e_2$ of $\partial\Delta_1$ with adjacent edges  $f_1w_1f_2$ of $\partial\Delta_2$, $v_1$ and $w_1$ are vertices,
		such that $\Phi(e_1)=a$ and $\Phi(e_2)=t$. 
		Then $M\not\in \mathcal{M}_2(W)$. In particular, if $M\in \mathcal{M}_3(W)$ then the poles 
		of $B_1$ and $B_2$ cannot be adjacent both on $\partial \Delta_1$ and on $\partial \Delta_2$.
	\end{lemma}
	
	\noindent {\bf Proof}
	
	\noindent
	Assume by way of contradiction that $M\in \mathcal{M}_2(W)$. Then by assumption, $M\in \mathcal{M}_3(W)$. 
	We have $\mathbb E^2\setminus (\Delta_1\cup B_1\cup B_2\cup\Delta_2)=P\cup P_\infty$, where P is the union of the bounded components and $P_\infty$ is the unbounded connected component. ($Int(P)$ is not necessarily connected but P is connected).
	Then $\partial P= \theta_1 v_1\theta^{-1}_2 w_1$,
	where $\theta_i$ are sides of $B_1$ and $B_2$, respectively. See Fig 42. 
	Let $U=\Phi(\theta_1)\Phi(\theta_2^{-1})$. We may assume that U is cyclically reduced, by carrying out diamond moves, if needed. 
	This leaves us in $\mathcal{M}_3(W)$. Also, it follows from Lemma 2.2.1 that $P\in\mathcal{M}_3(U)$. 
	Let $N=B_1\cup B_2\cup P$. Then $|N|<|M|$ and due to Lemma 2.2.1, $N\in \mathcal M_3(V),\ V=\Phi(\partial N)$, cyclically reduced. 
	Therefore by  Theorem T for $N$, $P$ is transferable by a sequence J of I-moves, such that ${M}^J\in\mathcal M_3(W)$, due to Lemma 2.2.1. 
	Hence in $M^J,B_1^J$ and $B_2^J$ are adjacent. $Reg (P^J)=\emptyset$. 
	But then by a sequence $J_1$ of I-moves, we can push $\Delta_1^J$ to $\Delta_2^J$ such that $e_1^{J_1}=f_1^{J_1}$ and $e_2^{J_1}=f_2^{J_1}$, 
	showing that $\Delta_1^{J\cdot J_1}$ and $\Delta_2^{J\cdot J_1}$ are Equivalent. 
	This however violates $M\in \mathcal{M}_2(W)$. 
	Hence $M\not\in \mathcal{M}_2(W)$, as stated.  \hfill $\Box$ 
	\begin{center}
		\begin{tikzpicture}[xscale=.7,yscale=0.5]
			\draw (0,0) ellipse (2 and 4);
			\draw (0,0) ellipse (3 and 5);
			
			\draw[fill=white] (0,4)--(-1,4.75)--(-1,5.5)--(-.5,6.5)--(.5,6.5)--(1,5.5)--(1,4.75)--(0,4) 
			node [label={[label distance=-2mm]below:$v_1$}] {$\bullet$};
			\draw[fill=white] (0,-4)--(-1,-4.75)--(-1,-5.5)--(-.5,-6.5)--(.5,-6.5)--(1,-5.5)--(1,-4.75)--(0,-4)
			node [label={[label distance=-2mm]above:$w_1$}] {$\bullet$};
			\draw (-3,0)--(-2,0) node [rotate=90,
			label={[label distance=-2mm]below:$\theta_1$}] {$>$};
			
			\draw (-2.9,-1)--(-1.9,-1);
			\draw (-2.8,-2)--(-1.7,-2);
			\draw (-2.4,-3)--(-1.3,-3);
			\draw (-2.9,1)--(-1.9,1);
			\draw (-2.8,2)--(-1.7,2);
			\draw (-2.4,3)--(-1.3,3);
			
			\draw (3,0)--(2,0) node [rotate=90,
			label={[label distance=-2mm] above:$\theta_2$}] {$>$};
			\draw (2.9,-1)--(1.9,-1);
			\draw (2.8,-2)--(1.7,-2);
			\draw (2.4,-3)--(1.3,-3);
			\draw (2.9,1)--(1.9,1);
			\draw (2.8,2)--(1.7,2);
			\draw (2.4,3)--(1.3,3);
			
			\node at (0,0) {$P$};
			\node at (0,-5.5) {$\Delta_2$};
			\node at (0,5.6667) {$\Delta_1$};
			\node[left] at (-3,0) {$B_1$};
			\node[right] at (3,0) {$B_2$};

			\draw (-.85,3.6)--(-1.55,4.3);
			\draw (-.85,-3.6)--(-1.55,-4.3);
			\draw (.85,3.6)--(1.55,4.3);
			\draw (.85,-3.6)--(1.55,-4.3);
			
			\node at (-0.9,4.15) {$e_1$};
			\node at (0.85,4.15) {$e_2$};
			
			\node at (-0.9,-4.3) {$f_1$};
			\node at (0.85,-4.3) {$f_2$};
			
		\end{tikzpicture}

		Figure 42
	\end{center}
	\begin{lemma}\ \\
	\addcontentsline{toc}{subsubsection}
	{\numberline{\thesubsubsection} 
		Lemma \quad (No vertices in $\widetilde{\mathbb M}^t$ with valency 2 and $\widetilde{\mathbb M}^t$ satisfies C(4))}
		Let $M\in\mathcal{M}_3(W)$ and assume
		$T(M)\neq\emptyset$. Let $t\in T(M)$. Then
		\begin{itemize}
			\item [(a)]
			$\widetilde{\mathbb M}^t$ has no inner vertices with valency 2, and
			\item [(b)]
			$\tilde{\mathbb{M}}^t$ satisfies the condition C(4).
		\end{itemize}
	\end{lemma}
	
	\noindent{\bf Proof }
	\begin{itemize}
		\item [(a)]
		Suppose by way of contradiction that $\widetilde{\mathbb M}^t$ contains an inner vertex $\widetilde{v}$ with valency 2.
		Then there are regions $\widetilde{\Delta}_1$ and $\widetilde{\Delta}_2$ in $\widetilde{\mathbb M}^t$ such that 
		$\widetilde{v}\in\partial\widetilde{\Delta}_1\cap\partial\widetilde{\Delta}_2=\tilde{\mu}_1\tilde{v}\tilde{\mu}_2$. 
		It follows by arguments similar to those used in the proof of Lemma 5.2.2 that $\Psi_t^{-1}(\widetilde{\Delta}_1\cup
		\widetilde{\Delta}_2\cup
		\widetilde{\mu}_1\cup
		\widetilde{\mu}_2\cup
		\widetilde{v}
		)$
		is connected and simply connected.  
		$\Psi_t^{-1}(\widetilde{\mu}_1\cup
		\widetilde{\mu}_2\cup
		\widetilde{\Delta}_1\cup
		\widetilde{\Delta}_2\cup
		\widetilde{v}
		)$
		is depicted in Fig. 43(a), where $B_1$ and $B_2$ are  t-bands connecting $\Delta_1$ and $\Delta_2$, 
		$\Delta_1=\Psi_t^{-1}(\widetilde{\Delta}_1)$,\ 
		$\Delta_2=\Psi_t^{-1}(\tilde{\Delta}_2)$,\ 
		$\partial\widetilde{\Delta}_1\cap\partial\widetilde{\Delta}_2=$
		$\tilde{\mu}'_1\tilde{\mu}''_1\tilde{v}\tilde{\mu}''_2\tilde{\mu}''_2$,and 
		$|\tilde{\mu}''_1|=|\tilde{\mu}'_2|=1$, 
		$P=\Psi_t^{-1}(\tilde v)$. 
		Let $U=\Phi(\partial P)$. 
		Then U is cyclically reduced because the sides of $B_1$ and $B_2$ are labelled by elements from $C_\Gamma(t)$, while $a\notin C_\Gamma(t)$, $b\notin C_\Gamma(t)$, 
		where $Supp(\Delta_1)=\{a,t\}$ and
		$Supp(\Delta_2)=\{b,t\}$.
		Also, it is easy to see that $P\in\mathcal{M}_3(U)$. Hence, Theorem B applies to P, implying that 
		$a=b$ and that P is a-abelian, since $||U||_a=2$. 
		$(\Phi(\alpha_i)\in\langle a\rangle \setminus \{1\}\, \alpha_i=\partial\Delta_i\cap\partial P_i\,i=1,2)$. 
		Hence, there is a band B connecting $\alpha'_1$ with $\alpha'_2$, $\alpha_1'\subseteq \alpha_1$, 
		$\alpha_2'\subseteq \alpha_2$,
		$|\alpha'_1|=|\alpha'_2|=1$,
		$\alpha'_1$ and $\alpha'_2$ adjacent to the poles of $B_1$
		on $\Delta_1$ and $\Delta_2$, respectively. See Fig. 43(b).
		But then there is a band $B_0$ with poles $\alpha'_1$ and $\alpha'_2$,
		adjacent to the poles of $B_1$ (or $B_2$), respectively,
		contradicting Lemma 5.4.1. 
		Hence $\mathbb M_t$ contains no inner vertices with valency 2.
		\begin{center}
			\begin{tabular}{cc}
				\begin{tikzpicture}[xscale=.6,yscale=0.4]
					\draw (0,0) ellipse (3 and 5);
					
					\draw[fill=white] (-0.5831,3.8334)--(-1,4.75)--(-1,5.5)--(-.5,6.5)--(.5,6.5)--(1,5.5)--(1,4.75)--(0.4998,3.8334);
					\draw[fill=white] (-0.4998,-3.8334)--(-1,-4.75)--(-1,-5.5)--(-.5,-6.5)--(.5,-6.5)--(1,-5.5)--(1,-4.75)--(0.4996,-3.8334);
					\draw (0,0) ellipse (2 and 4);
					\draw (-3,0)--(-2,0);
					\draw (-2.9,-1)--(-1.9,-1);
					\draw (-2.6334,-2.4998)--(-1.7,-2);
					\draw (-2.0668,-3.5831)--(-1.3,-3);
					\draw (-2.9,1)--(-1.9,1);
					\draw (-2.8,2)--(-1.7,2);

					\draw (3,0)--(2,0);
					\draw (2.9,-1)--(1.9,-1);
					\draw (2.5501,-2.4998)--(1.7,-2);
					\draw (2.0668,-3.5831)--(1.3,-3);
					\draw (2.9,1)--(1.9,1);
					\draw (2.8,2)--(1.7,2);
					\draw (2.1501,3.4998)--(1.3,3);
					
					\node at (0,0) {$P$};
					\node at (0,-5.5) {$\Delta_2$};
					\node at (0,5.5) {$\Delta_1$};
					\node[left] at (-3,0) {$B_1$};
					\node[right] at (3,0) {$B_2$};

					\draw (-1.3334,3)--(-2.1666,3.5);

				\end{tikzpicture}\quad
				&
				\quad \begin{tikzpicture}[xscale=.6,yscale=0.4]
					
					\draw (0,0) ellipse (3 and 5);
					
					\draw[fill=white] (-1.0829,3.4169) node (v0) {}--(-1.3332,4.4168)--(-1.3332,5.1668)--(-.5,6.5)--(.5,6.5)--(1.2499,5.2501)--(1.2499,4.5001)--(0.9163,3.5002);
					\draw[fill=white] (-0.9996,-3.4169) node (v1) {} --(-1.3332,-4.5001)--(-1.3332,-5.2501)--(-.5,-6.5)--(.5,-6.5)--(1.3332,-5.1668)--(1.3332,-4.4168)--(1.0827,-3.3336);
					\draw (0,0) ellipse (2 and 4);
					\draw (-3,0)--(-2,0);
					\draw (-2.9,-1)--(-1.9,-1);
					\draw (-2.6334,-2.4998)--(-1.7,-2);
					\draw (-2.0668,-3.5831)--(-1.4666,-2.7501);
					\draw (-2.9,1)--(-1.9,1);
					\draw (-2.8,2)--(-1.7,2);

					\draw (3,0)--(2,0);
					\draw (2.9,-1)--(1.9,-1);
					\draw (2.5501,-2.4998)--(1.7,-2);
					\draw (2.0668,-3.5831)--(1.4666,-2.7501);
					\draw (2.9,1)--(1.9,1);
					\draw (2.8,2)--(1.7,2);
					\draw (2.1501,3.4998)--(1.3,3);

					\draw (0.3,-1.7)--(-0.7,-1.5);
					\draw (0.1,-2.7)--(-0.8,-2.4);
					\draw (2.0668,-3.5831)--(1.4666,-2.7501);
					
					\draw (0.2,1.6)--(-0.7,1.4);
					\draw (2.1501,3.4998)--(1.3,3);

					\node at (1,0) {$B_0$};
					\node at (0,-5.5) {$\Delta_2$};
					\node at (0,5.5) {$\Delta_1$};
					\node[left] at (-3,0) {$B_1$};
					\node[right] at (3,0) {$B_2$};

					\draw (-1.4167,2.7501)--(-2.1666,3.5);

					\draw (v0.center) .. controls (-0.5,1.5) and (-0.5,-1) .. (v1.center);
					\draw (-0.2501,3.95) .. controls (0.5,1.5) and (0.5,-1) .. (-0.0833,-4);
					\draw (-0.8,2.4) -- (0.1,2.7);
					\draw (-0.6,0.6) -- (0.3,0.7)  ;
					\draw (-0.6,-0.4) -- (0.3,-0.5) ;
					\node[xshift=2mm, yshift=3.5mm] at (v0){$\alpha'_1$};
					\node[xshift=2mm, yshift=-4mm] at (v1){$\alpha'_2$};
				\end{tikzpicture}\\
				(a) & (b)
			\end{tabular}
			
			Figure 43
		\end{center}
		\item [(b)]
		If $\widetilde{\Delta}_1\in Reg (\widetilde{\mathbb M}^t)$ 
		is a region, $\tilde\Delta_1=$ $\Psi_t(\Delta_1)$ and  
		$\partial \Delta_1=$ $\tau_1\eta_1\ldots\tau_k\eta_k$, $\Phi(\tau_i)=t^{f_i}$,\quad
		$f_i\neq 0$ and $\Phi (\eta_i)\in H_t$, 
		then $\partial\widetilde{\Delta}_1=\widetilde{\tau}_1\widetilde{v}_1\cdots\widetilde{\tau}_k\widetilde{v}_k$,\ 
		$\widetilde{v}_i$ vertices, $\tau_i$ edges. 
		If $\widetilde{\Delta}_2$ is  a regions of $\widetilde{\mathbb M}^t$ and $\widetilde{\mu}$ is a common boundary path of
		$\widetilde{\Delta}_1$ and $\widetilde{\Delta}_2$ of length at least 2 which contains $\widetilde{v}_i$ then $\widetilde{v}_i$ is an inner vertex with valency 2. 
		Now, by part(a), inner vertices cannot have valency 2. In particular this implies that
		$|\partial\widetilde{\Delta_1}\cap\partial\widetilde{\Delta_2}|=1$ and hence
		$\widetilde{\Delta}_1$ has at least k neighbours in $\widetilde{\mathbb M}^t$ (with multiplicity). By Lemma 1.5.1 $k\geq n(\Delta)\geq 4$. Hence ${\mathbb M}^t$ 
		satisfies the condition C(4).
	\end{itemize}
	\ \hfill $\Box$
	
	\subsubsection{Proof of Theorem A}
	If $\widetilde{\mathbb M}^t$ does not satisfy the condition T(4) then due to Lemma 5.4.2,  
	$\widetilde{\mathbb M}^t$ contains an inner vertex $\widetilde{v}$ with valency 3.
	Let $\widetilde{\Delta}_1,\widetilde{\Delta}_2$ and $\widetilde{\Delta}_3$ be the regions of $\widetilde{\mathbb M}^t$ 
	which contain $\widetilde{v}$ on their boundary. 
	Let $\Delta_i=\Psi_t^{-1}(\widetilde{\Delta}_i),\, i=1,2,3$ and let 
	$Q=\Psi_t^{-1}(\widetilde{v})$. It is not difficult to show by arguments similar to those used in the proof of Lemma 5.2.2 that $\psi_t^{-1}\left(\bigcup\limits_{i=1}^{i=3}{\widetilde{\Delta}}_i\right)$ 
	is connected and simply connected and $\partial Q=\alpha_1\theta_1\alpha_2\theta_2\alpha_3\theta_3,\, \alpha_i=\partial\Delta_i\cap\partial Q$ and $\theta_i$ are sides of 
	bands $B_i,\theta_i=\partial B_i\cap\partial Q, i=1,2,3$, and either $\alpha_i$ is a point or a path with label $a_i^{k_i}, k_i\neq 0, Supp(\Delta_i)=\{t,a_i\},i=1,2,3$. 
	Consider $U:=\Phi(\partial Q)$. Since $Supp(\Phi(\theta_i))\in C_\Gamma(t)$ while if $\alpha_i$ is not a point
	then $\alpha_i$ is labelled with $a_i^{k_i}\notin C_\Gamma(t)$, the only occurrences of letters in U which do not commute with t are the non-trivial labels of $\alpha_i$. In particular
	\begin{equation*}
		\tag {*}
		\left.
		\begin{minipage}{10cm}
			there are at  most 3 disjoint occurrences of letters in U which do not commute with t
		\end{minipage}
		\right\}
	\end{equation*}
	
	\noindent We claim that if $\alpha_i$ and $\alpha_j$ are not points then $a_i=a_j$. Suppose not.
	By Lemma 1.1.c there are at least 2 disjoint occurrences of $a_i$ and 2 disjoint occurrences of $a_j$ in U, hence U has at least 4 such occurrences violating (*). Hence $a_i=a_j$. Denote $a=a_i=a_j$. 
	it is easy to see that U is cyclically reduced and $Q\in\mathcal{M}_3(U)$. In particular, Theorem B applies to Q. ($|Q|<|M|$). 
	Hence, if D is a region of Q with $a\in Supp(D)$ then $n(D)\leq 3$, due to (*). Since $n(D)\neq 3$ hence $n(D)=2$. 
	Hence Q is a-abelian. 
	Consequently, if $\alpha_i$ is not a point then from each edge of $\alpha_i$ there emanates a band L in Q which has one of its poles on $\alpha_i$ and the other on $\alpha_j$. 
	Observe that L and $B_i$ have adjacent poles on $\Delta_i$ and on $\Delta_j$. But this violates Lemma 5.4.1. 
	Hence all $\alpha_i$ are points. 
	But then since 3 is odd, the 3 changes of orientations of the poles of $\beta_i,i=1,2,3$ forces $t=t^{-1}$.
	This violates Lemma 1.1(a). 
	Hence $d(\widetilde{v})\neq 3$. 
	It follows, due to Lemma 5.4.2(b) that $\widetilde{\mathbb M}^t$ satisfies the condition C(4)\&T(4). \hfill $\Box$
	
	\subsection{Preparatory results for the proofs of Theorems S,L,T and B.}
	\begin {lemma}
	\addcontentsline{toc}{subsubsection}{\numberline{\thesubsubsection} Lemma\quad (The case $\widetilde{\mathbb M}^t=\{\widetilde{\Delta}\}$)}
	\ \\ Let $M\in\mathcal{M}_3(W)$ be a connected, simply connected $\mathcal R$-diagram with connected interior.
	Suppose $\partial M= u_1\omega_1 u_2\omega_2^{-1}$ and $T(M)\neq\emptyset$, where $u_1$ and $u_2$ are vertices.
	
	\noindent Let $t\in T(M)$ and consider $\widetilde{\mathbb M}^t$.
	Assume $Reg (\widetilde{M}^t)=\{\widetilde\Delta\}$.
	Let $\Delta=\Psi_t^{-1}(\widetilde{\Delta})$, 
	$Supp(\Delta) = \{t,a\}$ and define $Opt_{\omega_i}(\Delta):=\Psi^{-1}_t(\widetilde{\omega}_i), i=1,2$.
	Assume Theorem S holds for M. Then $\sigma(\omega_i)=1$ for $i=1$ or $i=2$. and $\|\omega_i\|_t\geq 2$.
\end{lemma}
\noindent {\bf Proof}\\
Let $\omega^*_i=\partial\Delta\cap_B\omega_i,\ i=1,2$
and let $u^*_1$ and $u^*_2$ be the endpoint of $\omega^*_1$ and
$\omega^*_2$. Then $\omega^*_1 u^*_1 {\omega^*_2}^{-1} u^*_2$
is a boundary cycle of $\Delta$.
We propose to show that $\Delta\in\mathbb D(\omega^*_i)$ for $i=1$ or for $i=2$.
Then the result follows by part (d) of Theorem S.
Now $|\omega^*_1|+|\omega^*_2|=|\partial\Delta|$. 
Either $|\omega^*_1|\leq |\omega^*_2|$ or $|\omega^*_2|\leq |\omega^*_1|$.
In the first case $|\omega^*_2|\geq n(\Delta)$ due to Lemma 1.3.2 hence $\Delta\in \mathbb{D}(\omega^*_2)$ and in the second case 
$|\omega^*_1|\geq n(\Delta)$, hence $\Delta\in \mathbb D(\omega^*_1)$, and in both cases $\|\omega_i\|_t\geq 2$, as required. \hfill $\Box$

\begin{lemma}
	\addcontentsline{toc}{subsubsection}{\numberline{\thesubsubsection} Lemma\quad (The case $\widetilde{\Delta}\in 
		\mathcal{D}_1(\tilde{\mu})\cup \mathcal{D}_2(\tilde{\mu})$)}
	Let $M\in\mathcal{M}_3(W)$ with boundary cycle $\mu u\nu^{-1} v$, 
	where $u$ and $v$ are vertices.
	Assume that $T(M)\ne\emptyset$ and let $t\in T(M)$.
	Let $\widetilde{\Delta}\in Reg(\widetilde{\mathbb M}^t)$ and let 
	$\Delta=\Psi_t^{-1}(\widetilde{\Delta})$.
	Let $\eta=(\partial\Delta\cap \partial\widehat{Opt}(\Delta))\cap_B \mu$ and let
	$\eta_0=(\partial\Delta\cap \partial{Opt}(\Delta))\cap_B \mu$.
	Assume $\eta$ is connected. Then in each of the following cases
	$\Delta\in \mathbb D(\eta)$ and $\|\Phi(\mu)\|_t\geq 2$
	\begin{itemize}
		\item [(a)] $\widetilde{\Delta}\in \mathcal{D}_1(\widetilde{\mu})$
		\item [(b)] $\widetilde{\Delta}\in \mathcal{D}_2(\widetilde{\mu})$ and there is an $H_t$-band $L$ emanating from $\eta \setminus \eta_0$ and ending on $\mu$.
	\end{itemize}
	Moreover, in Case (a)  $\|\Phi (\mu)\|_t\geq 3$
\end{lemma}

\noindent {\bf Proof}\\
We  have to show each of the following:
\begin{itemize}
	\item [i)] $\partial\Delta$ is a simple closed curve,
	\item [ii)] If $\Delta_1$ and $\Delta_2$ are B-connected regions in $\mathbb M$ then $\partial\Delta_1\cap_B\partial\Delta_2$ and
	$\partial\Delta_2\cap_B\partial\Delta_1$ are connected and
	$\|\partial\Delta_1\cap_B\partial\Delta_2\|=
	\|\partial\Delta_2\cap_B\partial\Delta_1\|=1
	$.
	Also, $\partial\Delta\cap\eta$ is connected,
	\item [iii)] $\|\partial\Delta\cap\eta \|\geq n(\Delta)$,
	\item [iv)] $|\Phi(\xi)|\leq |\Phi (\partial \Delta\cap
	\eta)|$, where $\xi$ is the complement of $\partial\Delta\cap\eta$ on $\partial\Delta$,
	\item [v)] $\|\Phi(\mu)\|_t\geq 2$.
\end{itemize}
We prove these in turn, (Notice that $\partial\Delta\cap\eta=\eta$)
\begin{itemize}
	\item [(i)] This follows by Proposition 2.1.7.
	
	\item [(ii)] 
	Suppose by way of contradiction that $\|\partial\Delta_1\cap_B\partial\Delta_2\|\geq 2$.
	Then there are bands $L_1$ and $L_2$, $L_1$ a t-band and $L_2$ an $H_t$-band such that one of the poles of $L_i$ is on $\partial \Delta_1$, and the other on $\partial\Delta_2, i=1,2$.
	Now, the poles of $L_1$ and $L_2$ are adjacent on $\partial \Delta_1$ and on $\partial\Delta_2$:
	If not then consider $N:=\Delta_1\cup \Delta_2\cup L_1\cup L_2$.
	We have $\mathbb E^2\setminus N=Q_1\cup \cdots\cup Q_k\cup Q_\infty$,\ $k\geq 1$,\ $Q_i$ are bounded, connected, simply connected $i=1,\ldots,k$,\ $Q_\infty$ the unbounded component.
	We may assume that $\Phi(\partial Q_i) $ is cyclically
	reduced and $Q_i\in \mathcal{M}_3(U_i), U_i=\partial Q_i,
	i=1,\ldots,k$.
	Hence Theorem A applies to $Q_i,\ i=1,\ldots, k$. 
	Hence, if $Q_i$ 
	is not t-abelian then $\widetilde{Q}_i$ satisfies the condition C(4)\& T(4).
	It follows, due to Lemma 1.2.2, that $\widetilde{Q}_i^t$ has a Greendlinger region $\widetilde{K}$ such that $K$ is B-connected to 
	$\partial\widetilde{Q}_i^t, i=1$ or $i=2$, and $\|\partial K\cap_B\partial\Delta_i\|\geq 2$. 
	Consequently, there are bands $B_1$ and $B_2$ with adjacent poles on $\partial\Delta_i$ and $\partial K$, violating Lemma 5.4.1.
	Hence $Q_i$ is t-abelian.
	Consequently, from every t-edge of $\partial Q_i\cap\partial\Delta_1$ there emanates a t-band L which ends on
	$\partial \Delta_2$.
	In particular we may assume from the outset that $L_1$ and $L_2$
	are  the closest pair of a t-band and an $H_t$-band which connect $\Delta_1$ with $\Delta_2$ in the sense that their poles on $\partial\Delta_1$ is closest, then L violates the assumption.
	Hence we may assume that $L_1$ and $L_2$ have adjacent poles, both on $\partial\Delta_1$ and on $\partial\Delta_2$.
	But this violates Lemma 5.4.1. 
	Hence 
	\begin{equation*}
		\tag* {(1}\|\partial\Delta_1\cap_B\partial\Delta_2\|=1.
	\end{equation*}		
	Consequently
	\begin{equation*}\tag*{(2}
		\|\xi\|_t\leq i(\widetilde{\Delta}) 
	\end{equation*}
	and
	\begin{equation*}\tag*{(3}
		\|\eta_0\|_t\geq n(\Delta)-i(\widetilde{\Delta})
	\end{equation*}
	
	Observe that if $V\in H_t*\langle t\rangle,||V||\geq 3$ and V starts and ends with $t^{\pm 1}$, then 
	$||V||\geq ||V||_t+(||V||_t-1)$, i.e. 
	\begin{equation*}\tag {*}
		||V||\geq 2||V||_t-1. 
	\end{equation*}
	Hence,
	\begin{equation*}\tag*{(4} 
		\|\eta_0\|\geq 2 \|\eta_0\|_t-1\geq 2n-2i(\widetilde{\Delta})-1
	\end{equation*}
	Consequently , if we write $n$ for $n(\Delta)$ then
	\begin{align}
		\|\eta_0\|\geq 2n-3\geq n+1,  \text{\quad if \quad}  
		i(\widetilde{\Delta})=1 \tag* {(5a}\\
		\|\eta_0\|\geq 2n-5\geq n-1,  \text{\quad if \quad} 
		i(\widetilde{\Delta})=2 \tag* {(5b}
	\end{align} 
	Hence (ii) holds for the case $i(\widetilde{\Delta})=1$ and also (iii) holds for the case $i(\tilde{\Delta})=2$ 
	provided that there is an $H_t$-band L emanating from $\eta\setminus\eta_0$ and ending on $\mu$, because then
	$\|\eta\|\geq\|\eta_0\|+1\geq n$.
	\item [(iii)] This follows by the assumption.
	\item [(iv)] Follows by Lemma 1.5.1, (2, (5a and (5b.
	\item[(v)] Follows by part (iii) and the fact that 
	$\Phi(\eta_0)$ contains a subword of t-length 
	$n-1\geq 4-1=3$.
\end{itemize}
\ \hfill $\Box$

\noindent {\bf Remark}\ Lemma 5.5.2 remains true if we replace $\mu$ by $\nu$.

\begin{lemma}
	\addcontentsline{toc}{subsubsection}{\numberline{\thesubsubsection} Lemma\quad ($\tilde{\mathbb M}^t\in O.L.$)}
	Let $M\in\mathcal{M}_3(W)$ be a connected, simply connected R-diagram over F with connected interior.
	Let $\partial M=v_1\omega_1 v_2\omega_2^{-1}$, $v_1$ and $v_2$ vertices and $\omega_1$ and $\omega_2$ boundary paths.
	Let $t\in T(M)$ and assume that $|\widetilde{\mathbb M}^t|\geq 2$ and $\widetilde{\mathbb M}^t\in O.L.(\widetilde{\omega_1},\widetilde{\omega_2}),\ \widetilde{\omega_i}=\Psi_t(\omega_i),\ 
	\widetilde{\mathbb M}^t=\langle \widetilde{\Delta}_1,\ldots,\widetilde{\Delta}_k),\ k\geq 2$.
	If $\sigma(\omega_1)=0$ and Theorem S holds for M then each of the following holds
	\begin{itemize}
		\item [(a)] $\sigma(\partial Opt (\Delta_i)\cap_B \omega_2)=1,\ 
		i=1,k$
		\item [(b)] ${\|\omega_2\|}_t\geq 3$
	\end{itemize}
\end{lemma}

\noindent {\bf Proof}\\
\begin{itemize}
	\item [(a)] 
	Let $\gamma_3=\partial\Delta_1\cap_B\partial\Delta_2$,
	$\gamma_2=\partial\Delta_1\cap_B \omega_2$ and
	$\gamma_1=\partial\Delta_1\cap_B \omega_1$.\\
	Then $\gamma_1 \gamma_2 \gamma_3$ is a boundary cycle of $\Delta_1$.
	Since $\widetilde{\Delta}_1\neq \widetilde{\Delta}_2$, 
	$\|\partial\Delta_1\cap_B\partial\Delta_2\|=1$,
	due to part (iii) in Lemma 5.5.2.
	Hence $\Phi(\gamma_3)=t^e,\ e\neq 0$.
	Consequently, Proposition 1.5.2 applies to $\Delta_1$.
	Hence $N(\gamma_i)$ contains a subdiagram $S_1\in 
	\mathcal{D}'(\omega_2)$, since $\sigma(\omega_1)=0$.
	The same argument applies to $\Delta_k$ yielding $S_k$.
	Hence (a) holds true.
	\item [(b)] Follow the notation of part (a) of the Lemma.
	Due to Theorem S part (a) $\|\omega_2\|_t\geq \|\omega_2^*\|_t$ 
	where $\omega_2^*=\partial S_1\cap \partial Opt(S_1)$.
	Here $Opt(S_1)$ is defined as $Opt(\Delta)$,\ $\Delta$ replaced by $S_1$.
	Since $S_1\in \mathcal{D}'(\omega_2)$, 
	$\|\partial S_1\cap\partial Opt (S_1)\|\geq n(S_1)$.
	Hence  $\|\omega_2\|_t\geq \|\omega_2^*\|_t\geq \frac {n(S_1)} 2\geq \frac 4 2=2$.
	Notice that $n(S_1)=n(\Delta_1)$.
	By the same argument for $S_k\subseteq\Delta_k$, 
	$\  \|\partial S_2\cap \partial Opt(S_k)\|\geq n(\Delta_k)$.
	It follows that $\|\omega_2\|_t\geq 2+2-1=3$.\\ \text{ }\hfill $\Box$
\end{itemize}

\begin{lemma}
	\addcontentsline{toc}{subsubsection}{\numberline{\thesubsubsection} Lemma\quad ($\sigma(\omega_1)=0$ implies $Q$ is a-abelian)}
	\ \\ Let notation and assumptions be as in Lemma 5.2.2. 
	Let $\eta_i=\partial \Delta_i\cap_B\omega_i,\ i=1,i=2$. 
	Let $\omega_1=\gamma$ and $\omega_2=\theta_1\alpha_1\theta_3\alpha_2\theta_2$.
	(Thus $|Q|<|M|$ and  $\omega_1 u\omega_2^{-1}$ is a boundary cycle of Q.)
	If $\sigma(\omega_1)=0$ then Q is a-abelian.
\end{lemma}

\noindent{\bf Proof} By induction on $|Q|$, the case $|Q|=1$ being clear in view of Case 1 below.
Assume by way of contradiction that $Q$ is not a-abelian.
Then $T(\widetilde{Q}^a)\neq\emptyset$.
Without loss of generality we may assume that $\Phi(\partial Q)$ is cyclically reduced and due to Lemma 2.2.1, $Q\in\mathcal{M}_3(U)$, where
$U=\Phi(\partial Q)$.
Hence Theorem A applies to Q, implying that $\widetilde{Q}^a$ satisfies
the condition C(4)\& T(4). Consequently, by Lemma 1.3.2, one of the following holds:
\begin{description}
	\item [Case 1] $\widetilde{Q}^a=\{\widetilde{\Delta_0}\}$
	\item [Case 2] $|\widetilde{Q}^a|\geq 2$ and $|\mathcal{D}(\widetilde{Q}^a)|=2$
	\item [Case 3] $|\widetilde{Q}^a|\geq 3$ and $|\mathcal{D}(\widetilde{Q}^a)|=3$
	\item [Case 4] $|\widetilde{Q}^a|\geq 4$ and $|\mathcal{D}(\widetilde{Q}^a)|\geq 4$
\end{description}
We derive a contradiction for each of these cases in turn.

\begin{description}
	\item[Case 1] Since $\omega_2=\theta_1\alpha_1\theta_3\alpha_2\theta_2$, (See 5.2.2 ) where $\Phi(\theta_i)\in C_\Gamma(t)$ and $\Phi(\alpha_i)=a^{f_i},\ f_i\neq 0$, 
	we have $\|\omega_2\|_a=2$.
	Let $Supp (\Delta_0)=\{a,b\}$.
	Let $\eta_i=\partial\Delta_0\cap_B\omega_i,\ i=1$ or $i=2$.
	Then $\|\eta_2\|=\|\eta_2\|_{\langle a,b \rangle}$, hence 
	$ 3\leq \|\eta_2\|_{\langle a,b \rangle}\leq 5$. (See (*) in the proof of Lemma 5.5.2).
	If $\|\eta_2\|\leq 4$, then $\|\eta_1\|\geq 2n(\Delta_0)-4\geq n(\Delta_0)$, i.e. $\|\eta_1\|\geq n(\Delta_0)$.
	It follows, due to Lemma 1.5.1 that $\Delta_0\in\mathbb D(\eta_1)$, hence, by part (d) of Theorem S applied to $Q$, $\sigma(\omega_1)=1$.
	This, however, violates the assumption that $\sigma(\omega_1)=0$.
	We claim that $\|\eta_2\|=5$ is not possible.
	Thus, assume $\|\eta_2\|=5$.
	Then due to the alternating nature of $\Phi(\partial\Delta_1)$ there is a b-band emanating from $\partial\Delta_0$ and ending on $\theta_i,\ i=1,2,3$.
	In particular $b\in Supp(\Phi(\theta_1))\subseteq C_\Gamma(t)$.
	Now, it follows by Theorem S part (a) applied to $Q$, that there is an a-band $L_0$ emanating from $\eta_1$ (and ending on $\omega_1$) and a t-band-bundle B
	such that $L_1\subseteq B$ ($L_1$ is the t-band which contains $\theta_1$ on its boundary)
	and $L_0$ has a pole adjacent to the pole of B on $\partial\Delta_1$.
	Hence either $b\in C_\Gamma(a)$ or $\mathcal{D}_2^2(\omega_1)\neq \emptyset$. See Fig. 44, where $\{E_1,E_2\}\subseteq\mathcal{D}_2^2(\omega_1)$.
	Both cases lead to contradictions to $\{a,b\}\subseteq T(Q)$ and $\sigma(\omega_1)=1$, respectively. 
	
	\begin{center}
		\begin{tabular}{cc}
			\begin{tikzpicture}[scale=.4]
				
				\draw (0,0) circle (2);
				\draw (9,0) circle (2);
				
				\draw (2,0)--(7,0);
				\draw (1.7,-1)--(7.2,-1);
				\draw (5.5,0)--(5.5,-1);
				\draw (3.5,0)--(3.5,-1);
				\draw (4.5,0)--(4.5,-1);

				\draw (-1.125,1.625)--(-1.125,7)--(1,7)--(1,1.7);
				\draw (9,2)--(9,7)--(10,7)--(10,1.7);
				
				\draw (-1.125,3)--(1,3);
				\draw (-1.125,4.5)--(1,4.5);
				\draw (-1.125,6)--(1,6);
				
				\draw (9,3)--(10,3);
				\draw (9,4.5)--(10,4.5);
				\draw (9,6)--(10,6);

				\draw (1,7)--(10,7);
				
				\node at (0,0) {$\Delta_1$};
				\node at (9,0) {$\Delta_2$};

				\node at (0.025,1) {$a$};
				\node[left] at (7.5,1.5) {$\alpha_2$};
				
				\node[above] at (5,3.5) {$Q$};
				\node[above] at (0,7.8) {$L_1$};

				\draw (0,7) -- (0,2);
				\node at (1.75,4.25) {$b$};
				\node[right] at (1.5,1.5) {$\alpha_1$};
				\node[above] at (-0.5,7) {$\alpha$};
				\node[above] at (0.5,7) {$t$};
			\end{tikzpicture} 
			&
			\begin{tikzpicture}[scale=.4]
				
				\draw (0,0) circle (2);
				\draw (9,0) circle (2);
				
				\draw (2,0)--(7,0);
				\draw (1.7,-1)--(7.2,-1);
				\draw (5.5,0)--(5.5,-1);
				\draw (3.5,0)--(3.5,-1);
				\draw (4.5,0)--(4.5,-1);

				\draw (0,2)--(0,7)--(1,7)--(1,1.7);
				\draw (0,3)--(1,3);
				\draw (0,4.5)--(1,4.5);
				\draw (0,6)--(1,6);

				\draw (-1.5,1.375)--(-3.125,5.875)--(-2,6.5)--(-0.5,2);
				\draw (-2,2.5)--(-0.75,3);
				\draw (-2.25,3.625)--(-1.25,4.125);
				\draw (-2.75,4.75)--(-1.625,5.375);

				\draw (9,2)--(9,7)--(10,7)--(10,1.7);
				\draw (9,3)--(10,3);
				\draw (9,4.5)--(10,4.5);
				\draw (9,6)--(10,6);

				\draw (1,7)--(10,7);
				
				\node at (0,1) {$a$};
				\node at (2,1.625) {$\alpha_1$};
				\node[left] at (7.5,1.5) {$\alpha_2$};
				
				\node[above] at (5,3.5) {$Q$};
				
				\node at (-3,7) {$t$};
				\node at (0.5,7.75) {$t$};
				\draw (0,0) node {$\Delta_1$} circle (2);
				\draw (9,0) node {$\Delta_2$} circle (2);			
				\node at (-2.35,5.65) {$E_1$};
				\node at (0.5,6.5) {$E_2$};
			\end{tikzpicture}
			\\
			(a) & (b)
		\end{tabular}
		
		Figure 44
	\end{center}

	\item[Case 2] $\widetilde{Q}^a\in O.L. (\omega_1,\omega_2)$,
	by Lemma 1.3.2. 
	Hence by Lemma 5.5.3, $\|\omega_2\|_a\geq 3$.
	This  contradicts $\|\omega_2\|_a=2$, which we showed in Case 1.
	Hence this case cannot occur.
	\item[Case 3] It follows by Lemma 5.5.2 and Lemma 1.3.2 that $\|\omega_2\|_a\geq 3$, 
	again contradicting $\|\omega_2\|_a=2$. Hence this case cannot occur.
	\item [Case 4] By Lemma 1.3.2 and the P.H.P one of the following holds
	\begin{itemize}
		\item [i)] $\mathcal{D}(\widetilde{\omega}_1)\neq \emptyset$
		and $\mathcal{D}(\widetilde{\omega}_2)\neq \emptyset$
		\item [ii)] $|\mathcal{D}(\widetilde{\omega}_2)|\geq 2$
	\end{itemize}
	
	\noindent In the first case, by Theorem L and the induction hypothesis for Q, $\sigma(\omega_1)=1$, violating $\sigma(\omega_1)=0$.
	In the second case it is easy to see that if $\widetilde{\Delta}_1$ and $\widetilde{\Delta}_2$ are regions in $\mathcal{D}(\widetilde{\omega}_2)$ then
	$\|\eta_1\|_a\geq 2$, hence $\|\omega_2\|_a\geq 2+2-1=3$.
	This contradicts $\|\omega_2\|_a=2$. 
\end{description}
\ \hfill $\Box$

\noindent {\bf Corollary}\quad 
{\it Let notation and assumptions be as in the Lemma.
	If $\widetilde{\Delta}\in\mathcal{D}_2(\widetilde{\mu})$ then 
	$\Delta\in \mathbb D(\mu)$.}

\noindent {\bf Proof\quad} Due to the Lemma, condition (b) of Lemma 5.5.2 (b) is satisfied.
Hence the result follows by part (b) of Lemma 5.5.2.\hfill $\Box$

\vspace{10pt} \noindent
\subsection{Proof of Theorem S}
\begin{itemize}
	\item [(a)] Let $\tau_1 w_1 \alpha w_2\tau _2$ be a subpath of $\partial\Delta\cap Opt(\Delta)$, $w_1$ and $w_2$ vertices, 
	$\Phi(\tau_1)=t^{\pm 1}$, $\Phi(\tau_2)=t^{\pm 1}$,
	$\Phi(\alpha)=a^f, f\neq 0$.
	It follows from the definition of $Opt(\Delta)$ that there are t-bands $L_1$ and $L_2$ emanating from $\tau_1$ and $\tau_2$ respectively and ending on $\partial M$.
	Let $\theta_1$ be the side of $L_1$ which has $w_1$ as an endpoint and let $\theta_2$ be the side of $L_2$ which has $w_2$ as an endpoint.
	Then $Supp (\Phi(\theta_i))\subseteq C_\Gamma(t)$,\ $i=1$ and $i=2$. Let $z_1$ and $z_2$ be the endpoints of $\theta_1$ and $\theta_2$ respectively, on $\partial M$ 
	and let $\gamma$ be the subpath of $\partial M\cap\partial Opt(\Delta)$ which starts  at $z_1$ and ends at $z_2$.
	Then $\sigma(\gamma)=0$ as $\gamma\subseteq \omega_1$ and  $\sigma(\omega_1)=0$. 
	Also, $U:=\Phi (\theta_1\alpha\theta_2\gamma^{-1})$ 
	is cyclically reduced, since $a\not\in C_\Gamma(t)$. 
	Let Q be the bounded component of $\mathbb E^2\setminus (\theta_1\alpha\theta_2\gamma^{-1})$.
	Assume first that $\partial Q 
	(=\theta_1\alpha\theta_2\gamma^{-1})$ is a simple closed path. 
	Hence $Int(Q)$ is connected.
	It follows from Lemma 2.2.1 that $Q\in\mathcal{M}_3(U)$.
	Hence Theorem A applies to $Q$ and it follows by the proof of Lemma 5.5.4 that if Q is not a-abelian then one of the 4 cases mentioned in Lemma 5.5.4 occurs.
	From this point on we follow the proof of Lemma 5.5.4, with 
	$\|\omega_2\|_a=1$ (here) instead of $\|\omega_2\|_a=2$, showing that each case leads to a contradiction, hence $Q$ is a-abelian.
	More precisely, in Case 1, by Lemma 5.5.1, in case 2, by Lemma 5.5.3,
	in Case 3, by Lemma 5.5.2 and in Case 4, by Theorem L.
	Therefore if $\beta=\beta_1\ldots\beta_k$,\ 
	$\Phi(\beta_i)=a^{\pm 1}$ then from every $\beta_i$ there emanates an a-band which ends on $\gamma$ , since $a\not\in C_\Gamma(t)$.
	(Notice that it cannot end on $\beta$.)
	If $Int (Q)$ is not connected then a similar argument applies to each component.
	\item [(b)] Let $B_1$ and $B_2$ be bands emanating from 
	$\partial\Delta\cap\partial Opt(\Delta)$ and ending on $\omega_1$.
	Assume that their poles on $\partial\Delta$ are adjacent.
	Let $\xi_1$ and $\xi_2$ be the sides of $B_1$ and $B_2$ respectively,
	with $w:=\xi_1\cap \xi_2\cap \partial\Delta\neq\emptyset$ a point.  
	We claim that $\xi_1\cap\xi_2$ is connected.
	Suppose not and consider $\mathbb{E}^2\setminus(\xi_1\cup\xi_2)$.
	Let P be a bounded connected component of it.
	Then $\partial P=(\partial P\cap \xi_1)\cup (\partial P\cap\xi_2)$.
	Let $V:=\Phi(\partial P)$.
	Then applying diamond moves at the endpoints of $\partial P\cap
	\xi_1$ if needed, we may assume that V is cyclically reduced.
	It follows from Lemma 2.2.1 that $P\in \mathcal{M}_3(V)$.
	Hence Theorem T applies to P ($|P|<|M|$).
	Hence by a sequence $J_0$ of I-moves, P is transferable
	beyond $B_1$ and $B_2$, such that $P^J$ contains no regions.
	Doing so independently for every connected component $P_i$, 
	by $J_i,\ i=1,\ldots,r$ we get that $(B_1\cap B_2)^J$ is connected, where $J=J_1\circ J_2\circ \cdots\circ J_r$.
	
	\noindent Consider now $\mathbb E^2\setminus\{L_1,\ldots,L_k\}$.
	Let $Q_1,\ldots,Q_m$ be its bounded connected components. \\
	Assume $k=3$ and $\partial L_1\cap\partial L_2$ is connected.
	Assume that for some $i,\ 1\leq i\leq m$, $Q_i$ occurs between $L_2$ and $L_3$ such that $\partial Q_i=$
	$(\partial Q_i\cap\partial L_2)\cup(\partial Q_i\cap\partial L_3)$.
	Suppose that we would like to transfer a region D from $Q_i$ via $L_2$.
	In principle it is possible that a previous step of transferring regions D occurred between $L_1$ and $L_2$ and was transferred to $Q_1$.
	Hence now we transfer it back and by this create a loop in the process.
	The point is that because $\partial L_1\cap \partial L_2$ is connected, 
	we may regard $L_1\cup L_2$ as a single band and this enables us to transfer D further beyond $L_1$, out of $\widehat{Opt}(\Delta)$ and 
	by this to "empty" (by $\mathcal{I}$-moves J) each $Q_i$ so to get $L_1^J\cap L_2^J$ and $L_2^J\cap L_3^J$ simultaneously connected.
	From this it is clear how to prove for general $k$.
	\item[(c)] (i) is clear from the construction of $Z$ 
	and the definition of $Q_i$. \\
	(ii) and (iii) It follows from Lemma 2.2.1 that $Q_i\in \mathcal{M}_3(U_i)$.
	It remains to show that either $Q_i=\{D\}$ with $n(D)=2$ or $Q_i$ is transferable.
	Thus, Theorem C is applicable to $Q_i$, implying that $\sigma(\rho_i\xi_i)=1$.
	(Recall that $\sigma(\gamma)=0$).
	Consequently, if $Q_i$ is not abelian then 
	$N(\rho_i\xi_i)$ contains a transferable region.
	Hence, by applying a finite number of I-moves we can transfer  all the regions D with $n(D)\geq 4$ beyond $L_i$
	and/or $L_{i+1}$ and either remaining with an abelian diagram 
	$Q'_i$ or with a diagram with $\mathcal{D}^2_2(\omega'_1)=\emptyset$, 
	where $\omega'_1$ is obtained from $\omega_1$ by the appropriate I-move.
	Applying Lemma 4.1.2, the result follows.
	Finally, since either every region D in Q is transferable via
	$L_j\cup L_{j+1}$ or $Q=\{D\}$, $Supp(\partial D)\subseteq C_\Gamma(a)\cup C_\Gamma(t)$. 
	But $\{a,t\} \cap(C_\Gamma(a)\cup C_\Gamma(t))=\emptyset$. 
	Hence $\{a,t\}\cap Supp(D)=\emptyset$.
	\item [(d)] We follow notation of part (c).\\
	(i) By assumption, $\Delta\in\mathbb D(\eta)$.
	Hence, $\|\partial\Delta\cap\eta\|\geq \|n(\Delta)\|\geq 4$.
	Consider $\widehat{Opt}(\Delta)$.
	Let $\{L_1,\ldots,L_k\}$ be the bands which emanate from $\eta$ and end on $\omega_1$.
	Consider $L_1\cup Q_1\cup L_2$. If $Q_1\neq \{D\}$ with
	$n(D)=2$ then we can transfer $Q_1$ beyond $L_1$ and $L_2$
	via a suitable sequence J of I-moves such that in $M^J$ 
	configurations (1),(2),(3) and (4), in Fig. 41, occur.
	In Conf. (1), $\mathcal{D}^2_2(\omega_1) \geq 2$, in configurations
	(2) and (3) $\mathcal{D}^2_2(\omega_1) \geq 1$
	and in configuration (4), $\mathcal{D}^2_2(\omega_1) \geq 0$.
	Similar arguments apply to $L_{k-1}\cup Q_{k-1}\cup L_k$.
	Notice that if $L_1\cup Q_1\cup L_2$ and 
	$L_{k-1}\cup Q_{k-1}\cup L_k$ are in configurations 1,2 or 3 then
	$|\mathcal{D}^2_2(\omega_1)| \geq 2$.
	This is part 2) in the definition of $\sigma(\omega_1)=1$.
	Hence, we have to show that if $\mathcal{D}^2_2(\omega_1) \leq 1$ then
	$\sigma(\omega_1)=1$, which contradicts $\sigma(\omega_1)=0$.
	
	\noindent Let $L_1\cup L_2$ be in configuration i and 
	$L_{k-1}\cup L_k$ in configuration j, $1\leq i\leq 4$ and $1\leq j\leq 4$.
	If $\|\mathcal{D}_2^2(\omega_1)|\leq 1$ then either $i=4$ or $j=4$ or both.
	Assume $L_1\cup L_2$ is in configuration 4 and 
	$L_{k-1}\cup L_k$ in configuration $j\leq 3$.
	Then $L_1\cup L_2$ is a homogeneous band-bundle.
	Let $\mathbb B_0$ be the maximal band-bundle that contains $L_1\cup L_2$.
	Then we claim that $\mathbb B_0=L_1\cup\cdots\cup L_{k-1}$.
	Thus, suppose $\mathbb B_0=L_1\cup\cdots\cup L_{l}$,\ 
	$l<k-1$ and consider $L_{l+1}$.
	If $Reg (Q_l)\neq \emptyset$ then we transfer $Q_l$
	via $L_l\cup L_{l+2}$ by a sequence J of I-moves, such that 
	$Reg (Q^J_l)=\emptyset$.
	Hence $L_l\cup L_{L+1}$ is in configurations 1,2 or 3,
	since $|\mathcal{D}_2^2(\omega_1)|=1$ 	implies that $l=k-1$.
	Observe that we can push up $\Delta$ by a sequence $J_1$ of 
	I-moves along $\mathbb B_0$
	and $\|\partial\Delta^J\cap\omega_1\|\geq n(\Delta)-2\geq 2$.
	This is part 3 in the definition of $\sigma(\omega_1)=1$.
	So it remains to deal with the case when both $L_1\cup L_2$ and 
	$L_{k-1}\cup L_k$ are in configuration 4.
	Let $\mathbb B_1$ be the maximal band-bundle which contains $L_{k-1}\cup L_k$.
	Then $\mathbb B_0\cup\mathbb B_1= Opt(\Delta)$.
	if $\mathbb B_0=\mathbb B_1$ then $Opt(\Delta)$ is a band-bundle and by pushing up $\Delta$ via $Opt(\Delta)$ 
	we get $\Delta^{J_2}\in \mathcal{D}(\omega_1)$.
	But then by part (4) of the definition of $\sigma(\omega_1)$,
	$\sigma(\omega_1)=1$.
	Finally, if $\mathbb B_0=\mathbb B_1$ then either 
	$\|\eta\cap \partial\mathbb{B}_0\|\geq 2$ or 
	$\|\eta\cap \partial\mathbb{B}_1\|\geq 2$
	and we fulfil part (3) of the definition of $\sigma(\omega_1)$,
	showing that $\sigma(\omega_1)=1$. If one of the $Q_i$ consist of a  single region $D$ with $n(D)=2$ then
	$|\mathcal{D}_2^2(\omega_1)|\geq 1$ and the above arguments apply to this case.\\
	(ii) The case when $\widetilde{\mathbb M}^t\in O.L.(\omega_1,\omega_2)$
	is dealt with similarly.
	We omit details.
\end{itemize}
\ \hfill $\Box$


\subsection{Proof of Theorem L}
\begin{itemize}
	\item [(a)] 
	If $|\widetilde{\mathbb M}^t| = 1$ then by Lemma 5.5.1, $\Delta\in \mathbb D(\eta)$.
	Assume $|\widetilde{\mathbb M}^t| \geq 2$. 
	By theorem S part (d), it is enough to show that 
	$\Delta\in \mathbb D(\eta)$. 
	By Lemma $i(\widetilde{\Delta})\leq 2$. 
	If $i(\widetilde{\Delta})=1$ then by Lemma 5.5.2 (a) $\Delta\in\mathbb D(\eta)$. 
	If $i(\widetilde{\Delta})=2$ then by Lemma 1.3.4,
	$\partial\widetilde{\Delta}\cap\widetilde\omega_1$ has an endpoint $\tilde v$ with valency 3. 
	Let $Q=\Psi^{-1}_t(\tilde{v})$. 
	If $\sigma(\omega_1) =0$, then $Q$ is a a-abelian, by Lemma 5.5.4, 
	hence there is an a-band-bundle emanating from 
	$\alpha_1$ and ending at $\omega_1$. 
	(We follow the notation of Lemma 5.5.4) 
	Consequently $\Delta \in\mathbb {D}(\eta)$ ,by Lemma 5.5.2 (b), as required. 
	Since $n(\Delta)\geq 4$, the rest follow. 
	\item [(b)] Let $\widetilde{\mathbb M}^t=\langle\widetilde{\Delta}_1,\ldots,
	\widetilde{\Delta}_k\rangle, k \geq 2$ and consider $\Delta_1$. 
	By Proposition 1.5.2, $\partial S \cap \partial \Delta_1$ is connected,
	$\|\partial S\cap\partial\Delta_1\|\geq n(\Delta_1)$  and 
	$|\xi|\leq|\partial S\cap\partial\Delta_1|$,  where $\xi$ is the complement of $\partial S\cap\partial\Delta_1$,  on $\partial\Delta_1$.
	Consequently, $S\in\mathbb D^v(\eta)$. 
	Therefore, the result follows by Theorem S part (d). 
	The rest follows from the fact that $\|\partial S\cap\partial\Delta_1\|\geq n(\Delta_1)\geq 4$, hence 
	$\|\partial S\cap\partial\Delta_1\|_a\geq 2$ 
	and $\|\partial S\cap\partial\Delta_1\|_t\geq 2$, due to the alternating nature (in a and t) of the boundary label of $\Delta_1$.
\end{itemize} 
\ \hfill $\Box$
\vspace{10pt}
\subsection{Proof of Theorem T}
By induction on $|P|$.
Let $\omega_r=\partial P\cap\partial B_r$ and let 
$\omega_l=\partial P\cap\partial B_l$.
Then $u\omega_r v\omega_l^{-1}=\partial P$,
$u$ and $v$ boundary vertices.
\noindent We start with all the transferable regions from $Reg_2(P)$ by applying successively a sequence $J_1$ of I-moves.
Let $P_1=P^{J_1}$ be the resulting diagram. 
Assume $Reg(P_1)\neq \emptyset$.
If $P_1\not\in\mathcal{M}_3(U_1)$ then we can replace $P_1$ by 
$P'_1\in\mathcal{M}_3(U_1)$ and make the same, following argument. 
Hence, we may assume that
$P_1\in\mathcal{M}_3(U_1)$, where $U_1$ is a boundary label of $P_1$, which we may assume to be cyclically reduced, 
Since we removed only regions of $Reg_2(M)$ and this cannot increase
$Reg_{4^+}(M)$ and $Reg_{4^+}(\mathbb M)$,
Theorem C applies to $P_1$ 
and implies that $\sigma(\omega^*_1) = 1$ or $\sigma(\omega^*_2) = 1$, where 
$\omega^*_i=\omega_i^{J_1},\ i=1,2$. 
Consequently, there is a region $\Delta$ with $n(\Delta)\geq 4$ or a submap $S\in \mathbb{D}^v(\omega^*_i)$ which can be transferred from $P_1$. 
See Remark following Theorem S.
Let the resulting diagram obtained from $P_1$ denoted by $P_2$.
Notice that $|P_2|<|P|$. 
Then we may continue transferring from P all the regions $\Delta$ with 
$n(\Delta)\geq 4$, ending up either with a diagram without regions proving the theorem, or ending up with an abelian diagram. 
In this case Lemma 4.1.2 completes the proof. \hfill $\Box$

\subsection{The Proofs of Propositions 2,3 and 4 relying on Theorems A,B and C}
\subsubsection* {Proof of Proposition 2} 
Define $N_t = \widetilde{\mathbb M}^t$. Then parts (a) and (b) follow from the definition of ${\mathbb M}^t$. Part (c) follows from the definition of $\psi_t$ and part (d) follows from Theorem A.
\hfill $\Box$ 
\subsubsection*{Proof of Proposition 3}
Part (a) follows from Theorem A. Parts (b) and (c) follow from part (a) and Proposition 1. \hfill $\Box$ 
\subsubsection*{Proof of Proposition 4} 
First we show that bands are adequate.
\noindent By carrying out diamond move \diamondmove{w} at $w$ on Fig. 45, disconnects $Int(K)$ into $Int (K_1) \cup Int (K_2)$ such that $K_1$ is an annulus. 
Since all the regions involved in \diamondmove{w} are from $Reg_2 (M)$,  the resulting diagram $M'$ remains in $\mathcal{M}_3(W)$. 
Let $M''$ be the diagram obtained from $M'$ by cutting out $K_1$ and sewing back $K_2$ with Q. See fig. 45. 
Since we cannot subdivide an 
$[\ ]_t$ and an $[\ ]_{H_t}$ equivalence class by the removal of $K_1, M'' \in \mathcal{M}_2(W)$. But then $|M''|_2<|M|_2$, where $|M|_2= |Reg_2(M)|$. A contradiction to $M$ being Minimal. 
Hence $\partial K$ is a simple closed curve. Since $M$ is Minimal, $\Phi(\partial K)$ is cyclically reduced. Next, assume that $B_1$ and $B_2$ are bands that intersect in two regions $\Delta_1$ and $\Delta_2$.
See Fig. 42. By theorem T there is a sequence J of $I$-moves such that $B_1^J\cap B_2^J$ is connected and $|Reg_2(M^J)|\leq |Reg_2(M)|$.
Now, by a sequence $J_1$, of I-moves we push $\Delta_1$ to $\Delta_2$ such that $|\partial \Delta_1^{J_1}\cap \partial \Delta_2^{J_2}|\geq 2$.
But then $\Delta_1^{J_1}$ cancels out $\Delta_2^J$, reducing 
$|M^J|_2$, a contradiction to Minimality.
Hence $Str (M)$ is adequate. See Fig. 37.
\hfill $\Box$
\begin{center}
	\begin{tabular}{cp{0.3cm}cp{0.3cm}cp{0.4cm}c}
		\begin{tikzpicture}[scale=0.45]
			
			
			\draw  (110:3) coordinate (v0) arc (110:430:3)
			coordinate [pos=0.1] (v1)
			coordinate [pos=0.2] (v2)
			coordinate [pos=0.3] (v3)
			coordinate [pos=0.4] (v4)
			coordinate [pos=0.5] (v5)
			coordinate [pos=0.6] (v6)
			coordinate[pos=0.7] (v7)
			coordinate[pos=0.8] (v8)
			coordinate [pos=0.9] (v9)
			coordinate (ve)
			;
			\draw  (90:2) coordinate (u0) arc (90:450:2)
			coordinate [pos=0.1] (u1)
			coordinate [pos=0.2] (u2)
			coordinate [pos=0.3] (u3)
			coordinate [pos=0.4] (u4)
			coordinate [pos=0.5] (u5)
			coordinate [pos=0.6] (u6)
			coordinate[pos=0.7] (u7)
			coordinate[pos=0.8] (u8)
			coordinate [pos=0.9] (u9)
			coordinate (ue)
			;
			\draw (v0)--(u0) node [pos=0.5,sloped,scale=0.7] {$<$};
			\draw (v1)--(u1);
			\draw (v2)--(u2);
			\draw (v3)--(u3);
			\draw (v4)--(u4);
			\draw (v5)--(u5);
			\draw (v6)--(u6);
			\draw (v7)--(u7);
			\draw (v8)--(u8);
			\draw (v9)--(u9);
			\draw (ve)--(ue) node [pos=0.5,sloped,scale=0.7] {$>$};;
			
			\node[below left]  at (v3) {$K$};
			\coordinate (v10) at (-2,5) ;
			\coordinate (v11) at (-1,5) ;
			\draw  (v10) -- (v0)
			coordinate [pos=0.25] (v1)
			coordinate [pos=0.5] (v2)
			coordinate [pos=0.75] (v3)
			;
			\draw (v11) -- (u0)
			coordinate [pos=0.25] (u1)
			coordinate [pos=0.5] (u2)
			coordinate [pos=0.75] (u3)
			;
			
			\draw (v1)--(u1);
			\draw (v2)--(u2);
			\draw (v3)--(u3);
			\draw (v10)--(v11);
			
			\coordinate (v10) at (2,5) ;
			\coordinate (v11) at (1,5) ;
			\draw  (v10) -- (ve)
			coordinate [pos=0.25] (v1)
			coordinate [pos=0.5] (v2)
			coordinate [pos=0.75] (v3)
			;
			\draw (v11) -- (ue)
			coordinate [pos=0.25] (u1)
			coordinate [pos=0.5] (u2)
			coordinate [pos=0.75] (u3)
			;
			
			\draw (v1)--(u1);
			\draw (v2)--(u2);
			\draw (v3)--(u3);
			\draw (v10)--(v11);

			\node at (0,0) {$Q$};
			\node[above=5mm] at (u0) {$w$};
			
		\end{tikzpicture}
		\hspace {-5mm} 
		& { \raisebox{1.9cm} {\huge $\rightsquigarrow$}} & 
		\begin{tikzpicture}[scale=0.45]
			
			\coordinate  (origin) at  (0,0);
			\coordinate (w) at (0,3.5);
			\node[above=10] at (w) {$K_2$};
			\draw  (110:3) coordinate (v0) arc (110:430:3)
			coordinate [pos=0.1] (v1)
			coordinate [pos=0.2] (v2)
			coordinate [pos=0.3] (v3)
			coordinate [pos=0.4] (v4)
			coordinate [pos=0.5] (v5)
			coordinate [pos=0.6] (v6)
			coordinate[pos=0.7] (v7)
			coordinate[pos=0.8] (v8)
			coordinate [pos=0.9] (v9)
			coordinate (ve)
			;
			\draw  (90:2) coordinate (u0) arc (90:450:2)
			coordinate [pos=0.1] (u1)
			coordinate [pos=0.2] (u2)
			coordinate [pos=0.3] (u3)
			coordinate [pos=0.4] (u4)
			coordinate [pos=0.5] (u5)
			coordinate [pos=0.6] (u6)
			coordinate[pos=0.7] (u7)
			coordinate[pos=0.8] (u8)
			coordinate [pos=0.9] (u9)
			coordinate (ue)
			;
			\draw (v0)--(u0) node [pos=0.5,sloped,scale=0.7] {$<$};
			\draw (v1)--(u1);
			\draw (v2)--(u2);
			\draw (v3)--(u3);
			\draw (v4)--(u4);
			\draw (v5)--(u5);
			\draw (v6)--(u6);
			\draw (v7)--(u7);
			\draw (v8)--(u8);
			\draw (v9)--(u9);
			\draw (ve)--(ue) node [pos=0.5,sloped,scale=0.7] {$>$};;
			
			\coordinate (v10) at (-2,5) ;
			\coordinate (v11) at (-1,5) ;
			\draw  (v10) -- (v0)
			coordinate [pos=0.25] (v1)
			coordinate [pos=0.5] (v2)
			coordinate [pos=0.75] (v3)
			;
			\draw (v11) -- (w)
			coordinate [pos=0.25] (u1)
			coordinate [pos=0.5] (u2)
			coordinate [pos=0.75] (u3)
			;
			\draw (v0)--(w) node [pos=0.5,sloped,scale=0.7]  {$<$};
			\draw (v1)--(u1);
			\draw (v2)--(u2);
			\draw (v3)--(u3);
			\draw (v10)--(v11);
			
			\coordinate (v10) at (2,5) ;
			\coordinate (v11) at (1,5) ;
			\draw  (v10) -- (ve)
			coordinate [pos=0.25] (v1)
			coordinate [pos=0.5] (v2)
			coordinate [pos=0.75] (v3)
			;
			\draw (v11) -- (w)
			coordinate [pos=0.25] (u1)
			coordinate [pos=0.5] (u2)
			coordinate [pos=0.75] (u3)
			;
			\draw (w)--(ve) node [pos=0.5,sloped,scale=0.7]  {$>$};
			\draw (v1)--(u1);
			\draw (v2)--(u2);
			\draw (v3)--(u3);
			\draw (v4)--(u4);
			\draw (v10)--(v11);

			\node at (0,0) {$Q$};
			
			\node at (-2.5,-.25) {$K_1$};
		\end{tikzpicture}
		%
		\hspace {-5mm}
		& {\centering \raisebox{1.9cm} {\huge $\rightsquigarrow$}} & 
		\begin{tikzpicture}[scale=0.45]
			
			
			\draw  (90:3) coordinate (v0) arc (90:450:3)
			coordinate [pos=0.1] (v1)
			coordinate [pos=0.2] (v2)
			coordinate [pos=0.3] (v3)
			coordinate [pos=0.4] (v4)
			coordinate [pos=0.5] (v5)
			coordinate [pos=0.6] (v6)
			coordinate[pos=0.7] (v7)
			coordinate[pos=0.8] (v8)
			coordinate [pos=0.9] (v9)
			coordinate (ve)
			;
			\draw  (90:2) coordinate (u0) arc (90:450:2)
			coordinate [pos=0.1] (u1)
			coordinate [pos=0.2] (u2)
			coordinate [pos=0.3] (u3)
			coordinate [pos=0.4] (u4)
			coordinate [pos=0.5] (u5)
			coordinate [pos=0.6] (u6)
			coordinate[pos=0.7] (u7)
			coordinate[pos=0.8] (u8)
			coordinate [pos=0.9] (u9)
			coordinate (ue)
			;
			\draw (v0)--(u0) ;
			\draw (v1)--(u1);
			\draw (v2)--(u2);
			\draw (v3)--(u3);
			\draw (v4)--(u4);
			\draw (v5)--(u5);
			\draw (v6)--(u6);
			\draw (v7)--(u7);
			\draw (v8)--(u8);
			\draw (v9)--(u9);
			\draw (ve)--(ue) node [pos=0.5,sloped,scale=0.7] {$<$};;
			\coordinate (w) at (0,4);
			\coordinate (v10) at (-2,5) ;
			\coordinate (v11) at (-1,5) ;
			\draw  (v10) -- (v0)
			coordinate [pos=0.25] (v1)
			coordinate [pos=0.5] (v2)
			coordinate [pos=0.75] (v3)
			;
			\draw (v11) -- (w)
			coordinate [pos=0.25] (u1)
			coordinate [pos=0.5] (u2)
			coordinate [pos=0.75] (u3)
			;
			\draw (v0)--(w) node [pos=0.5,sloped,scale=0.7] {$<$};
			\draw (v1)--(u1);
			\draw (v2)--(u2);
			\draw (v3)--(u3);
			\draw (v10)--(v11);
			
			\coordinate (v10) at (2,5) ;
			\coordinate (v11) at (1,5) ;
			\draw  (v10) -- (ve)
			coordinate [pos=0.25] (v1)
			coordinate [pos=0.5] (v2)
			coordinate [pos=0.75] (v3)
			;
			\draw (v11) -- (w)
			coordinate [pos=0.25] (u1)
			coordinate [pos=0.5] (u2)
			coordinate [pos=0.75] (u3)
			;
			
			\draw (v1)--(u1);
			\draw (v2)--(u2);
			\draw (v3)--(u3);
			\draw (v10)--(v11);

			\node at (0,0) {$Q$};
			\node at (-3,-2.5) {$K_1$};
			
			\node at (1.5,3.5) {$K_2$};
		\end{tikzpicture}
		\hspace {-5mm}
		& {\centering \raisebox{1.9cm} {\huge $\rightsquigarrow$}} & 
		\begin{tikzpicture}[scale=0.45]

			\draw  (90:3) coordinate (t) arc (90:450:3)
			coordinate [pos=0.1] (l)
			coordinate [pos=0.9] (r)
			
			;

			\coordinate (v10) at (-2,5) ;
			\coordinate (v11) at (-1,5) ;
			\draw  (v10) -- (l)
			coordinate [pos=0.25] (v1)
			coordinate [pos=0.5] (v2)
			coordinate [pos=0.75] (v3)
			;
			\draw (v11) -- (t)
			coordinate [pos=0.25] (u1)
			coordinate [pos=0.5] (u2)
			coordinate [pos=0.75] (u3)
			;
			
			\draw (v1)--(u1);
			\draw (v2)--(u2);
			\draw (v3)--(u3);
			\draw (v10)--(v11);
			
			\coordinate (v10) at (2,5) ;
			\coordinate (v11) at (1,5) ;
			\draw  (v10) -- (r)
			coordinate [pos=0.25] (v1)
			coordinate [pos=0.5] (v2)
			coordinate [pos=0.75] (v3)
			;
			\draw (v11) -- (t)
			coordinate [pos=0.25] (u1)
			coordinate [pos=0.5] (u2)
			coordinate [pos=0.75] (u3)
			;
			
			\draw (v1)--(u1);
			\draw (v2)--(u2);
			\draw (v3)--(u3);
			\draw (v10)--(v11);

			\node at (0,0) {$Q$};
			\node[above=3mm] at (0,3) {$K_2$};
		\end{tikzpicture}
		\\
		(i)\hspace {-5mm} & & (ii)\hspace {-5mm} & & (iii)\hspace {-5mm} & & (iv)
	\end{tabular}
	
	Figure 45
\end{center}
\subsection {Proofs of Theorems B and C}
\subsubsection*{Proof of Theorem B}
\begin{description}
	\item[Case 1] \quad $T(M)=\emptyset$. \\
	Then $Reg_{4^+}(M)=\emptyset$, hence $Reg(M)=Reg_2(M)$.
	Hence, if $Supp(D)=\{a,b\}$,\ $D\in Reg(M)$
	then $\|\partial D\|_a=\|\partial D\|_b=2$.
	Now, by Lemma 1.1.1(c), $\|\partial M\|_a\geq 2$.
	Hence $\|\partial M\|_a\geq \|\partial D\|_a$ and
	$\|\partial M\|_b\geq \|\partial D\|_b$, as required.
	\item[Case 2] \quad $T(M)\neq \emptyset$. \\
	Let $t\in T$ and consider $\widetilde{\mathbb M}^t$.
	Recall that for every region $\Delta$ in $\mathbb M_t$,
	$\|\partial\Delta\|_t=\|\partial\tilde{\Delta}\|$ and 
	$n(\Delta)=n(\tilde{\Delta})$, where $\tilde{\Delta}=\Psi_t(\Delta)$.
	Without loss of generality we may assume that 
	$\widetilde{\mathbb M}^t$ has connected interior.
	Hence due to Theorem A and Lemma 1.3.2, either
	$\widetilde{\mathbb M}^t=\{\widetilde{\Delta}\}$ or
	$|\mathcal{D}(\widetilde{\mathbb M}^t)|\geq 2$.
	In the first case, the result is clear.
	Assume, therefore, that 
	$\mathcal{D}(\widetilde{\mathbb M}^t)=\{\widetilde{\Delta}_1,
	\ldots,\widetilde{\Delta}_k\}$,\ $k\geq 2$.
	We show that
	\begin{equation*}\tag{1}
		\|\partial \widetilde{\mathbb M}^t\|\geq n(\Delta_j),\ j=1,\ldots,k
	\end{equation*}
	From this the result follows.
	We show this for $j=1$.
	The proof is the same for $j=2,\ldots,k$.
	Recall the notation of Lemma 5.5.2.
	By 3) in the Proof of Lemma 5.5.2, $\|\eta_0\|_t\geq n(\Delta_j)-i(\widetilde{\Delta}_j)$,\ 
	$j=1,\ldots, k$.
	Hence, by Theorem A and due to Lemma 1.3.2, the individual contribution of $\Delta_j$ to $\|\partial M\|_t$ is at least $n(\Delta_j)-2$.
	However, due to possible consolidation of labels of adjacent $\Delta_j$, the actual contribution of $\Delta_j$ to
	$\|\partial M\|_t$ is at least $n(\Delta_j)-3$.
	Observe that
	\begin{equation*}\tag{2}
		\|\partial {\mathbb M}\|_t\geq \|\eta_0\|_t
	\end{equation*}
	Consequently, due to (2),
	\begin{equation*}\tag{3}
		\|\partial {\mathbb M}\|_t\geq \sum_{j=1}^k [n(\Delta_j)-3]
	\end{equation*}
	Now, $\sum_{j=1}^k [n(\Delta_j)-3]=
	n(\Delta_1)+$
	$\left(\sum_{j=2}^k [n(\Delta_j)-3]\right)-3\geq$
	$n(\Delta_1)+\left(\sum_{j=2}^k [4-3]\right)-3=$ $n(\Delta_1)+(k-1)-3=n(\Delta_1)+k-4$, i.e.
	\begin{equation*}\tag{4}
		\|\partial {\mathbb M}\|_t\geq n(\Delta_1)+k-4\ 
		\text{(Recall that $n(\Delta_j)\geq 4$)}
	\end{equation*}
	Hence, if $k\geq 4$ then  $\|\partial {\mathbb M}\|_t\geq n(\Delta_1)$, as required.
	If $k=3$, then due to Lemma 1.3.2, $i(\widetilde{\Delta}_j)=1,\ j=1,2,3$ hence,
	$\|\partial\mathbb M\|_t\geq (n(\Delta_1)-2)+(n(\Delta_2)-2)+
	(n(\Delta_3)-2)\geq n(\Delta_1)+(4-2)+(4-2)-2>n(\Delta_1)$.
	Finally, if $k=2$ then by Lemma 1.3.2, $i(\Delta_j)=1,\ j=1,2$ hence $(n(\Delta_1)-2)+(n(\Delta_2)-2)=
	n(\Delta_1)+[(n(\Delta_2)-2)-2]\geq n(\Delta_1)$.
	i.e. due to 3), $\|\partial {\mathbb M}\|_t\geq n(\Delta_1)$, as required. \hfill $\Box$
	
\end{description}
\subsubsection*{Proof of Theorem C}
We have to show that either $\sigma(\omega_1)=1$ or 
$\sigma(\omega_2)=1$.
Let $t\in T(M)$ and consider $\widetilde{\mathbb M}^t$.
By Lemma 1.3.2, either $\widetilde{\mathbb M}^t=\{\widetilde{\Delta}\}$, or
$|\mathcal{D}(\widetilde{\mathbb M}^t)|\geq 2$.
If $|\mathcal{D}(\widetilde{\mathbb M}^t)|\geq 3$ then by the P.H.P either $|\mathcal{D}(\widetilde{\omega}_1)|\neq \emptyset$ or 
$|\mathcal{D}(\widetilde{\omega}_2)|\neq \emptyset$.
Assume $|\mathcal{D}(\widetilde{\omega}_2)|=\emptyset$.
Then $|\mathcal{D}(\widetilde{\omega}_1)|\neq \emptyset$, 
hence $\sigma(\omega_1)=1$, by Theorem L, as required.
If $|\mathcal{D}(\widetilde{\mathbb M}^t)|=2$ then by Lemma 1.3.2,
$\widetilde{\mathbb M}^t\in O.L.$.
Hence, again, by Theorem L, $\sigma(\omega_1)=1$.
Finally, if $\widetilde{\mathbb M}^t=\{\widetilde{\Delta}\}$ then by Lemma 5.4.3 $\sigma(\omega_1)=1$. \hfill $\Box$





\end{document}